\documentclass[10pt]{article}
\usepackage{amsmath}
\usepackage{amsthm}
\usepackage{amscd}
\usepackage{amssymb}
\usepackage{color}

\usepackage[all]{xy}

\definecolor{darkgreen}{cmyk}{1,0,1,.2}

\newcounter{enumitemp}
\newenvironment{enumeratecontinue}{
  \setcounter{enumitemp}{\value{enumi}}
  \begin{enumerate}
  \setcounter{enumi}{\value{enumitemp}}
}
{
  \end{enumerate}
}

\numberwithin{equation}{section}

\newcommand{\nb}[1]{#1\nobreakdash-}

\newcommand\ds\displaystyle

\theoremstyle{definition}

\newtheorem*{remark}{Remark}

\newtheorem{problem}{Problem}
\newtheorem{question}[problem]{Question}

\theoremstyle{plain}
\newtheorem*{theorem*}{Theorem}
\newtheorem{theorem}{Theorem}[section]
\newtheorem{proposition}[theorem]{Proposition}
\newtheorem{lemma}[theorem]{Lemma}

\newtheorem{corollary}[theorem]{Corollary}

\newtheorem{fact}[theorem]{Fact}

\newcounter{remarks}
{\paragraph*{Remarks}\smallskip
     \begin{list}{\arabic{remarks}. }{\usecounter{remarks}%
          \setlength{\leftmargin}{0in}%
          \setlength{\rightmargin}{0in}%
          \setlength{\labelsep}{0pt}%
          \setlength{\labelwidth}{0pt}%
          \setlength{\listparindent}{0pt}%
     }
}
{
\end{list}
}

\DeclareMathOperator{\Out}{Out}
\DeclareMathOperator{\Inn}{Inn}
\DeclareMathOperator{\Aut}{Aut}

\DeclareMathOperator{\Length}{Len}
\DeclareMathOperator{\Isom}{Isom}
\DeclareMathOperator\Min{Min}
\DeclareMathOperator\diam{diam}
\DeclareMathOperator\image{image}

\DeclareMathOperator\Fix{Fix}
\DeclareMathOperator\Per{Per}
\DeclareMathOperator{\Lip}{Lip}
\DeclareMathOperator{\valence}{valence}

\newcommand\reals{{\mathbf R}}
\newcommand\R\reals
\newcommand\real\reals
\newcommand\hyp{{\mathbf H}}

\newcommand\Z{{\mathbf Z}}
\newcommand\B{\mathcal B}
\newcommand\inject{\hookrightarrow}
\newcommand\homeo{\approx}
\newcommand\infinity{\infty}
\newcommand{\bdy}{\partial}
\newcommand{\from}{\colon}
\newcommand\composed{\circ}
\newcommand\suchthat{\bigm|}
\newcommand\inv{{-1}}
\newcommand\union{\cup}

\newcommand\abs[1]{\left| #1 \right|}

\newcommand\Id{\text{Id}}
\newcommand\intersect{\cap}
\newcommand\restrict{\bigm|}
\newcommand\Teichmuller{Teichm\"uller}
\newcommand\Poincare{Poincar\'e}

\newcommand\cross{\times}

\newcommand\A{\mathcal A}
\newcommand\C{\mathcal C}
\newcommand\F{\mathcal F}
\newcommand\G{\mathcal G}
\renewcommand\L{\mathcal L}

\renewcommand\P{\mathcal P}
\newcommand\Teich{\mathcal T}
\newcommand\<\langle
\renewcommand\>\rangle
\newcommand\wt\widetilde
\newcommand\wh\widehat

\DeclareMathOperator\interior{int}

\newcommand\act\circlearrowright

\DeclareMathOperator\GFL{\mathcal{GF}}
\DeclareMathOperator\GL{\mathcal{GL}}
\DeclareMathOperator\BCC{BCC}
\DeclareMathOperator{\inj}{inj}

\newcommand\FF{\mathcal{FF}}

\newcommand\X{\mathcal{X}}

\newcommand\CV\X

\newcommand\BookZero{\cite{BFH:laminations}}
\newcommand\BookOne{\cite{BFH:TitsOne}}

\newcommand\OuterLimits{\cite{BestvinaFeighn:OuterLimits}}

\newcommand\ti {\tilde}
\newcommand\gG{g \from G \to G}
\DeclareMathOperator{\PA}{PA}
\newcommand{\trtr}{train track}   

\newcommand\iNp{indivisible Nielsen path}
\newcommand\ipNp{indivisible periodic Nielsen path}

\newcommand\M{\mathcal{M}}

\newcommand\whclosure[1]{\widehat{\overline #1}}

\DeclareMathOperator\Domain{Dom}
\DeclareMathOperator\Range{Ran}
\newcommand\fol{\mathcal{F}}

\newcommand\T{\mathcal{T}}
\newcommand\PMF{\mathcal{PMF}}
\DeclareMathOperator\TT{TT}

\newcommand\vpp{{\vphantom\prime}}

\newcommand\TphiPlus{T_+^\phi}
\newcommand\TPlusphi\TphiPlus
\newcommand\TphiMinus{T_-^\phi}
\newcommand\TMinusphi\TphiMinus
\newcommand\TPlus{T_+}
\newcommand\TMinus{T_-}

\newcommand\LambdaMinusphi{\Lambda_-^\phi}
\newcommand\LambdaPlus{\Lambda_+}
\newcommand\LambdaMinus{\Lambda_-}
\newcommand\TPM{T_\pm}
\newcommand\TPMphi{T_\pm^\phi}
\newcommand\TphiPM\TPMphi
\newcommand\LambdaPM{\Lambda_\pm}
\newcommand\LambdaPMphi{\Lambda_\pm^\phi}

\title{Axes in Outer Space}
\author{Michael Handel and Lee Mosher}

\begin{document}
\maketitle

\begin{abstract}
We develop a notion of axis in the Culler--Vogtmann outer space $\X_r$ of a finite rank free group $F_r$, with respect to the action of a nongeometric, fully irreducible outer automorphism $\phi$. Unlike the situation of a loxodromic isometry acting on hyperbolic space, or a pseudo-Anosov mapping class acting on \Teichmuller\ space, $\X_r$ has no natural metric, and $\phi$ seems not to have a single natural axis. Instead our axes for $\phi$, while not unique, fit into an ``axis bundle'' $\A_\phi$ with nice topological properties: $\A_\phi$ is a closed subset of $\X_r$ proper homotopy equivalent to a line, it is invariant under $\phi$, the two ends of $\A_\phi$ limit on the repeller and attractor of the source--sink action of $\phi$ on compactified outer space, and $\A_\phi$ depends naturally on the repeller and attractor. We propose various definitions for $\A_\phi$, each motivated in different ways by train track theory or by properties of axes in \Teichmuller\ space, and we prove their equivalence.
\end{abstract}

\tableofcontents

\section{Introduction}
\label{SectionIntro}

There are many interesting and fruitful analogies between the group of isometries of hyperbolic $n$-space $\hyp^n$, the mapping class group of a finite type surface acting on \Teichmuller\ space, and the outer automorphism group $\Out(F_r)$ of a rank $r$ free group $F_r$ acting on outer space $\X_r$. 

In all cases the space has a natural compactification, obtained by adding a boundary at infinity to which the group action extends continuously. In all cases there is a class of ``hyperbolic'' elements whose extended action has source--sink dynamics, with a repeller and an attractor at infinity: loxodromic isometries of~$\hyp^n$; pseudo-Anosov mapping classes of a surface~\cite{McCarthy:Tits}; and outer automorphisms of $F_r$ which are \emph{fully irreducible}~\cite{LevittLustig:NorthSouth}, meaning that no proper, nontrivial free factor of $F_r$ has a conjugacy class which is periodic under the action of the outer automorphism.

In $\hyp^n$ and in \Teichmuller\ space, a ``hyperbolic'' element $\phi$ has a unique axis in the usual metric sense, a properly embedded, oriented, $\phi$-invariant geodesic line which $\phi$ translates in the forward direction. The axis of $\phi$ depends only on the ordered pair $(\text{repeller},\text{attractor})$, and this dependence is natural with respect to the group action: the axis is the unique geodesic line whose negative end converges to the repeller and whose positive end converges to the attractor.

In $\X_r$, however, it is not reasonable to expect such a nice notion of an axis for a fully irreducible $\phi \in \Out(F_r)$. There does not seem to be a natural candidate for a metric on $\X_r$, and certainly not one that has the unique geodesic properties of $\hyp^n$ or \Teichmuller\ space. So among the many properly embedded lines in $\X_r$ whose two ends converge to the repeller and the attractor of $\phi$, it seems hard to single out one of those lines in a natural way and call it \emph{the} axis of $\phi$. Nevertheless, among all such lines we shall select in a natural manner a particular collection of them, the union of which we call the ``axis bundle'' of $\phi$, a subset $\A_\phi\subset\X_r$ with the following properties:
\begin{description}
\item[Properness:] $\A_\phi$ is a closed subset of $\X_r$ proper homotopy equivalent to $\reals$. The negative end of $\A_\phi$ limits in compactified outer space on the repeller, and the positive end limits on the attractor.
\item[Naturality:] $\A_\phi$ depends only on the ordered pair $(\text{repeller},\text{attractor})$ in a manner which is natural with respect to the action of $\Out(F_r)$.
\end{description}
We shall do this only in the case where $\phi$ is \emph{nongeometric}, meaning that it does not arise from a homeomorphism of a compact surface with boundary. The geometric case, while conceptually much simpler and more well understood (see \cite{BestvinaHandel:tt}), has some peculiarities whose inclusion in our theory would overburden an already well laden paper.

Naturality means that for any nongeometric, fully irreducible $\phi,\phi' \in \Out(F_r)$ and any $\psi \in \Out(F_r)$, if $\psi$ takes the repeller--attractor pair of $\phi$ to the repeller--attractor pair of~$\phi'$, then $\psi$ takes $\A_\phi$ to $\A_{\phi'}$. One consequence of naturality is that $\A_\phi$ is invariant under~$\phi$, because the repeller--attractor pair of $\phi$ is invariant under~$\phi$. Another is that axes are indifferent to iterates (positive powers) of~$\phi$: $\phi^n$ and $\phi$ have the same axis for all $n>0$, because $\phi^n$ has the same repeller and the same attractor as~$\phi$. 

Certain differences between axes in hyperbolic space or \Teichmuller\ space on the one hand, and axis bundles in outer space on the other hand, may already be apparent from the above features. First, the topological relation between the axis bundle and the real line is not homeomorphism but instead proper homotopy equivalence. Second, axes in $\hyp^n$ or \Teichmuller\ space are indifferent to inverses --- the axis of $\phi^\inv$ is the orientation reversal of the axis of $\phi$. But as we have argued in \cite{HandelMosher:parageometric}, it seems unnatural to insist that axes in $\X_r$ be indifferent to inverses, due to the phenomenon that a fully irreducible $\phi \in \Out(F_r)$ can have an expansion factor which is different than the expansion factor of its inverse. Although our paper \cite{HandelMosher:logratio} offers a bare hint of an ``indifference to inverses'' property, by bounding the ratio of the logarithms of the expansion factors of $\phi$ and $\phi^\inv$, nevertheless we do not attempt here to construct a notion of axis which is indifferent to inverses.

The main results of this work are a theorem which proves equivalence of several characterizations of the axis bundle --- from some of which the naturality feature of the axis bundle is evident --- and a theorem which proves the properness feature of the axis bundle. In the remainder of this introduction we shall describe the three characterizations of the axis bundle, and the methods of proof. 

\subsection{Characterizations of the axis bundle.} We shall characterize the axis bundle by presenting three different definitions, the equivalence of which is one of our major theorems. These definitions are motivated by various aspects of the theory of outer automorphisms, and by analogies with various properties of \Teichmuller\ geodesics. 

\paragraph{Fold lines.} Fold lines form a natural class of paths in outer space arising from several concepts: train track maps due to Bestvina and Handel \cite{BestvinaHandel:tt}; Stallings folds \cite{Stallings:folding}; and semiflow concepts behind Skora's method for investigating homotopy properties of spaces of $\reals$-trees (\cite{Skora:deformations}, and see also \cite{White:FixedPoints}, \cite{Clay:contractibility}, \cite{GuirardelLevitt:outer}).

Recall that elements of compactified outer space $\overline\X_r = \X_r \union \bdy\X_r$ are represented by minimal, very small actions of $F_r$ on $\reals$-trees, what we shall call ``$F_r$-trees''. Two $F_r$-trees represent the same element of $\overline\X_r$ if and only if they are related by a homothety (a bijection that multiplies the metric by a constant) that is $F_r$-equivariant. An element of $\X_r$ itself is represented by a simplicial $F_r$-tree on which the action of $F_r$ is free. Passing to the quotient modulo $F_r$ of a free, simplicial $F_r$-tree, an element of outer space is also represented by a ``marked graph'', a metric graph $G$ without valence~1 vertices whose fundamental group is marked by an isomorphism $F_r \approx \pi_1(G)$; while we allow $G$ to have valence~2 vertices, their presence or absence does not effect the element of $\X_r$ that $G$ represents. The group $\Out(F_r)$ acts on $\overline\X_r$ from the right: the image under $\phi \in \Out(F_r)$ of a marked graph $G$ is a marked graph $G \cdot \phi$ defined by precomposing the marking $F_r \approx \pi_1(G)$ with an automorphism of $F_r$ that represents~$\phi$. 

Given a fully irreducible $\phi \in \X_r$, the repeller and attractor for $\phi$ are elements of $\bdy\X_r$ represented by $F_r$-trees that we denote $\TMinus=\TphiMinus$ and $\TPlus=\TphiPlus$.

A \emph{fold line} in outer space is a continuous, injective, proper function $\reals \mapsto \X_r$ defined by a continuous 1-parameter family of marked graphs $t \mapsto G_t$ for which there exists a family of homotopy equivalences $h_{ts} \from G_s \to G_t$ defined for $s \le t \in \reals$, each of which preserves marking, such that the following hold:

\break

\begin{itemize}
\item[] \textbf{Train track property:} $h_{ts}$ is a local isometry on each edge for all $s \le t \in \reals$.
\item[] \textbf{Semiflow property:} $h_{ut} \composed h_{ts} = h_{us}$ for all $s \le t \le u\in \reals$, and $h_{ss} \from G_s \to G_s$ is the identity for all $s \in \reals$.
\end{itemize}
Our first definition of the axis bundle is:
\begin{description}
\item[Definition (1): Fold lines.] \quad $\A_\phi$ is the union of the images of all fold lines 
$\gamma \from \reals \to \X_r$ such that $\gamma(t)$ converges in $\overline\X_r$ to $\TphiMinus$ as $t \to -\infinity$ and to $\TphiPlus$ as $t \to +\infinity$.
\end{description}
Notice that, just from the format of this definition, the axis bundle $\A_\phi$ has the desired naturality feature, because the collection of fold lines in $\X_r$ is invariant under the action of $\Out(F_r)$, and the points $\TphiMinus, \TphiPlus \in \bdy\X_r$ depend naturally on $\phi$.

The construction of $\phi$-periodic fold lines is given in detail in the beginning of Section~\ref{SectionFoldLines}, drawing on well known ideas reviewed in Section~\ref{SectionPaths}. Here is an outline of the construction. 

Recall that every fully irreducible $\phi \in \Out(F_r)$ is represented by an \emph{affine train track map}, a homotopy equivalence $g \from G \to G$ of a marked graph such that $g$ takes vertices to vertices, $g$ induces the outer automorphism $\phi$ on $\F_r \approx \pi_1(G)$, and there exists $\lambda>1$ so that the restriction of each iterate $g^i$ to the interior of each edge locally stretches the metric by the factor $\lambda^i$ (\cite{BestvinaHandel:tt}, and see Section~\ref{SectionTrainTracks}). The number $\lambda=\lambda(\phi)$ is independent of $g$ and is called the \emph{expansion factor} of $\phi$. Let $G_i = \frac{1}{\lambda^i} G \cdot \phi^i$ be renormalized marked graphs representing the orbit of $G$ under powers of $\phi$. For integers $i \le j$ the iterate $g^{j-i} \from G \to G$ of the map $g$ induces a marking preserving homotopy equivalence $g_{ji} \from G_i \to G_j$, and these maps satisfy the semiflow property, which is just a fancy way of saying that the following diagram commutes:
$$\xymatrix{
\cdots \ar[r] & G_i \ar[r]^{g_{i+1,i}} \ar@/_2pc/[rrrr]^{g_{ji}} & G_{i+1} \ar[r] & \cdots \ar[r] & G_{j-1} \ar[r]^{g_{j,j-1}} & G_j \ar[r] & \cdots
}
$$
Moreover, the definition of train track maps and the renormalization of $G_i$ together imply that $g_{ji}$ satisfies the above train track property. Now one uses Stallings method of folds to factor the map $g_{1,0} \from G_0 \to G_1$ as a product of fold maps, and then one interpolates each fold map continuously to obtain a path in $\X_r$ from $G_0$ to $G_1$. Concatenating the $\phi^i$ translates of this path for all $i \in \Z$, one obtains a $\phi$-periodic fold line in $\X_r$. Source--sink dynamics on compact subsets of outer space \BookZero\ implies that in $\overline\X_r$ this fold line converges in the negative direction to $\TMinus$ and in the positive direction to $\TPlus$. This shows that the point of $\X_r$ represented by $G$ satisfies Definition~(1) of $\A_\phi$.

Notice that we do \emph{not} claim that $\A_\phi$ is the union of all $\phi$-periodic fold lines. In fact, we do not even claim that $\A_\phi$ is the union of all $\phi^i$-periodic fold lines over all $i \ge 1$. See Section~\ref{SectionTTversusWTT} for an example of $G \in \A_\phi$ that does lie on any $\phi^i$-periodic fold line for any iterate~$\phi^i$. In particular, such a $G$ is not the domain of a train track representative of any iterate~$\phi^i$, nor is $G$ the domain of a \emph{nonsimplicial train track representative}\footnote{also called an ``efficient representative'' in an unpublished preprint of T.~White and in \cite{LosLustig:axes}.} of any iterate, a map that satisfies the definition of a train track representative except that it is not required to take vertices to vertices.

\subparagraph{An analogy with \Teichmuller\ geodesics.} Fold lines can also be motivated by an analogy with \Teichmuller\ geodesics, as follows. Recall that if $S$ is a closed oriented surface of genus $\ge 2$, each oriented geodesic $\gamma$ in the \Teichmuller\ space $\Teich(S)$ has the following essentially unique description. Choosing a unit speed parameterization $t \mapsto \sigma_t$ of $\gamma$, there is a transverse pair of measured foliations $\fol_-,\fol_+$ with transverse measures denoted $\abs{dy}$, $\abs{dx}$, respectively, such that for each $t \in \reals$ the singular Euclidean metric $ds_t = \sqrt{(e^t \abs{dx})^2 + (e^{-t} \abs{dy})^2}$ agrees with the conformal structure~$\sigma_t$. If $\abs{dx}$, $\abs{dy}$ are normalized so that $\int_S \abs{dx} \cdot \abs{dy}=1$ then the pair $\fol_-,\fol_+$ depends up to joint isotopy only on the parameterization of $\gamma$. Also, in Thurston's boundary $\bdy\Teich(S)=\PMF(S)$, the projective measured foliation space, the pair of points that is represented by $\fol_-,\fol_+$ depends only on $\gamma$ as an oriented geodesic in $\Teich(S)$. Notice that the identity map on $S$ underlies a system of maps $\sigma_s \mapsto \sigma_t$ that satisfies the semiflow property. Consider now the renormalized metric $ds'_t=e^{-t} ds_t = \sqrt{(\abs{dx})^2 + (e^{-2t} \abs{dy})^2}$, also in agreement with the conformal structure~$\sigma_t$. For each leaf $L$ of $\fol_-$ and each $t \le u \in \reals$, the restrictions to $L$ of the metric $ds'_t$ and the metric $ds'_u$ are in agreement, because each agrees with the restriction of the transverse measure $\abs{dx}$. To put it another way, each of the maps $\sigma_t \mapsto \sigma_u$ restricts to a local isometry on each leaf $L$ of $\fol_-$, with respect to the metrics $ds'_t$ and $ds'_u$. Thus we have demonstrated an analogy with the maps $G_s \mapsto G_t$ in the definition of a fold line, each of which restricts to a local isometry on each edge of $G_s$.

An even stronger analogy between \Teichmuller\ geodesics and fold lines will be used below to motivate Definition~(2).

\paragraph{Laminations, trees, and $\LambdaMinus$ isometries.} Thurston compactified \Teichmuller\ space with projective classes of measured laminations \cite{FLP}. Skora proved a duality theorem between $\reals$-trees and measured foliations, which provides an alternative description of the boundary of \Teichmuller\ space in terms of $\reals$-trees \cite{Skora:splittings}. The boundary of outer space $\bdy\X_r$ is described solely in terms of $\reals$-trees, no general theory of measured laminations having so far arisen. 

However, in \BookZero\ Bestvina, Feighn, and Handel define the \emph{expanding lamination} of a fully irreducible $\phi \in \Out(F_r)$ in terms of train track representatives of $\phi$, and they prove a duality relation between the expanding lamination and $\TMinus$ which is a natural analogue of Skora's duality. This relation is described by saying that the expanding lamination of $\phi$ is the ``length zero lamination'' of $\TMinus$ (see Section~3 of \BookZero, and Lemma~\ref{LemmaMinimalLamination} below). Because of this duality we shall denote the expanding lamination of $\phi$ as $\LambdaMinus = \LambdaMinusphi$. The lamination $\LambdaMinus$ can be regarded abstractly as an $F_r$-equivariant subset of $\bdy F_r \cross \bdy F_r$, or it may be realized as a family of leaves --- bi-infinite geodesics --- in any free, simplicial $F_r$-tree~$T$; see Sections~\ref{SectionLaminations} and~\ref{SectionExpandingLamination} below. 

Given a free, simplicial $F_r$-tree $T$, an $F_r$-equivariant map $f_T \from T \to \TPlus$ is called a \emph{$\LambdaMinus$ isometry} if, for each leaf $L$ of $\LambdaMinus$ realized in $T$, the restriction of $f_T$ to $L$ is an isometry onto a bi-infinite geodesic in $ \TPlus$.

Here is our next definition of the axis bundle:
\begin{description}
\item[Definition (2): Supports of $\LambdaMinus$ isometries to $\TPlus$.] \quad $\A_\phi$ is the set of all elements of $\X_r$ represented by a free simplicial $F_r$ tree $T$ for which there exists a $\LambdaMinus$ isometry $f_T \from T \to \TPlus$.
\end{description}

For example, given any train track representative $g \from G \to G$ of $\phi$, the universal cover $T = \wt G$, regarded as a free, simplicial $F_r$-tree, satisfies Definition (2). The details of why this is true are given in Corollary~\ref{CorollaryLeavesInTPlus}, but here is the idea. By lifting to the universal cover $T_i$ of each marked graph $G_i$ in the diagram above, we obtain $F_r$-equivariant maps $\ti g_{ji} \from T_i \to T_j$ that satisfy the semiflow property for integer values $i \le j$, that is, the following diagram commutes for all $i \le j$:
$$\xymatrix{
\cdots \ar[r] & T_i \ar[r]^{\ti g_{i+1,i}} \ar@/_2pc/[rrrr]^{\ti g_{ji}} & T_{i+1} \ar[r] & \cdots \ar[r] & T_{j-1} \ar[r]^{\ti g_{j,j-1}} & T_j \ar[r] & \cdots
}
$$
Each map $\ti g_{ji}$ is an $F_r$-equivariant lift of the train track map $g^{j-i} \from G \to G$, and the restriction of $\ti g_{ji}$ to each edge of $T_i$ is an isometry. In this situation, the trees $T_i$ and the maps $\ti g_{ji}$ form a direct system whose direct limit is an $F_r$ tree representing $\TPlus$, and from the direct limit construction one obtains a map $f_i \from T_i \to \TPlus$ for each~$i$. Moreover, from the construction of the expanding lamination $\LambdaMinus$ it follows that for each leaf $L$ of $\LambdaMinus$ realized in each $T_i$, the restriction of $g_{ji}$ to $L$ is an isometry. Passing to the direct limit, it follows that the restriction of $f_i$ to $L$ is an isometry, and so the map $f_i$ demonstrates that the point of $\X_r$ represented by $T_i$ or by $G_i$ satisfies Definition~(2) of $\A_\phi$.

\subparagraph{Another analogy with \Teichmuller\ geodesics.}
Further motivation for Definition~(2) comes again by an analogy with geodesics in \Teichmuller\ space $\T(S)$. Consider such a geodesic $\gamma$ with all the associated data as above: $\sigma_t$, $\fol_-$, $\fol_+$, $\abs{dy}$, $\abs{dx}$, $ds_t$, $ds'_t$. Lifting to the universal cover $\wt S$, the measured foliation $\wt \fol_+$ with transverse measure $\wt{\abs{dx}}$ is dual to an $\reals$-tree $T_+$ on which $\pi_1(S)$ acts: as a point set $T_+$ is obtained from $\wt S$ by collapsing each leaf of $\wt\fol_+$ to a point; as a metric space $T_+$ is obtained by thinking of $\wt{\abs{dx}}$ as a pseudo-metric on $\wt S$ and passing to the associated metric space obtained by collapsing to a point each maximal subset of pseudo-metric diameter zero. Let $f \from \wt S \to T_+$ be the collapse map. For each $t$, the conformal structure $\sigma_t$ lifts to a $\pi_1(S)$ equivariant conformal structure $\wt\sigma_t$, and we write the map $f$ in the form $f_t \from \wt\sigma_t \to T_+$. Under these conditions, the restriction of $f_t$ to each leaf $L$ of $\wt\fol_-$ realized in the singular Euclidean structure $ds'_t$ is an isometry, because the metric $ds'_t$ restricted to $L$ equals the transverse measure $\wt{\abs{dx}}$ restricted to $L$. Thus we can say that the map $f_t$ is an $\fol_-$-isometry (with respect to the metric $ds'_t$), and so each point $\sigma_t$ on $\gamma$ satisfies a property that is an analogue of Definition~(2). To state this property precisely, given a geodesic $\gamma$ determined by the data above, we say that \emph{$\sigma \in \Teich(S)$ satisfies property} ($2'$) if and only if there exists a singular Euclidean structure $ds$ agreeing with $\sigma$, and there is a map $f \from \wt S \to T_+$, such that locally $f$ is modelled on the projection of a Euclidean rectangle to one of its sides, and for each leaf $L$ of $\fol_-$ realized as a $\wt{ds}$ geodesic in $\wt S$, the restriction of $f$ to $L$ is an isometry onto a bi-infinite geodesic in $T_+$. We have shown that each point on $\gamma$ satisfies ($2'$), and the converse is easily seen to hold as well. Points of $\gamma$ are thus characterized by property ($2'$), just as points in the axis bundle $\A_\phi$ are characterized by~(2).

\paragraph{Train tracks.} Here is our final definition of the axis bundle, which could have been stated first but has been relegated to the last. Given a fully irreducible $\phi \in \Out(F_r)$, while ``train track representatives'' of $\phi$ are well established objects of study, ``train tracks'' themselves are not. We define a \emph{train track} for $\phi$ to be a point of $\X_r$ that is represented by a marked graph $G$ such that there exists an affine train track representative $g \from G \to G$ of $\phi$. Let $\TT(\phi) \subset \X_r$ be the set of train tracks of $\phi$.
\begin{description}
\item[Definition (3): Closure of train tracks.] $\A_\phi$ is the closure of $\displaystyle \bigcup_{i=1}^\infinity \TT(\phi^i)$ in $\X_r$.
\end{description}

\subparagraph{Remark.} In Section~\ref{SectionWTTexamples} we produce an element of $\A_\phi$ which is not a train track. We also believe there are easy examples showing that the closure of $\TT(\phi)$ itself need not equal $\A_\phi$, although we have not produced such an example. This raises the question: does there exist some $N \ge 1$ such that $\A_\phi$ is the closure of $\displaystyle \bigcup_{i=1}^N \TT(\phi^i)$?

\subsection{The main theorems.} 

Having proposed three definitions of the axis bundle, our first main result shows that these three definitions are equivalent:

\begin{theorem}
\label{TheoremEquivalence}
For each nongeometric, fully irreducible $\phi \in \Out(F_r)$ with repelling and attracting points $\TMinus,\TPlus \in \bdy\X_r$ and length zero laminations $\LambdaMinus,\LambdaPlus$, definitions~(1), (2), and~(3) of $\A_\phi$ are equivalent. In other words, for any $G \in \X_r$ the following are equivalent:
\begin{description}
\item[(1)] $G$ lies on some fold line $t \mapsto G_t$, $t \in \reals$, such that $G_t \to \TMinus$ as $t \to -\infinity$ and $G_t \to \TPlus$ as $t \to +\infinity$. 
\item[(2)] There exists a $\LambdaMinus$ isometry from the universal cover of $G$ to~$\TPlus$.
\item[(3)] $G$ is in the closure of the set of train tracks of iterates of $\phi$.
\end{description}
\end{theorem}

The set $A_\phi$ defined by this theorem is called the \emph{axis bundle} of $\phi$. While Definitions~(1) and~(3) may seem the most intuitive, Definition~(2) occupies an important middle ground between intuition and rigor, being easier to work with when proving certain properties, in particular properness. For much of the paper we focus on Definition~(2), and in Section~\ref{SectionWeakTT} we introduce the terminology ``weak train track'' to refer to $G$ satisfying Definition~(2). Of course, once Theorem~\ref{TheoremEquivalence} is proved, then a weak train track is any element of the axis bundle~$\A_\phi$. 

Note particularly that the natural dependence of $\A_\phi$ on the ordered pair $(\TMinus,\TPlus) \in \bdy\X_r \cross \bdy\X_r$ is an immediate consequence of definition~(1) and also of definition~(2). 

Our next main result establishes the properness feature of the axis bundle. Recall that a proper map $f \from A \to B$ is a \emph{proper homotopy equivalence} if there are proper maps $\bar f \from B \to A$, $H \from A \cross [0,1] \to B$, and $\bar H \from B \cross [0,1] \to A$ such that $H$ is a homotopy between $\bar f \composed f$ and $\Id_A$, and $\bar H$ is a homotopy between $f \composed \bar f$ and $\Id_B$.

\begin{theorem}
\label{TheoremProperness}
If $\phi \in \Out(F_r)$ is nongeometric and fully irreducible, and if $\gamma \from \reals \to \X_r$ is a fold line such that $\gamma(t) \to \TphiPM$ as $t \to \pm\infinity$, then the inclusion map $\gamma(\reals) \inject \A_\phi$ is a proper homotopy equivalence. Moreover, the end of $A_\phi$ corresponding to the negative end of $L$ converges in $\overline\X_r$ to $\TMinus$, and the end of $\A_\phi$ corresponding to the positive end of $L$ converges in $\overline\X_r$ to $\TPlus$.
\end{theorem}

\subsection{A question of Vogtmann.} In \cite{Vogtmann:OuterSpaceSurvey} Vogtmann asked: 
\begin{itemize}
\item[] ``Is the set of train tracks for an irreducible automorphism contractible?'' 
\end{itemize}
We address this question only for the case of a nongeometric, fully irreducible $\phi \in \Out(F_r)$. The question requires some interpretation, because the notion of ``train tracks'' for $\phi$ is not well established in the literature. The notion of train tracks that we have defined above --- points in outer space that support affine, irreducible train track representatives of $\phi$ --- is not appropriate, because the set of train tracks of $\phi$ has only countably many orbits under the action of the cyclic group $\<\phi\>$ on $\X_r$. 

By interpreting Vogtmann's question instead using weak train tracks --- elements of the axis bundle of $\phi$ --- Theorem~\ref{TheoremProperness} gives a strong positive answer to her question, establishing not just the homotopy type but the property homotopy type of the set of weak train tracks~$A_\phi$. 

\paragraph{Relations to work of Los and Lustig.} Another approach to Vogtman's question is given in a preprint of Los and Lustig \cite{LosLustig:axes}, who define a space of nonsimplicial train track representatives which they prove is contractible. One of our key technical results, Proposition~\ref{PropWeakTrainTrack}, is very similar to Lemma~3.16 of \cite{LosLustig:axes} which also plays a key role in their results. Our results differ from theirs in several aspects. They study a space of maps between marked graphs representing elements of $\X_r$, whereas our axis bundle is a subset of $\X_r$; see also the remarks at the end of Section~\ref{SectionGHConvergence}. Their proof uses an interesting semiflow which is necessarily different from the semiflow used by Skora in \cite{Skora:deformations}, whereas we prove proper homotopy equivalence to the line by plugging directly into Skora's theorem; see Section~\ref{SectionSkora}. Finally, in Section~\ref{SectionWTTexamples} we give an example of a weak train track that is not the domain of a nonsimplicial train track representative; hence, although the subset of $\X_r$ consisting of the domains of nonsimplicial train track representative can be shown to be a dense, contractible subset of the axis bundle, this subset is not equal to the whole axis bundle, and so it is not closed in $\X_r$, and its proper homotopy type is not evident.

\subsection{Contents and proofs.}

First we describe the contents of the individual sections, then we prove the two main theorems by referring to results found in various sections.

\paragraph{Section \ref{SectionPreliminaries}, Preliminaries.} We give very detailed preliminaries, some of which is a compilation of previous work. Some details give new geometric perspectives of familiar objects such as train track representatives (Section~\ref{SectionTrainTracks}) and laminations (Sections~\ref{SectionLaminations} and~\ref{SectionExpandingLamination}). Other details arise in making explicit what was only implicit in previous works, and in providing detailed proofs of results that, while previously evident, have not previously been recorded (see in particular Sections~\ref{SectionAttractingTree}, Section~\ref{SectionExpandingLamination}, and Section~\ref{SectionRelativeLeavesToTrees}). 

\paragraph{Section \ref{SectionIdealWhiteheadGraph}, The ideal Whitehead graph.} 
In this section and the next we introduce the study of the \emph{ideal Whitehead graph} of a nongeometric, fully irreducible $\phi \in \Out(F_r)$. Intuitively this graph describes topological and combinatorial aspects of the singularities of~$\phi$. The graph has one edge for each singular leaf of the expanding lamination of~$\phi$, with two singular half-leaves impinging on the same vertex whenever those half-leaves are asymptotic. The analogous construction carried out for a singularity of the unstable foliation of a pseudo-Anosov surface homeomorphism produces a circle with three or more vertices. In $\Out(F_r)$ the topology of a component of the ideal Whitehead graph can be different from a circle, as we show by an example in Section~\ref{SectionIWGExample}.

The ideal Whitehead graph of $\phi$ is defined in Section~\ref{SectionIWGDefinition} in terms of the actions of automorphisms representing $\phi$ on the boundary of $F_r$. In terms of a train track representative $g \from G \to G$, we show that the ideal Whitehead graph of $\phi$ decomposes into the local stable Whitehead graphs of $g$, in a pattern dictated by the Nielsen paths of $g$. Sections \ref{SectionAsymptoticLeaves} and~\ref{SectionIWGandTPlus} give other descriptions of the ideal Whitehead graph, in terms of: asymptotic relations among leaves of $\Lambda_-$; and branch points of the attracting $F_r$-tree $T_+$.

\paragraph{Section \ref{SectionCutAndPaste}, Cutting and pasting local stable Whitehead graphs.} Given a nongeometric, fully irreducible $\phi \in \Out(F_r)$ and a train track representative $g \from G \to G$, we study how the decomposition of the ideal Whitehead graph of $\phi$ into local stable Whitehead graphs of $g$ varies under some standard methods of varying $g$ that were first described in \cite{BestvinaHandel:tt}. This variation takes the form of cutting and pasting local stable Whitehead graphs, according to the formation or collapse of Nielsen paths of $g$. The main result of this section, Proposition~\ref{no cut vertices}, describes in a natural way those train track representatives which give the finest possible decomposition of the ideal Whitehead graph.

\paragraph{Section \ref{SectionWeakTT}, Weak train tracks.} In this section we take up in earnest the study of the axis bundle, by proving several results about weak train tracks. 
We show that many structural features of train tracks, e.g.\ the induced local decomposition of the ideal Whitehead graph, have analogues for weak train tracks. 

A key technical result, Proposition~\ref{PropWeakTrainTrack}, says roughly that if one fixes a train track $G$ then for every weak train track $H$, one can obtain $H$ from $G$ by a sequence of Stallings folds provided two conditions hold: the local decomposition of $H$ is no finer than the local decomposition of $G$; and the length of $H$ is sufficiently small in terms of the geometry of $G$. As mentioned above, this result is closely related to a key result of \cite{LosLustig:axes}.

Another important result in this section, Proposition~\ref{characterize train track}, shows how to tell when a weak train track is a train track. As a corollary we obtain one of the implications of Theorem~\ref{TheoremEquivalence} by proving that a dense set of weak train tracks are train tracks.

Also, in Section~\ref{SectionIrrigidity} we characterize when a weak train track $H$ is the domain of a unique $\Lambda_-$ isometry, and more generally we describe the space of all $\Lambda_-$ isometries with domain $H$.

\paragraph{Section \ref{SectionTopology}, Topology of the axis bundle.} This section proves several of our central results. 

Theorem~\ref{TheoremLengthPHE} describes a proper homotopy equivalence, the \emph{length map}, from the set of weak train tracks to the positive real axis $(0,\infinity)$. This theorem represents most of the work in proving the Properness Theorem~\ref{TheoremProperness}. The proof is an application of Skora's semiflow methods \cite{Skora:deformations}, which have also played a central role in other investigations of the homotopy types of spaces of $\reals$-trees  \cite{White:FixedPoints}, \cite{Clay:contractibility}, \cite{GuirardelLevitt:outer}. We shall directly apply the main theorem of \cite{Skora:deformations}, for which we must set up the appropriate maps and verify the hypotheses of that theorem. This is accomplished by applying the main results of the earlier sections, in particular Proposition~\ref{PropWeakTrainTrack} and the study of rigidity and irrigidity carried out in Section~\ref{SectionIrrigidity}.

Also proved in this section is the fact that the set of weak train tracks is a closed subset of outer space, a step in the proof of Theorem~\ref{TheoremEquivalence} (see below).

\paragraph{Section \ref{SectionFoldLines}, Fold Lines.} In this section we initiate a detailed study of fold lines. The main result proved is the equivalence of Definition~(1) of the axis bundle with the other definitions of the axis bundle. The hard work here is the proof that definition~(1) implies definition~(2): every point on a fold line is a weak train track. To do this, given a fold line $t \mapsto G_t$ in $\X_r$, using only the assumption that $\lim_{t \to \pm\infinity} G_t = T_\pm$ in $\overline\X_r$ we show that both the expanding lamination $\Lambda_-$ and the attracting tree $T_+$ can be reconstructed by natural geometric processes stated in terms of the marked graphs $G_t$ and the edge isometries $h_{ts} \from G_s \to G_t$. These processes allow us to construct a $\Lambda_-$ isometry from each $G_t$ to $T_+$, proving that $G_t$ is a weak train track.  

Along the way we also prove one other important piece in the proof of the Properness Theorem~\ref{TheoremProperness}, that the length map restricted to each fold line is a homeomorphism between that fold line and the positive real axis.

\bigskip

Throughout the paper we have probably erred on the side of more detail rather than less, particularly when the extra detail seems structurally interesting and can be illustrated with particular examples of train track maps; see for example Sections~\ref{SectionIWGExample} and~\ref{SectionIrrigidity}. We have tried to develop the theory of the axis bundle comprehensively, incorporating three different natural definitions into the statement of Theorem~\ref{TheoremEquivalence}. However, we have isolated all the material regarding Definition~(1) in Section~\ref{SectionFoldLines}, and so the reader who is willing to completely ignore Definition~(1), and consider only those portions of the theory concerned with Definitions~(2) and~(3) will find the paper somewhat shortened.

\begin{proof}[The Proof of Theorem \ref{TheoremEquivalence}] Let $\phi \in \Out(F_r)$ be nongeometric and fully irreducible, and denote $\TPM = \TPMphi$, $\LambdaPM = \LambdaPMphi$. We use the terminology \emph{weak train track of $\phi$} to refer to an element of $\X_r$ satisfying Definition~(2), being represented by a free, simplicial $F_r$-tree $T$ for which there exists a $\Lambda_-$ isometry $T \mapsto T_+$.

The equivalence of definitions~(1), (2), and~(3) breaks into several implications.

\emph{Definition~(2) implies Definition~(3):} A dense set of weak train tracks are train tracks, as proved in Corollary~\ref{CorollaryWeakIsDense}.

\emph{Definition~(3) implies Definition~(2):} Every train track for $\phi$ is a weak train track, as proved in Lemma~\ref{CorollaryLeavesInTPlus}. The set of weak train tracks is closed, as proved in Lemma~\ref{LemmaLengthFacts}~(\ref{NormalizedWTTClosed}). Together these prove that the closure of the set of train tracks is contained in the set of weak train tracks. 

\emph{Definition~(2) implies Definition~(1):} Every weak train track of $\phi$ lies on a fold line for $\phi$, as proved by a construction given in Proposition~\ref{PropWTTOnFoldLine}.

\emph{Definition~(1) implies Definition~(2):} Every point on a fold line for $\phi$ is a weak train track for $\phi$, as proved in Theorem~\ref{TheoremFoldLines}, which is stated in Section~\ref{SectionFoldLines}, and whose proof takes up the bulk of Section~\ref{SectionFoldLines}.
\end{proof}

\begin{proof}[The Proof of Theorem \ref{TheoremProperness}] Let $\A$ be the axis bundle of a nongeometric, fully irreducible $\phi \in \Out(F_r)$. In the beginning of Section~\ref{SectionTopology} we define the \emph{length map} $\Length \from \A \to (0,+\infinity)$,
by normalizing each weak train track in a natural manner and taking the length of this normalization, and we prove that the length map satisfies the equivariance condition
$$\Length(T \cdot \phi^i) = \lambda(\phi)^{-i} \Length(T) \quad\text{for any}\quad i \in \Z
$$ 
Theorem~\ref{TheoremLengthPHE} says that $\Length \from \A \to (0,+\infinity)$ is a proper homotopy equivalence. 

Consider a fold line $\gamma \from \reals \to \X_r$ such that $\gamma(t)$ converges in $\overline\X_r$ to $\TMinusphi$ as $t \to -\infinity$ and to $\TPlusphi$ as $t \to +\infinity$. From Theorem~\ref{TheoremEquivalence} it follows that $\gamma(\reals) \subset \A$. Corollary~\ref{CorollaryLengthOfFoldPath} of Theorem~\ref{TheoremLengthPHE} says that the restriction of the length map $\Length \from \A \to (0,+\infinity)$ to $\gamma(\reals)$ is a homeomorphism. It follows that the injection $\gamma(\reals) \inject \A$ is a proper homotopy equivalence. 

By combining the equivariance condition with the fact that $\Length \from \A \to (0,+\infinity)$ is a proper homotopy equivalence, it follows that the set $D=\Length^\inv [1,\lambda(\phi)]$ is a compact fundamental domain for the action of the cyclic group $\<\phi\>$ on $\A$. Moreover, the end of $\A$ that corresponds to the positive end of the fold line $\gamma(\reals)$ is represented by the subset $\union_{i=1}^{+\infinity} \phi^i(D)$. The unique limit point of this set in $\overline\X_r$ is $T_+$, because of the result from \BookOne\ that the action of $\phi^i$ on $\overline\X_r$ converges uniformly on compact subsets of $\X_r$ to the constant map with value $T_+$, as $i \to +\infinity$. Similarly, the end of $\A$ that corresponds to the negative end of $\gamma(\reals)$ converges to~$T_-$.
\end{proof}

\subsection{Problems and questions}

\subparagraph{Combinatorial structure of the axis bundle.}  Outer space has a combinatorial structure, a subdivision into closed subsets called ``cells'', each parameterized by a simplex minus some subcomplex of that simplex, so that the intersection of any two cells is a common face; see \cite{CullerVogtmann:moduli} or Section~\ref{SectionOuterSpace} below.

\begin{question}
Given a nongeometric, fully irreducible $\phi \in \Out(F_r)$, do the intersections of $\A_\phi$ with the cells of outer space define a cell decomposition of $\A_\phi$?
\end{question}

The Properness Theorem~\ref{TheoremProperness} provides a proper homotopy equivalence from $\A_\phi$ to $(0,\infinity)$ called the \emph{length function}. If the previous question has an affirmative answer, perhaps level sets of the length function can be studied using the Morse theory ideas of Bestvina and Brady \cite{BestvinaBrady:MorseTheory}, in order to address the following:

\begin{question}
Are the level sets of the length function contractible?
\end{question}

\subparagraph{Geometry of the axis bundles of $\phi$ and of $\phi^\inv$.} The Properness Theorem~\ref{TheoremProperness} implies that $\A_\phi$ has two ends, and that there exist subsets of $\A_\phi$ that are compact and that separate the two ends; we call such a subset a \emph{cross section} of $\A_\phi$. Our proof of the Properness Theorem carried out in Section~\ref{SectionTopology} provides a construction of a cross section, using a Stallings fold argument. Fold ideas are also central to the techniques of \cite{HandelMosher:logratio} relating the expansion factors of $\phi$ and of $\phi^\inv$, and we believe that these techniques may extend to relating their axis bundles as well. We propose here a question that would serve to explore this relation. 

To quantify the size of a cross section, fix an appropriate measurement of (pseudo)distance on $\X_r$, for example: an $\Out(F_r)$-equivariant geodesic metric for which closed balls are compact; or a combinatorial measurement where the distance between two points $x,y \in \X_r$ is the minimal length of a chain of cells $\sigma_0,\sigma_1,\ldots,\sigma_L$ with $x \in \sigma_0$, $y \in \sigma_L$, and $\sigma_{i-1} \intersect \sigma_{i} \ne \emptyset$ for all $i=1,\ldots,L$. Given a point $x \in \A_\phi$ define the \emph{girth of $\A_\phi$ at $x$} to be the infimal diameter of a cross section of $\A_\phi$ passing through $x$, and define the \emph{girth of $\A_\phi$} to be the supremal girth over all points $x \in \A_\phi$. Given $x \in \A_\phi \union \A_{\phi^\inv}$, define the \emph{joint girth at $x$ of $\A_\phi$ and $\A_{\phi^\inv}$} to be the infimal diameter of a compact subset of $\X_r$ which contains the point $x$ and which contains cross-sections of $\A_\phi$ and $\A_{\phi^\inv}$. Define the \emph{joint girth} of $\A_\phi$ to be the supremum over all $x \in \A_\phi \union \A_{\phi^\inv}$ of the joint girth at $x$.

There is no uniform bound for girth depending only on rank. The example in Section~\ref{SectionFoldExamples} of a nonconvergent split ray can be extended to a fold line for some reducible element $\phi \in \Out(F_r)$, and with some thought one sees that the union of fold lines for $\phi$ contains a plane properly embedded into $\X_r$. By carefully interjecting a few more folds into the fold description of $\phi$, one can produce nongeometric, fully irreducible elements of $\Out(F_r)$ with axis bundles of arbitrarily large girth.

Given a nongeometric, fully irreducible $\phi \in \Out(F_r)$, the ideas of \cite{HandelMosher:logratio} give a correlation between the expansion factors of $\phi$ and of $\phi^\inv$, and they do so in a manner which, we believe, also gives a correlation between large girth of $\A_\phi$ and some kind of quantitative measurement that tells us ``how irreducible'' $\phi$ is. We would therefore expect to obtain some control over the joint girths of the axis bundles of $\phi$ and $\phi^\inv$, as long as one of them is ``very irreducible''.

\begin{problem} Fixing the rank $r$, give conditions on a nongeometric, fully irreducible $\phi \in \Out(F_r)$ which imply a bound on the girth of $\A_\phi$. Also, give conditions that would imply a bound on the joint girth of $\A_\phi$ and $\A_{\phi^\inv}$.
\end{problem}

\subparagraph{Generalizations of axes.} In an attempt to find a broader generalization of \Teichmuller\ geodesics, it would be interesting to obtain results not just about source sink pairs $(T_-,T_+)$ for fully irreducible outer automorphisms, but for more arbitrary pairs of points in $\bdy\X_r$.

\begin{question} For what pairs of points $(T_0,T_1)$ in $\bdy \X_r$ is there a good theory of axes, mimicking our results for source sink pairs $(T_-,T_+)$?
\end{question}

\subparagraph{Singularity types.} As mentioned above, in Section~\ref{SectionIdealWhiteheadGraph} we present the ideal Whitehead graph of a nongeometric, fully irreducible outer automorphism as a description of the singularity structure of that outer automorphism.

First recall the singularity structure of a pseudo-Anosov homeomorphism $f \from S \to S$ of a closed, oriented surface $S$. The \emph{index} of a singularity $s$ of $f$ is defined to be $i(s)=1-\frac{k}{2}$ where $s$ has $k$ unstable directions and $k$ stable directions. The Euler-\Poincare\ formula gives an index equation $\sum_s i(s) = \chi(S)$, summed over all singularities of $f$. 

Consider now a nongeometric, fully irreducible $\phi \in \Out(F_r)$. If $C$ is a component of the ideal Whitehead graph, define its \emph{index} to be the number
$$i(C) = 1 - \frac{\#\text{vertices}(C)}{2} 
$$
The \emph{index} of $\phi$ is defined to be $i(\phi) = \sum_C i(C)$ summed over all components of the ideal Whitehead graph (modulo the action of $F_r$). In contrast to the situation on a surface, rather than an index equation there is in general only an index inequality $i(\phi) \ge 1-r$: in \cite{GJLL:index} the $F_r$-tree index $i(\TphiPlus)$ is defined and is proved to satisfy $i(\TphiPlus) \ge 1-r$, and the results of Section~\ref{SectionIdealWhiteheadGraph} imply that $i(\phi)=i(\TphiPlus)$. In \OuterLimits\ it is proved that the inequality $i(\phi) \ge 1-r$ is an equality if and only if $\TphiPlus$ is a geometric $F_r$-tree, which in our present context occurs if and only if $\phi$ is a parageometric outer automorphism.

\begin{question} What values of the index deficit $i(\phi) - (1-r)$ are possible? Does the maximum index deficit go to $+\infinity$ as $r \to +\infinity$?
\end{question}

As a finer invariant, define the \emph{index type} of $\phi$ to be the list of indices of the components of the ideal Whitehead graph of $\phi$, written in increasing order. For instance when $r=3$ then the possible index types which sum to $1-r=-2$ are $(-\frac{1}{2}, -\frac{1}{2}, -\frac{1}{2}, -\frac{1}{2})$, $(-\frac{1}{2}, -\frac{1}{2}, -1)$, $(-\frac{1}{2}, -\frac{3}{2})$, and $(-1,-1)$, and the possible index types whose sum is strictly greater than $-2$ are $(-\frac{1}{2}, -\frac{1}{2}, -\frac{1}{2})$, $(-\frac{1}{2}, -\frac{1}{2})$, $(-\frac{1}{2})$, and $(-1)$. See Section~\ref{SectionIWGExample} for an example in $F_3$ where the index type is $(-1)$. Masur and Smillie proved that on a finite type surface, all possible index types, which sum to the Euler characteristic, are acheived by pseudo-Anosov homeomorphisms, with a few exceptions in low genus \cite{MasurSmillie:singularities}. 

\begin{question} What possible index types, which satisfy the index inequality $i(\phi) \ge 1-r$, are acheived by nongeometric, fully irreducible elements of $\Out(F_r)$?
\end{question}

Even finer than the index type is the isomorphism type of the ideal Whitehead graph of $\phi$. In the analogous situation of a pseudo-Anosov surface homeomorphism $f \from S \to S$, a singularity $s$ of $f$ has an ideal Whitehead graph whose isomorphism type is determined by the index: the unstable directions at $s$ are arranged in a circle, and so if the index equals $i$ then the ideal Whitehead graph of $s$ is just a circle with $-1-2i$ vertices and edges. In contrast, the example given in Section~\ref{SectionIWGExample}  is of a nongeometric, fully irreducible $\phi \in \Out(F_3)$ having a connected ideal Whitehead graph homeomorphic to the triod, more precisely having 5~vertices and 4~edges in the topology of a tree with one vertex of valence~3, three of valence~1, and one of valence~2.  

\begin{question} What isomorphism types of graphs actually occur among the ideal Whitehead graphs of fully irreducible outer automorphisms of $F_r$?
\end{question}

\newcommand\vp{{\vphantom{\prime}}}

\section{Preliminaries} 
\label{SectionPreliminaries}

\subsection{Outer space and outer automorphisms}
\label{SectionOuterSpace}

\paragraph{Outer automorphisms.} 
Once and for all fix the integer $r \ge 2$. Let $F_r = \<x_1,\ldots,x_r\>$ be the free group of rank $r$ with free basis $x_1,\ldots,x_r$. Let $\Aut(F_r)$ be the automorphism group.  Let $i_{\gamma} \from F_r \to F_r$ denote the inner automorphism $i_{\gamma}(\delta) = \gamma \delta \gamma^{-1}$.  Let  $\Inn(F_r)=\{i_{\gamma}: \gamma \in F_r\}$ be the inner automorphism group. Let 
$$\Out(F_r) = \Aut(F_r) / \Inn(F_r)
$$
be the outer automorphism group, which acts on conjugacy classes of elements and of subgroups of $F_r$.

An outer automorphism $\phi \in \Out(F_r)$ is \emph{reducible} if there exists a nontrivial free decomposition $F_r = A_1 * \ldots * A_r * B$ such that $\phi$ permutes the conjugacy classes of $A_1,\ldots,A_r$. If $\phi$ is not reducible then it is \emph{irreducible}. If $\phi^i$ is irreducible for all integers $i \ge 1$ then $\phi$ is \emph{fully irreducible}; this is also known in the literature as ``IWIP'', an acronym for ``irreducible with irreducible powers''. $\phi$ is \emph{geometric} if there exists a homeomorphism $f \from S \to S$ of a compact surface with boundary $S$, and an isomorphism $F_r \approx \pi_1(S)$ which conjugates $\phi$ the outer automorphism of $\pi_1(S)$ induced by $f$. If $\phi$ is not geometric then it is \emph{nongeometric}.  

Fix a \emph{standard rose} $R_r$, a graph with one vertex $v$ and directed edges $e_1,\ldots,e_r$, and identify $F_r = \pi_1(R_r,v)$ where $x_i$ is the class of the loop $e_i$. This induces an isomorphism of $\Aut(F_r)$ with the group of homotopy classes of homotopy equivalences of the pair $(R_r,v)$, and an isomorphism of $\Out(F_r)$ with the group of homotopy classes of homotopy equivalences of $R_r$. We use these isomorphisms without comment henceforth.

\paragraph{Outer space.}
A \emph{graph} is a finite \nb{1}complex $G$ in which every vertex $v$ has $\valence(v) \ge 2$. The \emph{rank} of $G$ is the rank of $\pi_1 G$. The valence condition implies that every homotopy equivalence between graphs is surjective. In a graph of rank~$r$, an Euler characteristic argument shows that $\sum_v \frac{1}{2}(\valence(v)-2) = r-1$, and so there are at most $2r-2$ vertices of valence $\ge 3$. Also,  $\valence(v) \le 2r$ for each $v$, and if in addition $\valence(v) \ge 3$ for each vertex $v$ then $G$ has at most $3r-3$ edges.

A \emph{metric graph} is a graph $G$ equipped with a geodesic metric.  A
\emph{similarity} or \emph{homothety} between two metric graphs $G,G'$ is a homeomorphism $h \from G \to G'$ such that $d(hx,hy) = \lambda\, d(x,y)$ for some constant $\lambda > 0$ called  the \emph{stretch factor} of~$h$. We use the notation $\Length(\cdot)$ to represent length of objects in a metric graph $G$. For example: $\Length G$ is the sum of the length of the edges of $G$; for a path $p \from [a,b] \to G$ we use $\Length p$ for the path length; and if $c$ represents a homotopy class of maps $S^1 \to G$ then we use $\Length c$ to represent the minimal length of a representative of $c$. Sometimes to clarify the context we use $\Length_G(\cdot)$ to mean length
\emph{in the graph $G$}.

A \emph{marked graph} is a metric graph $G$ equipped with a homotopy equivalence, usually denoted $\rho=\rho_G \from R_r \to G$, called a \emph{marking} of $G$. Each marking $\rho_G$ induces an isomorphism $F_r \approx \pi_1(G, \rho_G(v))$, and this sets up a bijection between the set of markings of $G$ modulo homotopy and the set of isomorphisms $F_r \approx \pi_1(G)$ modulo inner automorphism; we use this correspondence henceforth without comment. Given two marked graphs $G,G'$, a homotopy equivalence $f\from G \to G'$ \emph{preserves marking} if $f \composed \rho_G$ is homotopic to $\rho_{G'}$. Two marked graphs $G,G'$ are \emph{(isometrically) equivalent} if there exists an isometry $G \mapsto G'$ that preserves marking, and $G,G'$ are \emph{projectively equivalent} if there exists a similarity $G \mapsto G'$ that preserves marking. 

Culler and Vogtmann's \emph{outer space} $\CV_r$ is the set of projective equivalence classes of marked graphs \cite{CullerVogtmann:moduli}. Define $\wh\CV_r$ to be the set of (isometric) equivalence classes of marked graphs. By our conventions, a marked graph may contain valence~2 vertices, and may be subdivided or unsubdivided at such points without altering the equivalence class in $\CV_r$ or in $\wh\CV_r$. With this flexibility, our definition of $\CV_r$ is equivalent to that of \cite{CullerVogtmann:moduli} where a marked graph has no valence~2 vertices.

We shall review below the topologies on $\X_r$ and $\wh\X_r$, after discussing the interpretation of these spaces in terms of $\reals$-trees.

\subparagraph{Notation:} Often we blur the distinction between a marked graph $G$, its isometric equivalence class in $\wh \CV_r$, and even its projective equivalence class in $\CV_r$. When we need more formality, we use $(\rho_G,G)$ as notation for a marked graph, and $[\rho_G,G]$ for its equivalence class, projective or nonprojective as the context makes clear.

\bigskip

The group $\Out(F_r)$ acts on the right of $\CV_r$ and $\wh\CV_r$ by changing the marking: given a marked graph $(\rho_G,G)$, and given $\phi \in \Out(F_r)$ corresponding to a homotopy equivalence $h \from R_r \to R_r$, the action is defined by
$$[\rho_G,G] \cdot \phi = [\rho_G \composed h,G] 
$$

Any homotopy equivalence $g \from G \to G$ of a marked graph $(\rho_G,G)$ determines a homotopy equivalence $h = \bar\rho_G \composed g \composed \rho_G \from R_r \to R_r$ well-defined up to homotopy independent of the choices of $\rho_G$ and its homotopy inverse $\bar\rho_G$ within their homotopy classes, and so $g$ determines an outer automorphism $\phi$ of $F_r$; we say that \emph{$g$ represents $\phi$}. We will always assume that $g$ is an immersion on edges.  If in addition $g$ takes vertices to vertices then we say that that $g$ is a \emph{topological representative} of $\phi$.  
Note that
$$[\rho_G,G] \cdot \phi = [\rho_G \composed h,G] = [g \composed \rho_G,G]
$$

\paragraph{Actions on trees.} We review group actions on $\reals$-trees, focussing on the group~$F_r$. 

An \emph{$\R$-tree} is a metric space $T$ such that for any $x\ne y\in T$ there exists a unique embedded topological arc, denoted $[x,y] \subset T$, with endpoints $x,y$, and this arc is isometric to the interval $[0,d(x,y)] \subset \reals$. Given $x \in T$, two arcs $[x,y]$, $[x,y']$ determine the \emph{same germ} or the \emph{same direction} at $x$ if $[x,y]\intersect [x,y'] \ne \{x\}$. This is an equivalence relation on arcs with endpoint $x$, and the set of equivalence classes is called the set of \emph{germs} or \emph{directions} at $x$. A \emph{triod} based at a point $x \in T$ is a union of three arcs with endpoints at $x$ no two of which have the same germ at $x$. The \emph{valence} of $x$ is the cardinality of the set of directions, equivalently, the number of components of $T-x$. The valence is $\ge 3$ if and only if there exists a triod based at $x$. We say that $T$ is \emph{simplicial} if the points of valence $\ne 2$ form a discrete set. If $T$ is simplicial and valence is everywhere finite then $T$ is homeomorphic to a CW-complex.

Every isometry $f \from T \to T$ of an $\reals$-tree is either \emph{elliptic} meaning that $f$ has a fixed point, or \emph{hyperbolic} meaning that $f$ has an \emph{axis} $\ell$, an $f$-invariant isometrically embedded copy of the real line on which the \emph{translation length} $\inf_{x \in T} d(x,f(x))$ is minimized. In the elliptic case the translation length equals zero. 

Consider an isometric action of $F_r$ on an $\reals$-tree $T$, denoted $\A \from F_r \to \Isom(T)$. Given $\gamma \in F_r$ and $x \in T$ we use various notations including $t_\gamma$ for $\A(\gamma)$ and $\gamma \cdot x$ for $t_\gamma(x)=\A(\gamma)(x)$. For a subset $S \subset T$ we denote $\gamma \cdot S = \{\gamma \cdot x \suchthat x  \in S\}$. The action $\A$ is \emph{minimal} if there does not exist a nonempty, proper, $F_r$-invariant subtree; this implies that $T$ has no valence~1 points, because the valence~1 points are $F_r$-invariant and their complement is a nonempty $F_r$-invariant subtree. The action $\A$ is \emph{very small} if the stabilizer of every nondegenerate arc is either trivial or a cyclic group generated by a primitive element of $F_r$, and the stabilizer of every triod is trivial. We shall refer to a very small, minimal action of $F_r$ on an $\reals$-tree $T$ as an \emph{$F_r$-tree}.

To each $F_r$-tree $T$ we associate the \emph{length function} $\ell_T \in \reals^{\C}$, where $\ell_T[\gamma]$ is the translation length of the action of $t_\gamma$ on $T$.

Consider two $F_r$-trees $T,T'$. A map $f \from T \to T'$ is \emph{equivariant} if $\gamma \cdot f(x) = f(\gamma\cdot x)$ for $x\in T$. We say that $T,T'$ are \emph{equivariantly isometric} or \emph{isometrically conjugate} if there is an equivariant isometry $T \mapsto T'$, and they are \emph{topologically conjugate} if there exists an equivariant homeomorphism. A result of Culler and Morgan \cite{CullerMorgan:Rtrees} says that $T,T'$ are isometrically conjugate if and only if $\ell_T = \ell_{T'}$ in $\reals^\C$. 

\paragraph{Outer space in terms of trees.} Given a marked graph $(\rho_G,G)$ with universal covering map $pr \from \wt G \to G$, we always assume that a basepoint in the universal cover $\wt G$ has been chosen which projects to $\rho_G(v)$, thereby inducing an action of $F_r$ on $\wt G$ by covering translations. This makes $\wt G$ a free, simplicial $F_r$-tree whose isometric conjugacy class is well-defined independent of choices. Conversely, every free, simplicial $F_r$-tree $T$ determines a quotient marked graph $G_T = T / F_r$, whose marking has homotopy class depending only on the isometric conjugacy class of $T$. Denoting $\wt G = T_G$, clearly we have $T_{G_T}=T$ and $G_{T_G}=G$, and so we obtain a natural bijection between $\hat\CV_r$ and the isometric conjugacy classes of free, simplicial $F_r$-trees. This bijection agrees with the length map to $\reals^\C$: if $T = \wt G$ then $\ell_T = \ell_G$ in $\reals^\C$. Under this bijection the right action of $\Out(F_r)$ on $\hat\CV_r$ can be described as follows: given a free, simplicial $F_r$ tree $T$ with action $\A \from F_r \to \Isom(T)$, and given $\phi \in \Out(F_r)$, choose $\Phi \in \Aut(F_r)$ representing $\phi$, and let $T\cdot \phi$ be the $F_r$ tree with underlying $\reals$-tree $T$ and with action $\A \composed \Phi \from F_r \to \Isom(T)$; this action is well-defined up to conjugacy independent of the choice of $\Phi$.

We can therefore describe $\widehat\CV_r$ \emph{either} as marked graphs up to marked isometry \emph{or} as free, simplicial $F_r$-trees up to isometric conjugacy; outer space $\CV_r$ has a similar pair of descriptions. Given a free, simplicial $F_r$-tree $T$, depending on context we use $[T]$ for the point in $\CV_r$ or $\widehat\CV_r$ represented by $T$, but we very often abuse notation by deleting the brackets and letting $T$ stand for the $F_r$-tree as well as the corresponding point in $\CV_r$ or $\widehat\CV_r$.

\paragraph{Ends, automorphisms, and lifts.} The Gromov boundary of the free group $F_r$ is a Cantor set denoted $\bdy F_r$. The action of $\Aut(F_r)$ on $F_r$ extends to an action on $\bdy F_r$; the extension of $\Phi \in \Aut(F_r)$ is denoted $\wh\Phi \from \bdy F_r \to \bdy F_r$. 

Given a free simplicial $F_r$-tree $T$, the space of ends of $T$ is denoted $\bdy T$, and there is a natural homeomorphic identification $\bdy F_r \approx \bdy T$ which we will use without further comment. Given $\gamma \in F_r$, the map $t_\gamma \from T \to T$ extends to a map of ends denoted $\hat t_\gamma \from \bdy T \to \bdy T$, whose fixed point set is equal to the pair of endpoints $t_{\gamma}^{\pm}$ of the axis of $t_{\gamma}$. 

Suppose that $G$ is marked graph and $g \from G \to G$ is a homotopy equivalence representing $\phi \in \Out(F_r)$. A lift of $g$ to the $F_r$-tree $\wt G$ will be denoted $\breve g \from \wt G \to \wt G$, where the ``\verb+\breve+'' notation $\breve{\,}$ is meant to be a reminder that the lifted map is \emph{not} $F_r$-equivariant (unless $\phi$ is the identity and $g$ preserves marking). The map $\breve g$ extends to a map of ends denoted $\hat g \from \bdy\wt G \to \bdy\wt G$. There is a bijection between the set of lifts $\breve g : \wt G\to \wt G$ of $g$ and the set of automorphisms $\Phi : F_r \to F_r$ representing~$\phi$. This bijection is characterized by the \emph{twisted equivariance equation}
\begin{equation}\label{EquationTwisted}
\breve g \composed t_\gamma = t_{\Phi(\gamma)} \composed \breve g \quad\text{for all $\gamma \in F_r$.}
\end{equation}
The bijection is also characterized by the condition that $\breve g$ and $\Phi$ extend to the same homeomorphism $\hat g = \hat \Phi : \partial F_r \to \partial F_r$. A simple computation using Equation~\ref{EquationTwisted} shows that if $\breve g, \breve g' \from \wt G \to \wt G$ are the lifts of $g$ with corresponding automorphisms $\Phi,\Phi' \in \Aut(F_r)$ that represent $\phi$, and if $\gamma \in F_r$ is the unique element such that $\Phi' = i_\gamma \composed\Phi$, then $\breve g' = t_\gamma \composed \breve g$. 

\subparagraph{Notation:} The notations $\hat g$ and $\hat \Phi$ are used interchangeably, depending on context.

\paragraph{The topology of outer space.} The topology on $\X_r$ has three equivalent descriptions which we review: the cellular topology; the length or axes topology; and the  Gromov topology. We use all three of these descriptions at various places in the paper.

\subparagraph{The cellular topology \cite{CullerVogtmann:moduli}.} Consider a marked graph $G$. With respect to the natural vertex set consisting of points of valence~$\ge 3$, denote the edge set by $E(G)$. The \emph{open cone} of $G$ is the set of elements in $\wh\X_r$ represented by marked graphs $G'$ for which there exists a homeomorphism $G \mapsto G'$ that preserves marking; the \emph{open cell} of $G$ in $\X_r$ is similarly defined. Note that if $g \from G \to G$ is a homeomorphism of a marked graph that preserves marking then $g$ takes each edge to itself. It follows that if $G'$ is in the open cone of $\wh\X_r$ containing $G$ then any two homeomorphisms $G \mapsto G'$ that preserve marking induce the same bijection $E(G) \leftrightarrow E(G')$, and therefore the lengths of edges on $G'$ determine a well-defined vector $L(G') \in [0,\infinity)^{E(G)}$. The map $L$ is a bijection between the open cone containing $G$ and the set $(0,\infinity)^{E(G)}$. By projectivizing, we obtain a bijection between the open cell of $G$ in $\X_r$ and the interior of the simplex $\P[0,\infinity)^{E(G)}$. 

The \emph{closed cone} in $\wh\X_r$ containing $G$ is the set of elements represented by marked graphs $G'$ for which there exists a quotient map $p \from G \to G'$ that preserves marking and that collapses to a point each component of a subforest $F \subset G$; the \emph{closed cell} in $\X_r$ containing $G$ is similarly defined. The map $p$ induces a well-defined bijection $E(G)-E(F) \leftrightarrow E(G')$, and therefore there is a well-defined vector $L_{G'} \in [0,\infinity)^{E(G)}$ which assigns zero to edges in $F$ and assigns the length of the corresponding edge of $G'$ to each edge in $G-F$. Again $L$ is a bijection between the closed cone containing $G$ and a certain subset of $[0,\infinity)^{E(G)}$, namely those vectors which assign coordinate zero to a subset of $E(G)$ that forms a subforest of $G$. The closed cone of $G$ is thus identified with a subset of $[0,\infinity)^{E(G)}$ which is the union of the interior and certain of the faces. By projectivizing, we obtain a bijection between the closed cell of $G$ in $\X_r$ and a union of faces of the simplex $\P[0,\infinity)^{E(G)}$. Note that if $F \subset F'$ are faces of the simplex and $F$ is a face of the closed cone of $G$ then $F'$ is also a face of the closed cone of $G$.

The \emph{cellular topology} on $\wh\X_r$ is the weakest topology such that for each marked graph $G$, the closed cone containing $G$ is a closed subset of $\wh\X_r$, and the map $L$ from this cone to $[0,\infinity)^{E(G)}$ is a homeomorphism onto its image. The cellular topology on $\X_r$ is the quotient by projectivization of the cellular topology on $\wh\X_r$.

We need some facts about the cellular topology. First, the closed cone (simplex) of $G'$ is a face of the closed cone (simplex) of $G$ if and only if there is a quotient map $G \mapsto G'$ that preserves marking and that collapses to a point each component of some subforest of $G$. Second, the closed cells (simplices) are locally finite: for any $G'$ there are only finitely many closed cells (simplices) that contain $G'$, and hence by the weak topology there is a neighborhood of $G'$ that intersects only those finitely many cells (simplices) that contain $G'$.

Given a marked graph $G$ define its \emph{length} $\ell(G)$ to be the sum of the lengths of its edges. This gives a well-defined function on $\wh\X_r$ which, on each closed cell, is just the sum of the coordinates. It follows that, in terms of the cellular topology:

\begin{proposition}
The map $\ell \from \wh\X_r \to (0,\infinity)$ is continuous.
\qed\end{proposition}

\subparagraph{The length or axes topology \cite{CullerMorgan:Rtrees}.} Let $\C$ be the set of nontrivial conjugacy classes in $F_r$, give $\R^\C$ the product topology and let $P\R^\C$ be the associated projective space. For each marked graph $G$, each $c \in \C$ corresponds bijectively, via the marking $\rho_G$, with an immersion $S^1 \to G$ also denoted $c$, which is well-defined up to a homeomorphism of the domain, and whose length $\Length_G(c)$ is well-defined. The map $G \to \Length_G$ defines an injection $\wh\CV_r \inject \R^\C$ and an induced injection $\CV_r \inject P\R^\C$. These injections may also be defined in terms of translations lengths of free, simplicial $F_r$-trees representing elements of $\X_r$. We give $\CV_r$ and $\wh\CV_r$ the topologies that make these injections homeomorphisms onto their images; we shall occasionally refer to these as the \emph{length topologies} on $\CV_r$ and $\wh\CV_r$. It follows that the natural action of $\reals_+$ on $\wh\CV_r$, where $t \cdot G$ is obtained from $G$ by multiplying the metric by the constant factor $t$, is an action by homeomorphisms with quotient space $\CV_r$. 

\subparagraph{The Gromov topology \cite{Paulin:GromovTopology}.} This topology on $\wh\CV_r$ is defined by the following basis. Given an $F_r$-tree $T \in \wh\CV_r$, $\epsilon>0$, a finite subset $K \subset T$, and a finite subset $P \subset F_r$, define the basis element $V(T;\epsilon,K,P)$ to be the set of $F_r$-trees $T' \in \wh\CV_r$ for which there exists a finite subset $K' \subset T'$ and a relation $R \subset K \cross K'$ such that $R$ projects onto $K$ and onto $K'$, and the following hold:
\begin{description}
\item[$R$ is an $\epsilon$-almost isometry:] For all $x,y \in K$ and $x',y' \in K'$, if $(x,x'), (y,y') \in R$ then $\abs{d(x,y) - d(x',y')} < \epsilon$.
\item[$R$ is $P$-equivariant:] For all $x \in K$, $x' \in K'$, and $\alpha \in P$, if $\alpha \cdot x \in K$ and $(x,x') \in R$ then $\alpha \cdot x' \in K'$ and $(\alpha \cdot x, \alpha \cdot x') \in R$.
\end{description}
In \cite{Paulin:GromovTopology} the Gromov topology is defined by forming a basis as above but using compact subsets of $T$ and $T'$ instead of finite subsets; as shown in Proposition~4.1 of \cite{Paulin:GromovTopology} the definition given above is a basis for the same topology. 

The quotient topology on $\CV_r$ is also called the Gromov topology.

The following theorem is a combination of results of \cite{CullerVogtmann:moduli} and \cite{Paulin:GromovTopology}:

\begin{theorem}
The cellular topology, the length topology, and the Gromov topology are identical on $\wh\CV_r$ and on $\CV_r$.
\qed\end{theorem}

\paragraph{Compactification and boundary of outer space.} The closure of $\CV_r$ in $P\reals^\C$ is denoted $\overline\CV_r$. Its boundary is $\bdy\CV_r = \overline\CV_r - \CV_r$. Both $\overline\CV_r$ and $\bdy\CV_r$ are compact \cite{CullerMorgan:Rtrees}. Let $\whclosure\CV_r$ be the inverse image of $\overline\CV_r$ under the projection $\reals^\C-\{O\} \mapsto P\reals^\C$, equivalently, the closure of $\wh\CV_r$ in $\reals^\C - \{O\}$. Let $\bdy\wh\CV_r = \whclosure\CV_r - \wh\CV_r$, which is just the inverse image under projection of $\bdy\CV_r$.

We next relate compactified outer space to actions on trees.

\begin{theorem}  \label{TheoremCullerMorgan}(Culler--Morgan \cite{CullerMorgan:Rtrees}) 
For each $\ell \in \whclosure\CV_r$ there exists an $F_r$-tree $T$ such that
$\ell=\ell_T$, and $T$ depends uniquely on $\ell$ up to equivariant isometry. 
\end{theorem}

Extending our notation, when $\ell=\ell_T$ as in this theorem, we use $[T]$ to denote the corresponding point in $\whclosure\CV_r$ or $\overline\CV_r$ as the context makes clear.

\begin{lemma}
\label{LemmaDiscreteTrees} 
Given $\ell \in \whclosure\CV_r$ and given $T$ with $\ell_T=\ell$, the following are equivalent:
\begin{enumerate}
\item $\ell \in \wh\CV_r$.
\item $\ell(c)$ is bounded away from zero for $c \in \C$.
\item The action of $F_r$ on $T$ is free and properly discontinuous.
\end{enumerate}
\end{lemma}

\begin{proof}  It is obvious that the first item implies the second. Since the action is by isometries the second item implies the third.  Finally, the third item implies that the quotient map from $T$ to  $G :=T / F_r$ is a covering space.  Since  $T$ is an $F_r$-tree, $G$ is a finite graph and so $T$ is  simplicial. 
\end{proof}

\subsection{Paths, Circuits and Edge Paths}
\label{SectionPaths} 
Consider a marked graph $G$ and the universal covering projection $pr : \wt G \to G$.  A proper map $\ti \sigma : J \to \wt G$ with domain a (possibly infinite) interval $J$  will be called a \emph{path in $\wt G$} if it is an embedding or if $J$ is finite and the image is a single point;  in the latter case we say that $\ti \sigma$ is \emph{a trivial path}.   If $J$ is finite, then every map $\ti \sigma : J \to \wt G$ is homotopic rel endpoints to a unique (possibly trivial) path $[\ti \sigma]$; we say that \emph{$[\ti \sigma]$ is obtained from $\ti \sigma$ by tightening}. 
 
We will not distinguish between paths in $\wt G$ that differ only by an  orientation preserving change of parametrization. Thus we are interested in the oriented image of $\ti \sigma$ and not $\ti \sigma$ itself. If the domain of $\ti \sigma$ is finite, then the image of $\ti \sigma$ has a natural  decomposition as a concatenation $\wt E_1 \cdot \wt E_2\cdot \ldots \cdot \wt E_k$, where $\wt E_i$ for $1 < i < k$ is an edge of  $\wt G$, $\wt E_1$ is the terminal segment of some edge, and $\wt E_k$ is the initial segment of some edge. If the endpoints of the image of $\ti  \sigma$ are vertices, then $\wt E_1$ and $\wt E_k $ are full edges. The sequence $\wt E_1\cdot \wt E_2\cdot\ldots \cdot\wt E_k$ is called \emph{the edge path associated to  $\ti \sigma$}, and the number $k$ is called the \emph{combinatorial length} of this edge path. When the interval of domain is half-infinite then an edge path has the form $\wt E_1\cdot \wt E_2\cdot \ldots$ or $\ldots \cdot \wt E_{-2}\cdot \wt E_{-1}$, and in the bi-infinite case an edge path has the form $\ldots \cdot\wt E_{-1}\cdot \wt E_0\cdot \wt E_1\cdot \wt E_2\cdot\ldots$.

A \emph{path in $G$} is  the composition of the projection map $pr$ with a path in $\wt G$. Thus a map $\sigma : J \to G$ with domain a (possibly infinite) interval will be called a path if it is an immersion or if $J$ is finite and the image is a single point; paths of  the latter type are said to be trivial. If $J$ is finite, then every map $ \sigma : J \to G$ is homotopic rel endpoints to a unique  (possibly trivial) path $[\sigma]$; we say that \emph{$[ \sigma]$ is obtained from $ \sigma$ by tightening}. For any lift $\ti \sigma : J \to \Gamma$  of $\sigma$, $[\sigma] = pr[\ti \sigma]$.  We do not distinguish between paths in $G$ that differ by an orientation preserving change of parametrization.  The \emph{edge path associated to $\sigma$} is the projected image of the edge path associated to a lift $\ti \sigma$.  Thus the edge path associated to a path with finite domain has the form $E_1 \cdot E_2 \cdot \ldots \cdot E_{k-1} \cdot E_k$ where $E_i$, $1 < i < k$, is an edge of $G$, $E_1$ is the terminal segment of an edge, $E_k$ is the initial segment of an edge; the \emph{combinatorial length} again equals $k$.

We reserve the word \emph{circuit} for an immersion $\sigma : S^1 \to G$. Any homotopically non-trivial  map $\sigma : S^1 \to G$ is homotopic to a unique circuit $[\sigma]$. As was the case with paths, we do not distinguish between circuits that differ only by an orientation preserving change in parametrization and we identify a circuit $\sigma$ with a \emph{cyclically ordered edge path} $E_1E_2\dots E_k$, and again the \emph{combinatorial length} equals $k$.

If $h \from H \to G$ is a homotopy equivalence of marked graphs and if $\ti h : \wt H \to \wt G$ is a lift then for each path $\ti\sigma$ in $\wt H$ we denote $[\ti h(\ti \sigma)]$  by $\ti h_\#(\ti \sigma)$, and for each path or circuit $\sigma$ in $H$ we denote $[h(\sigma)]$ by $h_\#(\sigma)$. Given a self homotopy equivalence $g \from G \to G$ of a marked graph and a path or circuit $\sigma$ in $G$, a decomposition $\sigma = \ldots \cdot \sigma_i \cdot \sigma_{i+1}\cdot \ldots$ is called a \emph{splitting} if $g^k_\#(\sigma)= \ldots \cdot g^k_\#(\sigma_i) \cdot g^k_\#(\sigma_{i+1})\cdot \ldots$ for all $k \ge 0$; in other words, to tighten $g^k(\sigma)$ one need only tighten the $g^k(\sigma_i)$. 

\vspace{.1in}

{ \it Throughout this paper we will identify paths and circuits with their associated edge paths.}
 
\vspace{.1in}

A path or circuit \emph{crosses} or \emph{contains} an edge if that edge occurs in the associated edge path.
For any path $\sigma$ in $G$ define $\bar \sigma$ to  be \lq $\sigma$ with its
orientation reversed\rq.   For notational simplicity, we sometimes refer to the inverse of $\ti \sigma$ by $\ti \sigma^{-1}$.

\paragraph{Bounded cancellation.} The bounded cancellation lemma was introduced in \cite{Cooper:automorphisms} and generalized to the form we need in \BookZero. Let $\Lip(f)$ be the lipschitz constant of a map $f$.

\begin{lemma}\label{bcl} [Bounded Cancellation, Lemma 3.1 of \BookZero] Let $T_0$ be a free, simplicial $F_r$-tree, let $T$ be an $F_r$-tree representing a point in $\overline\X_r$, and let $f \from T_0 \to T$ be an $F_r$-equivariant map. There exists a constant $\BCC(f) \le \Lip(f) \cdot \Length(T_0/F_r)$ such that for each geodesic segment $\alpha \subset T_0$, the image $f(\alpha)$ is contained in the $\BCC(f)$ neighborhood of the straightened image $f_\#(\alpha)$.
\end{lemma}

\subsection{Folds} 
\label{SectionFolds}

In this subsection we define morphisms, edge isometries, and fold paths. We also define fold maps and Stallings fold sequences, and we describe relations between these objects and fold paths.

\paragraph{Turns and gates.} In any graph a \emph{turn} is an unordered pair of directions with a common base point. The turn is  \emph{nondegenerate} if it is defined by two distinct directions and is \emph{degenerate} otherwise. 

Consider two graphs $T,T'$ and a map $H \from T \to T'$ which is locally injective on each edge. For each $x \in T$ there is an induced map $D_x H$ from the set of directions of $T$ at $x$ to the set of directions of $T'$ at $H(x)$: given a direction $d$ at $x$ represented by an oriented path $\alpha$ with initial endpoint $x$, $D_x H(d)$ is defined to be the direction at $H(x)$ represented by $H(\alpha)$; the map $D_x H$ is well-defined, independent of the choice of $\alpha$ representing $d$. We use $DH$ either as an abbreviation for $D_x H$ when $x$ is understood, or as the union of all of the maps $D_x H$ over all $x \in T$. 

Given an $F_r$-equivariant map of $F_r$-trees $H \from T \to T'$, we define the \emph{gates} at a point $x \in T$ with respect to $H$ to be the equivalence classes of directions at $x$ where $d,d'$ are equivalent if $D_x H(d)=D_x H(d')$. In applications where the base point $x$ of a direction is a vertex we can identify directions with oriented edges $E$ in the link of $x$, and if in addition $H$ takes vertices to vertices then $D_x H(E)$ is the first edge in the edge path associated to $H(E)$. We have here extended the terminology of a ``gate'' which in the literature has so far been applied only to train track maps; see Section~\ref{SectionTrainTracks}.

\paragraph{Piecewise isometries and edge isometries.}
An $F_r$-equivariant map $f \from T \to T'$ from a free, simplicial $F_r$-tree $T$ to an $F_r$ tree $T'$ is a \emph{piecewise isometry}\footnote{A ``piecewise isometry'' is the same thing as a ``morphism'', as commonly defined in the literature, but restricted so that the domain is free and simplicial.} if each edge of $T$ may be subdivided into finitely many segments on each of which $f$ is an isometry. We usually use the even stronger property that $f$ restricted to each edge of $T$ is an isometry, to which we refer by saying that $f$ is an \emph{edge isometry} or that $f$ is \emph{edge isometric}. For a piecewise isometry to be an edge isometry is equivalent to it being locally injective on each open edge. When $T'$ is also free and simplicial then we also define a piecewise isometry or edge isometry of marked graphs $G=T/F_r \mapsto T'/F_r=G'$, meaning that the lifted map $T \to T'$ is a piecewise isometry or edge isometry. Equivalently, an edge isometry between marked graphs is any homotopy equivalence that preserves marking and restricts to a local isometry on each open edge.

\paragraph{Fold paths.} A \emph{fold path} in outer space is a continuous, proper injection $\gamma \from I \to \CV_r$, where $I \subset \reals$ is connected, so that there exists a continuous lift $\wh\gamma \from I \to\wh\CV_r$, a marked graph $G_r$ representing $\wh\gamma(r)$ for each $r \in I$, and maps $g_{sr}\from G_r \to G_s$ for each $r \le s$ in $I$, satisfying the following:
\begin{enumerate}
\item \textbf{Train track property:}
\label{ItemTrainTrack} 
Each $g_{sr}$ is an edge isometry.
\item \textbf{Semiflow property:}
\label{ItemFold} $g_{tr} = g_{ts} \composed g_{sr}$ for each $r \le s \le t$ in $I$, and $g_{ss} \from G_s \to G_s$ is the identity for each $s \in I$.
\end{enumerate}
A fold path $\gamma$ whose domain $I$ is noncompact on both ends is called a \emph{fold line}. 

Fold paths are investigated in detail in Section~\ref{SectionFoldLines}. For now we turn to well known constructions that lead to simple examples of fold paths.

\paragraph{Fold maps.} Consider a marked graph $G$ and two distinct oriented edges $e,e'$ of $G$ that have the same initial vertex $v$. Parameterize $e,e'$ by arc length from their initial vertex, using the notation $e(\cdot), e'(\cdot)$. Choose a number $S$ such that $e([0,S]) \intersect e'([0,S]) = \{v\}$. For $s \in [0,S]$ define $G_s$ to be the graph obtained from $G$ by identifying $e(t)$ to $e'(t)$ for each $t \in [0,s]$. The quotient map $G \mapsto G_s$ is a homotopy equivalence, because of the fact that $e(s) \ne e'(s)$. Push forward the metric on edges of $G$ by the quotient map $G \mapsto G_s$, defining a metric on edges of $G_s$ which extends naturally to a geodesic metric on all of $G_s$. Postcompose the marking of $G$ by the homotopy equivalence $G \mapsto G_s$, defining a marking of $G_s$. It follows that the quotient map $G \mapsto G_s$ preserves the marking and is an edge isometry; when $s=0$ this is just the identity map on $G$. When $s>0$ we say that $G_s$ is obtained from $G$ by a \emph{length $s$ fold of $e$ and $e'$}, and the quotient map $G \mapsto G_s$ is called a \emph{fold map} from $G$ to $G_s$. 

Consider the map $\wh\gamma \from [0,S] \to \wh\X_r$ defined by $\wh\gamma(s)=G_s$. First, $\wh\gamma$ is a continuous path in $\wh\X_r$ from $G$ to $G_S$, and so it descends to a continuous path $\gamma \from [0,S] \to \X_r$ which, its domain being compact, is proper. Second, for $0 \le s \le t \le S$ the induced map $f_{ts} \from G_s \mapsto G_{t}$ is a fold map. Third, the family of maps $f_{ts}$ satisfies the train track property and the semiflow property, by construction. It follows that $\gamma$ is a fold path. Either of the paths $\gamma \from [0,S] \to \X_r$ or $\wh\gamma \from [0,S] \to \wh\X_r$, or any of their orientation preserving reparameterizations, are called \emph{interpolating paths} or just \emph{interpolations} of the fold map $G \mapsto G_S$. 

\paragraph{Stallings fold sequences.}
Given an edge isometry $h \from G \to H$ between marked graphs, a \emph{Stallings fold sequence} for $h$ is a finite sequence of edge isometric fold maps 
\begin{equation}
\label{EqFolds}
\xymatrix{
G=G_0 \ar[r]^{f_{1,0}} \ar@/_1pc/[rrr]_h & G_1 \ar[r]^{f_{2,1}} & \cdots \ar[r]^{f_{J,J-1}} & G_J=H
}
\end{equation}
whose composition equals $h$. If we define $f_{ji}$ for integer values $0 \le i \le j \le J$ to be the composition $f_{j,j-1} \composed \cdots \composed f_{i+1,i} \from G_i \to G_j$, then the maps $f_{ji}$ are edge isometries satisfying the semiflow property $f_{ki}=f_{kj}\composed f_{ji}$ for integer values $0 \le i \le j \le k \le J$. For $j=1,\ldots,J$ let $[j-1,j] \mapsto \wh\X_r$ be an interpolating path for the fold map $f_{j,j-1} \from G_{j-1} \to G_j$, and extend the map $f_{j,j-1}$ to a family of edge isometries $f_{ts} \from G_s \to G_t$, $j-1 \le s \le t \le j$ that satisfies the semiflow property. Concatenating the paths $[j-1,j] \mapsto \wh\X_r$ for $j=1,\ldots,J$ defines a continuous path $\wh\gamma \from [0,J] \mapsto \wh\X_r$, and the entire family of edge isometries $f_{ts}$ defined for $0 \le s \le t \le J$ satisfies the semiflow property. Thus $\wh\gamma$ descends to a continuous path $\gamma \from [0,J] \to \X_r$ which is a fold path, called the \emph{Stallings fold path} of the Stallings fold sequence~(\ref{EqFolds}). 

To complete the construction of Stallings fold paths, we review Stallings original topological construction of fold sequences \cite{Stallings:folding}, and adapt this construction to our present metric setting where fold maps are required to be edge isometric. At the same time we shall use this construction to bound the length of the fold sequence.

In a Stallings fold sequence~(\ref{EqFolds}), the fold $f_{j+1,j} \from G_j \to G_{j+1}$ is said to be \emph{complete} if it satisfies the following: assuming $f_{j+1,j}$ to be a length $S$ fold of oriented edges $e,e'$ of $G_j$, then $S$ is equal to the maximum value of $s$ such that restrictions of the map $f_{J,j} \from G_j \to G_J=H$ to $e([0,s])$ and to $e'([0,s])$ are the same path in $H$. In other words, $f_{j+1,j}$ is complete if the maximal possible initial segments of $e,e'$ are folded immediately.

Define the \emph{combinatorial length} of an edge isometry $h \from G \to H$ as follows. While at first regarding $G$ and $H$ with their natural vertex sets consisting solely of vertices of valence $\ge 3$, subdivide $H$ at the images of the vertices of $G$, which adds at most $2r-2$ new valence~2 vertices to $H$ and at most $2r-2$ new edges, bringing the number of edges in $H$ up to at most $5r-5$. For each edge $e$ of $G$ we may regard $h(e)$ as an edge path in the subdivided graph $H$, and we define the \emph{combinatorial length} of $h$ to be the sum of the combinatorial lengths of $h(e)$ over all edges $e$ of $G$. Equivalently, having subdivided $H$, if we then subdivide $G$ at the inverse image of the vertex set of $H$, and if we refer to the edges of this subdivision as the \emph{edgelets} of $G$ with respect to~$h$, then the combinatorial length of $h$ is the number of edgelets of $G$ with respect to $h$.

\begin{lemma}
\label{LemmaFoldBound}
Each edge isometry $h \from G \to H$ which is not an isometry has a sequence of complete Stallings folds. The length of any such fold sequence is bounded above by the combinatorial length of $h$.
\end{lemma}

\begin{proof} Subdivide $H$ at the images of the vertices of $G$.

Since the map $h$ is a local injection on each edge and the graphs $G,H$ have no valence~1 vertices, if $h$ were locally injective at each vertex then there would be an extension of $h$ to a covering map $h' \from G' \to H$ where the graph $G'$ deformation retracts to $G$ (see \cite{Stallings:folding}). Since $h$ is a homotopy equivalence it would follow that $h=h'$ is a homeomorphism and therefore $h$ is an isometry, contrary to the hypothesis. 

There exist therefore two oriented edges $e,e'$ of $G$ with the same initial vertex $v$, and there exists a positive number $R \le \min\{\Length(e),\Length(e')\}$, such that $h$ maps $e([0,R])$ and $e'([0,R])$ to the same oriented path in $H$. The maximal such $R$ has the property that $h(e(R))=h(e'(R))$ is a vertex of $H$ and $e([0,R]) \intersect e'([0,R]) = \{v\}$. Define $G \xrightarrow{h'} G'$ to be the length $R$ fold of $e$ and $e'$, and so $h$ factors as $G \xrightarrow{h'} G' \xrightarrow{h''} H$. Noting that each of $e([0,R])$ and $e'([0,R])$ is a union of edgelets of $G$ with respect to $h$, identifying $e([0,R])$ with $e'([0,R])$ has the effect of identifying their edgelets in pairs, and so the number of edgelets of $G'$ with respect to the edge isometry $h''$ --- the combinatorial length of $h''$ --- is strictly less than the number of edgelets of $G$ with respect to $h$ --- the combinatorial length of $h$. The existence of a sequence of complete Stallings folds now follows by induction on combinatorial length. 

The same argument shows that for \emph{any} sequence of complete Stallings folds~(\ref{EqFolds}), the combinatorial length of $f_{J,j} \from G_j \to G_J=H$ is a strictly decreasing function of $j$, proving the desired bound.
\end{proof}

\subsection{Train Track Maps}
\label{SectionTrainTracks}
A topological representative $\gG$ of an outer automorphism $\phi$ is \emph{an irreducible train track map} if $g^k \restrict E$ is an immersion  for each edge $E$ and each $k > 0$ and if for every pair of edges $E$ and $E'$ there exists $l > 0$ such that $g^l(E)$ crosses $E'$. Theorem 1.7 of \cite{BestvinaHandel:tt} states that every irreducible $\phi \in \Out(F_r)$ is represented by an irreducible train track map. Moreover, there exists $\lambda > 1$, called the \emph{expansion factor for $\phi$}, there exists a metric on $G$, and there exists a homotopy of $g$ so that, still using $g$ for the homotoped map, for all edges $E$ the restriction $g \restrict E$ uniformly expands length by the factor $\lambda$. We will always assume that the graph $G$ is equipped with such a metric, and we occasionally emphasize this point by referring to $g$ as an \emph{affine train track map}.

A nondegenerate turn in $G$ is \emph{illegal} with respect to $\gG$ if its image under some iterate of $Dg$ is degenerate, and otherwise the turn is \emph{legal}. A path $\sigma=E_1\cdot E_2 \cdot \ldots \cdot E_{k-1} \cdot E_k$ \emph{contains} or \emph{takes} each of the turns $\{\bar E_i,E_{i+1}\}$ for $0 \le i \le k-1$.  A \emph{path $\sigma \subset G$ is legal} if it contains only legal turns. It is immediate from the definitions that $Dg$ maps legal turns to legal turns and that $g$ maps legal paths to legal paths. Note that $g$ expands the length of any legal subpath by the factor~$\lambda$.

We define the \emph{gates} of a train track map $g \from G \to G$ as follows. The map $g$ induces a map $Dg$ on the set of directions at vertices of $g$. The \emph{gates} at a vertex $x$ with respect to $H$ are the equivalence classes of directions at $x$ where $d,d'$ are equivalent if there exists $i \ge 1$ such that $(Dg)^i(d)=(Dg)^i(d')$. It follows that two directions are in the same gate if and only if the turn they make is illegal. Note the distinction between gates of a train track map and gates of an $F_r$-equivariant map of $F_r$-trees which is defined in Section~\ref{SectionFolds}.

\paragraph{Local Whitehead Graphs.} \qquad
Consider an irreducible train track representative $g \from G \to G$ of an irreducible $\phi \in \Out(F_r)$. For each point $x \in G$, the \emph{local Whitehead graph at $x$}, denoted $W(x)$ or more explicitly $W(x;G)$, is the graph whose vertices are the directions $d_i$ at $x$ and with an edge connecting $d_i$ to $d_j$ if there exists an edge $E$ and $k>0$ such that the path $g^k(E)$ takes the turn $(d_i,d_j)$. Since $\phi$ is irreducible, $W(x)$ is connected (see section 2 of \cite{BFH:laminations}). For each $x \in G$, the map $D_x g$ extends to a simplicial map $W(x) \mapsto W(g(x))$ also denoted $D_x g$, because each turn at $x$ taken by some $g^k(E)$ maps to a turn at $g(x)$ taken by~$g^{k+1}(E)$. 

If $x \in \Per(g)$,  then the \emph{local stable Whitehead graph at $x$}, denoted $SW(x)$ or $SW(x;G)$, is the subgraph of $W(x)$ obtained by restricting to the periodic directions at $x$ and the edges between them. Since some iterate $g^k$ induces a continuous map $D_x g^k$ from $W(x)$ onto $SW(x)$, it follows that $SW(x)$ is connected. There are analogous definitions of $W(\ti x) = W(\ti x;\wt G)$ and  $SW(\ti x) = SW(\ti x;\wt G)$ for any vertex $\ti x$ of $\wt G$; these can also be defined as the lifts of  $W(x)$ and $SW(x)$. It follows that $W(\ti x)$ and $SW(\ti x)$ are connected.

\paragraph{Nielsen Paths.} \quad
A path $\sigma$ is a \emph{periodic Nielsen path}  if $g^p_\#(\sigma) = \sigma$ for some minimal $p \ge 1$ called the \emph{period} of $\sigma$.  A \emph{Nielsen path} is a periodic Nielsen path with period one.    A periodic  Nielsen path is \emph{indivisible} if    it can not be written as a non-trivial concatenation of periodic  Nielsen subpaths.  A path $\sigma$ is  a \emph{periodic pre-Nielsen  path} if it is not a periodic Nielsen path but some $g_\#^k(\sigma)$ is a periodic Nielsen path.  

We record information about Nielsen paths in the following lemma.   

\begin{lemma} \label{nielsen paths}  Assume that $\gG$ is an irreducible train track map.  
\begin{enumerate}
\item  There are only finitely many indivisible periodic or pre-periodic Nielsen paths for $g$.
\item  An indivisible periodic or pre-periodic Nielsen path $\sigma$ decomposes into subpaths $\sigma = \alpha \bar \beta$ where $\alpha$ and $\beta$ are legal and the turn $(\bar \alpha,  \beta)$  is illegal.
Moreover, if $\alpha_n$ and $\beta_n$ are the terminal segments of $\alpha$ and $\beta$ that are  identified by $g^n$, then the $\alpha_n \bar \beta_n$'s are an increasing sequence of subpaths of $\sigma$ whose union is the interior of $\sigma$.
\item For every path or circuit $\sigma$ there exists $K(\sigma) \ge 0$ such that for all $k \ge K(\sigma)$, there is a  splitting of $g^k_\#(\sigma)$ into  subpaths that are either legal or are \ipNp s. 
\item  There exists $K^*, L^* > 0$ so that if $\sigma$ is a path that decomposes as the concatenation of two legal subpaths each of length at least $L^*$  then $K(\sigma) \le K^*$. 
\item  If $\sigma$ is a bi-infinite path with exactly one illegal turn then one of the following holds for all $k \ge K^*$.
\begin{itemize}
\item  $g^k_\#(\sigma)$ is legal.   
\item $g^k_\#(\sigma)$ splits as $R_k\cdot \tau_k  \cdot R_k'$ where $R_k$ and $R_k'$ are legal rays and $\tau_k$ is an \ipNp. 
\end{itemize}
\end{enumerate}
\end{lemma}

\begin{proof}  The first three items follow from Lemma~4.2.5 and  Lemma~4.2.6 of \cite{BFH:TitsOne}. For the fourth item, let 
$$L^* = \Length(G) + \min\{\Length(e) \suchthat \text{$e$ is an edge of $G$}\}
$$
If $\sigma$ is legal then $K(\sigma)$ = 1 so we may assume that $\sigma$ has exactly one illegal turn. Decompose $\sigma$ as $\sigma = \sigma_1 \bar \sigma_2$  where $\sigma_1$ and $\sigma_2$ are legal and let $\sigma_{1,n}$ and $\sigma_{2,n}$ be the terminal segments of $\sigma_1$ and $\sigma_2$ that are identified by $g^n$. If there does not exist $n>0$ such that $\sigma_{i,k} = \sigma_{i,n}$ for all $k \ge n$ then the closure of $\cup_n\sigma_{1,n} \bar \sigma_{2,n}$ is an indivisible Nielsen or pre-Nielsen path by Lemma~4.2.5 of \cite{BFH:TitsOne} and $\sigma$ splits into an initial legal path, followed by $\cup_n\sigma_{1,n} \bar \sigma_{2,n}$ followed by a terminal legal path.  By item~(1), the pre-Nielsen paths become periodic Nielsen paths after a uniformly bounded number of iterations, and so $K(\sigma)$ is uniformly bounded. We may therefore assume that there is such an $n$. By Lemma~\ref{bcl} (see also Lemma~3.9 of \cite{BestvinaHandel:tt}) the length of $\sigma_{i,n}$ is bounded above by $\Length(G)$. Removing an initial segment of $\sigma_i$, $i=1,2$, that is disjoint from $\sigma_{i,n}$ does not change $K(\sigma)$. Since $\Length(\sigma_i) \ge L^*$ we may therefore remove such initial segments so that $\sigma$ has endpoints at vertices and has length at most $2\Length(G) + 2\max\{\Length(e) \suchthat \text{$e$ is an edge of $G$}\}$.  There are only finitely many such $\sigma$ and we take $K^*$ to be the maximum of their $K(\sigma)$'s.  This completes the proof of (4).

Lemma 4.2.2 of \cite{BFH:TitsOne} implies that a bi-infinite line with exactly one illegal turn splits as $R\cdot \tau \cdot R'$ where $R$ and $R'$ are legal rays and $\tau$ decomposes as the concatenation of two legal subpaths each of length at least $L^*$.  The fifth item therefore follows from the third and fourth items and the fact that $g_\#$ maps legal paths to legal paths and \ipNp s to \ipNp s.
\end{proof} 

\subparagraph{Assumption:} We assume from now on that for any irreducible train track map $\gG$, the endpoints of \ipNp s are vertices. This can always be arranged by a finite subdivision.

\paragraph{Principal vertices, principal lifts and rotationless train track maps.}  \hfill\break
We now recall from \cite{FeighnHandel:abelian} the concepts of principal automorphisms, principal lifts of train track representatives, and rotationless outer automorphisms. While we apply these concepts here only to fully irreducible outer automorphisms, more general definitions are given in \cite{FeighnHandel:abelian} which apply to arbitrary outer automorphisms and their relative train track representatives. 

Let $\phi \in \Out(F_r)$ be fully irreducible, and consider $\Phi \in \Aut(F_r)$ representing $\phi$ and its boundary extension $\hat\Phi \from \bdy F_r \to \bdy F_r$. Denote the fixed point set of $\hat\Phi$ by $\Fix(\hat\Phi)$, and denote the subset of nonrepelling fixed points of $\Fix(\hat \Phi)$ by $\Fix_N(\hat \Phi)$. We say that  $\Phi$ is a \emph{principal automorphism} if $\Fix_N(\hat \Phi)$ contains at least three points. The set of principal automorphisms representing $\phi$ is denoted $\PA(\phi)$. Two elements $\Phi, \Phi' \in PA(\phi)$ are \emph{equivalent}\footnote{This equivalence relation is called ``isogredience'' in \cite{Nielsen:InvestigationsIEnglish} and \cite{LevittLustig:SimpleDynamics}.} if there exists $\gamma \in F_n$ such that $\Phi'  = i_{\gamma} \composed\Phi \composed i_{\gamma}^{-1}$.  

If the fixed subgroup $\Fix(\Phi)$ is non-trivial, then $\Fix(\Phi) =\langle \gamma \rangle$  is infinite cyclic. This happens if and only if $\phi$ is geometric, in which case $\phi$ is represented by a pseudo-Anosov homeomorphism of a once punctured surface whose boundary is represented by $\gamma$ (see Theorem 4.1 of \cite{BestvinaHandel:tt}).  

The following lemma is contained in Lemma 2.3 and Lemma 2.4 of \cite{BFH:Solvable}  and in  Proposition 1.1 of \cite{GJLL:index}. 

\begin{lemma} \label{FixN} Assume that $\Phi \in \PA(\phi)$.
\begin{enumerate}
\item \label{ItemFiniteAttractors}
If $\Fix(\Phi)$ is trivial then $\Fix_N(\hat \Phi)$ is a finite set of attractors. 
\item \label{ItemAttractorOrbits}
If  $\Fix(\Phi) =\langle \gamma \rangle$ is infinite cyclic, then $\Fix_N(\hat \Phi)$ is the union of $t_{\gamma}^{\pm}$ with a finite set of $t_{\gamma}$-orbits  of attractors.
\item \label{ItemNoAxis}  
If $P \in \Fix_N(\hat \Phi)$ is an attractor then it is not $t_{\gamma}^+$ or $t_{\gamma}^-$ for any nontrivial covering translation $t_{\gamma}$.
\end{enumerate}
\end{lemma}
 
\begin{corollary} \label{distinct fixed sets} If $\Phi_1 \ne \Phi_2 \in \PA(\phi)$ then $\Fix(\hat \Phi_1)$ and $\Fix(\hat \Phi_2)$ are disjoint.  \end{corollary}

\begin{proof}  Suppose that $P \in \Fix(\hat\Phi_1) \intersect \Fix(\hat\Phi_2)$. Observe that $\Phi_1^\inv \Phi_2$ is a covering translation $t_\gamma$, for a nontrivial $\gamma \in F_r$. It follows that $t_\gamma(P)=P$ and hence that $\gamma \in \Fix(\Phi_1) \intersect \Fix(\Phi_2)$. In this case, $\phi$ is geometric and the corollary is well known; see for example Lemma~4.2 of \cite{HandelThurston:nielsen}. 
\end{proof}
  
Let $\phi \in \Out(F_r)$ be fully irreducible. We say that $\phi$ is \emph{rotationless} if for each $k > 0$ and each principal automorphism $\Phi_k$ representing $\phi^k$ there is a principal automorphism $\Phi$ representing $\phi$ such that $\Phi_k = \Phi^k$ and such that $\Fix_N(\hat \Phi) = \Fix_N(\hat \Phi_k)$. In other words, the set of principal automorphisms and their nonrepelling fixed points in $\partial F_r$ is stable under iteration.

Let $g \from G \to G$ be an irreducible train track representative of $\phi$. A lift $\breve g \from \wt G \to \wt G$ is a \emph{principal lift} if the corresponding automorphism is a principal automorphism, and so two principal lifts $\breve g$, $\breve g'$ are equivalent if and only if $\breve g'  = t_{\gamma} \composed \breve g \composed t_{\gamma}^{-1}$ for some $\gamma \in F_r$. Consider a periodic vertex $v \in G$. Since $v$ has at least one legal turn, $v$ has at least two periodic directions. We say that $v$ is a \emph{principal vertex} if there are at least three periodic directions at $v$ or if $v$ is the endpoint of an indivisible periodic Nielsen path.  The set of principal vertices is non-empty by Lemma 3.18 of \cite{FeighnHandel:abelian}.  If every principal vertex $v \in  G$ is fixed and if every periodic direction at a principal vertex is fixed then we say that $\gG$ is a \emph{rotationless train track map}.  

These definitions are connected by Proposition~3.24 of \cite{FeighnHandel:abelian}, which we recall.

\begin{proposition} Assume that $\gG$ is an irreducible train track map for a fully irreducible $\phi \in \Out(F_r)$.  Then $\phi$ is rotationless if and only if $\gG$ is rotationless.  
\end{proposition}

\subparagraph{Remark.}
It is clear that $g^p$, and hence $\phi^p$, is rotationless for some $p \ge 1$.  Some of our results are stated for rotationless $\phi$ and $\gG$; to apply these in general one replaces $\phi$ and $g$ by  suitable iterates.  

\bigskip

A pair of points $x,y \in \Fix(g)$ are \emph{Nielsen equivalent} if they are the endpoints of a Nielsen path in $G$.  The elements of a Nielsen equivalence class are usually referred to simply as  \emph{a Nielsen class}. It is well known, and easy to check, that $x$ and $y$ are Nielsen equivalent if and only if there is a lift $\breve g : \wt G \to \wt G$ and points $\ti x, \ti y \in \Fix(\breve g)$ that project to $x$ and $y$ respectively.    In particular, for any lift $\breve g$, the image of $\Fix(\breve g)$ in $G$ is either empty or is  an entire Nielsen class.  Moreover, if $\Fix(\breve g_1)$ and $\Fix(\breve g_2)$ are non-empty, then $\breve g_1$ and $\breve g_2$   are equivalent, as defined above in the context of principal lifts,  if and only if $\Fix(\breve g_1)$ and $\Fix(\breve g_2)$ project to the same Nielsen class in $\Fix(g)$. It is useful to note that if $\Fix(\breve g_1) \cap \Fix(\breve g_2) \ne \emptyset$ then $\breve g_1 = \breve g_2$.  This follows from the fact that $\breve g_1$ and $\breve g_2$ differ by a covering translation and the fact that every non-trivial covering translation is fixed point free. 

It will be convenient to refer to lifts $\ti v \in \wt G$ of principal vertices $v \in G$ as \emph{principal vertices of $\wt G$}.  

\begin{lemma} \label{principal essential}  Assume that $\gG$ is rotationless.  Then  $\PA(\phi) \ne \emptyset$.  Moreover,   
  $\breve g$ is a principal lift if and only if  $\Fix(\breve g) \ne \emptyset$ and some, and hence every, element of $\Fix(\breve g)$ is a principal vertex.  
\end{lemma}  

\begin{proof}  The second statement follows from Corollary~3.14 and Corollary~3.22 of \cite{FeighnHandel:abelian}.  The first follows from the second and  Lemma~3.18 of \cite{FeighnHandel:abelian}; see also Proposition I.5 of \cite{LevittLustig:PeriodicEnds}. 
\end{proof}

We also recall Lemma 3.23 of \cite{FeighnHandel:abelian}.  

\begin{lemma} \label{period one} If $g$ is rotationless then every periodic Nielsen path   has period one. \qed
\end{lemma}

\subsection{The attracting tree $\TPlus$}   
\label{SectionAttractingTree}
Given a fully irreducible $\phi \in \Out(F_r)$ there exists a unique point $\TPlus = \TphiPlus \in \bdy \X_r$ which is the attracting point for every forward orbit of $\phi$ in $\X_r$; this is proved in \BookOne. We review here the construction and properties of the $F_r$-tree $\TPlus$. Although we shall not need this fact here, we mention the result of  \cite{LevittLustig:NorthSouth} that the action of $\phi$ on $\overline X_r$ has source-sink dynamics with sink $\TphiPlus$ and source $\TphiMinus = \TPlus^{\phi^\inv}$. 

Let $\gG$ be an irreducible train track map that represents~$\phi$ and let $\rho_G \from R_r \to G$ be the marking on $G$.  For all $i\ge 0$, let $G_i$ be the marked graph obtained from $(g^i \rho_G, G)$ by multiplying the metric by $1/\lambda^i$, and so
$$\ell_{G_i}(c) = \ell_G(\phi^i(c))/\lambda^i \quad\text{for all $c \in \C$ and all $i \ge 0$.}
$$  
The map $g \from G \to G$ induces a homotopy equivalence denoted $g_{i+1,i} \from G_i \to G_{i+1}$ that preserves marking. By composition, for $0 \le i < j$, the map $g^{j-i} \from G \to G$ induces the marking preserving homotopy equivalence $g_{ji} = g_{j,j-1} \composed \cdots \composed g_{i+1,i} \from G_i \to G_j$ (see Figure~\ref{FigureNotation1}). We fix the base point of $G_i$ so that the marking $g^i \rho_G \from R_r \to G_i$ preserves base point. It follows that the maps $g_{ji} \from G_i \to G_j$ preserve base point.

\begin{figure}
$$\xymatrix{
	G_i		\ar[rr]^{g_{i+1,i}}
				\ar@/^3pc/[rrrrrr]^{g_{ji}}
	& & G_{i+1}	\ar[rr]^{g_{i+2,i+1}}
	& & \cdots 	\ar[rr]^{g_{j,j-1}}
	& & G_j
}$$
\caption{Notation associated to a train track map $g \from G \to G$.}
\label{FigureNotation1}
\end{figure}

The sequence $\ell_{G_i}$ converges  in $\whclosure\CV_r$ to an element $\ell_+ \in  \whclosure\CV_r$ whose projective class is $\phi$-invariant. To describe $\ell_+$ explicitly, consider $c \in \C$ represented by a circuit $\sigma \subset G$.   If $K(\sigma)$ is as in Lemma~\ref{nielsen paths} and  $k > K(\sigma)$ then $g^k_\#(\sigma)$ splits as a concatenation  of legal subpaths $\mu_{i,k}$ and periodic indivisible Nielsen paths $\nu_{i,k}$.  Excluding the $\nu_{i,k}$, define $\ell_+(c) = (\sum L_G(\mu_{i,k}))/\lambda^k$.   Since $g$ expands length of legal paths by the factor $\lambda$ and there are only finitely many periodic indivisible Nielsen paths,  $\ell_+(c)$ is independent of $k > K(\sigma)$ and is the limit of the sequence $\ell_{G_i}(c)$. It is immediate from the definition that  $\ell_+(\phi(c)) = \lambda \ell_+(c)$ and that if $c$ is represented by a legal circuit in $G$ then $\ell_+(c) = \ell_G(c)$. 

The $F_r$ tree $\TPlus$ that realizes $\ell_+$ is well defined up to isometric conjugacy.  

\paragraph{Direct limits.} It seems to be well known that $\TPlus$ is the $F_r$-equivariant Gromov limit of the trees $\wt G_i$; see e.g.\ \OuterLimits. We need some slightly stronger facts about this convergence.

Direct limits exist in the category of metric spaces and distance nonincreasing maps. That is, consider a direct system consisting of a partially ordered set $I$, an indexed collection of metric spaces $(A_i)_{i \in I}$ with metrics $d_i$, and an indexed collection of distance nonincreasing maps $(f_{ji}\from A_i\to A_j)_{i<j}$, such that $I$ is a \emph{directed set} meaning that  for all $i,j \in I$ there exists $k \in I$ with $i \le k$, $j \le k$, and the maps satisfy the compatibility property $f_{kj}\composed f_{ji} = f_{ki}$ for $i<j<k$ in $I$. Then there exists a metric space $A_\infinity$ and distance nonincreasing maps $f_{\infinity,i} \from A_i \to A_\infinity$ satisfying the following universality property:  given a metric space $B$ and a system of distance nonincreasing maps $h_i \from A_i \to B$ compatible with the maps $f_{ji}$, meaning that $h_j \composed f_{ji} = h_i$, there exists a unique distance nonincreasing map $h_\infinity \from A_\infinity \to B$ such that $h_\infinity \composed f_{\infinity,i} = h_i \from A_i \to B$ for all $i \in I$. To construct $A_\infinity$, first let $\hat A_\infinity$ be the direct limit of the system $(A_i)$, $(f_{ji})$ in the category of sets \cite{Spanier}: points of $A_\infinity$ are equivalence classes of the disjoint union of the $A_i$ where $x\in A_i$ is equivalent to $y \in A_j$ if there exists $k$ with $i,j\le k$ such that $f_{ki}(x) = f_{kj}(y)$. Then define a pseudometric on $\hat A_\infinity$ where for any two equivalence classes $[x], [y]$ we have
$$d_\infinity([x],[y]) = \inf_{k \ge i,j} d_k(f_{ki}(x),f_{kj}(y))
$$
The direct limit $A_\infinity$ is the natural quotient metric space of $\hat A_\infinity$, where $[x], [y] \in A_\infinity$ are equivalent if and only if $d_\infinity([x],[y])=0$. The map $f_{\infinity,i} \from A_i \to A_\infinity$ is the composition of the quotient maps $A_i \to \hat A_\infinity \to A_\infinity$. 

\begin{figure}
$$\xymatrix{
	\wt G=\wt G_0 		\ar[rr]^{\ti g_{1,0}} 
				\ar@/^5pc/[rrrrr]^(.25){f_g=f_{g_0}}
				\ar@(dl,dr)[]_{\breve g = \breve g_0}
	& & \wt G_1 	\ar[r]^{\ti g_{2,1}} 
				\ar@/^3pc/[rrr]^(.25){f_{g_1}}
				\ar@(dl,dr)[]_{\breve g_1}
	& \wt G_2 	\ar[r]^{\ti g_{3,2}} 
				\ar[r] \ar@/^2pc/[rr]^(.25){f_{g_2}}
				\ar@(dl,dr)[]_{\breve g_2}
	& \ldots 
	& \TPlus		\ar@(dl,dr)[]_{\Phi_+}
}$$
\caption{Notation associated to a train track map $g \from G \to G$ and automorphism $\Phi$ representing $\phi$. Maps denoted $\tilde g$ depend only on $g$ and are $F_r$-equivariant. Maps denoted $\breve g$ depend on $g$ and on $\Phi$ and satisfy the twisted equivariance equation~(\ref{EquationTwisted}). The maps $f_{g_i}$ are defined after Theorem~\ref{T+ is direct limit}, and $\Phi_+$ is defined in Theorem~\ref{TheoremPhiPlus}.}
\label{FigureNotation2}
\end{figure}
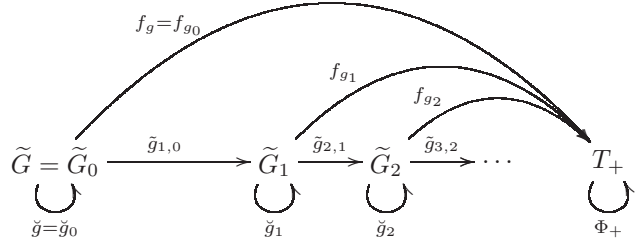

Consider the trees $\wt G_i$. We lift the base point of $G_i$ to $\wt G_i$, and we choose lifts $\wt g_{i+1,i}  \from \wt G_{i} \to \wt G_{i+1}$ of $g_{i+1,i} \from G_i \to G_{i+1}$ for all $i$, so that the maps $\wt g_{i+1,i}$ preserve base points, thereby making $\wt G_i$ into an $F_r$ tree, and making each map $\wt g_{i+1,i}$ into an $F_r$ equivariant map which is an isometry when restricted to each edge. We obtain a direct system with $\wt g_{ji} \from \wt G_i \to \wt G_j$ defined inductively by $\wt g_{ji} = \wt g_{j,j-1} \composed \wt g_{j-1,i}$. Whenever $\Phi \in \Aut(F_r)$ representing $\phi$ is implicitly chosen, the associated lifts of $g_i \from G_i \to G_i$ shall be denoted $\breve g_i \from \wt G_i \to \wt G_i$ for all $i \ge 0$ (see Figure~\ref{FigureNotation2}). These maps satisfy the following properties:
\begin{itemize}
\item $\breve g_i$ restricted to each edge is a homothety with stretch equal to $\lambda$. 
\item $\wt g_{ji}\breve g_i = \breve g_j \wt g_{ji}$. 
\item (Twisted equivariance, see (\ref{EquationTwisted})) $\breve g_i \composed t_\gamma = t_{\Phi(\gamma)} \composed \breve g_i$ for  each $\gamma \in F_r$.  
\end{itemize}
While we recall that the notation $\breve{g}$ is meant to suggest twisted equivariance (\ref{EquationTwisted}), on the other hand the notation $\wt g$ is meant to suggest $F_r$-equivariance.

The following theorem and corollary seem to be folklore. See e.g.\ \cite{HandelMosher:parageometric} for a proof of Theorem \ref{T+ is direct limit}, and see Section~\ref{SectionDirectLimits} below for a proof of a more general fact in the nongeometric case.

\begin{theorem}  \label{T+ is direct limit} The direct limit of the direct system $\wt g_{ij} :\wt G_{j} \to \wt G_i$, in the category of metric spaces with an $F_r$ action whose morphisms are $F_r$ equivariant, distance nonincreasing maps, is an $F_r$ tree in the same homothety class as $\TPlus$. \qed
\end{theorem}

In the context of this theorem, we shall normalize the metric on $\TPlus$ so that $\TPlus$ is identified with the direct limit.  From Theorem~\ref{T+ is direct limit} we obtain a direct limit map $f_g \from \wt G \to \TPlus$, and by construction $f_g$ is an isometry when restricted to each edge of $\wt G$. As the notation $f_g$ suggests, this map is well-defined up to isometric conjugacy of the $F_r$ action on $\TPlus$ (see Figure~\ref{FigureNotation2}), depending only on $g \from G \to G$, and \emph{not} depending on the other choices made in the definition of the direct system $\wt g_{ij} \from \wt G_j \to \wt G_i$. This well-definedness is a consequence of the direct limit description of $\TPlus$, from which it follows that $\TPlus$ is the metric space associated to a pseudometric on $\wt G_0$ where the distance between $\ti x,\ti y \in \wt G_0$ equals the limit of the nonincreasing sequence $\Length\left(g^k_\#(pr[\ti x,\ti y])\right) \bigm / \lambda^k$. Clearly the pseudometric depends only on $g$ and therefore so does $\TPlus$.

\begin{corollary} \label{first tree corollary} Assume that $\gG$ is rotationless and that $\wt G$ is a simplicial $F_r$-tree obtained from the universal cover of $G$ by choosing a marking and a lifting of the base point. Then there is a surjective equivariant map $f_g : \wt G \to \TPlus$ such that for all $\ti x, \ti y \in \wt G$, letting $\sigma = pr[\ti x, \ti y]$, the following are equivalent.
\begin{enumerate}
\item  $f_g(\ti x) = f_g(\ti y)$.  
\item  For some $k \ge 0$, the path $g^k_\#(\sigma)$ is either trivial or Nielsen.
\end{enumerate}
In particular,  $f_g$ restricts to an isometry on all legal paths. 
\end{corollary}

\begin{proof}   Assume the notation of Theorem~\ref{T+ is direct limit} with $\wt G = \wt G_0$ and let $f_g : \wt G \to \TPlus$ be  the direct limit map.  

     We first show that (2) implies (1). By construction, $g^k_\#(\sigma)$ is trivial if and only if $\wt g_{k,0}(\ti x) = \wt g_{k,0}(\ti y)$ and the latter implies that $f_g(\ti x) = f_g(\ti y)$. To prove (1) assuming that  $\sigma$ is a Nielsen or   pre-Nielsen path, it suffices to assume  that $\sigma$ is indivisible.    By Lemma~\ref{nielsen paths}-(2), there exist $\ti x_i \to \ti x$ and $\ti y_i \to \ti y$ such that $\wt g_{k,0}(\ti x_i) = \wt g_{k,0}(\ti y_i)$ for all $i$ and all sufficiently large $k$.  It follows that $f_g(\ti x_i) =  f_g(\ti y_i)$ for all $i$ and hence that  $f_g(\ti x) = f_g(\ti y)$.   

Suppose now that $f_g(\ti x) = f_g(\ti y)$. If each $g^k_\#(\sigma)$ is non-trivial, then Lemma~\ref{nielsen paths} implies that $g^k_\#(\sigma)$ splits as a concatenation of legal paths and \iNp s  for all sufficiently large $k$.  If $\mu$ is any legal path in this splitting, then the distance between $\wt g_{k,0}(\ti x)$ and $\wt g_{k,0}(\ti y)$ is uniformly bounded below by the length of $\mu$ in contradiction to the definition of $f_g$ and the assumption that $f_g(\ti x) = f_g(\ti y)$.  Thus  no such $\mu$ exists and $g^k_\#(\sigma)$ is Nielsen. 

We have now proven the equivalence of  (1) and (2).  Since a legal path does not contain subpaths $\sigma$ satisfying (2), the restriction of $\wt g_{ij}$ to any legal path is an embedding, and hence an isometry, and it follows that $f_g$ restricts to an isometry on all legal paths.
\end{proof} 

Next we review the well known fact that there is a kind of action on the $\reals$-tree $\TphiPlus$ by the set of automorphisms of $F_r$ representing a given fully irreducible outer automorphism $\phi$. This set is not a group, instead it is a coset of the normal subgroup $\Inn(F_r)$ in the group $\Aut(F_r)$ --- formally this coset is identified with the outer automorphism $\phi \in \Out(F_r) = \Aut(F_r) / \Inn(F_r)$. So while it does not make sense to ask for a ``group action'' of the coset, it does make sense to ask for an ``$\Inn(F_r)$-affine action'', which is what item~(\ref{ItemPhiPlusEquiv}) gives us:

\begin{theorem} \label{TheoremPhiPlus}
Let $\phi \in \Out(F_r)$ be fully irreducible, and denote $\TPlus = \TphiPlus$. To each $\Phi \in \Aut(F_r)$ representing $\phi$ there is associated a homothety $\Phi_+ \from \TPlus \to \TPlus$ with stretch factor $\lambda$ satisfying the following properties:
\begin{enumerate}
\item \label{ItemPhiPlusNatural}
For each irreducible train track representative $g \from G \to G$ of $\phi$, letting $\breve g \from \wt G \to \wt G$ be the lift corresponding to $\Phi$, we have
$$\Phi_+ \composed f_g = f_g \composed \breve g
$$
\item \label{ItemPhiPlusEquiv}
For each $\Phi,\Phi'$ representing $\phi$, if $\gamma \in F_r$ is the unique element such that $\Phi' = i_\gamma \composed\Phi$ then $\Phi'_+ = t_\gamma \composed\Phi_+$.
\end{enumerate}
\end{theorem}

\begin{proof} Choose any train track representative $g \from G \to G$ of $\phi$, and 
refer to Figure~\ref{FigureNotation1} for notations associated to $g$. Let $\breve g, \breve g' \from \wt G \to \wt G$ be the lifts associated to $\Phi,\Phi'$, respectively. Refer to Figure~\ref{FigureNotation2} for notations associated to $\breve g$. Refer again to Figure~\ref{FigureNotation2} for notations associated to $\breve g'$, but with the prime symbol ${}^\prime$ added to these notations. In particular we have identifications $\wt G = \wt G_0$,  $\breve g = \breve g_0, \breve g' =  \breve g'_0$.

Let $T^*_+$ be $\TPlus$ with the metric scaled down by the factor $\lambda$ and let $h : \TPlus \to T^*_+$ be the \lq identity\rq\ map; thus $h$ contracts length uniformly by $\lambda$.  By the universal property of direct limits applied to the maps $h\composed f_{g_i} \composed \breve g_i : \wt G_i \to T^*_+$ there is a map  $\Phi^*_+ : \TPlus \to T^*_+$ such that $h \composed f_{g_i} \composed \breve g_i = \Phi^*_+ \composed f_{g_i}$.  Let $\Phi_+ :\TPlus \to \TPlus$ be $h^{-1} \composed\Phi^*_+$. It follows that $f_{g_i} \composed \breve g_i = \Phi_+ \composed f_{g_i}$, which implies~(\ref{ItemPhiPlusNatural}) by taking $i=0$. Similarly $\Phi'_+$ is defined satisfying $f_{g_i} \composed \breve g'_i = \Phi'_+ \composed f_{g_i}$, which implies~(\ref{ItemPhiPlusEquiv}) together with the fact that $\breve g'_i = t_\gamma \composed \breve g_i$ and with the universal property of direct limits.  

Given $z_1$ and $z_2$ in $\TPlus$, choose $\ti x_j \in \wt G_0$ such that $f_0(\ti x_j) = z_j$. By the third item of Lemma~\ref{nielsen paths},  the projection into $G_i$ of the path in $\wt G_i$ connecting $\wt g_{i0}(x_1)$ to $\wt g_{i0}(x_2)$ splits into subpaths that are either legal or indivisible Nielsen paths for all sufficiently large $i$.  Corollary~\ref{first tree corollary} implies that the distance between $z_1$ and $z_2$ is the sum of the lengths of the legal subpaths of $G_i$ in this splitting.   Since $g$ expands length by a factor of $\lambda$ on all legal subpaths, (1) implies that $\Phi_+$ must expand the distance between $z_1$ and $z_2$ by this same factor and since $z_1$ and $z_2$ were arbitrary, it follows that  $\Phi_+$ expands length uniformly by $\lambda$. This proves that $\Phi_+$ is a homothety onto its image with stretch factor $\lambda$, but $\Phi_+$ is surjective by minimality of the $F_r$-tree $\TPlus$.
\end{proof}

Next we describe the connection between principal automorphisms representing a fully irreducible $\phi \in \Out(F_r)$ and the set of branch points of $\TphiPlus$.

\begin{lemma}\label{principal lifts and branch points}  Suppose that $\phi \in \Out(F_r)$ is fully irreducible and rotationless, $g \from G \to G$ is a train track representative,  $\Phi$ is a principal automorphism representing $\phi$ and $\breve g : \wt G \to \wt G$ is the principle lift corresponding to $\Phi$, and $\Phi_+ : \TPlus \to \TPlus$ is as in Theorem~\ref{TheoremPhiPlus}.  Then 
\begin{enumerate}
\item $f_g(\Fix(\breve g))$ is a branch point $b$ of $\TPlus$, and  $\Fix(\Phi_+)=\{b\}$.
\item Every direction based at $b$ is fixed by $D\Phi_+$ and has the form $Df_g(d)$ where $d$ is a fixed direction based at some $\ti v \in \Fix(\breve g)$.
\item If $f_g(\ti x) = b$  and $d$ is a  direction based at $\ti x$ then 
   \begin{itemize}
    \item For some $m>0$, $\breve g^m(\ti x) \in \Fix(\breve g)$.
    \item $f_g(\breve g^k(\ti x))$ and $Df_g(D\breve g^k(d))$ are independent of  $k\ge 0$.
    \end{itemize}
\item  The assignment $\Phi \mapsto   b= \Fix(\Phi_+)$ of (1) defines  a bijection between $\PA(\phi)$ and the set of branch points of $\TPlus$.
\end{enumerate}
\end{lemma}

\begin{proof}  First we prove~(1). Lemma~\ref{principal essential} implies that   $\Fix(\breve g)\ne \emptyset$ and that the elements of $\Fix(\breve g)$ are principal vertices.  By Corollary~\ref{first tree corollary} and Theorem~\ref{TheoremPhiPlus}, the point $b = f_g(\ti v)$ is independent of $\ti v \in \Fix(\breve g)$ and $\{b\}=\Fix(\Phi_+)$, and there is an induced map $Df_{g}$ from directions at $\ti v$ to directions at $b$ that is injective on the set of fixed directions. If there are at least three fixed directions at $\ti v $ then $b$ is a branch point. Otherwise, there is an indivisible Nielsen path $\sigma$ and a lift $\ti \sigma$ with one endpoint at $\ti v$ and the other at some principal $\ti w \in \Fix(\breve g)$. Let $d_1$ and $d_2$ be distinct fixed directions at $\ti v$ with $d_2$ being the initial direction of $\ti \sigma$ and let $d_3$ and $d_4$ be distinct fixed directions at $\ti w$ with $d_3$ being the initial direction of the inverse of $\ti \sigma$. Then $Df_{g}(d_1) \ne Df_{g}(d_2) = Df_{g}(d_3) \ne Df_{g}(d_4)$. To prove that $b$ is a branch point, it suffices to prove that $Df_{g}(d_1) \ne Df_{g}(d_4)$.  If this fails then there exist $\ti x$ close to $\ti v$ in the direction $d_1$ and $\ti y$ close to $\ti w$ in the direction $d_4$ such that $f_{g}(\ti x) = f_{g}(\ti y)$.  The second item of Lemma~\ref{nielsen paths} implies that the projection into $G$ of the path connecting $\ti x$ to $\ti y$ splits as a non-trivial legal path followed by the \iNp\ $\sigma$ followed by a non-trivial legal path. Corollary~\ref{first tree corollary} then implies that $\ti x$ and $\ti y$ have distinct $\ti f_{g}$-images and this contradiction  completes the proof of item~(1).

We now turn to (3).  Choose $\ti v \in \Fix(\breve g)$. Corollary~\ref{first tree corollary} implies that if $f_g(\ti x) = b$ then there exists $m \ge 0$ such that either $\breve g^m(\ti x)  =\ti v$ or the path connecting $\breve g^m(\ti x)$ to $\ti v$ is the lift of a Nielsen path.  In either case, $\breve g^m(\ti x) \in \Fix(\breve g)$, which implies that $f_g(\breve g^k(\ti x)) = b$ for all $k \ge 0$.    After increasing $m$, we may assume that $D \breve g^m(d) \in \Fix(D \breve g)$ for all directions $d$ based at $\ti x$.  From the fact that $f_g \circ \breve g = \Phi_+ \circ f_g$ we conclude that $ D\Phi_+^{k}\circ Df_g(d) = Df_g(D\breve g^k(d))$  is independent of $k \ge m$.  Since $\Phi_+$ is a homeomorphism, they are in fact independent of $k \ge 0$.   This completes the proof of item~(3) and also implies that $Df_g(d) \in \Fix(D\Phi_+)$.  Since every direction at $b$ occurs as $Df_g(d)$ for some $\ti x$ and $d$, we have also proved item~(2).
   
Next we prove injectivity of the map $\Phi \mapsto b$. Suppose that $\ti v_1$ and $\ti v_2$ are distinct principal vertices and that $f_{g}(\ti v_1) = f_{g}(\ti v_2)$.  Let $\breve g_1$ be the lift that fixes $\ti v_1$. Since principal vertices of $\wt G$ project into $\Fix(g)$, the restriction of $\breve g_1$ to the set of principal vertices of $\wt G$ is a bijection.  It follows that $\breve g^k(\ti v_1) \ne \breve g^k(\ti v_2)$ for all $k$.  The path $\ti \sigma$ connecting $\ti v_1$ to $\ti v_2$ cannot project to a pre-Nielsen path $\sigma$ which is not Nielsen.  Indeed if it did, then for some $k > 0$, the Nielsen path $\sigma'=f^k_\#(\sigma)$ has the same endpoints as $\sigma$ and $f_\#^k(\sigma') = f_\#^k(\sigma)$ which contradicts that $\sigma' \ne \sigma$ and the fact that $f$ is a homotopy equivalence.  Corollary~\ref{first tree corollary} therefore implies that $\ti \sigma$ projects to a Nielsen path and hence that $\ti v_2 \in \Fix(\breve g_1)$.  This proves that $\Phi \mapsto b$ is injective.

To see that $\Phi \mapsto b$ is onto it suffices to show that if $b$ is a branch point then $f_{g}(\ti c) = b$ for some principal vertex $\ti c$; one then chooses $\Phi$ to be the automorphism corresponding to the lift of $g$ that fixes $\ti c$.  For $1 \le l \le 3$ there are points $\ti a_l \in \wt G$ and directions $d_l$ at $\ti a_l$ such that $f_g(\ti a_l)=b$, and such that $Df_{g}(d_l)$ are three distinct directions at $b$.  Let $a_j \in G$ be the projected image of $\ti a_j$.  On the one hand, if $g^m(a_1) = g^m(a_2) = g^m(a_3)$ for some $m \ge 0$ then there are three directions at this point that have distinct $Df_g$-images and so must have distinct $Dg^k$-images for all $k \ge 0$.  Otherwise, Corollary~\ref{first tree corollary} implies that $g^m(a_j)$ is the endpoint of a Nielsen path for some $j$ and some $m \ge 0$. In either case $c := g^m(a_j)$ is a principal, and hence  fixed, vertex.  There is a path $\sigma$ connecting $a_j$ to $c$ so that $g^m_\#(\sigma)$ is trivial.  The lift $\ti \sigma$ with endpoint $\ti a_j$ terminates at a principal vertex $\ti c$ and $\breve g(\ti c) = \breve g(\ti a_j) = b$.  This 
completes the proof of (4) and so the proof of the lemma. 
\end{proof}

\subsection{Geodesic laminations in trees and marked graphs}
\label{SectionLaminations}
Given a topological space $S$ let $D^2 S$ denote the double space of $S$,
the space of distinct ordered pairs in $S$:
$$D^2 S = \{(\xi,\eta) \in S \cross S \suchthat \xi \ne \eta\}
$$
The group $\Z/2$ acts freely on $D^2 S$ by permuting coordinates.

Define the \emph{geodesic leaf space} of $F_r$, denoted $\G F_r$, to be
the set of distinct \emph{unordered} pairs of points in $\bdy F_r$, that
is, the quotient space of the action of $\Z/2$ on $D^2 \bdy F_r$:
$$\G F_r = (D^2 \bdy F_r) \biggm/ \Z/2
$$
The space $\G F_r$ is locally a Cantor set. We will usually abuse ordered
pair notation and write $(\xi,\eta)$ for a point in $\G F_r$. The action
of $F_r$ on $\bdy F_r$ induces an action on $\G F_r$ which has a dense
orbit. The quotient space $\G F_r / F_r$ is therefore a non-Hausdorff
topological space; this space is denoted $\B$ in \cite{BFH:TitsOne}. 

Given a free, simplicial $F_r$-tree $T$, we also define the
geodesic leaf space 
$$\G T = (D^2 \bdy T) \biggm/ (\Z/2)
$$
The natural equivariant homeomorphism $\bdy F_r \homeo \bdy T$ induces a
natural equivariant homeomorphism $\G F_r \homeo \G T$, and we will
henceforth identify $\G T$ with $\G F_r$. 

It will be convenient to think of $\G T$ not just in the abstract as a
set of unordered distinct pairs of boundary points, but more concretely as
a set whose elements are bi-infinite geodesics, and even more concretely
as a lamination in the sense of topology or dynamics, which is a
topological space decomposed into \nb{1}manifolds fitting together locally
as an interval crossed with a transversal. This approach is used, for
example, in \cite{CoornaertPapadopoulos:Hopf}. One advantage of this
approach is that it avoids conceptual difficulties of non-Hausdorff
spaces. We briefly review this point of view.

Define a \emph{parameterized geodesic} in $T$ to be an isometric embedding
$\gamma\from\reals \to T$, and define the \emph{geodesic flow} of $T$,
denoted $\GFL T$, to be the space of parameterized geodesics with the compact open
topology. The group $\reals$ acts homeomorphically on $\GFL T$ by
translating the domain: $\gamma\cdot r(t) = \gamma(t+r)$. An orbit of the
$\reals$ action on $\GFL T$ is called a \emph{geodesic flow line}, and
there are natural bijections between the set of geodesic flow lines, the
set of oriented bi-infinite geodesics, and the set $D^2 \bdy T$. The
geodesic flow lines decompose $\GFL T$, and they fit together locally
like the product of an interval crossed with a Cantor set, that is, a
\nb{1}dimensional lamination with Cantor cross sections. The group
$\Z/2$ acts on $\GFL T$ by reflecting the domain $\reals$ through the
origin. This action is free on the set of geodesic flow lines,
interchanging a flow line with its orientation reversal. The quotient
$\GL T = \GFL T / (\Z/2)$ is called the \emph{geodesic lamination} of
$T$. An element of $\GL T$ can be naturally identified with a ``pointed
geodesic'' in $T$, a bi-infinite geodesic equipped with a base point. The
decomposition of $\GFL T$ into flow lines descends to a decomposition of
$\GL T$ into sets called \emph{geodesic leaves}, where each leaf
corresponds to different choices of base points along some bi-infinite
geodesic in $T$; the geodesic leaves form a \nb{1}dimensional lamination
with Cantor cross sections. We obtain natural bijections between the set
of geodesic leaves, the set of unoriented bi-infinite geodesics, and the
set $\G T$. 

The action of $F_r$ on $\G T$ is free, properly discontinuous, and
cocompact, and the quotient space $\G T / F_r$ inherits the structure of
a \nb{1}dimensional lamination with Cantor cross sections. We can view
this more intrinsically as the \emph{geodesic lamination} of the quotient
marked graph $G = T / F_r$, as follows. We can define the geodesic flow
$\GFL G$ as above, as parameterized local geodesics $\gamma \from \reals
\to G$ with the compact open topology, equipped with the $\reals$ action
by translation of the domain. The group $\Z/2$ acts by reflecting the
domain through the origin, and we obtain the quotient $\GL G = \GFL G /
(\Z/2)$. Both $\GFL G$ and $\GL G$ are \nb{1}dimensional laminations with
Cantor cross sections, and there are natural identifications $\GFL G
\approx \GFL T / F_r$, $\GL G = \GL T / F_r$. 

Given $\lambda = (\xi,\eta) \in \G F_r$ let $\lambda_T$ denote the leaf of $\GL T$ with ideal endpoints $\xi,\eta$; we say that $\lambda_T$ is the \emph{realization} of $\lambda \in \G F_r$. This relation is an equivariant bijection between points of $\G F_r$ and leaves of $\GL T$. Moreover this relation respects the topologies in the following sense: a sequence $\lambda_i \in \G F_r$ converges to $\lambda \in \G F_r$ if and only if the sequence $(\lambda_i)_T$ converges to $\lambda_T$ in the Gromov-Hausdorff topology on closed subsets of $T$, meaning that the intersections with any closed ball in $T$ converge in the Hausdorff topology on closed subsets of that ball.

A \emph{sublamination} of $\G F_r$ is a closed, $F_r$-equivariant subset of $\G F_r$. Given a free, simplicial $F_r$-tree $T$, a \emph{sublamination} of $\GL T$ is a closed, $F_r$-equivariant union of leaves of $\GL T$. And, given a marked graph $G$, a \emph{sublamination} of $\GL G$ is a closed union of leaves of $\GL G$. For any $T$ and for $G = T / F_r$, and for any sublamination $\Lambda \subset \G F_r$, the \emph{realization} of $\Lambda$ in $T$ is the sublamination of $\GL T$ defined by  
$$\Lambda_T = \{\lambda_T \suchthat \lambda \in \Lambda\}
$$
Also, the \emph{realization} of $\Lambda$ in $G$ is the quotient $\Lambda_G = \Lambda_T / F_r$, a sublamination of $\GL G$. This gives bijections between sublaminations of $\G F_r$, sublaminations of $\GL T$, and sublaminations of $\GL G$. We make free use of these bijections without comment, and we freely use the symbol $\Lambda$ for sublaminations in any of these contexts. When $\Lambda$ stands for a sublamination of $\G F_r$ then we sometimes use $\Lambda(T)$ or $\Lambda(G)$ instead of $\Lambda_T$ or $\Lambda_G$. When $\Lambda$ stands for a sublamination of $\GL G$, the corresponding sublamination of $\GL T$ is usually denoted $\wt\Lambda$. 

A sublamination $\Lambda \subset \G F_r$ is \emph{minimal} if and only if every $F_r$-orbit in $\Lambda$ is dense in $\Lambda$; equivalently, for each free simplicial $F_r$-tree $T$, each orbit of leaves of $\Lambda_T$ is dense in $\Lambda_T$; equivalently, for each marked graph $G$, each leaf of $\Lambda_G$ is dense in $\Lambda_G$. 

The group $\Aut(F_r)$ acts on $\bdy F_r$ and therefore also on the geodesic leaf space $\G F_r$. Note that this action agrees, under restriction, with the given action of $F_r \approx \Inn(F_r)$ on $\G F_r$. Choosing a marked graph $G$ with universal cover $T$, $\G F_r$ is identified with the set of leaves of $\GL T$ and so $\Aut(F_r)$ acts on the latter, in agreement under restriction with the action of $F_r \approx \Inn(F_r)$. Modding out by $F_r \approx \Inn(F_r)$ we obtain an action of $\Out(F_r)$ on the set of leaves of $\GL G$. We will use these actions without reference in what follows.

\subsection{The expanding lamination $\Lambda_-$}
\label{SectionExpandingLamination}
The chief example of a sublamination of $\G F_r$ is the expanding lamination $\Lambda_-$ of a fully irreducible outer automorphism $\phi \in \Out(F_r)$, as constructed in \BookZero. We review here the construction, well-definedness, and properties of $\Lambda_-$. In particular we present key features of $\Lambda_-$ in Lemmas~\ref{first characterization} and~\ref{periodic leaves}, the second of which motivates our choice of notation. The construction of $\Lambda_-$ will be revisited in Section~\ref{SectionLegalLams} when we extend the construction to a more general context. Complete details on $\Lambda_-$ can be found in either section 1 of \cite{BFH:laminations} or section 3 of \cite{BFH:TitsOne}.  

\paragraph{Definition and basic properties of $\Lambda_-$.} Consider any train track representative $g \from G \to G$ of $\phi$. Choose a sequence of pairs $(E_i,x_i)$ so that each $E_i$ is an edge of $G$, $x_i \in \interior(E_i)$, and we assume that $\lambda^i \, d(x_i,\bdy E_i)  \to +\infinity$; this is true, for example, if $(x_i)$ is a periodic orbit. Then $(g^i(E_i),g^i(x_i))$ is a sequence of pointed geodesic arcs whose lengths go to $+\infinity$ so that the distance from the base point to the endpoints also goes to $+\infinity$. We may therefore choose a subsequence of $(E_i,x_i)$ so that $(g^i(E_i),g^i(x_i))$ converges to a bi-infinite, pointed geodesic in $G$. The union of such geodesics, over all the choices made, is clearly a sublamination of $\GL G$ that we denote $\Lambda_-(g)$. In \BookZero\ it is proved that the expanding lamination $\Lambda_- \subset \G F_r$ of $\phi$ is well-defined and that its realization in $G$, for each train track representative $g \from G \to G$ of $\phi$, is $\Lambda_-(g)$. If the train track map is clear from the context then we will simply denote $\Lambda_-(g)$ as $\Lambda_-(\phi)$ or just as $\Lambda_-$. 

The expanding lamination for $\phi^p$, $p > 0$, is evidently the same as the expanding lamination for $\phi$ so we may replace $\phi$ by an iterate when it is convenient.   

By irreducibility of $g \from G \to G$, every leaf or half-leaf of $\Lambda_-(g)$ covers every edge of $G$. Also, $\Lambda_-$ is carried by no proper free factor of $F_r$ (see Example 2.5~(1) in \BookZero), so the realization of $\Lambda_-$ in any marked graph has the property that every leaf covers every edge. 

\paragraph{Action of $\phi$ on $\Lambda_-$.} Under the action of $\Out(F_r)$ on the set of leaves of $\GL G$ described at the end of Section~\ref{SectionLaminations}, the outer automorphism $\phi$ leaves invariant the set of leaves of $\Lambda_-$. We denote this restricted action by $\ell \mapsto \phi_\#(\ell)$, for each leaf $\ell$ of $\Lambda_-$. This action is characterized by saying that for each train track representative $g \from G \to G$ of $\phi$ and each leaf $\ell$ of $\Lambda_-$, $g$ takes the realization of $\ell$ in $G$ to the realization of $\phi_\#(\ell)$ in $G$. For any $\Phi \in \Aut(F_r)$ that represents $\phi$, there is an action of $\Phi$ on the set of leaves of $\wt\Lambda_-$, denoted $\ti\ell \mapsto \Phi_\#(\ti\ell)$, which is well-defined by saying that for each train track representative $g \from G \to G$, letting $\breve g \from \wt G \to \wt G$ be the lift of $g$ corresponding to $\Phi$, the map $\breve g$ takes the realization of $\ti\ell$ in $\wt G$ to the realization of $\Phi_\#(\ti\ell)$ in $\wt G$. To justify the existence of these actions see the final comments of Section~\ref{SectionLaminations}.

A leaf $\ell$ of $\Lambda_-$ is \emph{periodic} if $\phi^i_\#(\ell)=\ell$ for some $i \ge 1$, and a leaf $\ti\ell$ of $\wt\Lambda_-$ is periodic if there exists $\Phi \in \Aut(F_r)$ representing an iterate of $\phi$ such that $\Phi_\#(\ti\ell)=\ti\ell$. Periodicity of leaves is preserved under projection from $\wt\Lambda_-$ to $\Lambda_-$.

\paragraph{$\Lambda_-$ and local Whitehead graphs.} From the definition of $\Lambda_-(g)$, for any $x \in G$ the local Whitehead graph $W(x;G)$ has a description in terms of $\Lambda_-(g)$, as follows. A nondegenerate turn $(d_1,d_2)$ at $x$ is taken by some leaf of $\Lambda_-(g)$ if and only if this turn is taken by $g^k(E)$ for some edge $E$ of $G$ and some $k \ge 1$, which holds if and only if $(d_1,d_2)$ are the endpoints of an edge of the local Whitehead graph $W(x;G)$.   

\paragraph{Singular leaves of $\Lambda_-$ and stable Whitehead graphs.} Next we introduce notation and terminology to refer to certain leaves of $\Lambda_-$, called ``singular leaves'' in analogy to singular leaves of the unstable foliation of a pseudo-Anosov surface homeomorphism. These singular leaves can be used to relate $\Lambda_-$ to stable Whitehead graphs.

Let $\gG$ be an irreducible train track map representing an iterate of $\phi$ such each principal or nonprincipal periodic vertex $v$ is fixed and each periodic direction at $v$ is fixed. Choose a lift $\ti v \in \wt G$ and a lift $\breve g : \wt G \to \wt G$ that fixes~$\ti v$. If~$d$ is a direction at $\ti v$ fixed by $D\breve g$ and $\wt E$ is the edge at $\ti v$ corresponding to $d$, then $\breve g^j(\wt E)$ is a proper initial subpath of $\breve g^{j+1}(\wt E)$ for each $j \ge 0$ and so $\wt R = \cup_{j=0}^{\infty} \breve g^j(\wt E)$ is a ray in $\wt G$ converging to some $P \in \Fix_N(\hat g)$. We say that $\wt R$ is the \emph{ray determined by $d$} (or \emph{by $\wt E$}), and if $v$ is principal then $\wt R$ is \emph{a singular ray}.

Let $\wt R_1$, $\wt R_2$ be rays at $\ti v$ determined by distinct fixed directions $d_1,d_2$, corresponding to edges $\wt E_1$, $\wt E_2$ at $\ti v$. We claim that if $d_1$ and $d_2$ are endpoints of an edge in $SW(\ti v;\wt G)$ (which is always true if $\ti v$ is not principal), then $\ti\ell = \wt R_1 \union \wt R_2$ is a leaf of $\wt \Lambda_-$ realized in $\wt G$. This follows because there exists an edge $\wt E$ of $\wt G$ and $l > 0$ such that $\breve g^l(\wt E)$ contains the turn $(d_1,d_2)$, and so as $i \to +\infinity$ the paths $\breve g^{l+i}(\wt E)$ converge (in the Gromov--Hausdorff topology on closed subsets) to $\ti\ell = \wt R_1 \union \wt R_2$. The leaf $\ti\ell$ is called \emph{the leaf at $\ti v$ determined by $d_1$ and $d_2$} (or {by $\wt E_1$ and $\wt E_2$}). If $\ti v$ is principal then we say that $\ti \ell$ is \emph{a singular leaf}. 

In the previous paragraph, the converse to the claim holds as well: if $\ti\ell = \wt R_1 \union \wt R_2$ is a leaf of $\wt\Lambda_-$ realized in $\wt G$ then $d_1,d_2$ bound an edge of $SW(\ti v;\wt G)$. This holds because, by definition of $\wt\Lambda_-$, each segment of $\ti\ell$ must be contained in $\breve g^l(\wt E)$ for some edge $\wt E$ and some $l>0$. It follows that $SW(\ti v;\wt G)$ can be identified with the graph having one vertex for each singular ray $\wt R$ based at $\ti v$ and having an edge connecting the vertices corresponding to a pair of singular rays $\wt R_1$, $\wt R_2$ at $\ti v$ if and only if $\wt R_1 \union \wt R_2$ is a singular leaf at $\ti v$.

\paragraph{Characterizing $\Lambda_-$.} The next lemma gives a useful characterization of $\Lambda_-$, and its corollary characterizes the nonsingular $\phi$-periodic leaves of $\Lambda_-$.

\begin{lemma}
\label{first characterization} 
$\Lambda_-$ is the unique minimal $\phi$-invariant lamination whose leaves are legal with respect to every train track representative $\gG$ of $\phi$.  Every leaf $\ti \ell$ constructed as above, and in particular every singular leaf, is contained in  $\Lambda_-$.
\end{lemma} 

\begin{proof}   With notation as above, consider the set $\cal L$ of  minimal $\phi$-invariant laminations whose leaves are legal with respect to every train track representative $\gG$ of $\phi$.  That  $\Lambda_- \in \cal L$ follows from the construction of $\Lambda_-$ in section~1 of \cite{BFH:laminations}, as reviewed above.   We will show that if  $\Lambda' \in \cal L$  and $\ti \ell$ is any leaf constructed as above, then $\ti \ell$ is a leaf of $\Lambda'$.  Since there can only be one minimal lamination that contains $\ti \ell$ it follows $\Lambda_-$ is the only element of $\cal L$.  

Pass to a power and choose an appropriate lift $\breve g \from \wt G \to \wt G$ as above, and choose $\ti v \in \wt G$ and distinct edges $\wt E_1$ and $\wt E_2$ at $\ti v$ all of which are fixed by $\breve g$, so that $\wt E_1, \wt E_2$ determine an edge in $SW(\ti v;\wt G)$.  There is an edge $\wt E$ and $l > 0$ such that $\breve g^l(\wt E)$ contains the turn $(\wt E_1, \wt E_2)$.  Since $g$ is irreducible and leaves of $\Lambda'$ are legal in $G$, there is a leaf $\ti \ell'$ of $\wt \Lambda'$ that contains $\wt E$. Each path in the nested increasing sequence  $\wt E_1^{-1}\wt E_2 \subset \breve g(\wt E_1^{-1} \wt E_2) \subset \ti  g^{2}(\wt E_1^{-1} \ti  E_2) \subset\cdots $  occurs as a subpath of a leaf of $\wt \Lambda'$. It follows that their union $\ti \ell$  is  a leaf of $\wt \Lambda'$.   
\end{proof}  

\begin{corollary} \label{periodic leaves} If $\Phi \in \Aut(F_n)$ represents an iterate of $\phi$ and $\Per_N(\hat \Phi) = \Fix_N(\hat \Phi) = \{P_1,P_2\}$ then the line $\ti \ell$ bounded by $P_1$ and $P_2$ is a periodic nonsingular leaf of $\wt \Lambda$. Conversely, if $\ti\ell$ is a nonsingular periodic leaf of $\wt\Lambda$ then there exists $\Phi\in\Aut(F_n)$ representing an iterate of $\phi$ such that $\Per_N(\hat\Phi)=\Fix_N(\hat\Phi)=\bdy\ti\ell$.
\end{corollary}

\begin{proof} Choose an irreducible train track map  $\gG$ representing $\phi$ and let $\breve g_k$ be the lift of an iterate of $g$ that corresponds to $\Phi$.  Since $P_1$ and $P_2$ are attractors, $\Fix(\breve g_k) \ne \emptyset$; see for example Lemma 3.9 and Lemma~3.13 of \cite{FeighnHandel:abelian}.  After replacing $\breve g_k$ by an iterate we may choose $\ti x \in \Fix(\breve g_k)$ and a pair of fixed directions $d_1$ and $d_2$ at $\ti x$.  After subdividing $G$, we may assume that $\ti x$ is a vertex. If there were any other periodic direction at $\ti x$ then the corresponding singular ray would give another point in $\Per_N(\hat\Phi)$ in addition to $P_1,P_2$, a contradiction. It follows that the stable Whitehead graph $SW(\ti x;\wt G)$ has just the two vertices $d_1,d_2$, and moreover the rays $R_1,R_2$ determined by $d_1,d_2$ must end at points of $\Per_N(\hat\Phi)$, and up to re-indexing these points must be $P_1,P_2$. Since $SW(\ti x;\wt G)$ is connected (see Section~\ref{SectionTrainTracks}), there must be an edge of $SW(\ti x;\wt G)$ connecting $d_1$ to $d_2$. An application of Lemma~\ref{first characterization} shows that $\ti \ell = R_1 \union R_2$ is in $\wt\Lambda$. 

For the converse, let $\ti\ell$ be a nonsingular periodic leaf of $\wt\Lambda$. Choose $\Phi \in \Aut(F_r)$ representing an iterate of $\phi$ so that $\Phi\#(\ti\ell)=\ti\ell$. Choose an irreducible train track map $g \from G \to G$ representing $\phi$ and let $\breve g_k \from \wt G \to \wt G$ be the lift of an iterate of $g$ that corresponds to $\Phi$. Using $\ti\ell$ to denote its own realization in $\wt G$, we have $\breve g_k(\ti\ell)=\ti\ell$. Since $\breve g_k$ expands $\ti\ell$ by the uniform factor $\lambda$, there is a point $x \in \ti\ell$ fixed by $\breve g_k$. Subdivide $G$ so that $x$ is a vertex. Up to passing to the square, $D\breve g_k$ fixes the two directions $d_1,d_2$ of $\ti\ell$ at $x$, and the rays determined by these two directions are the two rays of $\ti\ell$ based at $x$. It follows that $\bdy\ti\ell \subset \Fix_N(\hat\Phi) \subset \Per_N(\hat\Phi)$. If there were another point in $\Per_N(\hat\Phi)$ then, passing to a power so that this point is in $\Fix_N(\hat\Phi)$, it would follow that $\ti\ell$ is singular, contrary to hypothesis.
\end{proof}

\subsection{Relating $\Lambda_-$ to $T_-$ and to $\TPlus$}
\label{SectionRelativeLeavesToTrees}
Consider an $\reals$-tree $T$ representing a point in $\bdy\X_r$, and a nonempty minimal sublamination $\Lambda$ of $\G F_r$. Following \BookZero\ section 3 we say that $\Lambda$ \emph{has length zero} in $T$ if for every $G \in \X_r$ and every $F_r$-equivariant morphism $h \from \wt G \to T$ there exists $C \ge 0$ such that for every leaf $\ell \subset \wt G$ of $\wt\Lambda_G$ we have $\diam_{T_-}(h(\ell)) \le C$. 

\begin{lemma}
\label{LemmaMinimalLamination} For each fully irreducible $\phi \in \Out(F_r)$,
$\Lambda_-=\Lambda_-(\phi)$ is the unique minimal lamination that has length zero in $T_-=T_-(\phi)$.
\end{lemma}

\begin{proof}  This is proved in Lemma 3.5~(4) of \BookZero, but restricted to the so-called ``irreducible laminations'', meaning the set of expanding laminations of all fully irreducible outer automorphisms in $F_r$. However, the proof applies to arbitrary minimal sublaminations, for as is noted several times in \BookZero, the only property of irreducibility of $\Lambda$ that is really used is quasiperiodicity of leaves of $\Lambda$, which is equivalent to minimality of $\Lambda$. See for example the proof of \BookZero\ Proposition~2.12 which is used in the proof of Lemma 3.5~(4). 
\end{proof}

\paragraph{Remark} Lemma~\ref{LemmaMinimalLamination} shows that the expanding lamination of a fully irreducible outer automorphism $\phi$ is naturally related to the repelling tree $T_-$ rather than the attracting tree $\TPlus$.   This is why we have chosen to denote the expanding lamination with a subscript ``$-$'', altering the notation from the subscript ``$+$'' used in \BookZero. 

\bigskip

We next relate $\Lambda_-$ to $\TPlus$.

\begin{lemma}  \label{CorollaryLeavesInTPlus}
For any train track representative $g \from G \to G$ of $\phi$ and any leaf $\ti\ell$ of $\wt\Lambda_-$, the restriction of the map $f_g \from \wt G \to \TPlus$ to the realization of $\ti\ell$ in $\wt G$ is an isometry. Moreover, the bi-infinite line in $\TPlus$ which is the image of this isometry is independent of the choice of $g$, depending only on $\ti\ell$; this image is called the \emph{realization of $\ti\ell$ in $\TPlus$}.
\qed \end{lemma}

Motivated by the first part of this lemma, given a free, simplicial $F_r$-tree $T$, an $F_r$-equivariant map $f \from T \to \TPlus$ is said to be a \emph{$\Lambda_-$ isometry}, or to be \emph{$\Lambda_-$ isometric}, if the restriction of $f$ to any leaf of $\wt\Lambda_-$ realized in $T$ is an isometry. The first part of the lemma therefore says that the map $f_g$ is a $\Lambda_-$ isometry.

\begin{remark} \label{f is unique} The property that $T$ has a $\Lambda_-$ isometry $f \from T \to \TPlus$ does not, in and of itself, determine the map $f$.  In Section~\ref{SectionIrrigidity} we give an example of a marked graph $G$ such that there exist at least two distinct $\Lambda_-$ isometries $\wt G \to \TPlus$ --- in fact, that example will have a one-parameter family of distinct such maps. On the other hand, the arguments in Section~\ref{SectionIrrigidity} will show that a $\Lambda_-$ isometry is determined by where a single point of $\wt G$ goes in $\TPlus$.
\end{remark}

\begin{proof}[Proof of Lemma \ref{CorollaryLeavesInTPlus}] \quad The first part follows by combining Corollary~\ref{first tree corollary} and Lemma~\ref{first characterization}. 

To prove uniqueness of the realization in $\TPlus$ of a leaf $\ti\ell$ of $\wt\Lambda_-$, we must prove that $f_g(\ti\ell)$ is well-defined independent of $g$. Think of $\ti\ell$ abstractly as an element of $\G F_r$, and choose a sequence $\gamma_i \in F_r$ whose abstract axes $A(\gamma_i)$ in $\G F_r$ limit to $\ti\ell$. For any train track representative $g \from G \to G$ of $\phi$, the realizations $A(\gamma_i)_{\wt G}$ converge to the realization $\ti\ell_G$ in the Gromov--Hausdorff topology on closed subsets of~$\wt G$. It follows that the sets $p_i = \ti\ell_G \intersect A(\gamma_i)_{\wt G}$ form a longer and longer family of paths that exhausts $\ti\ell_G$. While $A(\gamma_i)_{\wt G}$ need not be legal in $\wt G$ and so may not embed in $\TPlus$ under $f_g$, its subpath $p_i$ is legal and does embed in $\TPlus$. By bounded cancellation, Lemma~\ref{bcl}, there exists $B$ independent of $i$ and depending only on $f_g$, such that if $p'_i$ is obtained from $p_i$ by cutting off segments of length $B$ at both ends, then $f_g(p'_i)$ embeds in $A(\gamma_i)_{\TPlus}$, the axis of $\gamma_i$ in $\TPlus$. But the sequence of paths $f_g(p'_i)$ form a longer and longer family of paths in $\TPlus$ that exhausts $f_g(\ti\ell_G)$. It follows that $A(\gamma_i)_{\TPlus}$ approaches $f_g(\ti\ell_G)$ in the Gromov-Hausdorff topology on closed subsets of $\TPlus$. Since the sequence $\gamma_i$ is independent of the choice of $g \from G \to G$, uniqueness of $f_g(\ti\ell)$ in $\TPlus$ follows.
\end{proof}

We conclude with a relation between $\Lambda_-$ and the homothety $\Phi_+ \from \TPlus \to \TPlus$.

\begin{corollary}\label{action on T+} For each $\Phi \in \Aut(F_n)$ representing $\phi$, and each leaf $\ti\ell$ of $\wt\Lambda_-$, the map $\Phi_+$ takes the realization of $\ti\ell$ in $\TPlus$ to the realization of $\Phi_\#(\ti\ell)$ in $\TPlus$.
\end{corollary}

\begin{proof}  Let $g \from G \to G$ be a train track representative of $\phi$ and $\breve g \from \wt G \to \wt G$ the lift corresponding to $\Phi$. Let $\ti\ell$ denote its own realization in $\wt G$, and so $f_g(\ti\ell)$ is its realization in $\TPlus$. Also, the realization of $\Phi_\#(\ti\ell)$ in $\wt G$ is $\breve g(\ti\ell)$, and its realization in $\TPlus$ is $f_g(\breve g(\ti\ell))$. But we know that 
$$
f_g \composed \breve g = \Phi_+ \composed f_g
$$
and so $\Phi_+$ takes the realization of $\ti\ell$ in $\TPlus$ to the realization of $\Phi_\#(\ti\ell)$ in $\TPlus$.

\end{proof}

\section{The ideal Whitehead graph} 
\label{SectionIdealWhiteheadGraph}

On a closed oriented surface $S$, the topological structure of the unstable foliation of a pseudo-Anosov homeomorphism $\phi \from S \to S$ near a singular point $s \in S$ is fully determined by the number $n \ge 3$ of unstable prongs at $s$. Associated to $s$ there is a principal region (that is, a complementary component) of the unstable geodesic lamination of $\phi$, this region is an ideal $n$-gon, and its $n$ sides are the singular leaves associated to $s$, which are circularly ordered by the asymptote relation. Also, in the $\reals$-tree dual to the unstable foliation, associated to each $n$-pronged singularity~$s$ there is an orbit of branch points of the $\reals$-tree, and the number of directions at each such branch point is exactly $n$. The singularity structure of $\phi$ is therefore described fully by a list of integers $n_i$ giving the number of unstable prongs at each singularity. The singularity structure is constrained by the Euler-\Poincare\ index equation $\sum (n_i - 2) = -2 \chi(S)$, and the paper \cite{MasurSmillie:singularities} proves that outside of a few fully understood exceptions that occur only for a few surfaces of low topological complexity, the index equation gives a complete description of all singularity structures that actually occur for pseudo-Anosov homeomorphisms of $S$.

Singularity structures of nongeometric, fully irreducible outer automorphisms $\phi \in F_r$ are less well understood. The most that is known so far is an index inequality of Euler-\Poincare\ type that is proved in \cite{GJLL:index}: there are finitely many $F_r$ orbits of branch points in $\TphiPlus$, each branch point has finitely many directions, and if $n_i$ is the number of directions at each branch point in the $i^{\text{th}}$ orbit then $\sum (n_i - 2) \le 2 (r-1)$. Whether this gives a complete description of all singularity structures that actually occur, along the lines of \cite{MasurSmillie:singularities}, has not been investigated. In the index inequality of \cite{GJLL:index}, the significant information about a singularity is, like in the surface case, a raw count of its number of ``unstable directions'', encoded as directions in $\TphiPlus$. 

In this section we study a finer invariant of the singularity structure of a nongeometric, fully irreducible $\phi \in \Out(F_r)$, the ``ideal Whitehead graph'' of $\phi$, which captures the asymptotic relations among singular leaves of $\Lambda_-(\phi)$. The ideal Whitehead graph is first defined in terms of principal automorphisms representing $\phi$ or, equivalently, in terms of principal lifts of a train track representative. The main results of this section give alternative descriptions of the ideal Whitehead graph: in terms of singular leaves of train track representatives (Lemma~\ref{singular leaves} and Corollary~\ref{ideal is union of local}); next in terms of the asymptote relation among singular leaves (Lemma~\ref{ideal whitehead});  and finally in terms of branch points in the $\reals$-tree $\TphiPlus$ (Lemma~\ref{TPlusWhitehead}). 

To describe the ideal Whitehead graph very roughly, to each ``singularity of $\phi$'' there is associated a component of the ideal Whitehead graph, a finite graph with singular leaves as edges, and with vertices describing asymptotic relations between singular leaves. By analogy, for a pseudo-Anosov surface singularity with $n$ unstable prongs, the corresponding component of the ideal Whitehead graph has $n$ vertices and $n$ edges forming a topological circle. For a general nongeometric, fully irreducible $\phi \in \Out(F_r)$, components of the ideal Whitehead graph need not be topological circles; indeed, they may have cut points, as we show by example in Section~\ref{SectionIWGExample}.

For present purposes we restrict to the case that $\phi$ is rotationless, but the reader will easily see how to proceed without this restriction. Also, we shall restrict attention to an $F_r$-equivariant version of the ideal Whitehead graph; each component having trivial stabilizer, the reader will easily see how to pass to a quotient under $F_r$ to obtain a nonequivariant version.

\bigskip

{\em For the rest of this section we fix a nongeometric, fully irreducible, rotationless $\phi \in \Out(F_r)$, and we denote $\Lambda_- = \Lambda_-(\phi)$ and $T_+ = \TphiPlus$. All automorphisms and all train track maps will be representatives of $\phi$}.

\subsection{Definition and structure of the ideal Whitehead graph}
\label{SectionIWGDefinition}

Recall that, abstractly, leaves of $\Lambda_-$ are just points $(P_1,P_2) \in D^2 \bdy F_r$.  

For every principal automorphism $\Phi$ let $L(\Phi)$ be the set of leaves $(P_1,P_2)$ of $\Lambda_-$ with $P_1,P_2 \in \Fix_N(\hat\Phi)$. Define the \emph{component $W(\Phi)$ of the ideal Whitehead graph determined by $\Phi$} to be the graph with one vertex for each point in $\Fix_N(\hat\Phi)$, and an edge connecting $P_1, P_2 \in \Fix_N(\hat \Phi)$ if there is a leaf $(P_1,P_2) \in L(\Phi)$. Corollary~\ref{distinct fixed sets} implies that for each point $P \in \partial F_r$ there is at most one principle automorphism $\Phi$ such that $P \in \Fix_N(\hat\Phi)$, and so $L(\Phi_1) \intersect L(\Phi_2) = \emptyset$ for distinct principle automorphisms $\Phi_1,\Phi_2$. It therefore makes sense to define the \emph{ideal Whitehead graph} $W(\phi)$ to be the disjoint union of its components $W(\Phi)$, one for each principal automorphism~$\Phi$ representing $\phi$.

The definition of $W(\phi)$ can be easily reformulated in terms of any train track representative $g \from G \to G$ of $\phi$. For every principal lift $\breve g \from \wt G \to \wt G$, let $L(\breve g)$ denote the set of leaves of $\ti \Lambda_-$ with both ideal endpoints in $\Fix_N(\hat g)$, define $W(\breve g)$ to be the graph with one vertex for each point in $\Fix_N(\hat g)$ and an edge connecting the vertices corresponding to $P_1, P_2 \in \Fix_N(\hat g)$ if there is a leaf $\ti \ell \in L(\breve g)$ with ideal endpoints $P_1$ and $P_2$, and define $W(g)$ to be the disjoint union of the $W(\breve g)$. Note that if $\Phi$ is the principal automorphism corresponding to $\breve g$ then there are natural isomorphisms $L(\breve g) \approx L(\Phi)$ and $W(\breve g) \approx W(\Phi)$, and there is a natural isomorphism $W(g) \approx W(\phi)$. Formally, it is convenient to let the vertex set of $W(\breve g)$ be $\Fix_N(\hat g)$ and to let the edge interiors of $W(\breve g)$ be the disjoint union of the leaves of $L(\breve g)$. 

For any principal lift $\breve g \from \wt G \to \wt G$ and $\ti v \in \Fix(\breve g)$, we identify the local stable Whitehead graph $SW_{\ti v}$ with a subgraph of $W(\breve g)$ as follows. Each vertex $d_i$ of $SW_{\ti v}$ is a direction at $\ti v$ fixed by $D\breve g$, $d_i$ is the initial direction of a unique singular ray $R_i$, and $d_i$ is identified with the vertex of $W(\breve g)$ which is the ideal endpoint $P_i$ of the ray~$R_i$. As shown in Section~\ref{SectionExpandingLamination}, two vertices $d_i,d_j$ of $SW_{\ti v}$ that are connected by an edge of $SW_{\ti v}$ correspond to a singular leaf $\ell = R_1 \union R_2$ with ideal endpoints $P_1,P_2$, and so $\ell$ corresponds to an edge of $W(\breve g)$ with endpoints $P_1,P_2$. 

\begin{lemma} \label{singular leaves} Let $\breve g$ be a principal lift of $g$.
\begin{enumerate}
\item Every $P \in \Fix_N(\hat g)$ is the endpoint of at least one element of $L(\breve g)$.  
\item Singular rays $R_1$ at $\ti v_1 \in \Fix(\breve g)$ and $R_2$ at $\ti v_2 \in \Fix(\breve g)$ terminate at the same point $P \in \Fix(\hat g)$ if and only if  the path $\ti \sigma$ connecting $\ti v_1$ to $\ti v_2$ projects to an \iNp\  $\sigma$ and determines the same direction at $\ti v_i$ as does $R_i$.
\item $L(\breve g)$ is the disjoint union, over all $\ti v \in \Fix(\breve g)$,  of the singular leaves at $\ti v$.
\end{enumerate}
\end{lemma}

\begin{proof}  Item (1) follows from Lemma~4.7 of \cite{FeighnHandel:abelian}. 

The path $\sigma$ of (2) is a Nielsen path. If $R_1$ and $R_2$ eventually coincide then $\ti \sigma = \ti \alpha \ti \beta^{-1}$ where $\ti \alpha \subset R_1$ and $\ti \beta \subset R_2$. This proves the only if part of (2). For the if part, suppose that $\sigma =\alpha \bar \beta$ is indivisible. Let $\ti \alpha$ and $\ti \beta$ be the lifts of $\alpha$ and $\beta$ that begin at $\ti v_1$ and $\ti v_2$ respectively. For all $n \ge 1$, there are paths $\ti \tau_n$ such that $\breve g^n(\ti \alpha) = \ti \alpha \ti \tau_n$ and $\breve g^n(\ti \beta) = \ti \beta \ti \tau_n$. The paths $\breve g^n(\ti \alpha)$ are an increasing sequence whose union is $R_1$ and the paths $\breve g^n(\ti \beta)$ are an increasing sequence whose union is $R_2$. This completes the proof of (2).

Each singular leaf based at an element of $\Fix(\breve g)$ is contained in $L(\breve g)$ by construction. Conversely, since the endpoints of any $\ti \ell \in L(\breve g)$ are attractors for $\hat g$, there is a necessarily unique fixed point $\ti x \in \ti l$; see for example Lemma~3.9 and Lemma~3.10 of \cite{FeighnHandel:abelian}. Lemma~\ref{principal essential} implies that $\ti x$ is  a principal vertex. The turn taken by $\ti \ell$ at $\ti x$ determines an edge of $SW_{\ti x}$ and hence a singular leaf $\ti \ell'$ at $\ti x$.  If $\ti \ell \ne \ti \ell'$ then there are rays $R$ and $R'$ based at $\ti x$ with fixed point free interiors, with a common initial edge and with distinct terminal endpoints $P$ and $P'$ in $\Fix_N(\hat g)$. The last point of $R \intersect R'$ is therefore mapped by $\breve g$ both to $R-R'$ and to $R'-R$, a contradiction. This proves that $\ell$ is a singular leaf at $\ti x$ and completes the proof of (3).  
\end{proof}

The following corollary says that each component of the ideal Whitehead graph is a union of local Whitehead graphs, pasted together along cut points:

\begin{corollary} \label{ideal is union of local}
Let $\breve g$ be a principal lift of $g$.
\begin{enumerate}
\item $W(\breve g)$ is connected.
\item $W(\breve g) = \cup SW(\ti v)$ where the union is taken over all $\ti v \in \Fix(\breve g)$.
\item For $i \ne j$, $SW(\ti v_i)$ and $SW(\ti v_j)$ intersect in at most one vertex. If they do intersect at a vertex $P$, then $P$ is a cut point of $W(\breve g)$, in fact $P$ separates $SW(\ti v_i)$ from $SW(\ti v_j)$ in $W(\breve g)$.
\end{enumerate}
\end{corollary}

\begin{proof} (2) follows from the first and third item of Lemma~\ref{singular leaves}. (1) follows from (2), item~2 of Lemma~\ref{singular leaves}, and the fact that any two elements of $\Fix(\breve g)$ are connected by a finite concatenation of \iNp s.  

To simplify notation, for proving (3) we assume $i=1$, $j=2$. Suppose that $P \in \Fix_N(\hat g)$ is a common vertex of  $SW(\ti v_1) \cap SW(\ti v_2)$.  There are singular rays $R_1$ and $R_2$ terminating at $P$ and initiating at $\ti v_1$ and $\ti v_2$ respectively. Since no leaf of $\ti \Lambda_-$ contains both  $\ti v_1$ and $\ti v_2$, it follows that $R_1 \cap R_2 \cap \{\ti v_1,\ti v_2\} = \emptyset$. If $Q \ne P$ is another common vertex of $SW(\ti v_1) \cap SW(\ti v_2)$ then there are singular rays $R'_1$ and $R'_2$ terminating at $Q$ and initiating at $\ti v_1$ and $\ti v_2$ respectively.  The initial direction of $R_1$  is distinct from that of $R_1'$ and similarly for $R_2$ and $R_2'$. It follows that $R_1 \cup R_2 \cup R_1' \cup R_2'$ contains an embedded circle, contradicting that $\wt G$ is a tree. This proves $SW(\ti v_1) \cap SW(\ti v_2) =\{P\}$.    

To prepare for proving the last part of (3), let $W^d(\breve g)$ denote $W(\breve g)$ with vertices deleted, namely, the disjoint union of all singular leaves based at all fixed points of $\breve g$; formally, one can think of $W^d(\breve g)$ as the set of ordered pairs $(l,x)$ where $l$ is a singular leaf based at a fixed point of $\breve g$ and $x \in l$. Let $r \from W^d(\breve g) \to \wt G$ be the natural map, taking each singular leaf to itself by the identity map; formally, $r(l,x)=x$. The set $r^\inv(\Fix(\breve g))$ consists of exactly one point in each edge of $W(\breve g)$, called the \emph{basepoint} of that edge, separating the edge into two \emph{open half-edges}. Suppose that $P$ does not separate $SW(\ti v_1)$ from $SW(\ti v_2)$ in $W(\breve g)$. Then there are edge basepoints $w_1,w_2 \in W^d(\breve g)$ with $r(w_i)=v_i$, and there is an embedded path $\gamma$ in $W(\breve g)$ from $w_1$ to $w_2$ disjoint from $P$. Let $\gamma^d = \gamma \intersect W^d(\breve g)$, that is, $\gamma$ with vertices of $W(\breve g)$ deleted. The image $r(\gamma^d)$ is a connected subset of $\wt G$, because near any deleted vertex $Q$ of $\gamma$ the two ends of $\gamma$ incident to $Q$ are mapped by $r$ to rays in $\wt G$ asymptotic to $Q$. Also, since $\gamma$ does not pass through $P$, and since $\gamma$ is embedded, it follows that $\gamma$ does not intersect any open half-edge incident to~$P$, from which it follows that the image $r(\gamma^d)$ is disjoint from $\interior(R_1) \union \interior(R_2)$. Therefore, $r(\gamma^d)$ contains an embedded path in $\wt G - \left( \interior(R_1) \union \interior(R_2) \right)$ from $\ti v_1$ to $\ti v_2$. It follows that $r(\gamma^d) \union R_1 \union R_2$ contains an embedded circle, contradicting that $\wt G$ is a tree.
\end{proof}

\subsection{Asymptotic leaves and the ideal Whitehead graph}
\label{SectionAsymptoticLeaves}

A pair of distinct leaves of $\ti \Lambda_-$  are {\em asymptotic} if they have a common infinite end.  This generates an equivalence relation on leaves of $\ti \Lambda_-$ whose non-trivial classes are called {\em asymptotic equivalence classes of $\ti \Lambda_-$}. 

Our next lemma  relates a component of the ideal Whitehead graph to asymptotic equivalence classes.

\begin{lemma} \label{ideal whitehead} The map $\breve g \mapsto L(\breve g)$ defines a bijection between the set of principal lifts of $g$ and the set of asymptotic equivalence classes of $\ti \Lambda_-$.  Two ends of edges of $W(\breve g)$ are incident to the same vertex if and only if the corresponding rays in leaves of $L(\breve g)$ are asymptotic.
\end{lemma}

\begin{proof} The second sentence is obvious.

Injectivity of the map $\breve g \mapsto L(\breve g)$ follows from the fact that $L(\breve g_1) \cap L(\breve g_2) = \emptyset$ for distinct $\breve g_1$ and $\breve g_2$.      Singular leaves at $\ti v \in \Fix(\breve g)$ belong to the same asymptotic equivalence class because $SW_{\ti v}$ is connected.  Since any two elements of $\Fix(\breve g)$ are connected by a finite concatenation of \iNp s, the second item of Lemma~\ref{singular leaves} implies that  all elements of $L(\breve g)$ are in the same asymptotic equivalence class.

       It remains to show that if $\ti \ell'$ and $\ti \ell''$  are asymptotic leaves in $\wt G$ then they are elements of    some $L(\breve g)$. It is convenient to work in $G$ rather than $\wt G$ so consider asymptotic leaves  $\ell',\ell'' \in\Lambda_-$.  For each $j \ge 0$, there are leaves  $\ell'_{j}$ and $\ell''_{j}$ such that $g^j(\ell'_{j}) =\ell'$ and  $g^j(\ell''_{j}) = \ell''$. The intersection of $\ell_j'$ and $\ell_j''$ is a ray $R_j$ based at a vertex $x_j$.  Denote the rays in $\ell_j'$ and $\ell_j''$ that are complementary to $R_j$ by $R_j'$ and $R_j''$ respectively.    Let $E_j,E_j'$ and $E_j''$ be the initial edges of  $R_j, R_j'$ and $R_j''$ respectively and let $d_j,d_j'$ and $d_j''$ be the directions that they determine.

As a first case suppose that $(d_j', d_j'')$ is legal for all $j\ge 0$.  Then $x_j$ is principal and $g^{i}(x_{j+i}) = x_j$ for all $i,j\ge 0$.   It follows that $x_j \in \Fix(g)$ is independent of $j$ and we  now denote $x_j$ by $v$.  The analogous argument shows that  $E_j,E_j'$ and $E_j''$ are independent of $j$ and we denote these by $E,E'$ and $E''$.       For all $j$, the paths $g^j(E)$ and $g^j( E')$ are respectively initial subpaths of $R_0$ and $R_0'$.    It follows that $\ell'$ is the singular leaf at $v$ determined by $E$ and $E'$.  The symmetric argument shows that $\ell''$ is the the singular leaf at $v$ determined by $E$ and $E''$.  The point $\ti v$ at which $\ti \ell'$ and $\ti \ell''$ diverge projects to $v$.  If  $\breve g$ is the lift that fixes $\ti v$, then $\ell'$ and $\ell''$ are singular leaves at $\ti v$ for $\breve g$.
This completes the proof in this first case.

Assume now that $(d_j', d_j'')$ is illegal for some, and hence every sufficiently large, $j$. Item 5 of Lemma~\ref{nielsen paths} and Lemma~\ref{period one} imply that there is an indivisible Nielsen path $\alpha_j \bar \beta_j$ such that $\bar \alpha_j$ is an initial path of $R'_j$, $\bar \beta_j$ is an initial path of $R''_j$, and $g^i_\#(\alpha_{i+j} \bar \beta_{i+j}) =  \alpha_j \bar \beta_j$ for all $i,j\ge 0$.  It follows that $\alpha_j,\beta_j,x_j$ and the initial endpoints $y_j$ and $z_j$ of $\alpha_j$ and $\beta_j$ respectively are independent of $j$; denote these by $\alpha,\beta,v,y$ and $z$.  Arguing as in the previous case we see that $\ell'$ and $\ell''$ are singular leaves at $y$ and $z$ respectively. The point $\ti v$ at which $\ti \ell'$ and $\ti \ell''$ diverge projects to the common terminal endpoint $v$ of $\alpha$ and $\beta$.  Let $\ti \alpha$ and $\ti \beta$ be the lifts that terminate at $\ti v$ and let $\breve g$ be the lift that fixes the initial endpoints $\ti y$ and $\ti z$ of $\ti \alpha$ and $\ti \beta$.  Then $\ell'$ and $\ell''$ are singular leaves for $\breve g$ at $\ti y$ and $\ti z$ respectively.  
\end{proof}

\subsection{$T_+$ and the ideal Whitehead graph}
\label{SectionIWGandTPlus}

Next we relate the ideal Whitehead graph to $T_+$. Recall Corollary~\ref{CorollaryLeavesInTPlus} which explains how the leaves of $\ti \Lambda_-$ are realized in $T_+$.   

Suppose that $b \in T_+$ is a branch point, $\Phi \in \PA(\phi)$ is the principal automorphism corresponding to $b$ by Lemma~\ref{principal lifts and branch points}, and $\breve g \from \wt G \to \wt G$ is the corresponding principal lift. Each singular leaf $\ell$ in $L(\breve g)$ passes through some $x \in \Fix(\breve g)$, and so $f_g(\ell)$, the realization of $\ell$ in $T_+$, passes through the point $f_g(x)=b$; we call these the \emph{singular leaves at $b$}. Each singular leaf at $b$ is divided by $b$ into {\em singular rays at $b$} which are the images under $f_{g}$ of the singular rays at some $x \in \Fix(\breve g)$.   

Define {\em the local Whitehead graph} at a point $z \in T_+$, denoted $W(z;T_+)$, to be the graph with one vertex for each direction $d_i$ at $z$ and with an edge connecting $d_i$ to $d_j$ if the turn $(d_i,d_j)$ is taken by the realization of some leaf of $\Lambda_-$ in $T_+$. For each $x \in \wt G$ the map $D_x f_g$ maps the directions of $\wt G$ at $x$ to the directions of $T_+$ at $f(b)$, and this map extends to a simplicial map $W(x;\wt G) \mapsto W(b;T_+)$, still denoted $D_x f_g$ (though we often drop $x$ from the notation $D_x f_g$.

With $b$, $\Phi$, $\breve g$ as above, for each $x \in \Fix(\breve g)$ the map $D_{x} f_g$ is injective on the set of fixed directions at $x$, and therefore $D_x f_g$ restricts to a simplicial embedding $SW(x;\wt G) \hookrightarrow W(b;T_+)$.

\begin{lemma} 
\label{TPlusWhitehead}
With $b$, $\Phi$, and $\breve g$ as above, as $x$ varies over $\Fix(\breve g)$ the embeddings 
$$D_x f_g \from SW(x;\wt G) \hookrightarrow W(b;T_+)
$$ 
induce an isomorphism of $W(\Phi)$ with $W(b;T_+)$.  More precisely,  
\begin{enumerate}
\item Every direction at $b$ is realized by a singular ray. 
\item Singular rays $R_1$ and $R_2$ in $\wt G$ with basepoints in $\Fix(\breve g)$ map to the same singular ray in $T_+$ if and only if they terminate at the same point in $\Fix_N(\hat g)$.     
\item A turn at $b$ is taken by some leaf  of $\ti \Lambda_-$ if and only if it is taken by a singular leaf at $b$. 
\end{enumerate}
 \end{lemma}

\begin{proof} Observe first that, by Corollary~\ref{ideal is union of local} and (2), the maps $D_x f_g \from SW(x;\wt G) \to W(b;T_+)$ for $x \in \Fix(\breve g)$ fit together to define a simplicial map $W(\Phi) \mapsto W(b;T_+)$, which by items (1)--(3) is an isomorphism. Item~(1) follows from the second item of Lemma~\ref{principal lifts and branch points} and the fact that every fixed direction at an element of $\Fix(\breve g)$ is realized by a singular ray.

The ``if'' part of (2) follows from the second item of Lemma~\ref{singular leaves} and Corollary~\ref{first tree corollary}. For the ``only if'' part, let $R_1$ and $R_2$ be singular rays at $\ti v_1$ and $\ti v_2$ such that  $f_g(R_1) = f_g(R_2)$. We prove that $R_1$ and $R_2$ converge to the same point in $\Fix(\hat g)$. Let $d_i$ be the direction determined by $R_i$ at $\ti v_i$. If $\ti v_1 = \ti v_2$, then $R_1 = R_2$ and the result is obvious.  Suppose then that $\ti v_1 \ne \ti v_2$.  The path connecting $\ti v_1$ to $\ti v_2$ splits as a concatenation $\ti \mu = \ti \mu_1 \cdot \ldots \cdot \ti \mu_l$ where each $\mu_i$ is an indivisible Nielsen path.  For $i=1,2$ choose $\ti x_i$ close to and in the direction $d_i$ from $\ti v_i$ such that $f_g(\ti x_1) = f_g(\ti x_2)$.  The path $\ti \mu'$ connecting $\ti x_1$ to $\ti x_2$ is obtained from $\ti \mu$ by making alterations at the beginning and end.  If $d_1$ is the direction determined by $\ti \mu_1$, then the path from $\ti v_1$ to $\ti x_1$ is removed from $\ti \mu_1$; otherwise the path from $\ti x_1$ to $\ti v_1$ is added. Similar changes are made at the terminal end.  Since there are only finitely many indivisible Nielsen or pre-Nielsen paths in $G$ we may assume that $\ti x_1$ and $\ti x_2$ have been choosen so that $\ti \mu'$ does not project to a Nielsen or pre-Nielsen path. Corollary~\ref{first tree corollary} implies that $\breve g^k(\ti x_1) = \breve g^k(\ti x_2)$ for some $k > 0$.  The second item of Lemma~\ref{nielsen paths} therefore implies that $\ti x_1$ and $\ti x_2$ belong to the interior of the same 
$\ti \mu_j$.  It follows that $\mu = \mu_1$ is indivisible and that $d_1$ and $d_2$ point into $\ti \mu$. The second item of Lemma~\ref{singular leaves} implies that $R_1$ and $R_2$ terminate at the same endpoint as desired.
 
The ``if'' part of (3) is obvious. Suppose the turn $(d_1^+, d^+_2)$ is taken by a leaf $\ti \ell$ of $\ti \Lambda_-$.  The realization of $\ti \ell$ in $\wt G$ contains some $\ti x \in f_g^{-1}(b)$.  By Lemma~\ref{principal lifts and branch points}(3), for some $m>0$ we have $\breve g^m(\ti x) \in \Fix(\breve g)$ and the two directions $d_1,d_2$ of $\breve g^m(\ti \ell)$ at $\breve g^m(\ti x)$ are fixed by $D\breve g$, so $f_g(g^m(\ti x)) = b$ and $Df_g(d_i)=d_i^+$. Also, since the turn $(d_1,d_2)$ is taken by some leaf, namely $\breve g^m(\ti\ell)$, it follows that this turn is taken by some singular leaf $\ti\ell'$ based at $\breve g^m(\ti x)$, and so $f_g(\ti\ell')$, the realization of $\ti\ell'$ in $T_+$, is a singular leaf based at $b$ taking the turn $(d_1^+,d_2^+)$.

\end{proof}

\subsection{An example of an ideal Whitehead graph}
\label{SectionIWGExample}

In this section we give an example of a train track map representing a nongeometric, fully irreducible element $\phi \in \Out(F_3)$, together with a complete computation of (the quotient by $F_3$ of) the ideal Whitehead graph of $\phi$. 

\paragraph{Example:} Consider $\phi \in \Out(F_3)$ defined by the automorphism 
\begin{equation*}
\Phi : \begin{aligned}
A &\mapsto AFGFGF\\
F &\mapsto FGF\\
G &\mapsto GFAFG
\end{aligned}
\end{equation*}
That $\Phi$ is indeed an automorphism is easily checked by computing its inverse. Direct calculation shows that the eigenvalues of the matrix of the abelianization of $\Phi$, accurate to the closest $2^{\text{nd}}$ decimal place, are $4.08$ and $.46 \pm .18 i$. No eigenvalue is a root of unity and so the following lemma, whose proof is found below, applies to each power of $\phi$ to show that $\phi$ is fully irreducible.

\begin{lemma} \label{eigenvalue 1} 
If $\phi \in \Out(F_3)$ is reducible then at least one eigenvalue of the abelianization of $\phi$ is a root of unity.  
\end{lemma}

The formula for $\Phi$ also defines an irreducible train track representative $g$ of $\phi$ acting on the rose with edges labelled $A,F,G$. We will see shortly that $g$ has no periodic Nielsen path. It follows by \cite{BestvinaHandel:tt} that $\phi$ is not geometric. It also follows, by Lemma~\ref{singular leaves}, that the local stable Whitehead graph at the valence~6 vertex of the rose is equal to the ideal Whitehead graph of $\phi$ (modulo the action of $F_3$).

The first of the following two graphs is the local Whitehead graph of $g$ at the valence~6 vertex. The second is the local stable Whitehead graph, as well as the ideal Whitehead graph of $\phi$ (modulo $F_3$), as we shall soon show:
$$
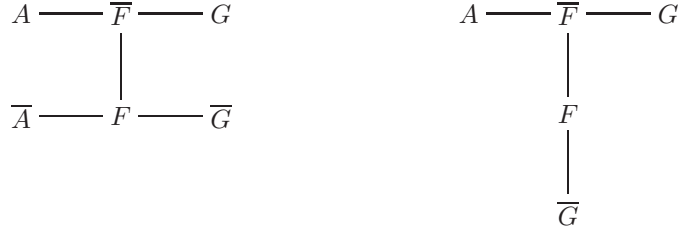
\xymatrix{
              A \ar@{-}[r]   & \overline F \ar@{-}[r] \ar@{-}[d]    &                G \\
\overline A \ar@{-}[r]  &                F \ar@{-}[r]                     & \overline G
} \qquad\qquad\qquad\qquad
\xymatrix{
              A \ar@{-}[r]   & \overline F \ar@{-}[r] \ar@{-}[d]    &   G \\
                                   &                F \ar@{-}[d]    &       \\
                                   & \overline G                                   &                                   
}
$$
Note that the edge $\{\overline F, F\}$ only appears in the local Whitehead graph after the second iterate of $g$, and that $\overline A$ is the only nonperiodic direction.

To show that the second graph above is the ideal Whitehead graph of $\phi$ (modulo $F_3$), it remains to demonstrate that $g$ has no periodic indivisible Nielsen path $\alpha * \bar \beta$. This demonstration amounts to a recursive construction of longer and longer candidates for terminal segments of $\alpha$ and $\beta$, throwing out impossible cases whenever encountered, until either $\alpha,\beta$ are constructed, or all cases are exhausted and $\alpha * \bar \beta$ is thereby proved not to exist. 

Suppose that there exists an indivisible periodic Nielsen path $\alpha * \bar \beta$. Since $\{\overline A, \overline F\}$ is the unique illegal turn of $g$, up to relabelling we have 

\medskip

\centerline{
\begin{tabular}{lrr}
(1) & $\alpha = \alpha_1 * A$  & $g \alpha = g \alpha_1 * AFGFGF$ \\
     & $\beta = \beta_1 * F$      & $g \beta  = g \beta_1 * FGF$
\end{tabular}
}

\medskip\noindent
The last letter of $g\beta_1$ must be $G$, because the last letters of $g\alpha$, $g\beta$ not in their common terminal segment must be $A,F$. The last letter of $\beta_1$ must therefore be $G$:
 
\medskip

\centerline{
\begin{tabular}{lrr}
(2) & $\alpha = \alpha_2 * A$  &      $g \alpha = g \alpha_2 * AFGFGF$ \\
     & $\beta = \beta_2 * GF$      & $g \beta  = g \beta_2 * GFAFGFGF$
\end{tabular}
}

\medskip\noindent
The recursion continues, with justifications to follow:
 
\medskip

\centerline{
\begin{tabular}{lrr}
(3)  & $\alpha = \alpha_3 * FA$    &     $g \alpha = g \alpha_3 * FGFAFGFGF$ \\
      & $\beta = \beta_3 * GF$        &     $g \beta = g \beta_3 *      GFAFGFGF$ \\ \hline
(4)   & $\alpha =  \alpha_4 *  FA$  &     $g \alpha = g\alpha_4 *   FGFAFGFGF$ \\
       & $\beta = \beta_4 *    FGF$  &     $g \beta = g \beta_4 * FGFGFAFGFGF$ \\ \hline
(5)   & $\alpha = \alpha_5 * GFA$ &      $g \alpha = g\alpha_5 * GFAFGFGFAFGFGF$ \\
       & $\beta = \beta_5 *     FGF$ &                 $g \beta = g\beta_5 * FGFGFAFGFGF$ 
\end{tabular}
}

\medskip\noindent
To justify (3), $g \alpha_2$ must end in $A$ or $F$, but no edge has image under $g$ ending in $A$ and so $g \alpha_2$ ends in $F$; this implies that $\alpha_2$ ends in $A$ or $F$, but $AA$ is not a leaf segment (because $\{\overline A, A\}$ is not an edge of the local Whitehead graph) and so $\alpha_2$ ends in $F$. The justification of (4) is similar to (3), and (5) is similar to (2). 

Next we conclude as in (3) that $g \beta_5$ ends in $F$, and so we branch into two cases: (5.1) $\beta_5$ ends in $A$; or (5.2) $\beta_5$ ends in $F$.

\medskip

\centerline{
\begin{tabular}{lrr}
(5.1)   & $\alpha = \alpha_{5} * GFA$ &     $g \alpha = g \alpha_{5}  * GFA.FGFGFAFGFGF$ \\
          & $\beta = \beta_{5.1} * AFGF$ & $g \beta = g \beta_{5.1} * AFGFGF.FGFGFAFGFGF$ \\ \hline
(5.2)   & $\alpha = \alpha_{5} *  GFA$ & $g \alpha = g\alpha_{5} *   GFA.FGFGFAFGFGF$ \\
        & $\beta = \beta_{5.2} *   FFGF$ & $g \beta = g\beta_{5.2} *  FGF.FGFGFAFGFGF$ \\ 
\end{tabular}
}

\medskip\noindent
In these examples the dot indicates the beginning of the longest common terminal segment. If $\alpha * \bar \beta$ were a fixed Nielsen path then it would follow that $\alpha_{5} * GFA = g \alpha_5 * GFA$ implying $\alpha_5 = g \alpha_5$ and we conclude that $\alpha_5$ is the empty word. But the same conclusion holds even under our presumption that $\alpha * \bar\beta$ is a periodic Nielsen path of some period $k \ge 1$, because repeating the analysis above allows us to conclude that for all $i \ge 1$, the edge path $g^i \alpha_5 * GFA$ is the initial segment of $g^i\alpha$ obtained by truncating the longest common terminal segment of $g^i\alpha$ and $g^i\beta$, and so $\alpha_5 = g^k \alpha_5$ implying that $\alpha_5$ is the empty word. 

We have shown that $\alpha = GFA$ is one-half of any periodic indivisible Nielsen path. Next we show that any periodic indivisible Nielsen path has period $k=1$. For suppose that $k>1$ were the period of $\alpha * \bar \beta_1$. Straightening $g(\alpha * \bar\beta_1)$ would yield $\alpha * \bar \beta_2$ with $\beta_1, \beta_2$ distinct from each other and from $\alpha$. This implies that the path obtained by straightening $\beta_1 * \bar \beta_2$ is a periodic indivisible Nielsen path, and so $\beta_1 = \gamma_1 * \delta$, $\beta_2 = \gamma_2 * \delta$ with $\gamma_1 * \bar \gamma_2$ a periodic indivisible Nielsen path. But it follows that $\gamma_i=\alpha=GFA$ for some $i=1,2$, which implies that $\Length(GFA)=\Length(\alpha)=\Length(\gamma_i) \le \Length(\beta_i) = \Length(\alpha)$ and so the lengths are equal and $\alpha = \gamma_i=\beta_i$, contradicting that $\alpha \ne \beta_i$.

We now continue analyzing cases (5.1) and (5.2), under the stronger condition that $\alpha * \bar\beta$ is a \emph{fixed} indivisible Nielsen path. The segment of $g \beta$ before the dot must agree with $\beta$, as it does in (5.2) but not in (5.1), and so (5.1) is ruled out. Continuing with (5.2), since $g\beta$ before the dot must agree with $\beta$, the last letter of $g\beta_{5.2}$ must be $F$, and so either (5.1.1) the last letter of $\beta_{5.2}$ is $A$, or (5.1.2) the last letter of $\beta_{5.2}$ is $F$:

\medskip

\centerline{
\begin{tabular}{lrr}
(5.2.1)   & $\alpha = GFA$  & $g \alpha = GFA.FGFGFAFGFGF$ \\
         & $\beta = \beta_{5.2.1} * AFFGF$      & $g \beta = g\beta_{5.2.1} *    AFGFGFFGF.FGFGFAFGFGF$ \\ \hline
(5.2.2)   & $\alpha = GFA$        & $g \alpha = GFA.FGFGFAFGFGF$ \\
         & $\beta = \beta_{5.2.2} * FFFGF$      & $g \beta = g\beta_{5.2.2} *  FGFFGF.FGFGFAFGFGF$ \\
\end{tabular}
}

\medskip\noindent
In neither (5.2.1) nor (5.2.2) does $\beta$ agree with $g \beta$ before the dot, and so both cases are ruled out. All cases being exhausted, it follows that $g$ has no periodic indivisible Nielsen path.

\bigskip

To complete the example we give:

\begin{proof}[Proof of Lemma \ref{eigenvalue 1}]
Since $\phi$ is reducible there is a free factor $F$ of rank one or of rank two whose conjugacy class is invariant under an iterate of $\phi$. After passing to a power of $\phi$,  the conjugacy class of $F$ becomes $\phi$-invariant. If $F$ has rank one then the conjugacy class of a generator of $F$ is $\phi$-invariant, and the image of this generator  in $\Z^3$ is an eigenvector with eigenvalue one.  We may therefore assume that $F$ has rank two.   Lemma~3.2 of \BookOne\ implies that, after passing to a power of $\phi$,  there exist $\Phi \in \Aut(F_3)$ representing $\phi$ and generators $a,b$ and $c$ of $F_3$ such that $a$ and $b$ generate $F$ and such that $\Phi(c) =ucv$ for some $u,v \in F$.   Letting $x$ be the abelianized image of $c$ and $Z^2$ the abelianized image of $F$, we have $\phi_\#(x)= x+s$ for some $s \in \Z^2$.  If  $1$ is an eigenvalue for the restriction of $\phi_\#$ to $\Z^2$, we are done.   Otherwise there exists $t \in \Z^2$ such that $\phi_\#(t) = t-s$ and $x+t$ is an eigenvector with eigenvalue one.  
\end{proof}







\section{Cutting and pasting local stable Whitehead graphs}
\label{SectionCutAndPaste}

Let $\phi \in \Out(F_n)$ be nongeometric and fully irreducible. Corollary~\ref{ideal is union of local} shows that choosing a train track representative $g \from G \to G$ of $\phi$ induces a decomposition of each component of the ideal Whitehead graph $W(\phi)$, expressing that component as the union of local stable Whitehead graphs pasted together along cut vertices. We refer to this as the \emph{local decomposition} of $W(\phi)$ induced by $g$.

In this section we study how the local decomposition changes under certain moves that alter the train track representative. Our main result, Corollary~\ref{no cut vertices}, explains how to construct a train track representative so that the local decomposition is as fine as possible, having no cut vertices in any local stable Whitehead graph. To put it another way, every finite graph has a unique decomposition into a union of maximal subgraphs without cut vertices, called the \emph{cut point decomposition} of the graph, and Corollary~\ref{no cut vertices} produces a train track representative that simultaneously realizes the cut point decomposition for all components of $W(\phi)$. 

The local decomposition of $W(\phi)$ is intimately connected with Nielsen classes of principal vertices of train track representatives. Given a train track representative $g \from G \to G$ of $\phi$, while the number of Nielsen classes of principal vertices of $g$ depends only on $\phi$, the number of principal vertices depends on the choice of $\gG$. As Lemma~\ref{singular leaves} and Corollary~\ref{ideal is union of local} show, for any train track representative $g \from G \to G$ of $\phi$, the decomposition of $W(\phi)$ into local stable Whitehead graphs of $g$ pasted together along cut vertices corresponds to the decomposition of Nielsen classes of principal vertices of $g$ into individual vertices joined by Nielsen paths. 

First, in Lemma~\ref{Nielsen collapse}, we describe a  condition which allows us to reduce the number of principal vertices, and a method for doing so --- this corresponds to coarsening the local decomposition of $W(\phi)$. Then, in Lemma~\ref{split Nielsen}, we give necessary and sufficient conditions for recognizing when a given irreducible train track map is the product of this method.  In this case the method can be reversed to increase the number of principal vertices, refining the local decomposition of $W(\phi)$. As an application, in Corollary~\ref{no cut vertices} we show that the finest decomposition of $W(\phi)$ is realized by maximizing the number of principal vertices of a train track representative; the maximum exists because the number of principal vertices is bounded above by the number of singular leaves of $\Lambda_-$.

\subsection{Pasting local stable Whitehead graphs} 

Close study of Lemma~\ref{singular leaves} and Corollary~\ref{ideal is union of local} indicates that we should be able to paste together local stable Whitehead graphs and coarsen the local decomposition of the ideal Whitehead graph, by decreasing the number of principal vertices, or equivalently by eliminating Nielsen paths. The following describes a process for carrying out this elimination; it is a variation on Lemma~3.8 of \cite{BestvinaHandel:tt}.

\begin{lemma} \label{Nielsen collapse}  Suppose that $\gG$ is an irreducible train track map, that $E_1$ and $E_2$ are edges in $G$ with a common terminal vertex and that $E_1 \bar E_2$ is an indivisible periodic Nielsen path that either has period one or is mapped by $f_\#$ to its inverse (and so has period two). Let $G'$ be the graph obtained from $G$ by isometrically identifying $E_1$ with $E_2$ and let $p:G \to G'$ be the quotient map. Then there is an induced irreducible train track map $g' : G' \to G'$ with the following properties:
\begin{enumerate}
\item $g' p = pg$. 
\item If $f_{g}:\wt G \to T_+$ and $f_{g'}:{\wt G'} \to T_+$ are as in Corollary~\ref{first tree corollary} then $p$ lifts to $\ti p :\wt G \to \wt G'$ such that $f_{g}= f_{g'}\ti p$. 
\item If $g$ is a train track representative of $\phi \in \Out(F_r)$ via a marking $\rho \from R_r \to G$, then $g'$ is also a train track representative of $\phi$ via the marking $p \composed \rho \from R_r \to G'$.
\item The function $p$ maps the principal vertices of $g$ onto the principal vertices of $g'$, preserving Nielsen classes. This map is one-to-one except that the endpoints of $E_1 \bar E_2$ are identified to a single point. 
\end{enumerate}
\end{lemma}

In this lemma, the local decomposition of $W(\phi)$ induced by $g$ is obtained from that induced by $g'$ by pasting the local stable Whitehead graphs of $g'$ at the initial endpoints of $E_1$, $E_2$, identifying the direction of $E_1$ with the direction of $E_2$.

We shall refer to the construction of this lemma as a \emph{Nielsen path collapse} performed on $E_1 \bar E_2$.

\proof The initial endpoints of $E_1$ and $E_2$ are distinct because $\phi$ does not fix the conjugacy class of any basis element.  There is a path $\tau$ such that $g(E_i) = E_{\sigma(i)} \tau$ where $i=1,2$ and $\sigma :\{1,2\} \to \{1,2\}$.  Denote $p(E_1) = p(E_2)$ by $E'$ and define $g'(E') = E' p(\tau)$.  Every other edge of $G'$ equals $p(E)$ for some edge $E$ of $G$ and we define $g'(p(E))= pg(E)$ which is an immersed path because $p$ restricts to an immersion on legal paths in $G$. For the same reason $(g')^k(p(E))= pg^k(E)$ for all $k > 0$ which implies that $g'$ is an irreducible train track map. Item~(1) is true by construction. Let $\ti p :\wt G \to \wt G'$ be the lift that fixes all edges that do not project to $E_1$ or $E_2$.  Then $f_{g}$ and $f_{g'}\ti p$ are equivariant, agree on vertices and restrict to isometries on all edges. It follows  that $f_{g}= f_{g'}\ti p$, proving (2). Item (3) is obvious.

To prove (4), let $w_i$ be the initial point of $E_i$ and let $w'=pw_1=pw_2$. The map $p$ is injective on vertices except that $pw_1=pw_2$. A vertex $v$ of $G$ is $g$-periodic if and only if $pv$ is $g'$-periodic; for each of the two cases $v \in \{w_1,w_2\}$ and $v \not\in \{w_1,w_2\}$ this is obvious, but for slightly different reasons.  We now assume that $v \in G$ is a principal vertex and prove, by consideration of cases, that $pv \in G'$ is principal. Case 1: if $v$ has at least three $Dg$-periodic directions then their images under $D_vp$ give at least three $Dg'$-periodic directions at $pv$. Case 2: if $v$ is the endpoint of a periodic Nielsen path $\gamma$ other than $E_1 \bar E_2$ or its inverse, then $p\gamma = p_\#\gamma$ is a periodic Nielsen path and has $pv$ as an endpoint. Case 3: if $v=w_1$ or $w_2$ then $Dg$ has at least two periodic directions at each of $w_1$ and $w_2$, and so $pv=w'$ has at least three periodic directions, namely $D_{w_1}p(E_1)=D_{w_2}p(E_2)$, one more in $\image(D_{w_1}p)$, and one more in $\image(D_{w_2}p)$. It follows that if $v$ is principal then so is $pv$. Conversely, every vertex of $G'$ with at least three periodic directions occurs as in Case 1 or Case 3, and every periodic Nielsen path of $g'$ occurs as in Case 2. This proves that $v$ is principal for $g$ if and only if $pv$ is principal for $g'$, and moreover that Nielsen classes are preserved. The statement about local decompositions is obvious by construction.
\endproof

A Nielsen path collapse is similar to the more familiar operation from \cite{BestvinaHandel:tt} of doing a fold on a turn of a train track map to produce a new train track map. We note some of the properties of the latter; the proof is similar to Lemma~\ref{Nielsen collapse} but easier.

\begin{lemma}
\label{LemmaStallingsFold}
Consider an irreducible train track map $g \from G \to G$ and oriented edges $E_1, E_2 \subset G$ with a common terminal vertex so that $g(E_1)=g(E_2)$. Let $p \from G \to G'$ be the \emph{Stallings fold}, the quotient map that identifies $E_1$ to $E_2$. Then there is an induced irreducible train track map $g' \from G' \to G'$ that satisfies items (1)--(3) of Lemma~\ref{Nielsen collapse}, and such that
\begin{enumerate}
\item[(4$'$)] The function $p$ restricts to a bijection between the principal vertices of $g$ and $g'$ preserving Nielsen classes. The local decompositions of the ideal Whitehead graph induced by $g$ and by $g'$ are identical.
\end{enumerate}
\qed\end{lemma}

We say that $g'$ is obtained from $g$ by a \emph{fold} of $E_1$ with $E_2$.

\subsection{Cutting local stable Whitehead graphs} 

A Nielsen path collapse is reversible, and one can ask how to recognize when a given irreducible train track map $g' : G' \to G'$ is the product of a Nielsen path collapse applied to some $\gG$ as in Lemma~\ref{Nielsen collapse}. Since the map $p:G \to G'$ restricts to an immersion on leaves of $\Lambda_-$ it induces a map of local Whitehead graphs $Dp:W_w \to W_{p(w)}$ for each vertex~$w$ of~$G$.  Similarly $g :G\to G$ induces a map $Dg:W_w \to W_{g(w)}$. If $w_i \in G$ is the initial vertex of $E_i$ and $w' =p(w_i) \in G'$ then $W_{w'}$ decomposes as the union of $Dp(W_{w_1})$ and $Dp(W_{w_2})$ and these two subgraphs of $W_{w'}$ intersect in the single vertex determined by $E' =p(E_i)$. This vertex is therefore a  {\em cut vertex} of $W_{w'}$, meaning that it separates $W_{w'}$. Note also that if $g'(v') =w'$ for some vertex $v' \ne w'$ then $Dg'(W_{v'})$ is contained in either $Dp(W_{w_1})$ or $Dp(W_{w_2})$.  Finally, $w'$ is a fixed point of $g'$ and the action of $Dg'$ on $W_{w'}$ preserves the decomposition $W_{w'} = Dp(W_{w_1}) \union Dp(W_{w_2})$, fixing or interchanging the two sets $Dp(W_{w_1})$, $Dp(W_{w_2})$ depending on whether $g_\#$ preserves or reverses $E_1 \bar E_2$. 

The following lemma states that these necessary conditions are also sufficient. The lemma is stated with sufficient breadth so that the construction may be used as well to reverse a fold.

\begin{lemma} \label{split Nielsen} Suppose that $g' : G'\to G'$ is an
irreducible train track map, that $w'$ is a vertex of $G'$ and that  
\begin{enumerate}
\item $W_{w'}$ is the union of two non-trivial subgraphs $X_1$ and $X_2$
that intersect in a single vertex $x$. 
\item If $g'(v') =w'$ for some vertex $v' \ne w'$ then $Dg'(W_{v'})$
is contained in either $X_1$ or $X_2$. 
\item If $g'(w') = w'$ then $Dg'(x) = x$ and there exists
$\sigma:\{1,2\} \to \{1,2\}$  such that $Dg'(X_i) = X_{\sigma(i)}$. 
\end{enumerate}
Then there exists an irreducible train track map $g : G \to G$ 
and a pair of edges $E_1$ and $E_2$ in $G$ with a common terminal vertex
and with distinct principal initial vertices $w_1$ and $w_2$ such that the
quotient map $p : G \to G'$ that identifies $E_1$ to $E_2$ is a
homotopy equivalence satisfying: 
\begin{enumeratecontinue}
\item $ pg = g' p$.
\item $ p(w_1)=  p(w_2) = w'$.
\item $ p(E_1) =  p(E_2) =E'$ is the edge in $W_{w'}$ corresponding to $x$.
\item $D p(W_{w_i}) =  X_i$ for $i=1,2$.
\end{enumeratecontinue}
Furthermore, if $g'(w')=w'$ then $g$ and the edges $E_1, E_2$ satisfy the hypotheses of Lemma~\ref{Nielsen collapse} and $g'$ is obtained from $g$ by a Nielsen path collapse along $E_1 \bar E_2$, whereas if $g'(w') \ne w'$ then $g$ and $E_1,E_2$ satisfy the hypotheses of Lemma~\ref{LemmaStallingsFold} and $g'$ is obtained from $g$ by a fold of $E_1$ with $E_2$.
\end{lemma}

\begin{proof}     For $i=1,2$ denote the edges of $G'$ that correspond to vertices of $X_i$ other than $x$ by $P'_i$. Let $u'$ be the terminal vertex of the edge $E'$ corresponding to $x$. Define a new graph $G$ from $G'$  by replacing $w'$ with a pair of vertices $w_i$, by replacing $E'$ with a pair of edges $E_i$ with initial vertex $w_i$ and terminal vertex $u'$ and by reattaching the edges in $P'_i$ to $w_i$.  

The map $p : G \to G'$  that identifies $E_1$ to $E_2$ is a homotopy equivalence.  We construct  $g : G \to G$ as follows.  If $g'(w')=w'$ then $g(w_i) := w_{\sigma(i)}$; otherwise $g(w_i) := p^{-1}g'(w')$. If 
$v$ is  a vertex of $G$ such that $g'(p(v)) \ne w'$, then  $g(v):=  p^{-1}(g'( p(v)))$.   If $v$ is neither $w_1$ nor $w_2$ and if $g'(p(v))= w'$ then $g(v) := w_i$ where $Dg'(W_v) \subset X_i$. If $E$ is an edge with endpoints $s$ and $t$ then define $g(E)$ to be the unique path with endpoints $g(s)$ and $g(t)$ such that $pg(E) =g'p(E)$. This unique lifting exists because the initial and terminal vertices have already  been lifted and the base vertex of each turn in $g'(p(E))$ lifts uniquely: at the only vertex $w'$ that does not lift uniquely to $G$ the turns of $g'(p(E))$ lift uniquely because they are contained in either $X_1$ or $X_2$.  This completes the definition of $g$. Properties (4)--(7) are clear from the construction.  

We must still show that $\gG$ is an irreducible train track map. The obvious induction argument shows that $pg^k=(g')^kp$ for all $k \ge 0$. If $E$ is an edge of $G$ then $p(E)$ is an edge of $G'$ and so $(g')^kp(E)$ is an immersed path for all $k \ge 0$.  It follows that $g^k(E)$ is immersed and hence that $g$ is a train track map.  
Each $(g')^k(pE)$ is a subpath of a leaf of the realization of $\Lambda_-$ in $G'$.  It follows that
$g^k(E)$ is a subpath of a leaf of the realization of $\Lambda_-$ in $G$.  Example 2.5(1) of 
\BookZero\ implies that $g^k(E)$ covers $G$ for all sufficiently large $k$ and hence that $g$ is 
irreducible.

It remains to prove the ``furthermore'' clause. If $g'(w')=w'$ then clearly $E_1 \bar E_2$ is a periodic Nielsen path mapped by $g_\#$ to itself or its inverse, so the Nielsen collapse of $E_1 \bar E_2$ is defined and its product is clearly $g'$. If $g'(w') \ne w'$ then $g(w_1)=g(w_2) = p^\inv(g'(w'))$, and clearly $g(E_1)=g(E_2)$ is the unique path in $G$ with initial point $p^\inv(g'(w'))$ that projects to $g'(E)$; the fold of $E_1$ with $E_2$ is therefore defined and its product is clearly $g'$.
\end{proof}

\subsection{The finest local decomposition} 

Here is our main result, in which we describe a process that starts with an arbitrary train track representative and produces one whose local decomposition is finer than any other. We do this only in the rotationless case; a more general result, while interesting, is not necessary for our present purposes, and the technical details would be a distraction.

\begin{proposition} \label{no cut vertices} Let $\phi \in \Out(F_r)$ be nongeometric, fully irreducible, and rotationless. If $g : G \to G$ is an irreducible rotationless train track map representing $\phi$ with the maximal number of principal vertices then $SW_v$ has no cut vertices for any principal vertex $v$ of $g$, and the local decomposition of $W(\phi)$ induced by $g$ is the cut point decomposition.
\end{proposition}

The proof gives a process that will produce such a $g$ starting from any train track representative of $\phi$. The steps in the process involve some folds, some inverse Nielsen path collapses, and some inverse folds.

\begin{proof} We begin with some preliminary folding that does not change the number of principal vertices. Recall that a pair of directions $d_1$ and $d_2$ at the same vertex $w$ are in the same gate at $w$ if they are identified by some iterate of $Dg$. Define the \emph{gate excess} at $w$ to be the number of directions minus the number of gates. Define the \emph{total gate excess} to be the sum of the gate excesses over all principal vertices and over all nonprincipal, and hence non-periodic, vertices with three or more gates. If $w$ has two directions $d_1,d_2$ that are identified by $Dg$ itself then, after possibly subdividing $G$, we can fold the edges corresponding to $d_1$ and $d_2$ to create a new nonprincipal vertex with two gates, and so the total gate excess is reduced as long as $w$ is not a nonprincipal vertex with two gates. Since $g$ is rotationless, this can be applied repeatedly to the principal vertices until each of them has zero gate excess. If the total gate excess is still positive then there exists a nonprincipal vertex $w$ with at least three gates and positive gate excess. It follows that some positive iterate of $g$ takes $w$ to a principal vertex, and so there exists a minimal $i \ge 1$ such that $g^i(w)$ has zero gate excess. Replacing $w$ by $g^{i-1}(w)$ we may assume that there are two directions $d_1,d_2$ at $w$ that are identified by $Dg$ itself, and so a subdivision and fold at $w$ reduces the total gate excess. Continuing by induction, we may assume that $g$ has zero gate excess. In particular, for every vertex $w$ which is either principal or has at least three gates, $Dg \restrict W_w : W_w \to W_{g(w)}$ is injective. Furthermore, $W_w = SW_w$ if $w$ is principal.

Assuming now that $g$ has a principal vertex $w$ such that $W_w = SW_w$ has a cut vertex, we describe a sequence of inverse folds followed by an inverse Nielsen collapse which refines the local decomposition of $W(\phi)$ induced by $g$. Let $X_1$ and $X_2$ be non-trivial subgraphs of $W_{w}$ whose union is $W_{w}$ and whose intersection is $\{x\}$.  Define $V$ to be the set of vertices $v \ne w$ for which there exists $i \ge 0$ so that $g^i(v) = w$ and so that $Dg^i(W_v)$ is not contained in either 
$X_1$ or $X_2$.   Thus $Dg^i(W_{v})$ contains $x$, at least one vertex in $X_1 \setminus\{x\}$ and at least one vertex in $X_2 \setminus\{x\}$.  In particular, $v$ has at  least three gates.  Let $X_1' \cup X_2'$ be the decomposition of $W_v$ defined by pulling back the decomposition $X_1 \cup X_2$ of $W_w$ via $Dg^i$.  It follows, by injectivity of $Dg \restrict W_v$, that $X_1' \cap X_2'$ is a cut vertex $x'$.   If $u$ is a vertex such that $g(u) = v$ and such that $Dg(W_u)$ is not contained in either $X_1'$ or $X_2'$ then $u \in V$.  

\subparagraph{Case 1: $V = \emptyset$.} In this case the three items in the hypothesis of Lemma~\ref{split Nielsen} are satisfied: (1) is obvious, (2) follows from the definition of $V$ and (3) follows from the assumption that $\phi$ is rotationless. Applying the lemma, an inverse Nielsen collapse produces an irreducible train track map with more principal vertices than $g$, refining the local decomposition of $W(\phi)$ induced by $g$.

\subparagraph{Case 2: $V \ne \emptyset$.} Since $g$ is rotationless, there exists  $v \in V$ that is not $g(u)$ for any $u \in V$.  Then Lemma~\ref{split Nielsen} applies to the decomposition $X_1'  \cup X_2'$ described above: (1) is obvious, (2) follows from the choice of $v$ and our discussion above, and (3) is vacuous.   Applying the lemma, an inverse fold produces a new irreducible train track map with the same number of principal vertices and the same local decomposition as $g: G \to G$ but with one fewer element in $V$. Continuing inductively until $V=\emptyset$, we reduce to Case 1.
\end{proof} 
%

%

\section{Weak train tracks}   
\label{SectionWeakTT}

Let $\phi \in \Out(F_r)$ be nongeometric and fully irreducible with expanding lamination $\Lambda_-=\Lambda_-(\phi)$ and attracting $F_r$ tree $T_+=T_+(\phi)$. Consider a marked graph $G$ or its universal cover $\wt G = T$. We say that $G$ or $T$ is a \emph{train track} for $\phi$ if there exists an irreducible rotationless train track map $g \from G \to G$ representing an iterate (a positive power) of $\phi$, and we say that $G$ or $T$ \emph{arises from $g \from G \to G$}. We say that $G$ or $T$ is a \emph{weak train track} for $\phi$ if, up to renormalizing $T$ and/or $T_+$, there exists a $\Lambda_-$ isometry $f_T \from T \to T_+$. The definitions of train tracks and weak train tracks are evidently invariant under the equivalence relation defining $\X_r$, and so we can speak about a point of $\X_r$ being a train track or a weak train track.

Since the definition of a weak train track for $\phi$ depends only on the lamination $\Lambda_-$ and the tree $T_+$, and since these two objects are insensitive to iterating $\phi$, a weak train track for $\phi$ is the same thing as a weak train track for any iterate of $\phi$. We will therefore replace $\phi$ with an iterate whenever it is convenient. 

Henceforth we shall always fix a normalization of $T_+$ a priori. Having done so, we would like to normalize each weak train track $T$ so that there exists a $\Lambda_-$ isometry $f_T \from T \mapsto T_+$.   As we will see in Section~\ref{SectionIrrigidity}, there are examples where $f_T$ is \emph{not} uniquely determined by $T$, although the failure of uniqueness is mild (if $\Lambda_-$ is nonorientable then $f_T$ is unique; whereas if $\Lambda_-$ is orientable then $f_T$ varies over a space of maps homeomorphic to a closed interval in the real line). It might seem that the normalization of $T$ depends on the choice of $f_T$. The following says that this does not happen:

\begin{lemma}
\label{LemmaNormalization}
Fix a normalization of $T_+$. For each weak train track $T \in \X_r$ there exists a unique normalization of $T$ in $\wh\X_r$ for which there exists a $\Lambda_-$ isometry from $T$ to $T_+$.
\end{lemma}

\begin{proof} Consider any normalization of $T$ for which there exists a $\Lambda_-$ isometry $f_T \from T \to T_+$. Consider the marked graph $G = T/F_r$. Each leaf $L$ of $\Lambda_-$ is realized by a surjective local isometry $\ell \from \reals \mapsto G$, and so there is an oriented edge $e$ of $G$ and points $r<s \in \reals$ such that $\ell(r)=\ell(s) \in \interior(e)$ and $\ell$ is orientation preserving near $r$ and near $s$. Surgery of $\ell$ produces a loop by restricting $\ell$ to $[r,s]$ and identifying $r$ and $s$. If $\gamma \in F_r$ represents this loop then the length of this loop is equal to the translation length of $\gamma$ acting on $T$. Moreover, by construction the axis of $\gamma$ in $T$ is legal with respect to the map $f_T$, and so is mapped by $f_T$ isometrically and $\gamma$-equivariantly onto the axis of $\gamma$ in $T_+$. The translation lengths of $\gamma$ acting on $T$ and on $T_+$ are therefore equal, thereby uniquely determining the normalization of $T$.
\end{proof}

By Lemma~\ref{CorollaryLeavesInTPlus}, every train track for $\phi$ is a weak train track for $\phi$. To be more explicit, for every train track representative $g \from G \to G$ of $\phi$, the train track $\wt G$ is a weak train track, and the map $f_g \from \wt G \to T_+$ is a $\Lambda_-$ isometry.   The main technical result of this section is Proposition~\ref{PropWeakTrainTrack} which explains how to obtain every sufficiently short weak train track by folding one particular train track. We also explain in Proposition~\ref{characterize train track} how to recognize when a weak train track is a train track. As an application we prove in Corollary~\ref{CorollaryWeakIsDense} that the set of weak train tracks for $\phi$ is dense in the set of train tracks for $\phi$, as subsets of $\CV_r$.

Much of what we do in this section is to generalize some of the structure of train tracks, such as the local stable Whitehead graphs  and the local decomposition of the ideal Whitehead graph $W(\Phi)$,  to the context of a weak train tracks.   These notions will depend not only on the weak train track $T$ but also  the choice of a $\Lambda_-$ isometry $f_T \from T \to T_+$.  We will sometimes emphasize the dependence on $f_T$ by referring, for example, to the local decomposition of $W(\Phi)$ induced by $T$ \emph{relative to $f_T$}. As we proceed we will begin to suppress this dependence when $f_T$ is clear from the context. On first pass, the reader may wish to assume uniqueness of $f_T$ (a perfectly safe assumption when $\Lambda_-$ is nonorientable).

\subsection{Local decomposition of the ideal Whitehead graph.} Let $T$ be a weak train track for $\phi$. Choose a $\Lambda_-$ isometry $f_T \from T \to T_+$. If $b$ is a branch point of $T_+$ and $\Phi \in \PA(\phi)$ corresponds to $b$ as in Lemma~\ref{principal lifts and branch points}, then the realization in $T$ of each leaf $\ti \ell \in L(\Phi)$ contains a unique point $w$ that is mapped by $f_T$ to $b$, called the \emph{base point} of $\ti\ell$ in $T$. Fixing $w \in T$, as $\ti\ell$ varies over the leaves in $L(\Phi)$ with base point $w$, the point $w$ divides the realization of $\ti \ell$ into two rays called {\em singular rays of $T$ at $w$}. Each singular ray of $T$ at $w$ is mapped by $f_T$ to a singular ray of $T_+$ at $b$. There are only finitely many $F_r$-orbits of basepoints $w$ so we may equivariantly subdivide $T$ to make each basepoint a vertex. A vertex of $T$ is {\em principal} if it is the basepoint of some singular leaf. Base points of singular leaves, singular rays, and principal vertices are all defined relative to $f_T$. 

For each point $v$ of $T$, the {\em local Whitehead graph} $W(v;T)$ is the graph whose vertices are directions at $v$ and with an edge connecting two vertices if the turn that they determine is taken by a leaf of $\ti \Lambda_-$. Local Whitehead graphs are defined \emph{absolutely}, not relative to $f_T$. 

For each principal vertex $w$ of $T$ relative to $f_T$, the {\em local stable Whitehead graph} $SW(w;T)$ is the graph whose vertices are initial directions of singular rays based at $w$ and with an edge connecting two vertices if the turn that they determine is taken by a singular leaf of $\ti \Lambda_-$ realized in $T$ with base point $w$. If $f_T(w) = b$ then $f_T$ induces an injection of $SW(w;T)$ into $W(b;T_+)$; we use this and Lemma~\ref{TPlusWhitehead} to identify $SW(w;T)$ with a subgraph of $W(b;T_+)$ and with a subgraph of the component $W(\Phi)$ of the ideal Whitehead graph $W(\phi)$ where  $\Phi \in \PA(\phi)$ corresponds to $b$. The local stable Whitehead graph $SW(w;T)$ is defined relative to $f_T$, so a more detailed notation which we prefer to avoid would be $SW(w;f_T)$.

The following lemma generalizes Lemma~\ref{ideal is union of local}.

\begin{lemma} \label{gen partition} Suppose that $T$ is a weak train track, $f_T \from T \to T_+$ is a $\Lambda_-$ isometry, $b$ is a branch point of $T_+$, $\Phi \in \PA(\phi)$ is the representative of $\phi$ corresponding to $b$, and $\{w_i\}$ are the principal vertices of $T$ that are mapped by $f_T$ to $b$.  Then 
\begin{enumerate} 
\item $\displaystyle W(\Phi) = \cup SW(w_i;T)$. 
\item For $i \ne j$, $SW(w_i;T) \cap SW(w_j;T)$ intersect in at most one vertex.  If they do intersect at a vertex $P$, then $P$ is a cut point of $W(\Phi)$, in fact $P$ separates $SW(w_i;T)$ from $S(w_j;T)$ in $W(\Phi)$.
\end{enumerate}
\end{lemma}

\begin{proof} (1) is immediate from the definitions and from our identification of $SW(w_i;T)$ with a subgraph of $W(\Phi)$.  The proof of (2) is the same as the proof of the third item of Corollary~\ref{ideal is union of local}. 
\end{proof}

As in the case of Corollary~\ref{ideal is union of local}, Lemma~\ref{gen partition} shows that the choice of a weak train track $T$ and a $\Lambda_-$ isometry $f_T \from T \to T_+$ induces a decomposition of the ideal Whitehead graph $W(\phi)$, expressing each component of $W(\phi)$ as a union of local stable Whitehead graphs of $T$ pasted together along cut vertices. We refer to this as the \emph{local decomposition} of $W(\phi)$ relative to $f_T$. 

Given two weak train tracks $S,T$ and $\Lambda_-$ isometries $f_S \from S \to T_+$, $f_T \from T \to T_+$, we say that $f_T$ \emph{splits at least as much as} $f_S$ if the local decomposition $W(\phi) = \cup SW(v_j;T)$ is at least as fine as the the local decomposition $W(\phi) = \cup SW(w_i;S)$, meaning that for each principal vertex $v_j$ of $T$ there is a principal vertex $w_i$ of $S$ such that $SW(v_j;T) \subset SW(w_i;S)$ where the inclusion takes place in $W(\phi)$.

To put it another way, a $\Lambda_-$ isometry $f_T \from T \to T_+$ imposes an equivalence relation on the set of singular leaves realized in $T$, where $l_1 \stackrel{f_T}{\sim} l_2$ if $l_1,l_2$ have the same base points relative to $f_T$. Obviously $l_1 \stackrel{f_T}{\sim} l_2$ if and only if their edges belong to the same component of the local decomposition of $W(\phi)$ relative to $f_T$. Thus, to say that $f_T$ splits at least as much as $f_S$ means that, for all singular leaves $l_1,l_2$, if the realizations of $l_1,l_2$ in $T$ have the same basepoint relative to $f_T$ then their realizations in $S$ have the same basepoint relative to $f_S$. 

\begin{lemma} \label{max split} There exists a rotationless irreducible train track map $\gG$ representing an iterate of $\phi$ such that for every weak train track $T$ and $\Lambda_-$ isometry $f_T \from T \to T_+$, the $\Lambda_-$ isometry $f_{g} \from
\wt G \to T_+$ splits at least as much as $f_T$. If $\phi$ is rotationless then we may choose $\gG$ to represent $\phi$.
\end{lemma}

\begin{proof} After replacing $\phi$ with an iterate we may assume that $\phi$ is rotationless. Apply Lemma~\ref{no cut vertices} to choose a rotationless train track map $\gG$ whose local stable Whitehead links have no cut vertices.  For each $\Phi \in \PA(\phi)$, there are decompositions $W(\Phi) = \cup SW(\ti v_i;\wt G)$ and $W(\Phi) = \cup SW(w_j;T)$ where $f_{\wt G}(\ti v_i) = f_T(w_j) = b$. Since $SW(\ti v_i;\wt G)$ has no cut vertices, it must be contained in some single $SW(w_j;T)$.  
\end{proof}

\subsection{Folding up to a weak train track}

We now come to the main technical result of this section. This result says that from any fixed train track we can fold to any sufficiently short weak train track which splits no more than the fixed train track. This is very similar to Lemma~3.16 of \cite{LosLustig:axes} which combines results from \cite{Lustig:conjugacy}. 

\begin{proposition}
\label{PropWeakTrainTrack}  Let $\phi \in \Out(F_n)$ be nongeometric and fully irreducible. For any train track representative $g \from G \to G$ for $\phi$ with associated $\Lambda_-$ isometry $f_g \from \wt G \to T_+$, there exists $\epsilon > 0$ so that if $T$ is any weak train track for $\phi$, if $f_T \from T \to T_+$ is any $\Lambda_-$ isometry, if $f_g \from \wt G \to T_+$ splits at least as much as $f_T$, and if $\Length(T) \le \epsilon$, then there exists a unique equivariant edge isometry $h \from \wt G \to T$ such that $f_g = f_T h$.  
This map $h$ is moreover a $\Lambda_-$ isometry.
\end{proposition}

In this statement, recall our convention of fixing a normalization of $T_+$, relative to which $\wt G$ and $T$ are normalized, and so $\Length(T)$ is well defined.

\begin{proof} Uniqueness of $h$, and the fact that $h$ must be a $\Lambda_-$ isometry, are proved as follows. For any leaf $L$ of $\Lambda_-$ realized in $\wt G$, since $f_g \restrict L$ is an isometry, $h \restrict L$ must be an injection, and since $h$ is an edge isometry the map $h \restrict L$ must be an isometry onto its image. Moreover, the image $h(L)$ must be the realization of $L$ in $T$, which we denote $L_T$, proving that $h$ is a $\Lambda_-$ isometry. Since $f_T \restrict L_T$ is an isometry, it follows that $h \restrict L = (f_T \restrict L_T)^\inv \composed (f_g \restrict L)$. This being true for any $L$, the map $h$ is uniquely characterized.

For the rest of the proof we turn to existence of $h$.

There is no loss in assuming that $\phi$ is rotationless.  After replacing $\gG$ by an iterate we may assume that for any vertex $v$ and any direction $d$ at $v \in G$, the vertex $g(v)$ is fixed and the direction $D_v g(d)$ is fixed. We now assume the notation of Figure 2, and we write $\wt G_0$ for $\wt G$ and $f_{g_0}$ for $f_g$.

The map $h:\wt G_0  \to T$ will be defined equivariantly, first on the vertices of $\wt G_0$, and then by linear extension over the edges.  With appropriate choice of $\epsilon$ we will then prove that $f_{g_0} = f_T h$, from which it will immediately follow that $h$ is edge isometric.  

If $\ti x \in \wt G_0$ projects to a fixed point of $g$ then there is a natural choice of $h(\ti x) \in T$.  To take care of a general vertex $\ti x \in \wt G_0$, we first push $\ti x$ forward via $\ti g_{1,0}$ from $\wt G_0$ to a point $\ti y \in \wt G_1$ which does project to a fixed point, and which now has a natural choice of image in $T$. In other words, for each vertex $\ti x$ of $\wt G_0$, we shall define $h(\ti x) = h_1 \ti g_{1,0}(\ti x)$ where $h_1$ is defined on $\ti y=\ti g_{1,0}(\ti x) \in \wt G_1$ as follows.

If $\ti y$ is principal then it is the basepoint of the realization of some singular leaf $\ti \ell$ in $\wt G_1$ and we define $h_1(\ti y)$ to be the basepoint of the realization of $\ti \ell$ in $T$. This is independent of the choice of $\ti l$, and therefore well defined, by the 
hypothesis that $\wt G_0$, and so also $\wt G_1$, splits at least as much as $T$. If $\ti y$ is not principal then there are exactly two fixed directions at $\ti y$ and they determine a leaf $\ti \ell$ of $\ti \Lambda_-$ by Lemma~\ref{first characterization}.  There is a lift of $g$ that fixes $\ti y$ and if $\Phi$ is the automorphism corresponding to that lift then $\Phi_\#(\ti \ell) = \ti \ell$.  The homeomorphism $\Phi_+ : T_+ \to T_+$ of Corollary~\ref{action on T+} preserves $\ti \ell$, expanding the metric by $\lambda(\phi)$, and so has a unique fixed point $c$ in the realization of $\ti \ell$ in $T_+$.  Define $h_1(\ti y)$ to be the unique point in the realization of $\ti \ell$ in $T$ that is mapped to $c$ by $f_T$.

Having  defined $h$ on the vertices of $\wt G_0$, we now extend it linearly to all of $\wt G_0$. Choose $\epsilon>0$ to be less than half the length of the shortest edge in $G_1$. We show next that for any edge $\wt E$ of $\wt G_0$, the path $h(\wt E)$ is legal with respect to $f_T$. 

Let $\ti x_0$ and $\ti x_1$ be the endpoints of $\wt E$. For $i=0,1$, let $\ti y_i = \ti g_{1,0}(\ti x_i)$. There is a lift of $g$ that fixes $\ti y_i$ and that fixes the direction at $\ti y_i$ determined by $\ti g_{1,0}(\wt E)$; let $\Phi_i$ be the automorphism corresponding to this lift of $g$. We may choose the leaf $\ti\ell_i$ used in the definition of $h_1(\ti y_i)$ to \lq contain\rq\  the direction at $\ti y_i$ determined by $\ti g_{1,0}(\wt E)$.  To be precise, $\ti\ell_i$ is the leaf through $\ti y_i$ whose realization in $\wt G_1$ is fixed by $(\Phi_i)_\#$. The leaf $\ti\ell_i$ contains a non-trivial segment of $\ti g_{1,0}(\wt E)$ incident to $\ti y_i$, and its realization in $T$ contains $h_1(\ti y_i) = h(\ti x_i)$. Let $\ti \ell_2$ be any leaf whose realization in $\wt G_0$ contains $\wt E$, and so the realization of $\ti \ell_2$ in $\wt G_1$ contains $\ti g_{1,0}(\wt E)$. 

Let $\ti z_0,\ti z_1$ be the midpoints of the initial and terminal edges of  $\ti g_{1,0}(\wt E)$ respectively.  For $i=0,1$, the realizations of $\ti \ell_i$ and $\ti \ell_2$ in $T_+$ contain an interval in $T_+$ of radius greater than $\epsilon$ centered on $f_{g_1}(\ti z_i)$.   Since $T$ has length at most $\epsilon$, the number $\epsilon$ is a cancellation constant (Lemma~\ref{bcl}) for  $f_T \from T \to T_+$.   It follows that the intersection of the realizations of $\ti \ell_i$ and $\ti \ell_2$ in $T$ contains an open arc about a point $\ti w_i$ that maps by $f_T$ to $f_{g_1}(\ti z_i)$. One can therefore travel along a path in $T$ that starts from $h(\ti x_0)=h_1(\ti y_0)$, goes along the realization of $\ti \ell_0$ to $\ti w_0$, then along the realization of $\ti \ell_2$ to $\ti w_1$, and then along the realization of $\ti \ell_1$ to $h_1(\ti y_1)=h(\ti x_1)$, without ever taking an illegal turn. This legal path is therefore the unique embedded arc in $T$ from $h(\ti x_0)$ to $h(\ti x_1)$, which by definition is $h(\wt E)$.

By construction, $f_{g_0}(x) = f_T h(\ti x)$ for all vertices $\ti x$ of $\wt G_0$.  Thus $f_{T}h \from \wt G_0 \to T_+$ is an equivariant map that agrees with $f_{g}$ on vertices and is linear on edges.  It follows that $f_{g_0} = f_T h$ as desired. Also, since $f_{g_0}$ and $f_T$ are both isometric on edges, it follows that $h$ is isometric on edges.
\end{proof}

\subsection{Comparing train tracks to weak train tracks}
\label{SectionTTversusWTT}
There are two moves one can perform on a train track map $\gG$ representing an iterate of $\phi$ to produce a weak train track that is not a \trtr.  The first move is to partially fold a pair of edges to create a vertex that is not preperiodic for the induced map. In this case, although there is an induced homotopy equivalence on the quotient space, this map does not take vertices to vertices no matter how the quotient space is subdivided as a CW complex. The second move is to collapse a pre-Nielsen path by a folding operation that is analogous to the one used in Lemma~\ref{Nielsen collapse} without also collapsing the Nielsen path to which it maps.  In this case there is no induced map on the resulting quotient space.   The next proposition shows that these are the only possibilities, and the following corollary uses this to prove that train tracks are dense in the set of weak train tracks.

We say that $c \in T_+$ is a {\em periodic point} if there is an automorphism $\Phi$ representing a rotationless iterate of $\phi$ such that $c \in \Fix(\Phi_+)$. Every branch point of $T_+$ is a periodic point by Lemma~\ref{principal lifts and branch points}. Also, if $g \from G \to G$ is a train track representative of some rotationless iterate of $\phi$ then Theorem~\ref{TheoremPhiPlus} guarantees that if $x \in \wt G$ projects to a periodic point of $g$ then $f_g(x)$ is a periodic point of $T_+$, proving the necessity of item~(\ref{ItemVerticesPeriodic}) in the following:

\begin{proposition} \label{characterize train track} 
Let $T=\wt G$ be a weak train track and $f_T \from T \to T_+$ a $\Lambda_-$ isometry. The following are necessary and sufficient conditions for the existence of a train track representative $g \from G \to G$ of $\phi^i$ for some $i$, such that $f_T=f_g$:
\begin{enumerate}
\item \label{ItemVerticesPeriodic}
For every vertex $w$ of $T$, $f_T(w) \in T_+$ is a periodic point.  
\item \label{ItemLWGconsistent}
For every vertex $y$ of $T$, if $f_T(y) = b$ is a branch point of $T_+$ then there exists a principal vertex $w$ of $T$ such that $D f_T(W(y;T)) \subset D f_T(SW(w;T)) \subset W(b;T_+)$.
\end{enumerate}
\end{proposition}

\begin{proof} For the necessity of item~(\ref{ItemLWGconsistent}), choose $i \ge 1$ so that $g^i$ maps every vertex to a fixed point of $g^i$ and $D g^i$ maps every direction at a vertex to a fixed direction. Referring to the notation of Figure~\ref{FigureNotation2}, it follows that there is a principal vertex $w$ of $T=\wt G = \wt G_0$ such that $y$ and $w$ have the same image $w_i \in \wt G_i$ under $\ti g_{i,0}=\ti g_{i,i-1} \composed\cdots\composed \ti g_{1,0}$, and $D g_{i,0}(W(y;T)) \subset SW(w_i;\wt G_i)$. Since $D g^i$ is the identity on the set of fixed directions at fixed vertices of $g^i$, the map $D \ti g_{i,0}$ induces an isomorphism between $SW(w;T))$ and $SW(w_i;\wt G_i)$. Applying the ``chain rule'' $D f_T = D f_{g_0} = D f_{g_i} \composed D f_{g_{i,0}}$, it follows that $D f_T(W(y;T)) \subset D f_T(SW(w;T)) \subset W(b;T_+)$.

For sufficiency, let $\cal R$ be the set of rotationless train track maps $\gG$ representing an iterate of $\phi$ for which there exists an equivariant map $h \from \wt G \to T$ such that $h$ is simplicial (vertices map to vertices and edges map to edges) and $f_{g} = f_T h$.  

We first show that $\cal R \ne \emptyset$.  By Lemma~\ref{max split} there exists a rotationless train track map $\gG$ representing an iterate of $\phi$ such that the associated \trtr\ $f_{g} : \wt G \to T_+$ splits as least as much as $f_T:T \to T_+$.  As there are only finitely many indivisible pre-Nielsen paths for $g$, we may assume, after replacing $g$ with an iterate,  that if $\mu$ is a pre-Nielsen path for $g$, then $g_\#(\mu)$ is a Nielsen path for $g$.  Choose $\epsilon$ as in Proposition~\ref{PropWeakTrainTrack} applied to $f_{g} : \wt G \to T_+$ .  Replacing $G$ with $\phi^{-n}G$ has the effect of replacing the constant $\epsilon$ with $\lambda^n \epsilon$.  Choosing $n$ sufficiently large, we may therefore assume that $\Length(T) \le \epsilon$ and hence that there is an equivariant map $h :\wt G \to T$ such that $f_{\wt G} = f_Th$.  Up to the action of $F_r$, there are only finitely many points of $\wt G$ that map via $h$ to vertices of $T$. Suppose that $w$ is a vertex of $T$, that $\ti x$ is a point of $\wt G$ and that $h(\ti x) = w$.  Let $x$ be the image of $\ti x$ in $G$.  After replacing $\phi$ with an iterate,   (1) and Theorem~\ref{TheoremPhiPlus}(1) imply the existence of a lift $\widetilde{g(x)}$ of $g(x)$ such that $f_g(\widetilde{g(x)}) = f_g(\ti x)$.  Corollary~\ref{first tree corollary} implies that for some $m \ge 0$ either $g^m(x)$ is the endpoint of Nielsen path or $g^{m+1}(x) = g^m(x)$.  In either case $g^m(x) \in \Fix(g)$.   We may therefore subdivide along the $g$-orbit of  $x$ to make $x$ and hence $\ti x$ a vertex.  After finitely many such subdivisions, the map $h$ is simplicial. We have now shown that $\cal R \ne \emptyset$.
 
For each $\gG$ in $\cal R$ the  associated map $h:\wt G \to T$ descends to a simplicial homotopy equivalence  $\bar h : G \to H$ where the marked graph $H$ is the quotient space of $T$ by the action of $F_r$, and $\bar h$ preserves marking.  If $\bar h$ is a homeomorphism then $T = \wt G$ and $f_T=f_g$, and so we are done.  Otherwise, there exist edges $E_1$ and $E_2$ of $G$, oriented so as to have a common terminal endpoint, such that $E_1, E_2$ are identified by the map $\bar h$. Arguing by cases that depend on the nature of $E_1,E_2$, we will show that there is an element $g' : G' \to G'$ of $\cal R$ such that $G'$ has fewer edges than $G$.  The proof then concludes by induction.

There are two easy cases.  In case~1, $E_1 \bar E_2$ is an indivisible Nielsen path and Lemma~\ref{Nielsen collapse} produces the desired $g' :G' \to G'$.  In case~2, there exists $k \ge 1$ such that $g^k_\#(E_1) = g^k_\#(E_2)$, and replacing $g$ by $g^k$, Lemma~\ref{LemmaStallingsFold} produces the desired $g' :G'  \to G'$.

If neither case 1 nor case 2 applies to the pair $E_1,E_2$, then Corollary~\ref{first tree corollary} implies, after replacing $g$ with an iterate if necessary, that $\sigma := g_\#(E_1 \bar E_2)$ is an indivisible Nielsen path.  Lift $\sigma$ to a path $\ti \sigma \subset \wt G$ with endpoints $\ti v_1$ and $\ti v_2$.  There are lifts $\wt E_i$ with a common terminal endpoints and with initial endpoints $\ti x_i$ such that (in the notation of Figure 2) $\breve g_0(\ti x_i) = \ti v_i$.  Let $y = h(\ti x_1) = h(\ti x_2)$ and $w_i = h(\ti v_i)$.  In particular, the endpoints $v_1$ and $v_2$ of $\sigma$ are principal.  Lemma~\ref{principal essential} and Lemma~\ref{principal lifts and branch points} imply that  $f_T(y)  = f_g(\ti x_i)=f_g(\ti v_i) = f_T(w_i) = b\in T_+$ is a branch point.  By construction, $Df_g(W(\ti x_i;\wt G)) \subset Df_T(W(y;T))$ and  $Df_g(W(\ti x_i;\wt G)) \subset Df_g(SW(\ti v_i;\wt G)) \subset Df_T(SW(w_i;T))$.  
For each $i=1,2$, the graph $Df_T(W(y;T))$ therefore intersects $Df_T(SW(w_i;T))$ in at least one edge. Applying (2) and Lemma~\ref{gen partition} it follows that $Df_T(W(y;T)) \subset Df_T(SW(w_i;T))$, and so $Df_T(SW(w_1;T))$ and $Df_T(SW(w_2;T))$ intersect each other in at least one edge. Another application of Lemma~\ref{gen partition} implies that $w_1 = w_2$. In particular, letting $(E_1',E_2')$ be the illegal turn in $\sigma$, we have that $\bar h(E_1') = \bar h(E_2')$.  Replacing $E_1$ and $E_2$ with $E_1'$ and $E_2'$ we are in case 1 if $\sigma = E_1'\bar E_2'$ and in case 2 otherwise.  
\end{proof}

\begin{corollary}  
\label{CorollaryWeakIsDense}
For any nongeometric, fully irreducible $\phi$, a dense set of weak train tracks are \trtr s.
\end{corollary}

\begin{proof}   A point in $T$ is said to be \emph{preperiodic} if its image under $f_T$ is a periodic point in $T_+$. 

We define a complexity $C(T)$ for  weak train train tracks $T$ as follows. First define the \lq valence complexity\rq\ $VC(T)$ to be the non-increasing sequence of integers with one entry $k_i$ for each $F_r$-orbit of vertices  $w_i$ of  valence $k_i \ge 3$. Then define the \lq aperiodic complexity\rq\ $AC(T)$ to be the number of orbits of vertices in $T$ that are not preperiodic. Finally, define $C(T) := (VC(T),AC(T))$ ordered lexicographically. Thus $C(T') < C(T)$ if  $VC(T') < VC(T)$ or if $VC(T') = VC(T)$ and $AC(T') < AC(T)$.

It suffices to show that if $T$ is not a \trtr\ then there is an arbitrarily close weak train track $T'$ with $C(T') < C(T)$.   

Let $\gG$ be a rotationless irreducible train track map representing an iterate of $\phi$. In the tree $\wt G$ there is a dense set of points that project to periodic points in $G$, and so for any leaf $ \ti \ell$ of $\wt \Lambda_-$ realized in $T_+$ the set of periodic points is dense in $\ti \ell$. The set of preperiodic points in $T$ is therefore dense. 

Suppose that $w$ is a vertex of $T$ that is not preperiodic. Then $f_T(w)$ is not a branch point, and so the local Whitehead graph $W(f_T(w);T_+)$ is a single edge. Consider the map $D f_T \from W(w;T) \to W(f_T(w);T_+)$. If $w$ has valence three then there are directions $d_1 \ne d_2$ at $w$ which have the same image under $D f_T$. Folding these directions up to nearby points whose $f_T$-image is periodic produces a new valence three vertex whose $f_T$-image is preperiodic and reduces $w$ to valence two; we then remove $w$ from the list of vertices of $T$. This does not change $VC(T)$ and lowers $AC(T)$. If the valence of $w$ is higher than three then perform a similar small fold to reduce $VC(T)$.

We may therefore assume that each vertex of $T$ is preperiodic. By this same folding argument we may assume that if the image $Df_T(W(w;T))$ is not a single edge of $W(f_T(w);T_+)$ then the restricted map $Df_T \from W(w;T) \to W(f_T(w);T_+)$ is an injection. 

Proposition~\ref{characterize train track} implies that if $T$ is not a train track, then there is a vertex $y$ of $T$ such that $f_T(y) = b$ is a branch point  of $T_+$ but $D f_T(W(y;T)) \subset W(b;T_+)$ is not contained in $D f_T(SW(w;T))$ for any principal vertex $w$ of $T$. In particular, $D f_T(W(y;T))$ contains at least three vertices of $W(b;T_+)$, so $D f_T(W(y;T))$ is not a single edge of $W(f_T(w);T_+)$, and so $Df_T$ is injective on $W(y;T)$. Lemma~\ref{gen partition} implies that $W(y;T)$ has a cut vertex $d$, and so $W(y;T)$ is the union of two subcomplexes $X_1$ and $X_2$, each containing $d$ properly, such that $X_1 \intersect X_2 = \{d\}$. We now perform a small splitting very similar to the one in the proof of Lemma~\ref{split Nielsen}.  Let $E$ be the edge in the link of $y$ that corresponds to $d$ and let $P_i$ be the edges in the link of $y$ that correspond to the vertices of $X_i \setminus\{d\}$. After subdividing $T$ at preperiodic points, we may assume that the length of $E$ is arbitrarily small.  Define a new tree $T'$ from $T$  by equivariantly replacing $y$ with a pair of vertices $y_i$,  by  equivariantly  replacing $E$ with a pair of edges $E_i$ with initial vertex $y_i$ and terminal vertex equal to the terminal vertex of $E$ and by equivariantly reattaching the edges in $P'_i$ to $y_i$. The map $f_{T'}$ defined in the obvious way makes $T'$ a weak train track of lesser complexity. 
\end{proof}

\subsection{Rigidity and irrigidity of $\Lambda_-$ isometries}
\label{SectionIrrigidity}

In this section we study uniqueness of $\Lambda_-$ isometries. Given a marked graph $G \in \wh\X_r$ with universal cover $T = \wt G$ let $I(T,T_+)$, denote the space of   $\Lambda_-$ isometries $T \mapsto T_+$ equipped with the topology of pointwise convergence.   Assuming that $T$ is a weak train track, in other words $I(T,T_+)$ is nonempty, we wish to know whether this space has more than one point, and if it does then we wish to have some description of its structure. The answer, as it turns out, depends on orientability properties of $\Lambda_-$: if $\Lambda_-$ is nonorientable then $I(T,T_+)$ is a single point; and if $\Lambda_-$ is orientable then $I(T,T_+)$ is a compact interval, possibly degenerating to a point. 


Recall that $D^2 \bdy F_r$ is the set of ordered pairs of elements of $\bdy F_r$, that $\G F_r$ is the space of unordered pairs of elements of $\bdy F_r$ and that   $\Lambda_-$ is a  closed $F_r$-equivariant subset of $\G F_r$.  An {\em orientation of $\Lambda_-$}   is a lift of the inclusion map $\Lambda_- \inject \G F_r$    to a continuous, $F_r$-equivariant map    $o \from \Lambda_- \to D^2 \bdy F_r$.   Equivalently, an orientation of $\Lambda_-$ is a choice of orientation on each leaf of $\Lambda_-$ with the following consistency condition: for some, and hence every, $F_r$ tree $T$, there is a constant $C$ so that if $L^\vpp_T$ and $L_T'$ are the realizations in $T$ of oriented leaves $L$ and $L'$ of $\Lambda_-$ and if $L^\vpp_T \cap L_T'$ has length at least $C$ then the orientations of $L^\vpp_T$ and $L_T'$ agree on $L^\vpp_T \cap L_T'$.  If $\Lambda_-$ has an orientation then it is {\em orientable}.  Since $\Lambda_-$ is minimal, there are covering translations $\gamma \in F_r$ such that $L_T' \cap \gamma \cdot L^\vpp_T$ is arbitrarily long.  It follows that an orientation on $\Lambda_-$ is determined by an orientation of any of its leaves and hence that  $\Lambda_-$ has either zero or two orientations.   

An {\em orientation of $\Lambda_-$ in $T$} is a choice of orientation on each leaf of $\Lambda_-$ with the feature that the orientations on $L^\vpp_T$ and $L_T'$ agree on $L^\vpp_T \cap L_T'$ whenever $L^\vpp_T \cap L_T'$ contains a nontrivial interval. An orientation of $\Lambda_-$ in $T$ therefore induces a well defined $F_r$-equivariant orientation on each edge of $T$. Recall from Lemma~\ref{CorollaryLeavesInTPlus} that each leaf of $\Lambda_-$ has a well defined realization in $T_+$.  Thus this definition works equally well to define an {\em orientation of $\Lambda_-$ in $T_+$}.

It is obvious that if $\Lambda_-$ is orientable in some $T$ then it is orientable.  We have the following partial converse.

\begin{lemma} If $\Lambda_-$ is orientable then it is orientable in $T_+$ and is orientable in any weak train track $T$.
\end{lemma}

\begin{proof}   Let $g : G \to G$ be a train track representative of $\phi$, let $\wt G$ be the universal cover of $G$ and let $o$ be an orientation on $\Lambda_-$. Choose a lift $\breve g \from \wt G \to \wt G $ of $g$. Since $\hat g :\partial F_r \to \partial F_r$ is continuous and $F_r$-equivariant, it acts on the set of orientations of $\Lambda_-$ and this action is independent of the choice of lift $\breve g$.  After replacing $g$ with $g^2$ if necessary, we may assume that $\hat g$ preserves $o$. Given $L$ and $L'$ so that $L_{\wt G}^\vpp \cap L_{\wt G}'$ is a nontrivial interval, choose $n > 0$ so that $\breve g^n (L_{\wt G}') \cap \breve g^n(L_{\wt G}^\vpp)$ has length greater than the constant $C$ associated to $o$ in the definition of an orientation. The  orientations of $\breve g^n(L_{\wt G}')$ and $\breve g^n(L_{\wt G}^\vpp)$ agree on $\breve g^n (L_{\wt G}') \cap \breve g^n(L_{\wt G}^\vpp)$ which implies that the orientations of $L_{\wt G }'$ and $L_{\wt G}^\vpp$ agree on $L_{\wt G}' \cap L_{\wt G}^\vpp$. This proves that  $\Lambda_-$ is orientable in $\wt G$.  

If  $L^\vpp_{T_+} \cap L'_{T_+}$  is a non-trivial interval then there is a train track representative $g : G \to G$ with universal cover $\wt G$ so that $L_{\wt G}^\vpp$ and $L_{\wt G}'$ intersect in a non-trivial interval.  By our previous argument, the orientations of $L_{\wt G}^\vpp$ and of $L_{\wt G}'$ agree on $L_{\wt G}^\vpp\cap L_{\wt G}'$.  Since $f_g : \wt G \to T_+$ preserves orientation on the leaves of $\Lambda_-$, it follows that  the orientations of $L^\vpp_{T_+}$ and $L'_{T_+}$  agree on $L^\vpp_{T_+} \cap L'_{T_+}$.  This proves that $\Lambda_-$ is orientable  in $T_+$.  If $T$ is a weak train track then leaves of $\Lambda_-$ that intersect in $T$ also intersect in $T_+$.  An orientation of $\Lambda_-$ in $T_+$ therefore pulls back to an orientation on $\Lambda_-$ in $T$. 
\end{proof}
  
The next theorem gives an abstract description of the set $I(T,T_+)$. An orientation $o$ of $\Lambda_-$ induces orientations of $\Lambda_-$ in $T$ and in $T_+$ and each element $f \in I(T,T_+)$ preserves the induced orientations. It follows that for each vertex $v$ of $S$, if $e_1,e_2$ are two edges at $v$ in the same gate with respect to $f$ then $e_1,e_2$ point in the same direction at $v$: either both point towards $v$ in which case that gate is said to be a \emph{negative} gate at $v$; or both point away from $v$ in which case the gate is said to be \emph{positive}.

\begin{theorem}
\label{TreeTheorem}
Fix a weak train track $G \in \wh\X_r$ with universal cover $T=\wt G$, and fix $f \in I(T,T_+)$.  
\begin{enumerate}
\item If $\Lambda_-$ is not orientable in $G$ then $f$ is the unique element of $I(T,T_+)$.
\item If $\Lambda_-$ is orientable in $G$ then:
\begin{enumerate}
\item 
\label{ItemCompactParams}
There exists a compact subinterval $\Delta$ of the real line, and a homeomorphism $\Delta \mapsto I(T,T_+)$. We denote this homeomorphism $t \mapsto f_t$.
\item 
\label{ItemLeafParams}
For any leaf $L$ of $\Lambda_-$ and for any $x \in L_{T_+}$, the map $\Delta \mapsto T_+$ given by $t \mapsto f_t(x)$ is an orientation preserving isometry onto a segment of the leaf $L_{T_+}$.
\item
Given $f \in I(T,T_+)$, the following hold:
\begin{enumerate}
\item 
\label{ItemInterior}
$f$ is an interior point of $I(T,T_+)$ if and only if every vertex of $T$ has exactly one positive and one negative gate with respect to $f$.
\item 
\label{ItemEndpoints}
$f$ is the positive (resp.\ negative) endpoint of $I(T,T_+)$ if and only if some vertex of $T$ has more than one positive (resp.\ negative) gate with respect to $f$.
\item 
\label{ItemUnique}
$f$ is the unique point of $I(T,T_+)$ if and only if there exist vertices $v,w \in T$ such that $v$ has at least two positive gates with respect to $f$ and $w$ has at least two negative gates with respect to $f$.
\end{enumerate}
\end{enumerate}
\end{enumerate}
\end{theorem}

Even when $\Lambda_-$ is orientable, this lemma leaves open the possibility that $f$ is unique, that is, the interval $\Delta$ degenerates to a point. See below for examples where $\Delta$ is nondegenerate and where it is degenerate.

In item~(\ref{ItemCompactParams}) of this theorem, we use the topology of pointwise convergence on $I(T,T_+)$, but as should be clear from the proof this is equivalent to the Gromov topology on maps (for the definition of this topology see Section~\ref{SectionSkora}, proof of step 3(d)).

\begin{proof} Consider two $\Lambda_-$ isometric maps $f,f' \from T \to T_+$. For each oriented leaf $L$ of $\Lambda_-$, let $L_{T} \subset T$ and $L_{T_+} \subset T_+$ denote the realizations of $L$ as oriented, bi-infinite geodesics. Each of the restricted maps $f, f' \from L_{T} \to L_{T_+}$ is an orientation preserving isometry, and so these two maps differ by a real number $\delta(f,f';L)$, meaning that the map $f' \from L_{T} \to L_{T_+}$ is obtained by postcomposing the map $f \from L_{T} \to L_{T_+}$ by the translation of $L_{T_+}$ that displaces each point by the amount $\delta(f,f';L)$ (remember that $L_{T_+}$ is oriented). 

Next we study how $\delta(f,f';L)$ depends on $L$. Consider another oriented leaf $L'$ of $\Lambda_-$. If $L_{T} \intersect L'_{T} = J$ is a nondegenerate segment then the restriction to $J$ of the map $f \from L_{T} \to L_{T_+}$ equals the restriction to $J$ of the map $f \from L'_{T} \to L'_{T_+}$, and the same is true for $f'$. It follows that if $L_{T}$ and $L'_{T}$ induce the same orientation on $J$ then $\delta(f,f';L)=\delta(f,f';L')$, whereas if $L_{T}$ and $L'_{T}$ induce opposite orientations on $J$ then $\delta(f,f';L) = - \delta(f,f';L')$. Note also that $\delta(f,f';L)$ is invariant under the action of $F_r$ on oriented leaves. 

Suppose now that $\Lambda_-$ is nonorientable in $T$. Then for any oriented leaf $L$ of $\Lambda_-$ there exists a covering transformation $g \in F_r$ such that $L_{T} \intersect g \cdot L_{T}$ is a nondegenerate segment $J$ and the orientations induced on $J$ by $L_{T}$ and $g \cdot L_{T}$ are opposite. It follows that 
$$\delta(f,f';L) = \delta(f,f';g \cdot L) = -\delta(f,f';L)
$$
and so
$$\delta(f,f';L) = 0
$$
implying that $f \restrict L_{T} = f' \restrict L_{T}$. Since this is true for every leaf $L$, and since realizations of leaves of $\Lambda_-$ cover $T$, it follows that $f=f'$.

Suppose now that $\Lambda_-$ is orientable in $T$, and fix an orientation. If $L,L'$ are any two oriented leaves then there exists $g \in F_r$ such that $g \cdot L_{T} \intersect L'_{T} = J$ is a nondegenerate segment and the orientations on $J$ induced by $g \cdot L_{T}$ and $L'_{T}$ agree. it follows that 
$$\delta(f,f';L') = \delta(f,f';g \cdot L) = \delta(f,f';L)
$$
This shows that $\delta(f,f';L)$ is a constant independent of $L$, which we denote $\delta(f,f')$. The map $f'$ is therefore determined by $f$ and the constant $\delta(f,f')$: for each $x \in T$, pick any leaf $L$ whose realization in $T$ contains $x$, and then $f'(x)$ is obtained by displacing $f(x)$ along $L_{T_+}$ by the amount $\delta(f,f')$. 

Fix $f \in I(T,T_+)$.  Let $\Delta$ be the set of all real numbers of the form $\delta(f,f')$ for some $f' \in I(T,T_+)$.  Each $d \in \Delta$ equals  $\delta(f,f')$ for a unique $f' \in I(T,T_+)$ and we define   $f_d = f'$. To complete the proof of (\ref{ItemCompactParams}) we must show that $\Delta$ is a closed, bounded interval; continuity follows by noticing that for any $d,d' \in \Delta$ and any $x \in T$, the distance in $T$ between $f_d(x)$ and $f_{d'}(x)$ is exactly $\abs{d-d'}$. Item~(\ref{ItemLeafParams}) follows from item~(\ref{ItemCompactParams}) and the definition of $\delta(f,f')$.

To show that $\Delta$ is connected, given a $\Lambda_-$ isometry $f' \from T \to T_+$, and given a number $d$ between $0$ and $\delta(f,f')$ it suffices to exhibit a $\Lambda_-$-isometry $f'' \from T \to T_+$ such that $\delta(f,f'')=d$. For $x \in T$ the segment $\overline{f(x) f'(x)}$ in $T_+$ has length $\delta(f,f')$, and this segment is oriented with initial point $f(x)$ when $\delta(f,f') >0$ and terminal point $f(x)$ when $\delta(f,f') < 0$. Now define $f''(x)$ to be the point on this segment at displacement $d$ from $f(x)$. 

Boundedness of $\Delta$ follows from the observation that there exist leaves $L,L' \subset \Lambda_-$ such that $L^\vpp_{T} \intersect L'_{T} \ne \emptyset$ and $L^\vpp_{T_+} \intersect L'_{T_+}$ is bounded. 

The fact that $\Delta$ is closed follows by application of the Ascoli-Arzela theorem together with the observation that the compact-open limit of a sequence of $\Lambda_-$ isometries $T \mapsto T_+$ is a $\Lambda_-$ isometry, and the fact that on $I(T,T_+)$ the compact open topology agrees with the pointwise convergence topology.

Item~\ref{ItemUnique} follows from items~\ref{ItemInterior} and~\ref{ItemEndpoints}, to whose proofs we now turn. Suppose first that every vertex $v \in T$ has a unique positive gate with respect to~$f$, so there is a positive length oriented segment $J_v$ in $T_+$ with initial endpoint $f(v)$ such that every edge in $T$ with initial endpoint $v$ has an initial subsegment mapping isometrically to $J_v$. There being only finitely many orbits of vertices, we can truncate to assume that all the segments $J_v$ have the same length $\epsilon>0$. We can then alter $f$ to a new $\Lambda_-$ isometry, by moving each vertex $v$ along $J_v$ to the positive endpoint of $J_v$, and extending by an isometry on each edge of $T$. Applying item~(\ref{ItemLeafParams}) it follows that $f$ is not the positive endpoint of $I(T,T_+)$. 

Conversely, suppose that there exists a vertex $v \in T$ having two or more positive gates, and so there exist edges $e,e'$ of $T$ each with initial endpoint $v$ such that $f(e)$, $f(e')$ have distinct germs at $f(v)$. Choose leaves $L,L' \subset \wt\Lambda_-$ whose realizations in $T$ contain $e,e'$, and hence each passes through $v$. The realizations of $L,L'$ in $T_+$, while passing through $f(v)$, are disjoint in the positive direction beyond $f(v)$. By item~(\ref{ItemLeafParams}) it follows that $f$ is the positive endpoint of $I(T,T_+)$. Combining this argument with that of the last paragraph, and repeating these arguments with orientation reversed, items~\ref{ItemInterior} and~\ref{ItemEndpoints} follow.
\end{proof}

\subsection{Examples of exceptional weak train tracks}
\label{SectionWTTexamples}
We describe first a weak train track $T$ for which the space of $\Lambda_-$ isometries $I(T,T_+)$ is not just a single point. Then we describe a weak train track $T$ that is not a train track.

\paragraph{Example: a weak train track $T$ with nontrivial $I(T,T_+)$.}  Let $a,b$ and $c$ be generators of $F_3$. Consider $\phi \in \Out(F_3)$ represented by the following positive automorphism $\Phi$, shown together with its inverse:
$$ \Phi \from \begin{cases}
a &\mapsto bacaaca \\
b &\mapsto baca \\ 
c &\mapsto caaca
\end{cases}
\qquad
\Phi^\inv \from \begin{cases}
a &\mapsto \bar b a \bar b a \bar c \\
b &\mapsto b \bar a b \\
c &\mapsto c \bar a b c \bar a b \bar a b
\end{cases}
$$
As realized on the three petalled rose, these formulas each define train track maps representing $\phi$ and $\phi^\inv$, respectively. The automorphism $\Phi$ being positive, its transition matrix equals the matrix of its abelianization. Direct calculation shows that no eigenvalue of the abelianization of $\phi$ is a root of unity, so $\phi$ is fully irreducible by Lemma~\ref{eigenvalue 1}. Direct calculation of the Perron-Frobenius eigenvalues of the transition matrices of $\Phi$ and $\Phi^\inv$ also show that $\lambda(\phi) > 6$ and $\lambda(\phi^\inv) < 5$, so $\phi$ is not geometric.

Let $g \from G \to G$ be the affine train track representative of $\phi$ determined by the above formula for $\Phi$. Let $\lambda$ be the expansion factor, $\Lambda_-$ the expanding lamination, and $T_+$ the attracting tree for $\phi$. When the three petalled rose $G$ is oriented so that its three oriented edges represent $a,b$ and $c$ then it induces an orientation on  $\Lambda_-$ in $T = \wt G$. It is clear that at the unique vertex of $G$ the map $g$ has exactly two positive gates $\{a,b\}$ and $\{c\}$ and one negative gate $\{\bar a, \bar b, \bar c\}$. Since each of these gates contains a fixed direction, these are also the gates for $f_g \from \wt T \to T_+$. Theorem~\ref{TreeTheorem} item~\ref{ItemEndpoints} implies that $f_g$ is the positive endpoint of  $I(T,T_+)$ but is not the negative endpoint.
 
Now we shall slide the map $g$ backwards to locate the negative endpoint of $I(T,T_+)$. Looking at the formula for $\Phi$ we see that the $g$ image of each edge ends with the path $aca$. Postcomposing $\Phi$ with the inner automorphism $w \mapsto (aca)w(aca)^\inv$ we obtain
$$\Phi' \from \begin{cases}
a &\mapsto acabaca \\
b &\mapsto acab \\ 
c &\mapsto acaca
\end{cases}
$$
This formula defines a train track representative $g' \from G \to G$ of $\phi$ supported on the \emph{same marked graph $G$} as the train track representative $g$.  Notice now that $g'$ has one positive gate $\{a,b,c\}$ and two negative gates $\{\bar a,\bar c\}$ and $\{\bar b\}$ at its unique vertex.  As before each gate contains a fixed direction so these are also the gates of  $f_{g'}$.   Thus $f_{g'}$ is the negative endpoint of $I(T,T_+)$.

\paragraph{Example: a weak train track that is not a train track.} We describe an example of a nongeometric, fully irreducible $\phi \in \Out(F_r)$ and a weak train track $H$ for $\phi$ that is not a train track for any iterate of $\phi$. We shall also sketch the proof that $H$ is not even the domain of a nonsimplicial train track representative of any iterate of $\phi$.

Start with the graph $G$ and train track map $g \from G \to G$ given by
\begin{equation*}
\xymatrix{
q \ar@(ul,dl)[]_E     \ar@/^2pc/[rr]^{A} \ar@/_2pc/[rr]_{D}  & & p \ar[ll]^F
}
\qquad\qquad g : 
\begin{aligned}
A &\mapsto AF \\
D &\mapsto DFAFDFE \\
E &\mapsto DFAFDFE \\
F &\mapsto DF
\end{aligned}
\end{equation*}
One easily shows by the method of Stallings folds that $g$ is a homotopy equivalence. 
The only illegal turn of $g$ being $\{E,D\}$, and the words $\overline E D$ and $\overline D E$ not occurring in the definition of $g$, it follows that $g$ is a train track map (of course, any map of a directed graph which preserves direction is a train track map).
  
Let $g$ represent $\phi \in \Out(F_3)$ with respect to any marking of $G$. Lemma~\ref{eigenvalue 1} applies to show that $\phi$ is fully irreducible, and by computing some train track $\bar g$ for $\phi^\inv$ one shows that the expansion factors of $g$ and $\bar g$ are not equal and hence $\phi$ is not geometric (in fact this outer automorphism is conjugate to the one in Section~\ref{SectionIWGExample}, as one can show by computing and collapsing a Nielsen path of $g$ to get the train track map in Section~\ref{SectionIWGExample}).

Next we compute the local Whitehead graphs at $p$ and $q$. The turns taken by the images of edges under $g$ are
$$\{\overline A, F\}, \{\overline D, F\}, \{\overline F, A\}, \{\overline F, D\}, \{\overline F, E\}
$$
The map $Dg$ acts as follows: it takes $\{\overline A, F\}$ and $\{\overline F,E\}$ to $\{\overline F, D\}$ and fixes the latter; it takes $\{\overline D, F \}$ to $\{\overline E, D\}$ and fixes the latter; and it fixes $\{\overline F, A\}$. 
These are all the turns taken by iterates of edges under $g$, so the local Whitehead graphs are
$$\xymatrix{
\text{At $p$:}     & \overline A \ar@{-}[r]    
                          & F \ar@{-}[r]    
                          & \overline D \\
\text{At $q$:}     & A \ar@{-}[r] 
                          & \overline F \ar@{-}[r] \ar@{-}[d] 
                          & D \ar@{-}[r] 
                          & \overline E \\
                          &  & E
}$$

\medskip\noindent
Note that $D$ is a cut vertex of the local Whitehead graph at $q$. We construct $H$ by a small splitting of $G$ at $q$ in the direction $D$, as was done in Lemma~\ref{split Nielsen} and Corollary~\ref{CorollaryWeakIsDense}. The effect of this is to split the vertex $q$ into two vertices $q_1,q_2$. An initial segment of the edge $D$ is split into two new edges $D_1,D_2$, with the tail of $D_i$ incident to $q_i$, and with the heads of $D_1,D_2$ incident to a new valence~3 vertex $r$. The remaining terminal segment of $D$ becomes an edge $D_3$ from $r$ to $p$. The other four ends $A$, $E$, $\overline E$, $\overline F$ of $G$ at $q$ are detached from $q$ and reattached as follows: $A$, $E$, and $\overline F$ are reattached to $q_1$, and $\overline E$ is reattached to~$q_2$.
$$
\xymatrix{
q_1 \ar@/_2pc/[d]_E \ar[dr]^{D_1} \ar@/^2pc/[rr]^A & & p \ar[ll]_F\\
q_2 \ar[r]_{D_2} & r \ar[ur]_{D_3} & 
}
$$
To be precise about metric and marking, the edges lengths of $A,E,F$ in $H$ are the same as in $G$; $\Length_H(D_1)=\Length_H(D_2)=\epsilon>0$; and $\Length_H(D_3)=\Length_G(D)-\epsilon$. Let $h \from H \to G$ be the edge isometry that is the ``identity'' in $A,E,F$, and that maps $D_i * D_3$ isometrically to $D$ for $i=1,2$. The marking of $G$ composed with a homotopy inverse of $h$ defines the marking of $H$.

Consider the $\Lambda_-$ isometry $f_g \from \wt G \to T_+$, and let $\ti h \from \wt H \to \wt G$ be an $F_3$ equivariant lift of $h$. To show that $H$ is a weak train track we show that $f_h=f_g \composed \ti h \from \wt H \to T_+$ is a $\Lambda_-$ isometry, and to do this we must check that each leaf $L$ of $\Lambda_-$ realized in $G$ lifts isometrically via $\ti h$ to the realization of $L$ in $\wt H$. This follows immediately from the observation that in the local Whitehead graph of $G$ at $q$, the four edges $A \overline F$, $E \overline F$, $D \overline F$, $D \overline E$ lifts to edges $A \overline F$, $E \overline F$, $D_1 \overline F$ at $q_1$ and $D_2 \overline E$ at $q_2$.

We prove that $f_h$ is the unique $\Lambda_-$ isometry from $\wt H \to T_+$, by applying Theorem~\ref{TreeTheorem}. Since $\Lambda_-$ is oriented we must check condition~\ref{ItemUnique} in the theorem. Note that the directions $\overline A$ and $\overline D$ of $G$ at $p$ are in distinct negative gates, because $\overline A$ and $\overline D$ are mapped by $Dg$ to $\overline E$ and $\overline F$ at $q$, and the latter are fixed by $Dg$; it follows that any lift to $\wt H$ of the vertex $p$ of $H$ has two negative gates, containing lifts of the directions $\overline A$ and $\overline D$ respectively. Note also that in $G$ the directions $A$ and $D$ at $q$ are in distinct positive gates, because they are fixed by $g$; it follows that any lift to $\wt H$ of $q_1$ has two positive gates, containing lifts of the directions $A$ and $D_1$ respectively.

To prove that $H$ is not a train track of any iterate of $\phi$, by uniqueness of the $\Lambda_-$ isometry $f_h$ it suffices to prove that $f_h$ does not satisfy item~(\ref{ItemLWGconsistent}) of Proposition~\ref{characterize train track}. Let $\ti p \in \wt H$ be a lift of $p \in H$, and so $\ti h(\ti p)$ is a lift of $p \in G$. We know that there is a lift $\wt q \in G$ of $q \in G$ such that $b=f_g(\ti h(\ti p)) = f_g(\ti q)$ is a branch point of $T_+$ and $D f_g$ maps the four stable directions of $\wt G$ at $\ti q$, named $\{A,\bar F,D,\bar E\}$ to the four directions of $T_+$ at $b$, which shall have the same names. We also know that $f_g$ maps $\ti h(\ti p)$ taking the three directions $\{\bar A,F,\bar D\}$ at $\ti h(\ti p)$ to the three directions $\{\bar F,D,\bar E\}$ at $b$. It follows that $D f_h$ maps the three directions $\{\bar A,F,\bar D\}$ at $\ti p$ to the three directions $\{\bar F,D,\bar E\}$ at $b$. The only principal vertices in $\wt H$ that map to $b$ via $f_h$ are the lifts $\ti q_1,\ti q_2$ of $q_1,q_2$ such that $\ti h(\ti q_1)=\ti h(\ti q_2) = \ti q$. However, the image under $D f_h$ of the set of stable directions of $\wt H$ at $\ti q_i$ is the set $\{A, \bar F, D\}$ at $b$ when $i=1$, and the set $\{D, \bar E\}$ when $i=2$, and neither of these two sets contains $\{\bar F, D, \bar E\}$, as is required by  item~(\ref{ItemLWGconsistent}) of Proposition~\ref{characterize train track}.

We also sketch the proof that $H$ is not the domain of any \emph{nonsimplicial train track representative} of any iterate of $\phi$, meaning a homotopy equivalence $g \from H \to H$ that represents $\phi$ such that for each edge $e$ of $H$ and each $i \ge 1$ the map $g \restrict \interior(e)$ is a local homothety with stretch factor $\lambda(\phi)^i$. To prove this requires generalizing several results about train track representatives to the context of nonsimplicial train track representative. Suppose that there is a nonsimplicial train track representative $g \from H \to H$ of some iterate of $\phi$. We can lift $g$ to define maps as denoted in Figure~\ref{FigureNotation2} (replacing the symbol $G$ by $H$), and then Theorem~\ref{T+ is direct limit} generalizes to produce $T_+$ as a direct limit, and to show that the unique $\Lambda_-$ isometry $f_h \from \wt H \to T_+$ is the direct limit map denoted by $f_g$ in Figure~\ref{FigureNotation2}. Now the argument given in Proposition~\ref{characterize train track} for the necessity of item~(\ref{ItemLWGconsistent}) generalizes in a straightforward manner. Finally, we have already proved that $f_h$ fails to satisfy item~(\ref{ItemLWGconsistent}).

\section{Topology of the axis bundle}
\label{SectionTopology}

For this section fix a nongeometric, fully irreducible $\phi \in \Out(F_r)$ with $\Lambda_-=\Lambda_-(\phi)$ and $T_+=T_+(\phi)$, and fix a normalization of $T_+$. Having so far proved the equivalence of two of the three proposed definitions of the axis bundle --- the set of weak train tracks (2), and the closure of the set of train tracks (3) --- in this section we work with definition~(2), in terms of which we prove a version of Theorem~\ref{TheoremProperness}. 

Let $\A \subset \X_r$ be the set of all weak train tracks of $\phi$. By Lemma~\ref{LemmaNormalization}, each $T \in \A$ has a normalization $\wh T \in \wh\X_r$ well-defined by requiring existence of a $\Lambda_-$ isometry $\wh T \mapsto T_+$, and so $\Length(T) = \Length(\wh T/F_r)\in (0,+\infinity)$ is well-defined. This map satisfies the equivariance condition 
$$\Length(T \cdot \phi^i) = \lambda(\phi)^{-i} \Length(T) \quad\text{for any}\quad i \in \Z
$$ 
because from any $\Lambda_-$ isometry $\wh T \mapsto T_+$ we obtain a $\Lambda_-$ isometry $\wh T \cdot \phi^i \mapsto T_+ \cdot \phi^i = \lambda(\phi)^i T_+$ which, after renormalizing, becomes a $\Lambda_-$ isometry $(\lambda(\phi)^{-i} \wh T) \cdot \phi^i \to T_+$.

Let $\wh\A \subset \wh\X_r$ be the set of normalized weak train tracks; projection $\wh\X_r \mapsto \X_r$ restricts to a bijection $\wh\A\mapsto\A$. Abusing notation we often use the same symbol to represent a weak train track and its normalization, although we will the symbol $\,\wh{}\,$ to maintain the distinction when needed.  

Here is the main result of this section, a major step in the proof of Theorem~\ref{TheoremProperness}:

\begin{theorem}
\label{TheoremLengthPHE}
The map $\Length \from \A \to (0,+\infinity)$ is a proper homotopy equivalence.
\end{theorem}

\paragraph{Outline.} Section~\ref{SectionLengthContinuous} proves basic topological facts about $\A$, including continuity of the normalization map $\A \mapsto \wh\A$ and the resulting continuity of the length map $\Length \from \A \to (0,+\infinity)$. Section~\ref{SectionGHConvergence} proves a technical result, strengthening the topological properties of $\wh\A$, that will be useful in the proof of Theorem~\ref{TheoremLengthPHE}. Section~\ref{SectionProperness} proves that $\Length \from \A \to (0,+\infinity)$ is a proper map, using folding methods combined with Proposition~\ref{PropWeakTrainTrack}. Section~\ref{SectionSkora} proves Theorem~\ref{TheoremLengthPHE} using Skora's method \cite{Skora:deformations} for studying homotopy types of spaces of trees (see also \cite{White:FixedPoints}, \cite{Clay:contractibility}, \cite{GuirardelLevitt:outer}).

\subsection{Continuity properties of the normalized axis bundle}
\label{SectionLengthContinuous}

We prove here several facts about the set $\wh\A$ and related facts about the length map.

\begin{lemma}
\label{LemmaLengthFacts} 
\quad
\begin{enumerate}
\item 
\label{NormalizedWTTClosed}
$\A$ is a closed subset of $\X_r$ and $\wh\A$ is a closed subset of $\wh\X_r$.
\item 
\label{ItemNormalizingContinuous}
The normalization map $\A \to \wh\A \subset \wh\X_r$ is a homeomorphism. 
\item 
\label{ItemLengthProps}
The map $\Length \from \A \to (0,+\infinity)$ is continuous, surjective, and $\phi$-equivariant, where $\phi$ acts on $(0,+\infinity)$ by multiplying by $\frac{1}{\lambda}$. 
\item 
\label{ItemLogLengthProps}
The map $-\log\Length \from \A \to \reals$ is also continuous, surjective and $\phi$-equivariant, where $\phi$ acts on $\reals$ by adding $\log\lambda$.
\end{enumerate}
\end{lemma}

\begin{proof} Note that (\ref{ItemLogLengthProps}) follows from (\ref{ItemLengthProps}). 
Also, in proving (\ref{ItemLengthProps}), continuity of $\Length$ follows from (\ref{ItemNormalizingContinuous}) and equivariance follows from Theorem~\ref{TheoremPhiPlus}. For surjectivity, choose a train track representative $g \from G \to G$, let $T_0 = \wt G$, let $T_1 = T_0 \cdot \phi$, lift $g$ to an $F_r$-equivariant $\Lambda_-$ isometry $f \from T_0 \to T_1$, and interpolate $f$ by a Stallings fold path from $T_0$ to $T_1$. For each $T$ on this path there are equivariant edge isometries $T_0 \xrightarrow{f'} T \xrightarrow{f''} T_1$ that factor $f$, and so $f''$ followed by a $\Lambda_-$ isometry $T_1 \to T_+$ is a $\Lambda_-$ isometry $T \to T_+$, proving that the fold path is in $\wh\A$. Length varies continuously and monotonically along this path from $T_0$ to $T_1$ with image $[\Length(G),\Length(G)/\lambda]$. Surjectivity of $\Length\from\wh\A\to(0,+\infinity)$ then follows from $\phi$-equivariance.

Consider a sequence $G_i \in \A$ which converges in $\X_r$ to some limit $G$. We first show that the normalized sequence $\wh G_i$ converges to some limit $\wh G \in \X_r$, and second that $\wh G \in \wh A$. This proves~(\ref{NormalizedWTTClosed}). It also proves that the normalization map $\A \to \wh\A$ is continuous, but its inverse is the restriction of the continuous projection map $\wh\X_r\mapsto\X_r$, proving~(\ref{ItemNormalizingContinuous}) as well. Let $T,T_i$ denote the universal covers of $\wh G,\wh G_i$.

Cells in $\X_r$ being locally finite, by restricting to each of the finitely many closed cells containing $G$ we may assume that all the $G_i$ are in a single open cell,  so there is a family of marked homeomorphisms $g_{ji} \from G_i \to G_j$ such that $g_{kj} \composed g_{ji} = g_{ki}$, well-defined by requiring each $g_{ji}$ to be affine on each edge of $G_i$.  Lift $g_{ji}$ to an $F_r$-equivariant homeomorphism $\ti g_{ji} \from T_i \to T_j$ which is affine on each edge.

Convergence in $\wh\A$ is defined by convergence of translation length functions $\ell_\gamma$ for all $\gamma \in F_r-\{\Id\}$, whereas convergence in $\A$ is defined by convergence of ratios $\ell_\gamma/\ell_{\gamma'}$ for all $\gamma,\gamma' \in F_r-\{\Id\}$. Since $G_i$ converges to $G$ in $\A$, to prove convergence of $\wh G_i$ it suffices to find a single $\gamma \in F_r - \{\Id\}$ such that $\ell_\gamma(T_i)$ is constant. Pick any leaf $L$ of $\Lambda_-$ and any $i$. There is a directed edge $e$ of $G_i$ that is traversed twice in the same direction by the realization of $L$ in $G_i$. Surgery of $L$ along $e$, as in the proof of Lemma~\ref{LemmaNormalization}, produces a legal loop in $G_i$, which maps to a legal loop in each $G_j$ under $g_{ji}$. Choose $\gamma \in F_r$ representing this loop. The axes of $\gamma$ in $T_i$ all correspond under the maps $\ti g_{ji}$, and these axes are all legal, and so $\ell_\gamma(T_i)=\ell_\gamma(T_+)$, completing the proof that $\wh G_i$ converges to some $\wh G \in \X_r$. 

Choose $\Lambda_-$ isometries $f^+_i \from T_i \to T_+$. We show that, in some sense, a subsequence of $f^+_i$ converges to a $\Lambda_-$ isometry $T \mapsto T_+$, proving that $\wh G \in \wh\A$. 

Since $G$ is in the closure of the open cell of $\X_r$ containing all of the $G_i$, for each~$i$ there exists a unique marked homotopy equivalence $p^i \from \wh G_i \to \wh G$ which collapses each component of a subforest to a vertex, and which is affine on each noncollapsed edge. This subforest is trivial if and only if $G$ is in the  same open cell as all of the $G_i$. Note that $p^j \composed g_{ji} = p^i$. 

Choose a finite set $A$ and an indexing $\{e^\alpha_i\}_{\alpha \in A}$ of the edges of $\wh G_i$ so that $g_{ji}(e^\alpha_i)=e^\alpha_j$. Let $B \subset A$ correspond to edges that are not collapsed by each $p^i$. Index the edges of $\wh G$ as $\{e^\beta\}_{\beta \in B}$ so that $p^i(e^\beta_i)=e^\beta$. The edge length function of $\wh G_i$, denoted $\ell_i \from A \to (0,+\infinity)$, converges to a function $\ell \from A \to [0,\infinity)$ which is nonzero precisely on $B$. By restricting to $B$, $\ell$ is identified with the edge length function for $\wh G$.

Choose a marked homotopy equivalence $q^i \from \wh G \to \wh G_i$ as follows. 
For each vertex $v \in \wh G$, consider the subtree $\tau_i(v) = (p^i)^\inv(v) \subset \wh G_i$, and choose $q^i(v) \in \tau_i(v)$ arbitrarily. For each edge $e^\beta$ of $\wh G$ with endpoints $v,w$, there is a unique continuous extension of $q^i$ to $e^\beta$ which is affine on $e^\beta$, whose image is contained in $e^\beta_i \union \tau_i(v) \union \tau_i(w)$, and such that $p^i \composed q^i \restrict e^\beta$ is homotopic rel endpoints to the identity map. This defines $q^i \from \wh G \to \wh G_i$ which is clearly a marked homotopy inverse for $p^i$. Denote by $\ti p^i, \ti q^i$ the lifts of $p^i,q^i$ to $F_r$-equivariant maps of universal covers.

Let $\epsilon_i$ denote the maximal diameter in $\wh G_i$ of the trees $\tau_i(v)$, and so $\epsilon_i \to 0$ as $i \to \infinity$. Also let $r_i \ge 1$ denote the maximum value of the ratios
$$\frac{\Length e^\beta_{\vphantom i}}{\Length e^\beta_i}, \frac{\Length e^\beta_i}{\Length e^\beta_{\vphantom i}}, 
\frac{\Length e^\beta_{\vphantom i}}{\Length e^\beta_i+2\epsilon_i}, 
\frac{\Length e^\beta_i+2\epsilon_i}{\Length e^\beta_{\vphantom i}}
$$
as $\beta$ varies over $B$, so $r_i \to 1$ as $i \to \infinity$, and $r_i$ is an upper bound for the amount by which either of the maps $p^i$, $q^i$ stretches or compresses an edge with index in $B$; of course $p^i$ compresses by $+\infinity$ any edge with index in $A-B$.

Next we show that for each vertex $v$ of $T$, the points $f^+_i \ti q^i(v)$ vary over a compact subset of $T_+$. Let $\{E_1,\ldots,E_K\}$ be the edges of the tree $\tau_i(v)$. For each $E_k$ choose leaves $L_k,L'_k$ of $\Lambda_-$ which run through $E_k$ and which are not asymptotic in either direction, and so the realizations of $L_k$ and $L'_k$ in $T_+$ intersect in a compact interval $I_k$ independent of $i$. It follows that $f^+_i \ti q^i(v) \in I_1 \union \cdots \union I_K$, the latter being compact.

Since there are only finitely many orbits of vertices we can pass to a subsequence so that $f^+_i \ti q^i(v)$ converges to a point denoted $f(v)$. Extend $f$ affinely to an $F_r$-equivariant map $f \from T \to T_+$. Since $r_i \to 1$ it follows that $f$ is an edge isometry. Given a leaf $L$ of $\Lambda_-$ realized in $T$, with realizations in $T_i, T_+$ denoted $L_i, L_+$, note that $\ti p^i(L_i)=L$, and that $\ti q^i(L)$ is contained in the $\epsilon_i$ neighborhood of $L_i$. Since $f^+_i \restrict L_i$ is an isometry onto $L_+$, the map $f^+_i \ti q^i \restrict L$ has image contained in the $\epsilon_i$ neighborhood of $L_+$, and so $f(L) = L_+$. The map $f \from L \to L_+$ restricts to an isometry on each edge of $T$ crossed by $L$, and so this map will be an isometry once we prove that it is injective. The map $f^+_i \ti q^i \restrict L$ can only identify a pair of points $x,y \in L$ if there is a vertex $v$ such that $\ti q^i\left([x,v] \union [v,y]\right) \subset \tau_i(v)$, and so
$$\sup_{x,y \in L, \,\, f^i \ti q^i(x)=f^i \ti q^i(y)} d_T(x,y) \le 2 r_i \epsilon_i \to 0 \quad\text{as}\quad i \to \infinity
$$
The map $f \from L \to L_+$ is therefore an isometry. Since $L$ is an arbitrary leaf of $\Lambda_-$, we have proved that $f$ is a $\Lambda_-$ isometry, and so $\wh G \in \wh\A$.
\end{proof}

\subsection{The Gromov topology on weak train tracks}
\label{SectionGHConvergence}

The proof of Lemma~\ref{LemmaLengthFacts} item~(\ref{ItemNormalizingContinuous}) can be interpreted in terms of the Gromov topology on $\wh\X_r$: for any $D \ge 0$ and $\epsilon>0$, if $i$ is sufficiently large then the restriction of $\ti q_i$ to any subtree of diameter $\le D$ is an $\epsilon$-almost isometry; we verify this below. In fact, this proves one of the easy directions of Paulin's Theorem \cite{Paulin:GromovTopology}, namely that if $\wh G_i$ converges to $\wh G$ in the cellular topology then $\wh G_i$ converges to $\wh G$ in the Gromov topology. For use in the proof of Theorem~\ref{TheoremLengthPHE}, in Lemma~\ref{LemmaRelation} we present a version of this argument tailored specifically to weak train tracks, couched in terms of a relation among weak train tracks called the ``$\Lambda_-$ relation''. This result will be applied in the verification of step~3(d) in Section~\ref{SectionSkora}. 

Given a weak train track $T \in \wh\A$, from among the set $I(T,T_+)$ of all $\Lambda_-$ isometries $T \mapsto T_+$ we pick out one, denoted $k^+_T$, as follows. If $\Lambda_-$ is nonorientable then Theorem~\ref{TreeTheorem} says that $I(T,T_+)$ has a unique element, taken to be $k^+_T$. If $\Lambda_-$ is orientable then, fixing an orientation, Theorem~\ref{TreeTheorem} says that $I(T,T_+)$ is a point or a compact oriented arc, and we take $k^+_T$ to be the positive endpoint of $I(T,T_+)$, characterized by $T$ having a vertex with two or more positive gates with respect to $k^+_T$. Choosing the positive rather than the negative endpoint of $I(T,T_+)$ is arbitrary, the point being to make a choice which is sufficiently canonical to force the appropriate continuity property --- namely, Lemma~\ref{LemmaRelation} --- to be true.

Given two weak train tracks $T,T' \in \wh\A$, the \emph{$\Lambda_-$ relation} between $T$ and $T'$ is the set $R \subset T \cross T'$ consisting of all $(x,x')$ with the following two properties:
\begin{enumerate}
\item $k^+_T(x) = k^+_{T'}(x')$
\item There exists a leaf $L$ of $\Lambda_-$ such that $x$ is contained in the realization of $L$ in $T$ and $x'$ is contained in the realization of $L$ in $T'$.
\end{enumerate}
The following lemma says that the $\Lambda_-$ relation can be used to exhibit closeness among weak train tracks in the Gromov topology:

\begin{lemma}
\label{LemmaRelation}
For any $T \in \wh\A$, $\epsilon>0$, any finite subset $t \subset T$, and any finite subset $P \subset F_r$, there exists a neighborhood $U \subset \wh\A$ of $T$ such that for any $T' \in U$, if $R$ denotes the $\Lambda_-$ relation between $T$ and $T'$, and if we define a finite subset 
$$t' = \{x' \in T' \suchthat \quad\text{there exists}\quad x \in t \quad\text{such that}\quad (x,x') \in R\}
$$
then the relation $R \intersect (t \cross t')$ is a $P$-equivariant $\epsilon$-almost isometry between $t$ and $t'$.  
\end{lemma}

\begin{proof} Choose a sequence $T_i$ converging to $T$. We use various notations, and we know various facts, from the proof of Lemma~\ref{LemmaLengthFacts} item~(\ref{ItemNormalizingContinuous}), in particular:  marked graphs $\wh G_i = T_i / F_r$, $\wh G = T / F_r$, marked homotopy equivalences $p_i \from \wh G_i \to \wh G$, $q_i \from \wh G \to \wh G_i$ lifting to $F_r$ equivariant maps $\ti p_i \from T_i \to T$, $\ti q_i \from T \to T_i$, etc. Sequences $r_i \ge 1$ converging to~$1$ and $\epsilon_i \ge 0$ converging to~$0$ were shown to have the property that for any $\eta>0$ there exists a neighborhood $U$ of $T$ such that if $T_i \in U$ then $\abs{r_i-1}, \abs{\epsilon_i} < \eta$. It therefore suffices to prove that if $r_i$ is sufficiently close to $1$ and $\epsilon_i$ is sufficiently close to $0$ then the conclusions of the lemma hold.

Consider now the $\Lambda_-$ isometries $k^+_i \from T_i \to T_+$. We know that $k^+_i \ti q_i \from T \to T_+$ converges to a $\Lambda_-$ isometry $f^+_T \from T \to T_+$. We wish to show that $f^+_T$ has a vertex with two or more positive gates, allowing us to conclude that $f^+_T=k^+_T$. Since there are only finitely many orbits of turns under the action of $F_r$, we may pass to a subsequence so that the gates in $T_i$ with respect to the maps $k^+_i$ all correspond under the $F_r$-equivariant homeomorphisms $\ti g_{ji}$. Pick a vertex $v \in T$ so that each $\tau_i(v)$ contains a vertex $w_i$ with two or more positive gates relative to $k^+_i$. The construction of the maps $\ti q_i$ gives us the freedom to choose $\ti q_i(v)$ in $\tau_i(v)$, and we choose $\ti q_i(v)=w_i$. It follows that there is a nondegenerate turn $\{e,e'\}$ at $v$ such that the directions $D_{\ti q_i}(e)$ and $D_{\ti q_i}(e')$ are contained in distinct gates of $k^+_i$ at $w_i$. Passing to the limit, this proves that $e,e'$ are contained in distinct gates of $f^+_T$ at $v$, proving that $f^+_T = k^+_T$.

Let $R_i$ denote the $\Lambda_-$ relation between $T$ and $T_i$: $x \in T$ and $y \in T_i$ are related by $R$ if $k^+_T(x) = k^+_i(y)$ and $x,y$ are contained in realizations of the same leaf of $\Lambda_-$ in $T,T_i$, respectively. Assuming this is so, we argue that the distance between $y$ and $\ti q^i(x)$ converges to zero uniformly as $i \to +\infinity$, independent of the choice of $(x,y) \in R$. Pick any leaf $L$ of $\Lambda_-$ realized in $T$, with realization $L_i \subset T_i$ and $L_+ \subset T_+$, so that $x \in L$ and $y \in L_i$. The arguments of Lemma~\ref{LemmaLengthFacts} item~(\ref{ItemNormalizingContinuous}) show that $\ti q^i(x)$ is within distance $\epsilon_i$ of some point $z_i \in L_i$, and so $k^+_i \ti q^i(x)$ is within distance $\epsilon_i$ of $k^+_i(z_i) \in L_+$. On the other hand $k^+_i \ti q^i(x)$ converges uniformly to $k^+_T(x) = k^+_i(y) \in L_+$, and so $k^+_i(z_i)$ converges to $k^+_i(y)$ in $L_+$, and thus the distance between $k^+_i(z_i)$ and $k^+_i(y)$ converges uniformly to zero. Since $k^+_i$ restricts to an isometry from $L_i$ to $L_+$, the distance between $y$ and $z_i$ in $L_i$ converges uniformly to zero. This proves that the distance between $\ti q^i(x)$ and $y$ converges uniformly to zero, as desired. 

Pick $\eta>0$, pick a finite subset $t \subset T$, and let $t_i \subset T_i$ be the finite set of points $y \in T_i$ related to some $x \in t$ by the relation $R_i$. The set $R_i \intersect (t \cross t_i)$ projects surjectively to $t$ and to $T_i$ because $R_i \subset T \cross T_i$ projects surjectively to $T$ and to $T_i$. Let $D$ be the diameter of $t$. We must show that for sufficiently large $i$ the relation $R_i \intersect (t \cross t_i)$ is an $\eta$-almost isometry (equivariance with respect to any finite subset of $F_r$ is easily verified),
in other words if $(x,y), (x',y') \in R_i \intersect (t \cross t_i)$ then 
$$\abs{d_T(x,x') - d_{T_i}(x',y')} < \eta
$$
Picking $i$ sufficiently large so that $d_{T_i}(\ti q^i(x),y) < \frac{\eta}{4}$ for any $(x,y) \in R_i$, it remains to show that the restriction $\ti q^i \restrict t$ is an $\frac{\eta}{2}$-almost isometry, that is, if $x,x' \in t$ then
\begin{equation}
\label{EqAlmostIsom}
d_T(x,x') - \frac{\eta}{2} \le d_{T_i}(\ti q^i(x),\ti q^i(x')) \le d_T(x,x') + \frac{\eta}{2}
\end{equation}
There is an integer $J \ge 0$ such that for any $x,x' \in t$ the arc $[x,x']$ makes at most $J$ turns. The restriction $\ti q^i \restrict [x,x']$ stretches or compresses each edge by at most $r_i \ge 1$, and makes a fold at each turn of length at most $\epsilon_i \ge 0$, as shown in the proof of Lemma~\ref{LemmaLengthFacts} item~\ref{ItemNormalizingContinuous}. It follows that
$$\frac{1}{r_i}  d_T(x,x') - J \epsilon_i   \le d_{T_i}(\ti q^i(x),\ti q^i(x')) \le r_i d_T(x,x')
$$
The right hand inequality of (\ref{EqAlmostIsom}) follows by choosing $i$ so large that $r_i < 1 + \frac{\eta}{2D}$, from which it follows that $r_i d_T(x,x') \le d_T(x,x') + \frac{\eta}{2}$. The left hand inequality follows by first choosing $i$ so large that $\frac{\eta}{2} - J \epsilon_i > 0$, and then choosing $i$ so large that $(1-\frac{1}{r_i})D \le \frac{\eta}{2} - J \epsilon_i$ which implies that $d_T(x,x') - \frac{\eta}{2} \le \frac{1}{r_i}  d_T(x,x') - J \epsilon_i$. 
\end{proof}

\subparagraph{Remarks.} In the course of this proof we have shown a little more regarding the structure of the space $I(\wh\A,T_+) = \union_{T \in \wh\A} I(T,T_+)$ in the Gromov topology (we refer the reader to \cite{Skora:deformations}, or to Section~\ref{SectionSkora}, verification of Step 3(d), for the definition of the Gromov topology on such spaces of maps). There is a continuous, open, surjective map $\pi \from I(\wh\A,T_+) \to \wh\A$, the ``domain map'', which maps each function $f \from T \to T_+$ in $I(\wh\A,T_+)$ to its domain $T$ in $\wh\A$. In the above proof, we have selected a particular element $k^+_T$ from each fiber $\pi^\inv(T) = I(T,T_+)$, thereby defining a section of the map~$\pi$. What we have proved amounts to the statement that this section is a continuous map $k^+ \from \wh\A \to I(\wh\A,T_+)$. In the case that $\Lambda_-$ is nonorientable one has in fact that $\pi$ is a homeomorphism with inverse $k^+$. But when $\Lambda_-$ is orientable and when there exist nontrivial fibers, one can also select the negative endpoint, and the same proof gives a continuous section $k^- \from \wh\A \to I(\wh\A,T_+)$. This shows that, in some sense, the intervals $I(T,T_+)$ all fit together in a continuous fashion, making the space $I(\wh\A,T_+)$ into an ``interval bundle'' over $\wh\A$ whose fiber over $T \in \wh\A$ is $I(T,T_+)$. To make sense out of this one would have to investigate local triviality, and one would have to deal with the fact that the fibers may not be topologically homogeneous, some being nondegenerate intervals and others degenerating to points. We shall not elaborate on this further, except insofar as to point out that the space of efficient representatives studied by Los and Lustig in \cite{LosLustig:axes} is naturally a subspace of the total space $I(\wh\A,T_+)$ of this interval bundle.

\subsection{Properness of the length map.} 
\label{SectionProperness}
Our present goal, a piece of the proof that $\Length \from \A \to (0,+\infinity)$ is a proper homotopy equivalence, is the following:

\begin{proposition}[Length is proper]
\label{PropAxisBundleCompact}
The map $\Length \from \A \to (0,+\infinity)$ is proper. It follows that a subset $C \subset \A$ is compact if and only if it is closed in $\A$ and the set $\Length(C)$ is contained in a compact subinterval of $(0,+\infinity)$.
\end{proposition}

The next lemma gives a compactness criterion that will be used several times in the proof of Proposition~\ref{PropAxisBundleCompact}. Recall our convention from Section~\ref{SectionPaths} that all fold maps are edge isometries. Given a subset $K \subset \wh\X_r$ define $F(K;1)$ to be the union of $K$ with the set of all $H \in \wh\X_r$ for which there exists $H_0 \in K$ and a fold map $f \from H_0 \to H$. By induction define $F(K;M)=F(F(K;M-1);1)$; equivalently this is the set of all $H \in \wh\X_r$ for which there exists $H_0 \in K$ and a piecewise isometry $f \from H_0 \to H$ (not required to be an edge isometry) which factors into $M$ marking preserving homotopy equivalences $H_0 \mapsto H_1 \mapsto \cdots \mapsto H_M=H$ each of which is either an isometry or a fold map. 

\begin{lemma}
\label{LemmaFoldsCompact}
For each compact $K \subset \wh\X_r$ and integer $M \ge 0$ the set $F(K,M)$ is compact.
\end{lemma}

\begin{proof} By induction it suffices to prove that $F(K;1)$ is compact. 

Let each closed cell $\sigma$ of $\X_r$ be labelled $\sigma(G)$ for some choice of marked graph $G \in \interior(\wh\sigma)$, so if $H \in \X_r$ then $H \in \interior(\sigma(G))$ if and only if there exists a marking preserving homeomorphism $G \approx H$; moreover, such a homeomorphism is unique up to isotopy rel vertices. Let $\wh\sigma(G)$ denote the preimage of $\sigma(G)$ in $\wh\X_r$. By local finiteness of the cell decomposition of $\X_r$, the compact set $K$ intersects only finitely many closed cells of $\wh\X_r$, and by intersecting each of these one at a time it suffices to fix a closed cell $\wh\sigma(G)$ and to assume that $K \subset \wh\sigma(G)$. 

For each ordered turn $(e,e')$ of $G$ let $F(K;e,e')$ denote the union of $K$ with the set of all $H \in \wh\X_r$ for which there exists $H_0 \in K$ and a fold map $f \from H_0 \mapsto H$ such that, with respect to the marked homeomorphism $G \approx H_0$, the map $f$ folds the turn $(e,e')$ and the length of $e$ in $H_0$ is greater than or equal to the length of $e'$ in $H_0$. By finiteness of the number of turns of $G$ it suffices to fix an ordered turn $(e,e')$ of $G$ and to prove that $F(K;e,e')$ is compact.

We consider two cases, depending on the valence of the common initial vertex $v$ of the edges $e,e'$.

Suppose first that $v$ has valence~$\ge 4$. Let $G^\#$ be the marked graph obtained from $G$ by folding proper initial segments of $e$ and of $e'$ to obtain an edge $e^\#$ of $G^\#$. The cell $\wh\sigma(G)$ is a codimension~1 face of $\wh\sigma(G^\#)$, defined by a marked homotopy equivalence $G^\# \to G$ that collapses $e^\#$ to a point and is otherwise a bijection. Under this collapse map, the edges $e,e'$ of $G$ correspond to two edges of $G^\#$ also denoted $e,e'$. Let $\Delta = \{(H_0,t) \in \wh\sigma(G) \cross [0,\infinity) \suchthat 0 \le t \le \Length_{H_0}(e') \le \Length_{H_0}(e)\}$. Define a map $\Delta \mapsto \wh\sigma(G^\#)$ taking $(H_0,t)$ to the element of $\wh\sigma(G^\#)$ which assigns length $t$ to $e^\#$, which subtracts $t$ from the lengths of the edges $e,e'$ in $H_0$, and which leaves the lengths of all other edges of $H_0$ unchanged. Notice that when $t=0$ the length assigned to $e^\#$ equals zero; whereas when $t=\Length_{H_0}(e')$ then the length assigned to $e'$ equals zero and, in the case  $\Length_{H_0}(e') = \Length_{H_0}(e)$, the length assigned to $e$ also equals zero; and so in either of these two cases the image of $(H_0,t)$ is in a face of the closed cell $\wh\sigma(\G^\#)$. The map $\Delta \mapsto \wh\sigma(G^\#)$ is clearly continuous, and $F(K;e,e')$ is the image under this map of a compact subset of $\Delta$, namely the set $\{(H_0,t) \in \Delta \suchthat H_0 \in K\}$, proving that $F(K;e,e')$ is compact.  

Suppose next that $v$ has valence~$3$. Let $e''$ be the third edge incident to $v$. Note that $e''$ is distinct from both $e$ and $e'$, for if say $e''=e$, consider any oriented leaf $L$ of $\Lambda_-$ travelling backwards along $e$; when $L$ hits the vertex $v$, it cannot take the illegal turn $(e,e')$ and so $L$ must continue into $e''=e$ again, and this must repeat forever, contradicting surjectivity of half-leaves of $\Lambda_-$. Now define $\Delta$ as above, and define the map $\Delta \mapsto \wh\sigma(G)$ taking $(H_0,t)$ to the element of $\wh\sigma(G)$ which subtracts $t$ from the lengths of $e,e'$ and adds $t$ to the length of $e''$. The proof is then completed as in the previous case.
\end{proof}

\begin{proof}[Proof of Proposition \ref{PropAxisBundleCompact}]  Since $\wh\A$ is a closed subset of $\wh\X_r$ and since $\Length$ is continuous, for each closed interval $[a,b] \subset (0,+\infinity)$ the set $\wh\A_{[a,b]} = \wh\A \intersect \Length^\inv[a,b]$ is closed, and it suffices to prove that $\wh\A_{[a,b]}$ is compact. Lemma~\ref{bounded combinatorial length} below implies that there exists a train track $G \in \wh\A$ and an integer $M \ge 0$ such that for each $H \in \wh\A_{[a,b]}$ there is a $\Lambda_-$ isometry $G \mapsto H$ which has a Stallings fold sequence consisting of at most $M$ folds. It follows that $\wh\A_{[a,b]} \subset F(\{G\};M)$, and the latter set is compact by Lemma~\ref{LemmaFoldsCompact}. Since $\wh\A_{[a,b]}$ is a closed subset of a compact set, it is compact. This completes the proof subject to our verifying Lemma~\ref{bounded combinatorial length}.
\end{proof}

\begin{lemma} \label{bounded combinatorial length} For any closed interval $[a,b] \subset (0,+\infinity)$ there exists a train track map $g  \from G \to G$ with $G$ normalized to be in $\wh\A$, and there exists an integer $M \ge 0$, so that if $H \in \wh\A_{[a,b]}$, and if $f_{\wt H}: \wt H \to T_+$ is a $\Lambda_-$-isometry, then there exists an edge isometry $h \from G \to H$ such that $f_g = f_{\wt H} \ti h$, and for any such edge isometry $h \from G \to H$, the combinatorial length of $h$ is at most $M$, and $h$ has a Stallings fold sequence of length at most $M$.
 \end{lemma}

\begin{proof} Pick a train track $G \in \A$ with the finest possible Nielsen structure.  Pick $\epsilon>0$ for $G$ by applying Proposition~\ref{PropWeakTrainTrack}. We may assume that $\epsilon > b$, for otherwise we replace $G$ by $\phi^{-n}(G)$ for sufficiently large $n$, which has the effect of replacing $\epsilon$ by $\epsilon \cdot \lambda(\phi)^n$. By Proposition~\ref{PropWeakTrainTrack}, for every $H \in\A_{[a,b]}$ and every $\Lambda_-$ isometry $f_{\wt H}: \wt H \to T_+$ there is an edge isometry $h \from G \to H$ so that $f_g = f_{\wt H} \ti h$ where $\ti h$ is a lift of $h \from G \to H$.  

Let $h \from G \to H$ be any edge isometry such that $f_g = f_{\wt H} \ti h$. Recall from the discussion preceding Lemma~\ref{LemmaFoldBound} that, after subdividing $H$ at the $h$-images of the vertices of $G$, the subdivided $H$ has at most $5r-5$ edges. Since $\Length(H) \ge a$, there exists an edge $e$ of $H$ with $\Length(e) \ge a/(5r-5)$. For any edge isometry $G' \to H$ whose vertices map disjointly from $e$, define an \emph{$e$-arc in $G'$} to be an open arc of $G'$ which maps isometrically onto the interior of~$e$. Two $e$-arcs in $G'$ are disjoint and each $e$-arc has length equal to $\Length(e)$, and so $G'$ has at most $\Length(G')(5r-5)/a$ distinct $e$-arcs.

Using the notation of (\ref{EqFolds}), apply Lemma~\ref{LemmaFoldBound} to obtain a Stallings fold sequence $G=G_0 \xrightarrow{g_1} G_1 \xrightarrow{g_2} \cdots \xrightarrow{g_J} G_J=H$ for $h$ whose length $J$ is bounded by the combinatorial length of $h$. For $j=0,\dots,J$ let $M_j$ be the number of $e$-arcs in $G_j$. The sequence $M_j$ is nonincreasing: each $e$-arc in $G_{j-1}$ maps isometrically onto an $e$-arc in $G_j$; the $e$-arcs in $G_{j-1}$ are precisely the components of the inverse images of the $e$-arcs of $G_j$; and so either $M_j = M_{j-1}$, or $M_j < M_{j-1}$ which occurs precisely when the fold map $g_j \from G_{j-1} \to G_j$ identifies at least one pair of $e$-arcs in $G_{j-1}$ to a single $e$-arc of $G_j$. There is therefore a sequence $0 \le j_1 < \cdots < j_N = J$ of some length $N$ such that the subsequence $M_j$, for $j_{n-1} < j \le j_n$, is constant, and such that $M_{j_n} > M_{j_n+1}$. It follows that $N \le M_0 \le \Length(G) (5r-5) / a$. 

The following lemma, applied with $K$ equal to the singleton $\{G\}$, and with $h$ equal to the composed map $G_{0} \mapsto G_{j_1-1}$, implies that $G_{j_1-1}$ belongs to a compact subset $K'$ of $\wh\X_r$ and that the map $G_{0} \mapsto G_{j_1-1}$ has uniformly bounded combinatorial length. Then apply Lemma~\ref{LemmaFoldsCompact} to the single fold $G_{j_1-1} \mapsto G_{j_1}$ to obtain a compact set $K''$ that contains $G_{j_1}$. Then reapply the following lemma with $K=K''$ and with $h$ equal to the composed map  $G_{j_{1}} \mapsto G_{j_2-1}$. Continue inductively in this manner to prove that each composed map $G_{j_{n-1}} \mapsto G_{j_n-1}$ has uniformly bounded combinatorial length. Summing up we obtain a uniform bound $M$ for the combinatorial length of $h$. Applying Lemma~\ref{LemmaFoldBound}, $h$ has a Stallings fold sequence of length at most $M$.
\end{proof}

\begin{lemma} 
\label{LemmaEdgeInjective}
For any compact set $K \subset \wh\X_r$ there exists a compact set $K' \subset \wh\X_r$ and an integer $L \ge 0$ so that for any $G \in K$, $G' \in \wh\X_r$, and edge isometry $h \from G \to G'$, if $\ell_\gamma(G') \ge \ell_\gamma(T^+)$ for all $\gamma \in F_r$, and if there exists an edge of $G'$ over the interior of which $h$ is injective, then $G' \in K'$ and $h$ has combinatorial length $\le L$.
\end{lemma}

In this lemma, edges of $G'$ are determined as usual by subdividing at the images of the vertices of $G$ under the map $h$, and so $G'$ has at most $5r-5$ edges.

\begin{proof}   A proper, connected, nontree subgraph of a marked graph $G$ determines a conjugacy class of proper, nontrivial free factors of $F_r$ by pulling back the fundamental group of the subgraph under the marking $R_r \to G$. Let $\FF(K)$ be the set of all free factor conjugacy classes associated to the proper, connected, nontree subgraphs of $G$, as $G$ varies over $K$. This set is finite, because the cells of the cellular topology on $\wh\X_r$ are locally finite (see Section~\ref{SectionOuterSpace}). By Theorem~5.4 of \BookZero, for each proper free factor conjugacy class $A$ in $F_r$, the infimum of $\ell_\gamma(T^+)$ for $\gamma \in A$ is a positive number denoted $\epsilon(A)$. Let $\epsilon>0$ be the minimum value of $\epsilon(A)$ as $A$ varies over the finitely many elements of $\FF(K)$. If $G,G'$ are as in the lemma then it follows that
for each $\gamma \in \F_r$, if the circuit in $G$ representing $\gamma$ does not surject onto $G$, then the circuit in $G'$ representing $\gamma$ has length $\ge \epsilon$. 

Let $e' \subset G'$ be an edge over whose interior $h$ is injective, let $X'$ be the complement in $G'$ of the interior of $e'$, and let $X = h^{-1}(X') \subset G$. Then $h(X)=X'$ and, since $h$ is a homotopy equivalence, $h$ restricts to a homotopy equivalence between $X$ and $X'$. It follows that $X,X'$ have the same number of components --- two or one depending on whether $e'$ is separating or nonseparating --- these components are connected, proper, subgraphs of $G,G'$, respectively, and neither $X$ nor $X'$ is a tree for otherwise $G$ or $G'$ would have a valence~1 vertex. Moreover, components of $X$, $X'$ that correspond under $h$ represent the same proper, nontrivial free factor conjugacy class of $F_r$. Thus every simple circuit in $X'$ has length at least $\epsilon$ and so crosses an edge of length at least $\epsilon/(5r-6)$, the subgraph $X'$ having at most $5r-6$ edges. Letting $E'$ be the union of $e'$ and the set of edges of $G'$ with length at least $\epsilon/(5r-6)$, it follows that every simple circuit in $G'$ intersects an edge in $E'$, and so the complement of $E'$ is a forest $U'$.  For each edge $e$ of $G$, since $\Length(e) \le \Length(G)$, the edge path $h(e)$ contains at most $1 + \Length(G) \cdot (5r-6)/\epsilon$ elements of $E'$. Removing these edges from the edge path $h(e)$ leaves at most $2+\Length(G)\cdot (5r-6)/\epsilon$ complementary subpaths each of which is contained in the forest $U'$ and so is embedded, crossing each edge of $U'$ at most once.  We conclude that there is a uniform bound to the combinatorial length of the edge path $h(e)$.  Since $e$ was arbitrary, the combinatorial length of $h$ is uniformly bounded. The map $h$ therefore has a Stallings fold sequence of uniformly bounded length. Applying Lemma~\ref{LemmaFoldsCompact}, $G'$ lies in a compact subset $K'$ of $\X_r$. 
\end{proof}

\subsection{Applying Skora's method to the Properness Theorem \ref{TheoremLengthPHE}}
\label{SectionSkora}
Let $\A \subset \X_r$ be the set of weak train tracks of a nongeometric, fully irreducible $\phi \in \Out(F_r)$, with a normalization $\A \inject \wh\X_r$ determined as usual by a choice of normalization of $T^+_\phi$. Rather than continue with the notation $\wh\A$ for the image of the normalization map, in this section we assume implicitly that all marked graphs and $F_r$-trees representing elements of $\A$ are normalized.

The major tool that we use in proving Theorem~\ref{TheoremLengthPHE}, that the map $\Length \from \A \mapsto (0,+\infinity)$ is a proper homotopy equivalence, is Skora's method \cite{Skora:deformations} for detecting the homotopy type of various spaces of trees. The preprint \cite{Skora:deformations} being unpublished and therefore difficult to obtain, other accounts to consult include \cite{White:FixedPoints}, \cite{Clay:contractibility}, and \cite{GuirardelLevitt:outer}; we shall quote a technical result from the latter which simplifies the method. While Skora applies his technique to very general spaces of group actions on trees, we shall consider only free, simplicial, minimal $F_r$-trees, that is, elements of outer space $\X_r$.

\paragraph{Spaces of edge isometries.} Skora's method requires defining topologies on various spaces of piecewise isometries between trees. We adopt here an addition to the method suggested by Definition~3.8 and Proposition 3.10 of \cite{GuirardelLevitt:outer} which guarantees that minimality of trees is true at all stages of Skora's method, and which allows us to focus solely on edge isometries among free, simplicial, minimal $F_r$ trees, that is, elements of $\wh\X_r$.

Given free, simplicial, minimal $F_r$ trees $S,T$, a piecewise isometry $f \from S \to T$ satisfies the \emph{minimality condition} if for each $x \in S$ there exists an open arc embedded in $S$ and containing $x$ such that the restriction of $f$ to this arc is an isometry; equivalently, $f$ is an edge isometry and each vertex of $S$ has at least two gates with respect to $S$. The motivation for the minimality condition is that if it failed, if there were a point $x$ having a single gate, then one could construct an equivariant Stallings fold sequence for $f$ by first folding all of the turns at $x$, producing a tree which has a valence~1 vertex and therefore is not minimal. Particularly important is the fact that $\Lambda_-$ isometries automatically satisfy the minimality condition.

Let $\M$ be the set of equivalence classes of piecewise isometries $f \from S \to T$ satisfying the minimality condition, where $S,T$ are free, simplicial, minimal $F_r$ trees, and where $f \from S \to T$ and $f' \from S' \to T'$ are equivalent if there exist $F_r$ equivariant isometries $h \from S \to S'$, $k \from T \to T'$ such that $kf=f'h$. Formally the element $f \from S \to T$ is to be regarded as an ordered triple $(f,S,T)$ which allows us to define two projections $\Domain,\Range \from \M \to \X_r$, the domain projection $\Domain(f,S,T)=S$ and the range projection $\Range(f,S,T) \mapsto T$. 

The group $\Out(F_r)$ acts on $\M$ from the right so that the domain and range projections are $\Out(F_r)$ equivariant. The action is defined as follows. Given $f \from S \to T$ in $\M$ and given $\phi \in \Out(F_r)$, choose $\Phi \in \Aut(F_r)$ representing $\phi$. Recall that the $F_r$-trees $S \cdot \phi$, $T \cdot \phi$ have underlying $\reals$-trees $S,T$ but with $F_r$ actions precomposed by $\Phi$. The $F_r$-equivariant map $f \from S \to T$ therefore remains $F_r$-equivariant when the actions are precomposed by $\Phi$, and we define this map to be the image of $f$ under $\phi$, well-defined in $\M$ independent of the choice of $\Phi$.

We put off until later the definition of the Gromov topology on $\M$, as given in the references above. For the moment, we note that the domain and range projections $\Domain,\Range \from \M \to \X_r$ are continuous maps with respect to this topology.

Notice that the reparameterized maps $-\log\Length = \log\left(\frac{1}{\Length}\right) \from \A \to \reals$ and $\frac{1}{\Length} \from \A \to (0,+\infinity)$ are  ``correctly oriented'' in the sense that, on an orbit $G_i = \phi^i(G_0)$, the values of these maps increase as $i$ increases. The map $\Length \from \A \to (0,+\infinity)$, on the other hand, is ``incorrectly oriented'', and the map $-\log\Length$ carries heavier notational baggage, so in the proof that follows we often think in terms of the map $\frac{1}{\Length}$.

\paragraph{Proof of the Properness Theorem \ref{TheoremLengthPHE}.}
We give the proof in the several steps. 

\subparagraph{Step 1: A bi-infinite fold path $\L$.} Our goal being to prove that $\frac{1}{\Length} \from \A \to (0,+\infinity)$ is a proper homotopy equivalence, we begin by choosing a subset $\L \subset \A$ such that $\frac{1}{\Length} \from \L \to (0,+\infinity)$ is a homeomorphism. This subset $\L$ will be a particular fold line, periodic under an appropriate power of $\phi$.

Passing to a positive power of $\phi$, we may assume $\phi$ is rotationless. Applying Proposition~\ref{no cut vertices}, choose $g \from G \to G$ to be a train track representative of $\phi$ whose local Whitehead decomposition is as fine as possible, meaning that the local stable Whitehead graphs of $g$ have no cut points. We assume that the train track $G$ is normalized just as $\A$ is, and so there is an associated $\Lambda_-$ isometry $f_g \from \wt G \to T_+$. Define $G_0 = G$ and $G_i = \lambda^{-i}\phi^i(G)$ for all $i \in \Z$. The map $g \from G \to G$ induces marked edge isometries $g_{i+1,i} \from G_i \to G_{i+1}$. Let $L_i = \Length(G_i) = \lambda^{-i} \Length(G_0)$. 

For later use recall that if we define $g_{ji} = g_{j,j-1} \composed g_{j-1,j-2} \composed \cdots \composed g_{i+1,i} \from G_i \to G_j$ then the universal covering trees $\wt G_i$ and the $F_r$-equivariant edge isometries $\ti g_{ji}$ form a direct system in the category of metric spaces and distance nonincreasing maps, the direct limit of this system is the $F_r$ tree $T_+$, we obtain direct limit maps $\ti g^+_{i} \from \wt G_i \to T_+$ such that the two maps $f_g, \ti g^+_{0}  \from \wt G = \wt G_0 \to T_+$ are equal to each other.

As described in Section~\ref{SectionPaths}, choose a Stallings fold sequence for $g_{1,0}$ and interpolate the fold maps to get a Stallings fold path for $g_{1,0}$ connecting $G_0$ to $G_1$. We parameterize this fold path by $\frac{1}{\Length}$, with parameter interval $[\frac{1}{L_0},\frac{1}{L_1} = \frac{\lambda}{L_0}]$. Translating this path by $\phi^i$ and reparameterizing, one obtains a Stallings fold path for $g_{i-1,i}$ connecting $G_{i-1}$ to $G_i$ and parameterized by $\frac{1}{\Length}$, with parameter interval $[\frac{\lambda^i}{L_0},\frac{\lambda^{i+1}}{L_0}]$. Concatenating these paths over all $i \in \Z$ we obtain a bi-infinite Stallings fold path $(0,+\infinity) \mapsto \A$. The image of this path is the desired set $\L$, completing Step~1. 

Note that $\frac{1}{\Length} \from \L \to (0,+\infinity)$ is a homeomorphism. The inverse of this homeomorphism, which we denote $\rho \from (0,+\infinity) \to \L$, will be our preferred parameterization for the fold path $\L$. 

\subparagraph{Step 2: The map $\A \to \L$.} We construct a proper, continuous map $\sigma \from \A \to \L$ whose restriction to $\L$, while not the identity on $\L$, is properly homotopic to the identity on $\L$. In later steps we shall prove that $\sigma$ is a proper homotopy equivalence.

Applying Proposition~\ref{PropWeakTrainTrack}, there exists a number $\eta>0$ such that if $T$ is any weak train track with $\Length(T) < \eta$ and if $f_T \from T \to T_+$ is a $\Lambda_-$ isometry, then there exists a $\Lambda_-$ isometry $\psi \from \wt G_0 \to T$ such that $f_T \composed \psi = f_g$. 

Define the map $\sigma \from \A \to \L$ by
$$\sigma(T) = \rho \left( \frac{\eta}{L_0 \lambda \Length(T)} \right) = \rho \left(  \frac{\eta}{L_0 \lambda} \rho^\inv(T) \right)
$$
This map is evidently continuous, and it is proper by Proposition~\ref{PropAxisBundleCompact}. The restriction to $\L$, when conjugated by $\rho$, is the map from $(0,+\infinity)$ to itself given by $r \mapsto \frac{\eta}{L_0 \lambda} r$. This map being an orientation preserving homeomorphism of $(0,+\infinity)$, the map $\sigma \restrict \L \from \L \to \L$ is properly homotopic to the identity.

\subparagraph{Remark.} For later purposes we note that $\eta < L_0 < \lambda L_0$ and so $\frac{\eta}{\lambda L_0} < 1$. In other words, the map $\sigma$ decreases the value of the parameter $\frac{1}{\Length}$, multiplying it by the number $\frac{\eta}{\lambda L_0}$, and so $\sigma$ moves each point of $\L$ strictly to the left.

\subparagraph{Step 3: The map to $\M$.} For each $T \in \A$ we define a $\Lambda_-$ isometry $\psi^T \from \sigma(T) \to T$ with the property that the map $(\sigma(T),T) \mapsto \psi^T$ is a continuous map from the set $\{(\sigma(T),T) \in \L \cross \A \suchthat T \in \A\}$ to the set $\M$. Letting $u \in (0,+\infinity)$ be such that $\sigma(T)=\rho(u)$, the map $\psi^T \from \sigma(T) \to T$ will be defined by factoring a particular $\Lambda_-$ isometry $f^+_u \from \rho(u)=\sigma(T) \to T_+$ through $T$, where the other factor is a particular $\Lambda_-$ isometry $k^+_T \from T \to T_+$. We construct these $\Lambda_-$ isometries, and we define $\psi^T$ and verify its properties, in several substeps.

\subparagraph{Step 3(a): The $\Lambda_-$ isometry $f^+_u$ from $\rho(u) \in \L$ to $T_+$.} Define the \emph{semiflow} of the flow line $\L$ as follows. This semiflow will be a family of $F_r$-equivariant $\Lambda_-$ isometries $f_{ts} \from \rho(t) \to \rho(s)$, $0 < s < t < +\infinity$, so that the semiflow identity $f_{ut} \composed f_{ts} = f_{us}$ is satisfied. For each integer $i$ the map $f_{1/L_{i},1/L_{i-1}} \from \rho(1/L_{i-1}) \to \rho(1/L_{i})$ is defined to be the map $\ti g_{i,i-1} \from \wt G_{i-1} \to \wt G_{i}$. For each integer $i$, the subfamily of maps $f_{ts}$ for $\frac{1}{L_{i-1}} \le s \le t \le \frac{1}{L_{i}}$ is (up to reparameterization) defined to be the semiflow of the chosen Stallings fold sequence for $g_{i,i-1} \from G_i \to G_{i-1}$. The definition of $f_{ts}$ extends uniquely to $0 < s \le t < \infinity$ in such a way as to satisfy the semiflow identity. 

The trees $\rho(t)$ and the maps $f_{ts}$ form a direct system in the category of $F_r$-equivariant metric spaces and distance nonincreasing maps. The direct limit of this direct system is the $F_r$ tree $T_+$, because the trees $\wt G_i$ and the maps $g_{ji}$ form a co-finite subsystem, whose direct limit is known to be $T_+$ by Theorem~\ref{T+ is direct limit}. We obtain direct limit maps $f^+_u \from \rho(u) \to T_+$ which satisfy the property that $f^+_u \composed f_{ut} = f^+_t$ for $0<t<u<+\infinity$. When $\rho(u)=\sigma(T)$ we also denote $f^+_u = f^+_{\sigma(T)} \from \sigma(T) \to T_+$.

\subparagraph{Step 3(b): The $\Lambda_-$ isometry $k^+_T$ from $T \in \A$ to $T_+$.} 
This was defined in Section~\ref{SectionGHConvergence} using Theorem~\ref{TreeTheorem}: if $\Lambda_-$ is nonorientable then $k^+_T$ is the unique $\Lambda_-$ isometry $T \mapsto T_+$; and when $\Lambda_-$ is oriented then $k^+_T$ is the positive endpoint of the compact oriented interval $I(T,T_+)$ of $\Lambda_-$ isometries.

\subparagraph{Step 3(c): The $\Lambda_-$ isometry $\psi^T \from \sigma(T) \to T$.} 
For each $T \in \A$ the map $\psi^T$ is characterized as the unique $\Lambda_-$ isometry from $\sigma(T)$ to $T$ with the property that $f^+_{\sigma(T)} = k^+_T \composed \psi^T$. Uniqueness of $\psi^T$ follows from the proof of uniqueness in Proposition~\ref{PropWeakTrainTrack}. However, we still have to establish that $\psi^T$ exists with this property.

The existence of $\psi^T$ being invariant under the action of $\phi$ on $\A$,
it suffices to prove existence after replacing $T$ by $\phi^k(T)$ for a particular power $k$. Since $\phi$ increases the parameter $\frac{1}{\Length}$ by multiplying it by the quantity $\lambda$, after acting by the appropriate power of $\phi$, we may assume that
$$\frac{1}{\eta} \le \frac{1}{\Length(T)} \le \frac{\lambda}{\eta}
$$
Recalling that $\frac{1}{\Length(\sigma(T))} = \frac{\eta}{L_0 \lambda \Length(T)}$ we have
$$\frac{1}{L_0 \lambda} \le \frac{1}{\Length(\sigma(T))} \le \frac{1}{L_0}
$$
Setting $t_0 = \frac{1}{\Length(\sigma(T))}$, the semiflow map $f_{\frac{1}{L_0},t_0} \from \sigma(T) \to \wt G$ is therefore defined. Since $\Length(T) \le \eta$, Proposition~\ref{PropWeakTrainTrack} gives us a $\Lambda_-$ isometry $h_T \from \wt G = \rho(1/L_0) \to T$ with the property that $f_g = k^+_T \composed h_T$. Define $\psi^T = h_T \composed f_{\frac{1}{L_0},t_0}$, and so
\begin{align*}
k^+_T \composed \psi^T &= k^+_T \composed h^{\vphantom{+}}_T \composed f^{\vphantom{+}}_{\frac{1}{L_0},t_0} \\
  &= f_g \composed  f_{\frac{1}{L_0},t_0}  \\
  &= f^+_{\frac{1}{L_0}} \composed f^{\vphantom{+}}_{\frac{1}{L_0},t_0}  \\
  &= f^+_{t_0} = f^+_{\sigma(T)}
\end{align*}
as desired.

\subparagraph{Remark:} What we have done is to define $\psi^T$ first for $T$ in a fundamental domain of the action of $\phi$ on $\A$, and then to extend over all of $\A$ by $\phi$-equivariance. This raises the issue of well-definedness of $\psi^T$ where one end of the fundamental domain, defined by $\frac{1}{\Length(T)} = \frac{1}{\eta}$, is translated by $\phi$ to the other end of the fundamental domain, defined by $\frac{1}{\Length(T)} = \frac{\lambda}{\eta}$. This issue is made moot by the fact that $\psi^T$ is already uniquely characterized: it is the unique $\Lambda_-$ isometry from $\sigma(T) \to T$ that satisfies the equation $f^+_t = k^+_T \composed \psi^T$.

\subparagraph{Step 3(d): Continuity of the map $(\sigma(T),T) \mapsto \psi^T$.}
We must prove that the map from the set $\{(\sigma(T),T) \suchthat T \in \A\} \subset \L \cross \A$ to the set $\M$, taking the pair $(\sigma(T),T)$ to the map $\psi^T \from \sigma(T) \to T$, is a continuous map. 

This is the most technical of the steps and requires us to look carefully at the topology of the mapping space $\M$. We put this issue off for a while more, in order to get through the steps in the proof of Theorem~\ref{TheoremLengthPHE}.

\subparagraph{Step 4:} Combining Steps~2 and~3, we obtain a continuous map $\A \mapsto \M$ defined by $T \mapsto \psi^T$, so that postcomposing with $\Domain \from \M \to \X_r$ we obtain the map $\sigma \from \A \mapsto \L$, and postcomposing with the map $\Range \from \M \to \X_r$ we obtain the identity map on $\A$.

Now apply the following theorem:

\begin{theorem}[\cite{Skora:deformations}, \cite{GuirardelLevitt:outer}]
\label{TheoremSkora}
Letting $\nabla = \{(s,t) \in [0,\infinity] \cross [0,\infinity] \suchthat s \le t\}$, there exists a continuous map $\theta \from \M \cross \nabla \to \M$, denoted formally as 
$$\theta(\psi,s,t) = (\psi_{ts}, T_s, T_t)
$$
with the following properties for each $\psi \from S \to T$ in $\M$:
\begin{enumerate}
\item The map $\psi_{\infinity,0}  \from T_0 \to T_\infinity$ equals $\psi \from S \to T$.
\item The map $\psi_{0,0} \from T_0 \to T_0$ equals the identity on $S$.
\item The map $\psi_{\infinity,\infinity}$ equals the identity on $T$.
\item (The semiflow property) $\psi_{ut} \composed \psi_{ts} = \psi_{us}$ for $0 \le s \le t \le u \le \infinity$.
\end{enumerate}
\end{theorem}

\paragraph{Remarks on Theorem \ref{TheoremSkora}.} Rather than the general $F_r$ trees allowed in \cite{Skora:deformations}, we are restricting elements $(\psi,S,T) \in \M$ by requiring $S,T$ to be free, simplicial, minimal $F_r$ trees, and by requiring $f$ to satisfy the minimality condition. To see that this is allowable, first note that if $\M$ is defined using ordinary piecewise isometries rather than piecewise isometries satisfying the minimality condition, then Theorem~\ref{TheoremSkora} is exactly the main result, Theorem 4.8, of \cite{Skora:deformations} except for the restrictions to free simplicial trees (and except for correcting the erroneous $\{(s,t) \in [0,1] \cross [0,1]\}$ occurring in the statement of Theorem 4.8). Then note that for a map $\theta(\psi,s,t)=\psi_{ts} \from T_s \to T_t$ obtained by the constructions of \cite{Skora:deformations} or of \cite{GuirardelLevitt:outer}, the tree $T_t$ is simplicial by Proposition~3.6 of \cite{GuirardelLevitt:outer}, and $T_t$ is free because $T$ is free, and $\psi_{ts}$ satisfies the minimality condition by Propositions~3.4, 3.10, and~3.11 of \cite{GuirardelLevitt:outer}. 

\bigskip

To complete Step~4, as a consequence of Theorem~\ref{TheoremSkora} the following chain of compositions is continuous:
$$\A \cross [0,\infinity] \xrightarrow{(T,s) \mapsto (T,s,\infinity)} \A \cross \nabla \xrightarrow{(T,s,t) \mapsto (\psi^T,s,t)} \M \cross \nabla \xrightarrow{\theta} \M \xrightarrow{\Domain} \X_r
$$
This homotopy, which we denote $h \from \A \cross [0,\infinity] \to \X_r$, has the property that $h(T,0) = \sigma(T)$ and $h(T,\infinity) = T$, in other words $h$ is a homotopy from the map $\sigma\from\A \mapsto \L$ to the identity map on $\A$. 

As a consequence of the semiflow property, for each $t \in [0,\infinity]$ we obtain a commutative diagram as follows, where we denote $T_t = h(T,t)$:
$$\xymatrix{
\sigma(T)=T_0 \ar[rr]^{\psi^T} \ar[dr]_{\psi^T_{t,0}} & & T_\infinity=T  \ar[r]^{k^+_T} & T_+ \\
 & T_t \ar[ur]_{\psi^T_{\infinity,t}}
}$$
We reiterate that the maps $\psi^T_{t,0}, \psi^T_{\infinity,t} \in \M$ are edge isometries with the minimality property.

\subparagraph{Step 5:} We must verify that $h$ takes values in $\A$, in other words, in the above diagram the tree $T_t$ is a weak train track. For any leaf $L$ of $\Lambda_-$, its realization in $\sigma(T)$ denoted $L_{\sigma(T)}$ maps injectively under $\psi^T$ and so also maps injectively under the edge isometry $\psi^T_{t,0}$, with image therefore equal to the realization of $L$ in $T_t$, denoted $L_t$. Note that $\psi^T_{\infinity,t}$ also restricts to an injection on $L_{t}$, and its image equals $\psi^T(L_{\sigma(T)})$, which is the realization of $L$ in $T$ denoted $L_T$. Since $\psi^T_{\infinity,t}$ is an edge isometry it follows that $\psi^T_{\infinity,t}$ restricts to an isometry from $L_t$ to $L_T$. The map $k^+_T \composed \psi^T_{\infinity,t}$ therefore restricts to an isometry from $L_t$ into $T_+$. This being true for any leaf $L$ of $\Lambda_-$, we have proved that $k^+_T \composed \psi^T_{\infinity,t} \from T_t \to T_+$ is a $\Lambda_-$ isometry, and so $T_t \in \A$ as desired.

\subparagraph{Step 6:} We must verify that the homotopy $h \from \A \cross [0,\infinity] \to \A$ is a proper map, that is, for each compact set $C \subset \A$ the set $h^\inv(\C)$ is compact. This will be a consequence of Theorem~\ref{PropAxisBundleCompact}. Since $C$ is compact, there is a closed interval $[a,b] \subset (0,\infinity)$ such that $\frac{1}{\Length} \restrict C$ takes values in $[a,b]$. For each $T \in \A$, the definition of $h$ and the remark at the end of the verification of Step~2 combine to show that as $t$ goes from $0$ to $\infinity$ the value of $\frac{1}{\Length(h(T,t))}$ increases monotonically from $\frac{1}{\Length(h(T,0))} = \frac{1}{\Length(\sigma(T))} = \frac{\eta}{\lambda L_0 \Length(T)}$ to $\frac{1}{\Length(h(T,\infinity))} = \frac{1}{\Length(T)}$. It follows that
$$\frac{1}{\Length(h^\inv(C))} \subset \left[ \frac{\eta}{\lambda L_0} a,b\right]
$$
that is, the function $\frac{1}{\Length}$ restricted to $h^\inv(C)$ takes  values in a closed subinterval of $(0,+\infinity)$. Since $h$ is continuous, $h^\inv(C)$ is also closed. Applying Proposition~\ref{PropAxisBundleCompact}, it follows that $h^\inv(C)$ is compact.

\bigskip

Steps~1 through~6 prove that the map $\sigma \from \A \to \L$ is a proper homotopy equivalence. Composing this map with the homeomorphism $\Length \from \L \to (0,+\infinity)$, and then composing with the self homeomorphism of $(0,+\infinity)$ which is multiplication by the number $\frac{\lambda L_0}{\eta}$, we obtain the map $\Length \from \A \to (0,+\infinity)$, proving that this map is a proper homotopy equivalence. This finishes the proof of Theorem~\ref{TheoremLengthPHE}, subject to completing the verification of Step 3(d).

\paragraph{Verifying Step 3(d): Continuity of the map $(\sigma(T),T) \mapsto \psi^T$.} \hfill\break
Let $\P=\{(\sigma(T),T) \suchthat T \in \A\}$, so $\P$ equals the subset of $\L \cross \A$ consisting of all $(S,T)$ such that $\Length(T)/\Length(S) = \frac{\eta}{\lambda L_0}$, equivalently, all $(\rho(u),T)$ such that $u \Length(T) = \frac{\eta}{\lambda L_0}$. We shall show that the map $\P \mapsto \M$ taking $(\rho(u),T)$ to $\psi^T$ is continuous.
We remind the reader that $\psi^T \from \rho(u) \to T$ is the unique $\Lambda_-$ isometry that factors the direct limit $\Lambda_-$ isometry $f^+_{\rho(u)} \from \rho(u) \to T_+$ through the $\Lambda_-$ isometry $k^+_T \from T \to T_+$. In this argument we shall denote $\psi^T=\psi^T_u$, the superscript being the range and the subscript denoting the domain $\rho(u)$.

Recall from Section~\ref{SectionOuterSpace} the  Gromov topology on $\X_r$, with basis element $V(T;\epsilon,t,P)$ determined by $T \in \wh\X_r$, $\epsilon>0$, a finite subset $t \subset T$, and a finite subset $P \subset F_r$, consisting of all $T' \in \wh\X_r$ for which there exists a finite subset $t' \subset T'$ and a $P$-equivariant $\epsilon$-almost isometry $\rho \subset t \cross t'$ that projects surjectively to $t$ and to $t'$.

Next we review from \cite{Skora:deformations} the Gromov topology on $\M$. A basis element $U(\psi;\epsilon,s,t,P)$ for this topology is determined for each $\psi \from S \to T$ in $\M$, $\epsilon>0$, finite subsets $s \subset S$, $t \subset T$, and finite set $P \subset F_r$. This basis element consists of all $\psi' \from S' \to T' \in \M$ for which there exist finite subsets $s' \subset S'$ and $t' \subset T'$, and $P$-equivariant $\epsilon$-almost isometries $R_0 \subset s \cross s'$, $R_1 \subset t \cross t'$ such that $(x,y) \in R_0$ implies $(\psi(x),\psi'(y)) \in R_1$. To put this more colloquially, $\psi' \from S' \to T'$ is a member of $U(\psi;\epsilon,s,t,P)$ if we can exhibit $S'$ as a member of $V(S;\epsilon,s,P)$ and $T'$ as a member of $V(T;\epsilon,t,P)$ in such a way that the two $P$-equivariant $\epsilon$-almost isometries that exhibit these memberships are respected by the maps $\psi$ and~$\psi'$. Note the evident continuity of the domain and range maps $\Domain, \Range \from \M \to \wh\X_r$.

We remark that the Gromov topology on $\M$ as defined in \cite{Skora:deformations} uses arbitrary compact subsets of metric spaces, not just finite subsets. However, just as Paulin does in Proposition~4.1 of~\cite{Paulin:GromovTopology}, one obtains a different basis for the same topology on $\M$ by restricting to finite subsets.

The relations that we shall use for $\epsilon$-almost isometries will be restrictions of the $\Lambda_-$ relations defined in Section~\ref{SectionGHConvergence}. 

Pick $(\rho(u),T) \in \P$, denote $S=\rho(u)$, and let $\psi = \psi^T_u \from S \to T$. Pick a basis element $U(\psi;\epsilon,s,t,P) \subset \M$ about $\psi$. Corresponding to this are basis elements $V(S;\epsilon,s,P) \subset \X_r$ about $S$ and $V(T;\epsilon,t,P) \subset \wh\X_r$ about $T$. We must produce $P$-equivariant $\epsilon$-almost isometries which exhibit membership in these two basis elements.

Applying Lemma~\ref{LemmaRelation}, we obtain a neighborhood $W \subset \A$ of $T$ so that for each $T' \in W$, there is a finite subset $t' \subset T'$ such that if $R_1$ denotes the restriction to $t \cross t'$ of the $\Lambda_-$ relation between $T$ and $T'$, then $R_1$ surjects to $t$ and to $t'$, and $R_1$ is a $P$-equivariant $\epsilon$-almost isometry, thereby exhibiting that $T' \in V(T;\epsilon,t,P)$.

We also need a relation $R_0$ that exhibits membership in $V(S;\epsilon,s,P)$ as long as $\abs{u-u'}$ is sufficiently small. We will define this relation in two cases, depending on whether $u \le u'$ or $u' < u$.

When $u \le u'$ then the semiflow associated to the Stallings fold path $\L$ gives a $\Lambda_-$ isometry $f_{u',u} \from S=\rho(u) \to \rho(u')=S'$. Set $s'=f_{u',u}(s)$. Regarding $f_{u',u}$ as a subset of $S \cross S'$, let $R_0 = f_{u',u} \intersect (s \cross s')$; to put it another way, $R_0$ is the restriction of $f_{u',u}$ to $s$. $P$-equivariance of $R_0$ is an immediate consequence of $F_r$-equivariance of $f_{u',u}$. There is an integer $J$ such that any two points of the finite set $s$ are connected by an edge path with no more than $J$ turns. We choose $\abs{u-u'}$ sufficiently small so that the map $f_{u',u}$ is a fold of length no more than $\epsilon/J$. So, if $x,y \in s$ then the segment between $x$ and $y$ undergoes at most $J$ folds each of length at most $\epsilon$, and so the distance $d_{S'}(f_{u',u}(x),f_{u',u}(y))$ is less than $d_{S}(x,y)$ by an amount that is no more than $J(\epsilon/J)=\epsilon$. This shows that $f_{u',u} \restrict s \cross s'$ is an $\epsilon$-almost isometry.

Next assume that $u' \le u$. The semiflow gives a $\Lambda_-$ isometry $f_{u,u'} \from S'=\rho(u') \to S=\rho(u)$. Let $s' = f_{u,u'}^\inv(s)$ and let $R_0 = (f_{u,u'} \restrict s')^\inv \subset s \cross s'$. $P$-equivariance is again immediate. Letting $H(s)$ be the convex hull of $s$ in $S$, we can think of $s'$ as obtained from $s$ by splittings at certain vertices of $H(s)$. To be precise, for any number $\eta > 0$ such that $\eta$ is less than half the length of the shortest edge of $H(s)$, if $S'$ is sufficiently close to $S$ in $\wh\X_r$ then $s'$ is obtained from $s$ in the following manner: for certain vertices $v$ in $H(s)$, choose a direction $d$ at $v$, choose a segment $[v,w]$ representing $d$ of length at most $\eta$, partition the directions at $v$ distinct from $d$ into two nonempty subsets $D_1 \union D_2$, and then split the segment $[v,w]$ into two segments $[v_1,w]$ and $[v_2,w]$, reattach the directions $D_1$ to $v_1$, and reattach the directions $D_2$ to $v_2$. From this construction it follows that there is a constant $J$ depending only on $s$ such that any two points $x,y \in s'$ are connected by an edge path with no more than $J$ turns. This implies that $d_S(f_{u',u}(x),f_{u',u}(y))$ is less than $d_{S'}(x,y)$ by an amount that is no more than $J \eta$, and so by choosing $\eta < \epsilon/J$ it follows that $f_{u,u'}^\inv \restrict s \cross s'$ is an $\epsilon$-almost isometry.

Finally, we must verify that $R_0$ and $R_1$ are respected by the maps $\psi^T_u$ and $\psi^{T'}_{u'}$. That is, we must show that if $(x,y) \in R_0$ then $(\psi^T_u(x),\psi^{T'}_{u'}(y)) \in R_1$. Focussing on the case that $u \le u'$ (the opposite case $u' \le u$ being handled similarly), we have $y=f_{u',u}(x)$, and we must verify that 
\begin{equation}
\label{EqLambdaRel}
(\psi^T_u(x),\psi^{T'}_{u'} \composed f_{u',u}(x)) = (\psi^T_u(x),\psi^{T'}_u(x)) \in R_1
\end{equation}
$R_1$ being the restriction of the $\Lambda_-$ isometry between $T$ and $T'$, there are two things to check. First, if $L$ is any leaf of $\Lambda_-$ realized in $\rho(u)$ such that $x \in L$, then since $\psi^T_u$ and $\psi^{T'}_u$ are both $\Lambda_-$ isometries it follows that $\psi^T_u(x)$ is in the realization of $L$ in $T$, and $\psi^{T'}_u(x)$ is in the realization of $L$ in $T'$. Second, we note that
$$k^+_T \composed \psi^T_u(x) = f^+_u(x) = k^+_{T'} \composed \psi^{T'}_u(x)
$$
completing the verification of (\ref{EqLambdaRel}), and thereby completing the verification of Step~3(d) and the proof of Theorem~\ref{TheoremLengthPHE}.

\newcommand\Turns{\text{Turns}}
\renewcommand\L{\mathcal{L}}
\newcommand\LL{\mathcal{LL}}

\section{Fold lines} 
\label{SectionFoldLines}
In this section we investigate Definition~(1) of the axis bundle of a nongeometric, fully irreducible $\phi \in \Out(F_r)$: the union of fold lines converging in $\overline\X_r$ to $\TMinusphi$ in the negative direction and to $\TPlusphi$ in the positive. Our main result, Theorem~\ref{TheoremFoldLines}, says that every point on such a fold line is a weak train track for $\phi$, thus proving that Definition~(1) of the axis bundle implies Definition~(2). We also prove the converse, that Definition~(2) implies Definition~(1); see Proposition~\ref{PropWTTOnFoldLine}.

Recall from Section~\ref{SectionFolds} that a fold path in $\X_r$ is an injective, continuous, proper map $\gamma \from I \to \X_r$ with domain a connected subset $I \subset \reals$, for which there is a continuous lift $\wh\gamma \from I \to \wh\X_r$, marked graphs $G_r$ representing $\wh\gamma(r)$ for each $r \in I$, and a family of edge isometries $f_{sr} \from G_r \to G_s$ defined for $s \le r$ in $I$ and satisfying the semiflow condition $f_{ts} \composed f_{sr}=f_{tr}$. Reparameterizing $\gamma$ if necessary, we ofttimes impose the following:
\begin{description}
\item[Normalization Condition:] $\Length(G_r) = e^{-r}$ for each $r \in I$.
\end{description}
If $\gamma$ is normalized then $G_r$ is uniquely determined up to an isometry that preserves marking, and so the lift $\wh\gamma$ is uniquely determined. A translation in the parameter amounts to a change in normalization: for each $T \in \reals$, the reparameterization $\gamma'(r') =\gamma(r'+T)$ is also a fold path, and $\wh\gamma'(r') = e^{T} \cdot \wh\gamma(r'+T)$.

The set of data $(G_r)_{r \in I}$, $(g_{sr} \from G_r \to G_s)_{r<s \in I}$ is called a \emph{connection} for the fold line $\gamma$. While $(G_r)$ is uniquely determined by $\gamma$ in the presence of the normalization condition, the maps $(g_{sr})$ are not determined in general, as we saw from the example in Section~\ref{SectionWTTexamples}.

A normalized fold path $\gamma$ whose domain $I$ is noncompact only on the
positive end then $\gamma$ is a \emph{fold ray}; and if $I$ is noncompact
only on the negative end then $\gamma$ is called a \emph{split ray}.
Denote 
$$I_- = \inf(I) \in \reals \union \{-\infinity\}, \qquad I_+ = \sup(I) \in \reals \union \{+\infinity\}
$$
As we shall see in the arguments below, if $I$ is noncompact on the negative end then $I_- = -\infinity$, in other words $\Length(G_r)$ diverges to $+\infinity$ as $r \to I_-$. But we shall shortly give an example in which $I$ is noncompact on the positive end, $\gamma$ is indeed proper on the positive end as required, and yet $I_+ \ne +\infinity$, that is, $\Length(G_r)$ may converge to a positive number as $r \to I_+$.

\subsection{Examples of fold paths}
\label{SectionFoldExamples}

In this section we give several examples and constructions of fold paths. For each nongeometric, fully irreducible $\phi \in \Out(F_r)$ we construct a periodic fold line through every train track for $\phi$, and we construct a fold line for $\phi$ passing through every weak train track for $\phi$; the latter construction proves that Definition~(1) of the axis bundle of $\phi$ implies Definition~(2). In order to motivate some of the later results, we shall also describe a fold ray whose length does not converge to zero, and a split ray which does not converge to a unique point in compactified outer space. 

\paragraph{Example: periodic fold lines.} Given a nongeometric, fully irreducible $\phi \in\Out(F_r)$ with expansion factor $\lambda=\lambda(\phi)$, we outlined in the introduction the construction of a $\phi$-periodic fold line. To give the construction in detail, consider a train track map $g\from G\to G$ representing $\phi$. Normalize $G \in \wh\X_r$ by requiring $\Length(G)=1$. Let $G_i = \frac{1}{\lambda^i} G \cdot \phi \in \wh\X_r$. For integers $i \le j$ let $g_{ji} \from G_i \to G_j$ be the marked homotopy equivalence which agrees with $g^{j-i} \from G \to G$. Following the method described in Section~\ref{SectionFolds}, choose a Stallings fold sequence for $g_{1,0} \from G_0 \to G_1$, choose an interpolating path for each fold map in that sequence, and concatenate these paths to get a Stallings fold path $\wh\gamma^0$ connecting $G_0$ to $G_1$. The path $\wh\gamma^0$, when normalized, becomes $\wh\gamma^0 \from [0,\log\lambda] \to \wh\X_r$ with $\wh\gamma^0(0)=G_0$ and $\wh\gamma^0(\log\lambda)=G_1$. For each integer $i$ define the path $\wh\gamma^i \from [i \log \lambda, (i+1) \log\lambda] \to \wh\X_r$ to be $\wh\gamma^i(t) = (\wh\gamma^0 \cdot \phi^i)(t - i \log\lambda)$. Define the path $\wh\gamma \from \reals \to \wh\X_r$ to be the concatenation of $\wh\gamma^i$ over all integers~$i$. Connection data for $\wh\gamma$ is easily produced starting from the maps $g_{ji}$ and connection data for $\wh\gamma^0$. Define the path $\gamma \from \reals \to \X_r$ by projectivizing $\wh\gamma$. To prove that $\gamma$ is proper, it suffices to prove that $\gamma(t)$ converges in $\overline\X_r$ to $\TMinusphi$ as $t \to -\infinity$ and to $\TPlusphi$ as $t \to +\infinity$: by construction the path $\gamma$ satisfies the $\phi$-periodicity condition $\gamma(t - i \log(\lambda)) = \gamma(t) \cdot \phi^i$, and the desired convergence properties follow from the uniform source-sink dynamics of the action of $\phi$ on compact subsets of $\X_r$ \cite{BFH:laminations}.

\paragraph{Example: Weak train tracks.} We prove here that Definition~(2) of the axis bundle implies Definition~(1), by constructing a fold line through every weak train track:


\begin{proposition}\label{PropWTTOnFoldLine}
Let $\phi \in \Out(F_r)$ be fully irreducible and let $\A_\phi \subset \X_r$ denote the set of weak train tracks of $\phi$. For each $H \in \A_\phi$ there is a fold line $\gamma \from \reals \to \X_r$ passing through $H$ such that $\gamma(t) \to T_\pm=T_\pm^\phi$ as $t \to \pm\infinity$. 
\end{proposition}

\begin{proof} Choose a normalization of $\TPlus$. Let $\wh\A_\phi \subset \wh\X_r$ be the corresponding normalization of $\A_\phi$, as discussed in the beginning of Section~\ref{SectionTopology}. Denote the normalizations of $H$ and $f$ by $H_0 \in \wh\A_\phi$ and by $f_0 \from \wt H_0 \to T_+$. 

As in the proof of Lemma~\ref{LemmaFoldBound}, let $h_1 \from H_0 \to H_1$ be a fold map whose lift $\ti h_1 \from \wt H_0 \to \wt H_1$ factors the map $f_0$ with induced map $f_1 \from \wt H_1 \to T_+$, and such that ($h_1$ and) $\ti h_1$ are complete, meaning that if $e, e' \subset \wt H_0$ are oriented edges being folded by $\ti h_1$ then the maximal possible initial segments of $e$ and $e'$ that are identified in $T_+$ are identified immediately in $\wt H_1$. Repeat this operation to define an infinite sequence of fold maps
$$
H_0 \xrightarrow{h_1} H_1 \xrightarrow{h_2} H_2 \xrightarrow{h_3} \cdots
$$
with $H_j \in \wh\X_r$, and maps $f_i \from \wt H_i \to T_+$, such that $f_{i-1} = f_{i} \composed \ti h_{i}$ for integers $i \ge 1$. By induction $f_i$ is a $\Lambda_-$ isometry and so $H_j \in \wh\A_\phi$. By construction, any finite segment of this infinite sequence is a sequence of complete Stallings folds in the sense of Section~\ref{SectionFolds}. Applying Lemma~\ref{LemmaFoldBound} it follows that the combinatorial length of the composition $h_i \composed \cdots \composed h_1 \from H_0 \to H_i$ goes to infinity as $i \to +\infinity$. Lemma~\ref{bounded combinatorial length} therefore implies that the $H_i$'s eventually leave every compact subset of $\A_\phi$. Proposition~\ref{PropAxisBundleCompact} implies that there is a compact fundamental domain for the action of $\phi$ on $\A_\phi$, and the source--sink action of $\phi$ on compact subsets of $\X_r$ \cite{BFH:laminations} then implies that the closure of $\A_\phi$ in $\overline\X_r$ is $\A_\phi \union \{\TMinus,\TPlus\}$. The sequence $H_i$ therefore accumulates on $\{\TMinus,\TPlus\}$ as $i \to +\infinity$. But $\Length(H_i)$ is monotonically decreasing, and so $H_i$ converges to $\TPlus$.

By interpolating the fold maps $h_{i}$ we produce a fold ray that begins at $H$ and converges to $T_+$.  By Lemma~\ref{bounded combinatorial length} there is $\Lambda_-$ isometry from some train track $G$ to $H$.  Factoring this as a sequence of folds and interpolating each fold produces a fold ray that begins with $G$, passes through $H$ and converges to $T_+$.  Finally, preconcatenating with a periodic split ray that ends at $G$ produces the desired fold line $\gamma$.
\end{proof}

\paragraph{Example: A fold ray with length $\not\to 0$.} To construct such an example observe that you can do a lot of folds without driving the length to zero, for example by only allowing folds supported in a proper subgraph. For a cheap example, consider any fold ray $\gamma \from [0,\infinity) \to \X_r$, with $G_s$ representing $\gamma(s)$. Choose $x \in G_0$. For each $s \ge 0$, let $G'_s$ be a marked graph of rank $r+1$ obtained from $G_s$ by attaching a loop of length $1$ based  at the point $g_{s0}(x) \in G_s$ a loop of length $1$. This yields a fold ray $\rho' \from [0,\infinity) \to \X_{r+1}$ with $\rho'(s)$ represented by $G'_s$, but in the positive direction $\Length(G'_s)$ converges to $1$, not to zero.  

Although we shall prove below that every fold ray in outer space converges to a unique limit point in the boundary of outer space, this example shows that the limiting $F_r$-tree may not have dense branch points.

\paragraph{Example: A nonconvergent split ray.} Let $\rho \from (-\infinity,0] \to \X_r$ be a split ray in rank $r$ outer space $\X_r$ that converges to a tree $S \in \bdy\X_r$, and let $\sigma \from (-\infinity,0] \to \X_{r'}$ be a split ray in rank $r'$ outer space $\X_{r'}$ that converges to $T \in \bdy\X_{r'}$. Choose connections: a rank $r$ marked graph $G_t$ representing $\rho(t)$ and connection maps $g_{st} \from G_{t} \to G_s$; and a rank $r'$ marked graph $H_t$ representing $\sigma(t)$ and connection maps $h_{st} \from H_{t} \to H_s$. Choose base points $x_s \in G_s$ so that $g_{st}(x_{t}) = x_s$ and $y_t \in H_t$ so that $h_{st}(y_{t}) = y_s$. Identify $F_r = \pi_1(G_s,x_s)$ and $F_{r'} = \pi_1(H_t,y_t)$. Let $G_s * H_t$ denote the graph obtained from the disjoint union of $G_s$ and $H_t$ by identifying $x_s$ to $y_t$ to form a base point $xy_{st}$, and so we have identified $\pi_1(G_s * H_t,xy_{st}) = F_{r+r'}$, thus making $G_s * H_t$ into a rank $r+r'$ marked graph. Consider the third quadrant $(-\infinity,0] \cross (-\infinity,0]$ in $s,t$-space. Any continuous path starting from $(s,t)=(0,0)$ for which $s+t$ decreases strictly monotonically to $-\infinity$ defines a split ray in $\X_{r+r'}$, using $\rho$ to split the $G$ factor and $\sigma$ to split the $H$ factor. Holding $t=0$ and letting $s \to -\infinity$, we obtain a split ray in $\X_{r+r'}$ converging in $\overline\X_{r+r'}$ to a tree $\hat S$ on which the $F_{r'}$ factor of $F_{r+r'}=F_r * F_{r'}$ has a fixed point, and so that the $F_r$ action on $\hat S$ is isometrically conjugate to $S$. Similarly, holding $s=0$ and letting $t \to -\infinity$, we obtain a split ray in $\X_{r+r'}$ converging in $\overline\X_{r+r'}$ to a tree $\hat T$ on which $F_r$ acts with a fixed point, and $F_{r'}$ acts isometrically conjugate to $T$. Now construct a split ray as follows: start from $G_0 * H_0$; hold $t$ constant and let $s$ decrease, splitting the $G$ factor, so as to approach very close to $\hat S$; next, hold $s$ constant and let $t$ decrease, splitting the $H$ factor, so as to approach very close to $\hat T$; alternate back and forth in this manner, coming ever closer to $\hat S$ during the times when $G$ is split, and ever closer to $\hat T$ when $H$ is split. The result is a split ray in $\X_{r+r'}$ that accumulates on both $\hat S$ and $\hat T$.

\subsection{Characterizing fold lines} 

We have shown in Proposition~\ref{PropWTTOnFoldLine} that Definition~(2) of the axis bundle of a fully irreducible $\phi \in \Out(F_r)$ implies Definition~(1): every weak train track for $\phi$ lies on some fold line going from $T_-$ to $T_+$. The main result of this section is the converse: each point on each fold line going from $T_-$ to $T_+$ is a weak train track for $\phi$, thus proving that Definition~(1) implies Definition~(2).

\begin{theorem}
\label{TheoremFoldLines}
Given a fully irreducible $\phi \in \Out(F_r)$ with attracting and repelling fixed points $T_+,T_- \in \bdy\X_r$, if $\gamma\from I\to \X_r$ is a normalized fold line such that 
$$\lim_{t \to I_+} \gamma(t) = T_+ \quad\text{and}\quad \lim_{t \to I_-}
\gamma(t) = T_- \quad\text{in}\quad \overline\X_r
$$
then $I=\reals$ and for each $t \in \reals$, $\gamma(t)$ is a weak train track for $\phi$.
\end{theorem}

Since $\Length(\gamma(t))=e^{-t}$ it immediately follows that: 

\begin{corollary}
\label{CorollaryLengthOfFoldPath}
The map $\Length \from \gamma(\reals) \to (0,\infinity)$ is a homeomorphism.
\qed\end{corollary}

The proof of Theorem~\ref{TheoremFoldLines} takes up the rest of Section~\ref{SectionFoldLines}. First we set up the proof, breaking it into three steps, which are carried out in the following subsections.

Choose connection data for $\gamma$: for each $t \in I$ a marked graph $G_t$ of length $e^{-t}$ representing $\gamma(t)$, and for each for $s<t \in I$ an edge isometry $g_{ts}\from G_s\to G_t$, so that $g_{ts}\composed g_{sr} = g_{tr}$. Let $T_t$ denote the universal cover of $G_t$ for each $t\in I$. Choose a lift $f_{ts} \from T_s\to T_t$ of each $g_{ts}$, so that $f_{ts}\composed f_{sr} = f_{tr}$. 

We must prove that for each $t \in I$ there is a map $f_t^+  \from T_t \to T_+$ which is injective when restricted to each leaf of $\Lambda_-$ realized in $T_t$. Along the way we also prove that $I=\reals$. There are three steps to the proof.

\paragraph{Step 1: Direct limits of fold rays (Section \ref{SectionDirectLimits}).} Recall the discussion of direct limits in Section~\ref{SectionAttractingTree}. Recall also Theorem~\ref{T+ is direct limit} which says that from a train track representative $g \from G \to G$ of $\phi$ one obtains a sequence of marked graphs $G_i = \frac{1}{\lambda(\phi)^i} G \cdot \phi^i$ whose universal covers fit into a direct sequence with direct limit $T_+$. In other words, $T_+$ is the direct limit of any $\phi$-periodic fold ray. 

In this section we prove the generalization that for \emph{any} fold ray which converges to $T_+$, this convergence is manifested as a direct limit. Regarding the fold ray as a direct system, we first prove that its direct limit is an $F_r$-tree, and then that $T_t$ converges in $\wh\X_r$ to its direct limit. Using the hypothesis that $T_t$ converges to $T_+$, we may identify $T_+$ with the direct limit. We obtain direct limit maps $f_t^+ \from T_t \to T_+$ satisfying the property that $f_t^+ \composed f_{ts} = f_s^+$ for all $s<t$. These maps $f_t^+$ will be edge isometries. Then, using properties of $T_+$, we will prove that $I_+=+\infinity$.

\paragraph{Step 2: Legal laminations of split rays (Section \ref{SectionLegalLams}).} Using the fact that $T_t$ converges to $T_-$ in $\overline\X_r$ at $t \to I_-$ we first prove that the injectivity radius of $G_t$ diverges to $+\infinity$ as $t \to I_-$. This immediately implies that $\Length(G_t) \to +\infinity$ and therefore $I_- = -\infinity$. Next we focus on the ``legal lamination'', a sublamination of the geodesic lamination of $F_r$, defined as the set of all leaves $\ell$ of the geodesic lamination such that for all $r<s$ in $I$, the realization of $\ell$ in $G_r$ is legal for the map $f_{sr}$. We prove that the legal lamination is nonempty, and so by the Hausdorff maximal principle has a minimal nonempty sublamination $\L$. We then prove that this sublamination has zero length in $T_-$ in the sense of \BookZero, and quoting a result of \BookZero\ (Lemma~\ref{LemmaMinimalLamination} above) it follows that $\L = \Lambda_-$.

\paragraph{Step 3: Fold lines through weak train tracks (Section \ref{SectionStep3}).} Putting Steps~1 and~2 together, we prove that the realization of $\Lambda_-$ in each $G_t$, when lifted to $T_t$, is legal for the map $f^+_t$, and so $G_t$ is a weak train track.

\subsection{Direct limits of fold rays} 
\label{SectionDirectLimits}
In this section we carry out Step~1 in the proof of Theorem~\ref{TheoremFoldLines}, that a fold ray converging to $T_+$ may be regarded as a direct system with direct limit $T_+$.

Consider a normalized fold path $\rho \from I \to \X_r$ which is noncompact in the positive direction. We shall make no other assumptions on $\rho$ until near the end of the argument when we bring in the hypothesis that $\rho(s) \to T_+$ as $s \to I_+$. 

Choose a normalized connection for $\rho$: marked graphs $G_t$ representing $\rho(t)$ for each $t\in I$ so that $\Length(G_t) = e^{-t}$; and edge isometries $g_{ts}\from G_s \to G_t$ so that $g_{ts} \composed g_{sr} = g_{tr}$. Denote the universal covers as $T_t = \wt G_t$. We can lift the map $g_{ts}$ to an $F_r$ equivariant map $f_{ts} \from T_s \to T_t$  so that $f_{ts} \composed f_{sr} = f_{tr}$. The maps $f_{ts}$ are all distance nonincreasing edge isometric maps, and so we have a direct system $(T_i)_{i \in I}$, $(f_{ts} \from T_s \to T_t)_{s<t \in I}$. Let $T_\infinity$ be the direct limit and let $f^\infinity_t \from T_t \to T_\infinity$ be the direct limit maps, distance nonincreasing maps such that $f^\infinity_s \composed f_{st} = f^\infinity_t$. Since each map $f_{ts}$ restricts to an isometry on each edge of $T_s$, the map $f^\infinity_s$ also restricts to an isometry on each edge. (Despite the notation we do \emph{not} assume at this stage that $I_+ = +\infinity$; this will be proved at the appropriate point). The $F_r$ actions on the trees $T_t$ yield an $F_r$ action on $T_\infinity$ such that each map $f^\infinity_t$ is $F_r$ equivariant: each $\gamma \in F_r$ acts isometrically on each $T_s$ compatible with the maps $f_{st}$, and the existence clause of the universality property gives the action of $\gamma$ on $T_\infinity$; the uniqueness clause gives the action equation $(\gamma\cdot\delta) \cdot x = \gamma \cdot (\delta\cdot x)$ for $\gamma,\delta \in F_r$, $x \in T_\infinity$. Since $\gamma$ and $\gamma^\inv$ are both distance nonincreasing on $T_\infinity$, and since $\gamma\gamma^\inv = \gamma^\inv\gamma = \Id$, it follows that $\gamma$ acts as an isometry on $T_\infinity$.

\paragraph{$T_\infinity$ is an $\reals$-tree.} First observe that the 0-hyperbolic inequality
$$d(w,x) + d(y,z) \le \Min\{d(w,z) + d(x,y), d(w,y) + d(x,z) \}
$$
holds in each tree $T_s$, and so evidently holds as well in $T_\infinity$ after taking infimums. It follows that $T_\infinity$ embeds isometrically in an $\reals$-tree (see e.g. \cite{GhysHarpe:afterGromov}). 

To prove that $T_\infinity$ is an $\reals$-tree, it remains to show that if $\xi,\eta \in T_\infinity$ and if $a,b > 0$ and $a+b=d(\xi,\eta)$ then there exists $\zeta \in T_\infinity$ such that $d(\xi,\zeta) = a$ and $d(\zeta,\eta) = b$.  Fix $s_0 \in I$, and choose $x_{s_0},y_{s_0} \in T_{s_0}$ such that $f^\infinity_0(x_0)= \xi$ and $f^\infinity_0(y_0)= \eta$.   For each $s \ge s_0$   let $x_s = f_{s0}(x_{s_0})$ and $y_s = f_{s0}(y_{s_0})$.   For all $s_0 \le s<t \in I$ we have $f_{ts}(x_s)=x_t$, $f_{ts}(y_s)=y_t$, $f^\infinity_s(x_s) = \xi$ and $f^\infinity_s(y_s) = \eta$.  For each $s \ge s_0$ choose $v^s \in \overline{x_{s_0} y_{s_0}}$ so that $f_{s, s_0}(v^s) \in \overline{x_s y_s}$ and so that $d(x_s,f_{s, s_0}(v^s)) \ge a$ and $d(f_{s, s_0}(v^s),y_s) \ge b$; this is possible because $\overline{x_s y_s} \subset f_{s,s_0}\left( \overline{x_{s_0} y_{s_0}} \right)$ and because $d(x_s,y_s) \ge d(\xi,\eta) = a+b$. Some subsequence of the sequence $v^s$ converges to a point $z_{s_0} \in \overline{x_{s_0} y_{s_0}}$. Set $z_s = f_{s, s_0}(z_{s_0})$ and $\zeta = f^\infinity_{s_0}(z_{s_0}) = f^\infinity_s(z_s)$, for $s \ge s_0$. We shall show that $\zeta$ is the desired point. 

Observe that for all $s \ge s_0$ we have $d(x_s,z_s) \ge a$. Indeed, if there existed $s \ge s_0$ such that $d(x_s,z_s) < a$ then for all sufficiently large $s' > s$ we would have $d(x_s,f_{s, s_0}(v^{s'})) < a$ which would imply that $d(x_{s'},f_{s', s_0}(v^{s'})) < a$, a contradiction.  By a symmetric argument we have $d(z_s,y_s) \ge b$. Taking infimums as $s \to \infinity$, it follows that $d(\xi,\zeta) \ge a$ and $d(\zeta,\eta) \ge b$.

Next observe that $d(x_s,z_s) + d(z_s,y_s) = d(x_s,y_s)$, and taking infimums as $s \to \infinity$ we obtain $d(\xi,\zeta) + d(\zeta,\eta) = d(\xi,\eta) = a+b$. Combined with the above, it follows that $d(\xi,\zeta)=a$ and $d(\zeta,\eta) = b$.

\paragraph{$T_\infinity$ is nontrivial and minimal.} 
Given $s<t \in I \union \{\infinity\}$, and given $x \in T_s$, recall the induced map $Df_{ts}$ taking directions  of $T_s$ at $x$ to directions of $T_t$  at $f_{ts}(x)$, and recall that two directions at $x$ are in the same gate relative to $f_{ts}$ if they are identified under $Df_{ts}$.  

\begin{lemma}  \label{at least two gates}   For all sufficiently large $s$, each point in $T_s$ has $\ge 2$ gates relative to $f^\infty_s$.
\end{lemma}

\begin{proof} If $d_1$ and $d_2$ are identified by $Df^\infty_s$ then they are identified by $Df_{ts}$ for some finite $t$. It therefore suffices to prove that if $s$ is sufficiently large then each point in $T_s$ has at least two gates relative to $f_{ts}$ for all finite $t>s$. The property of being in the same gate is invariant under the action of $F_r$ so it suffices to prove that if $s$ is sufficiently large then for each $x \in G_s$ and $t>s$ there are directions $d_1$ and $d_2$ at $x$ that are not identified by $Dg_{ts}$.

Choose $a \in I$. For each $x \in G_a$ consider the following property: 
\begin{description}
\item[$(*)$] There exists $c \ge a$ such that if $c \le s \le t$ then $g_{sa}(x)$ has $\ge 2$ gates relative to $g_{ts}$. 
\end{description}
If $x$ is not a vertex of $G_a$ then Property~$(*)$ is true with $c=a$, as a consequence of Property~(\ref{ItemTrainTrack}) in the definition of a fold path. The set of vertices being finite, to complete the proof it remains to verify $(*)$ for each vertex. If $x$ is a vertex, and if there exists $t \ge a$ such that $g_{ta}^\inv(g^{\vphantom\inv}_{ta}(x))$ contains a point that is not a vertex of $G_a$, then Property~$(*)$ is true with $c=t$. Consider the set $V_a$ of vertices $x$ of $G_a$ such that for all $t \ge a$, the set $g_{ta}^\inv(g^{\vphantom\inv}_{ta}(x))$ consists solely of vertices. As $t$ increases the sets $V_t=g_{ta}(V_a)$ have nonincreasing cardinality, and so for some $c \ge a$ the cardinality of $V_t$ is constant for $t \ge c$. Suppose that $c \le s \le t$, $x \in V_a$ and consider $y=g_{ca}(x) \in V_c$. We have shown that $g_{tc}^\inv(g^{\vphantom\inv}_{tc}(y))=\{y\}$. Since $g_{tc}$ maps $G_c$ onto $G_t$, it maps a neighborhood of $y$ in $G_c$ onto a neighborhood of $g_{tc}(y)$ in $G_t$, and so there must be at least two directions at $y$ that are not identified by $Dg_{tc}$. Property~(\ref{ItemFold}) in the definition of a fold path implies that the $Dg_{sc}$ images of these two directions are not identified by $Dg_{ts}$.
%
%
\end{proof}

For the remainder of the proof we fix $s$ satisfying Lemma~\ref{at least two gates} and we focus on  $f^\infinity_s$.   Given $x \in T_s$, a pair of directions  at $x$ form a \emph{legal turn} if they are not in the same gate relative to $f^\infinity_s$.   A bi-infinite geodesic $\ell\subset T_s$ is \emph{legal} if, for every $x \in \ell$, the turn made by $\ell$ at $x$ is legal; equivalently, the restriction of $f^\infinity_s$ to $\ell$ is an isometric embedding $\ell \to T_\infinity$. By Lemma~\ref{at least two gates},   every direction at every point $x \in T_s$ forms a legal turn with some other direction at $x$. This implies that every point in $T_s$ lies on some bi-infinite legal geodesic: starting with any edge in $T_s$, we may extend the edge at each endpoint legally, and continuing to extend by induction we obtain a bi-infinite legal geodesic.

To prove that the action of $F_r$ on $T_\infinity$ is nontrivial and minimal, consider $\xi\in T_{\infinity}$, and choose $x \in T_s$ such that $f^\infinity_s(x) = \xi$. Let $\ell \subset T_s$ be a bi-infinite legal geodesic passing through~$x$ (in this context we do not assume $\ell$ is a leaf of~$\Lambda_-$). Let $\ell_-,\ell_+$ be the components of $\ell-x$. Since there are only finitely many $F_r$-orbits of oriented edges in $T_s$, the ray $\ell_-$ must pass over two oriented edges $e\ne e'\subset T_s$ in the same $F_r$ orbit. Let $\gamma_- \in F_r$ take $e$ to $e'$. Construct a $\gamma_-$-invariant line by gluing together translates of the segment of $\ell$ connecting $e$ to $e'$.  This line must be the axis $A_s(\gamma_-)$ of $\gamma_-$ in $T_s$, and the interior of the segment of $\ell$ from $e$ to $e'$ contains a fundamental domain for $A_s(\gamma_-)$, implying that $A_s(\gamma_-)$ is legal. Doing the same in the other direction, we obtain $\gamma_+\in F_r$ whose axis $A_s(\gamma_+)$ is legal. It follows that $\gamma_\pm$ are hyperbolic on $T_\infinity$ and their axes $A_\infinity(\gamma_\pm)$ are the isometric images under $f^\infinity_s$ of the axes $A_s(\gamma_\pm)$; in particular, we now know that the action of $F_r$ on $T_\infinity$ is nontrivial. By a standard exercise in group actions on $\reals$-trees, the axes $A_\infinity(\gamma_\pm)$ are contained in the minimal subtree of $T_\infinity$. By construction, there is a segment of $T_s$ passing through $x$ with endpoints respectively on $A_s(\gamma_\pm)$, and so there is a segment of $T_\infinity$ passing through $\xi$ with endpoints respectively on $A_\infinity(\gamma_\pm)$, implying that $\xi$ is contained in the minimal subtree of $T_\infinity$. Since $\xi \in T_\infinity$ is arbitrary, we proved that the action of $F_r$ on $T_\infinity$ is minimal.

\paragraph{Convergence in $\reals^\C$.} Next we show that $T_s$ converges to $T_\infinity$ in $\reals^\C$, that is, for each $\gamma \in F_r$ the translation length of $\gamma$ acting on $T_s$ converges to the translation length of $\gamma$ acting on $T_\infinity$, as $s \to I_+$. 

Consider the axis $A_s(\gamma) \subset T_s$. For $s<t$ we have $f_{ts}(A_s(\gamma)) \supset A_t(\gamma)$. We may therefore choose fundamental domains $[x_s,y_s] \subset A_s(\gamma)$ for $s \in I$ so that $f_{ts}(x_s)=x_t$, $f_{ts}(y_s)=y_t$. Letting $\xi = f^\infinity_s(x_s)$, $\eta = f^\infinity_s(y_s)$, it follows that $d(\xi,\eta)$ equals the infimum of the translation lengths. But if $\xi \ne \eta$ then clearly $\union_{i \in \Z} \gamma^i [\xi,\eta]$ for $i \in \Z$ is an axis for $\gamma$ with translation distance $d(\xi,\eta)$; whereas if $\xi = \eta$ then $\gamma(\xi)=\eta$ and so $\gamma$ has translation distance zero.

\subparagraph{Remark.} So far we have not used the hypothesis that $T_s$ converges to $\TphiPlus$. Our arguments therefore give a completely general proof that every fold ray converges in $\overline\X_r$ to a unique point $T_\infinity \in \bdy\X_r$, and that limit may be regarded as a direct limit. We described earlier an example of a fold ray for which length did not converge to zero. For this fold ray it easily follows that limiting $F_r$-tree $T_\infinity$ has nonzero ``colength'', meaning that there is a positive lower bound to the length of any finite subtree $\tau \subset T_\infinity$ whose translates under $F_r$ cover $T_\infinity$. This implies in turn that $T_\infinity$ has a nondegenerate arc whose interior is an open subset of $T_\infinity$, and so the action of $F_r$ on $T_\infinity$ is not transitive. The action of $F_r$ on $T_\infinity$ therefore has a nontrivial decomposition, in the sense of \cite{Levitt:graphs}, as a finite graph of actions with dense branch points.

\paragraph{Identification of $T_\infinity$ with $T_+$.} We have shown that $T_s$ converges to $T_\infinity$ in $\reals^\C$ as $s \to I_+$. Projectivizing, it follows that (the projective class of) $T_\infinity$ is a point in $\overline\X_r$, and that $T_s$ converges to $T_\infinity$ in $\overline\X_r$ as $s \to I_+$. But now we use the hypothesis that $T_s$ converges to $T_+$ as $s \to I_+$, to conclude that $T_\infinity$ equals $T_+$ in $\X_r$, in other words, there exists an $F_r$-equivariant homothety $T_\infinity \approx T_+$. By scaling the metric on $T_+$ we may conclude that this homothety is an isometry.

We thus obtain edge-isometric morphisms $f^+_s = f^\infinity_s \from T_s \to T_+ = T_\infinity$ for $s \in I$, with $f^+_s = f_t^+ \composed f_{ts}$ for $s<t$ in $I$.

\paragraph{$I_+$ equals $+\infinity$.} From the definition of a normalized fold line combined with the definition of the direct limit, it follows that $e^{-I_+}$ equals the limit of the lengths of the marked graphs $G_s = T_s / F_r$ as $s \to I_+$. We show that this limit equals zero. If not, then there exists an open subarc $a_s \subset G_s$ for some $s$, so that for all $t \ge s$, if $x \in a_s$, $y \in G_s$, and $g_{ts}(x)=g_{ts}(y)$ then $x=y$. Letting $A_s \subset T_s$ be a connected lift of $a_s$, it follows that for all $t \ge s$, if $x \in A_s$, $y \in T_s$, and $f_{ts}(x)=f_{ts}(y)$ then $x=y$. The open arc $A_s$ therefore embeds as an open subarc of $T_+$ under the direct limit map $f^+_s \from T_s \to T_+$. However, each orbit of the action of $F_r$ on $T_+$ is dense and so no such open arc exists in $T_+$. 

\bigskip

This completes Step 1.

\subsection{Legal laminations of split rays}
\label{SectionLegalLams}  

In this section we carry out Step~2 in the proof of Theorem~\ref{TheoremFoldLines}: for any fold ray that converges to $T_-$, the lamination $\Lambda_-$ arises naturally as the ``legal lamination'' of that fold ray. One can think of this as a generalization of Lemma~\ref{first characterization}, in which it was shown that $\Lambda_-$ is the unique minimal lamination whose leaves are legal with respect to any train track representative of $\phi$.


We review the notation: $\rho \from I \to \X_r$ is a normalized fold path assumed to be noncompact in the negative direction; there is a normalized connection for $\rho$ consisting of marked graphs $G_t$ representing $\rho(t)$ for each $t\in I$ so that $\Length(G_t) = e^{-t}$, and edge isometric morphisms $g_{ts}\from G_s \to G_t$ so that $g_{ts} \composed g_{sr} = g_{tr}$. Denote the universal covers as $T_t = \wt G_t$, and lift the map $g_{ts}$ to an $F_r$ equivariant map $f_{ts} \from T_s \to T_t$ of $g_{ts}$ so that $f_{ts} \composed f_{sr} = f_{tr}$.

\paragraph{Injectivity radius.}
Define the \emph{injectivity radius} of $G_t$, denoted $\inj(G_t)$, to be the supremum of $\epsilon>0$ such that for all $x \in G_t$, the $\epsilon$-neighborhood $N_\epsilon(x)$ is a tree and the inclusion map $N_\epsilon(x) \inject G_t$ is an isometry with respect to the path length metric on $N_\epsilon(x)$. Note that $\inj(G_t)$ equals half the length of the shortest closed curve on~$G_t$.

For the remainder of this section, choose $\mu_t>0$ for each $t \in I$ so that
$$\lim_{t \to I_-} \mu_t T_t = T_- \quad\text{in}\quad \whclosure\X_r
$$ 

\begin{proposition}
\label{PropInjRad} As $t \to I_-$ we have: \quad (1)~$\Length(\mu_t G_t) \to 0$; \quad (2)~$\mu_t \to 0$; \quad (3)~$\inj(G_t) \to +\infinity$; \quad and (4)~$\Length(G_t) \to +\infinity$.
\end{proposition}


\begin{proof} Note that the $F_r$-tree $T_-$ has an arbitrarily short fundamental domain, meaning that for each $\epsilon>0$ there exists a point $p \in T_-$ and a free basis $g_1, \ldots, g_r$ of~$F_r$ such that the distances $d_{T_-}(p,f_i \cdot p)$ sum to less than $\epsilon$. This is true on general principles because $T_-$ has dense branch points, and it can be seen very easily in this case by starting with arbitrary $p,g_1,\ldots,g_r$, iterating the action on $T_-$ of any automorphism representing $\phi$, and using the fact that $\phi$ contracts distance on $T_-$ by a factor of $\lambda(\phi^\inv) > 1$.

Since $\mu_t T_t$ converges to $T_-$ as $t \to I_-$ in the Gromov--Hausdorff topology on closed subsets of $T_-$, it follows that for any $\epsilon>0$, if $t$ is sufficiently close to $I_-$ then there exists $p \in T_t$ and a free basis $g_1, \ldots, g_r$ of $F_r$ such that the distances $d_{\mu_t T_t}(p,f_i \cdot p)$ sum to less than $\epsilon$. From this it follows immediately that $\Length(\mu_t G_t) \to 0$ as $i \to I_-$. But $\Length(G_s) \ge \Length(G_t)$ for $s<t \in I$, and so $\lim_{t \to I_-} \Length(G_t)$ exists and is $>0$, possibly $+\infinity$. It follows that $\ds\lim_{t \to I_-} \mu_t = 0$.

Suppose that $\inj(G_t) \not\to +\infinity$ as $t \to I_-$. Since $\inj(G_t)$ is a nonincreasing function of $t$, it follows that $\inj(G_t)$ converges to some $R>0$ as $t \to I_-$. Recall that $\C$ is the set of nontrivial conjugacy classes in $F_r$. Let $\C_t = \{C \in \C \suchthat \Length(C;G_t) \le 2R\}$, so $\C_t$ is finite and nonempty for all $t \in I$. If $s<t \in I$ then $\C_s \subset \C_t$, because $g_{ts}$ is distance nonincreasing.  It follows that $\C_t$ stabilizes for $t$ sufficiently close to $I_-$, say $t \le t_0$. Let $\gamma$ be the minimal length closed curve in $G_{t_0}$ and let $C \in \C_{t_0}$ be the conjugacy class of $\gamma$, so $\Length(C;G_t) \le 2R$ for all $t \in I$. Note that $\gamma$ is simple and so each representative of $C$ is a generator of a free factor of $F_r$. Applying Theorem~5.4 of \BookZero, it follows that $\Length(C;T_-) > 0$.  
On the other hand,
\begin{align*}
\Length(C;T_-) &= \lim_{t \to I_-} \Length(C;\mu_t T_t) 
   = \lim_{t \to I_-} \mu_t \Length(C;T_t) \\
 &\le \lim_{t \to I_-} \mu_t 2 R 
   = 0
\end{align*}
a contradiction.

Since $\inj(G_t) \to +\infinity$ as $t \to I_-$, it follows as well that $\Length(G_t) \to +\infinity$. 
\end{proof}

Combining Proposition~\ref{PropInjRad} with the normalization condition on the fold ray, we have proved one of the pieces of Theorem~\ref{TheoremFoldLines}:

\begin{proposition}
\label{PropMinusInfinity}
$I_- = -\infinity$.
\qed\end{proposition}

\subparagraph{Remark.} \quad We have proved that $\lim_{t \to I_-} \Length(G_t) = +\infinity$ via a proof that $\inj(G_t) \to +\infinity$, which itself uses the hypothesis that the split ray converges to $T_-$. But $\inj(G_t) \not\to +\infinity$ for an arbitrary normalized split ray, as shown by the example described earlier of a split ray in which a certain positive rank subgraph denoted $G_0$ does not change, as the length of the rest of the subgraph goes to $+\infinity$. 

Nevertheless, here is a proof that $\lim_{t \to I_-} \Length(G_t) = +\infinity$ which is valid for an arbitrary split ray. If not then by monotonicity we have $\lim_{t \to I_-} \Length(G_t) = L<+\infinity$ and $\Length(G_t) \le L$ for all $t \in I$. Using this we show for each $C \in \C$ that $\Length(C;G_t)$ an upper bound independent of $t$, for then by monotonicity the translation length function  of $G_t$ has a limit in $\reals^\C$ whose coordinates are bounded away from zero, and by Lemma~\ref{LemmaDiscreteTrees} this corresponds to a point in $\X_r$, contradicting properness of split rays. To bound $\Length(C;G_t)$, fix $c \in I$, and for $t \le c \in I$ consider a segment $\alpha \subset T_t$ whose endpoints are identified by $f_{ct}$. We have $\Length(G_t) - \Length(G_c) \ge \frac{1}{2} \Length(\alpha)$, and so $\Length(\alpha) \le 2(L-\Length(G_c))$. So if $\gamma \in F_r$ represents $C$, and if the axis of $\gamma$ in $T_c$ has $K$ vertices in a fundamental domain, it follows that $\Length(C;G_t) = \Length(\gamma;T_t) \le \Length(\gamma;T_c) + 2K(L-\Length(G_c))$.

\paragraph{A minimal legal lamination.} Consider an element $\ell$ of the geodesic leaf space $\G F_r$. Let $\ell_t$ denote the realization of $\ell$ in $T_t$. We say that $\ell$ is a \emph{legal leaf} if for all $s<t \in I$ the restriction of $f_{ts} \from T_s \to T_t$ to 
the geodesic $\ell_s \subset T_s$ is an embedding, and hence is an isometry onto $\ell_t$. The set of legal leaves is a closed, $F_r$-equivariant subset of $\G F_r$, called the \emph{legal lamination} and denoted $\wt\LL$. For each $t \in I$ the realization of $\wt\LL$ in $T_t$ forms an $F_r$-equivariant sublamination $\wt\LL_t \subset \GL T_t$ such that $f_{ts}$ pushes $\wt\LL_s$ forward to $\wt\LL_t$. Let $\LL_t \subset \GL G_t$ be the downstairs image of $\wt\LL_t$, and so $g_{ts}$ pushes $\LL_s$ forward to $\wt\LL_t$.

\begin{fact} $\wt\LL$ is not empty.
\end{fact}

\begin{proof}
Let $E_s$ be the longest edge of $G_s$ and so applying Proposition~\ref{PropMinusInfinity} it follows that $\Length(E_s) \ge \frac{1}{3r-3} \Length(G_s) \to +\infinity$ as $s \to -\infinity$. Fixing $t \in I$, for each $s<t$ the path $f_{ts} \restrict E_s$ is a legal path in $G_t$ of length equal to $\Length(E_s)$, by item~(\ref{ItemTrainTrack}) in the definition of foldline.  Parameterizing this path by a symmetric interval $[-\Length(E_s)/2,\Length(E_s)/2]$ and extending the domain to $\reals$ so that the map is constant on each ray complementary to this interval, we may then pass to a subsequence converging in the compact open topology to a geodesic immersion $\ell_t \from \reals \to G_t$. Since $\ell_t$ is a union of subsegments of paths $f_{ts} \restrict E_s$, it follows by item~(\ref{ItemTrainTrack}) in the definition of foldline that $\ell_t$ is legal with respect to all maps $f_{ut} \from G_t \to G_u$, $u \ge t$.

Fixing a sequence $t_1 > t_2 > \cdots$ diverging to $-\infinity$, by a diagonalization argument we obtain $\ell_{t_i}$ so that $g_{t_1 t_i}(\ell_{t_i}) = \ell_{t_1}$. Lifting these geodesics to the universal covers $T_{t_i}$ we therefore obtain a legal leaf $\ell$.
\end{proof}

We now invoke the Hausdorff maximal principle to conclude that $\wt\LL$ has a minimal nonempty, closed, $F_r$-invariant subset which we denote $\wt\L$. This is seen more clearly in one of the compact spaces $\LL_t \subset \GL G_t$ for any fixed $t$, because the intersection of any nested sequence of nonempty closed sublaminations of $\LL_t$ is clearly a nonempty closed sublamination. The realization of $\wt\L$ in $\GL T_t$ is denoted $\wt\L_t$ and the image downstairs in $\GL G_t$ is denoted $\L_t$. 

Recall that $\L_t$ is said to be \emph{quasiperiodic} if for each compact leaf segment $\alpha \subset \ell$ there exists $k$ such that for each compact leaf segment $\alpha'$, if $\Length(\alpha') > k$ then $\alpha'$ has a subsegment whose corresponding path in $G_t$ is identical with the path corresponding to $\alpha$.

We list here the properties of the sublaminations $\L_t$:

\begin{lemma}
\label{LemmaLegalLamination}
For each $t \in I$, $\L_t$ is a nonempty, minimal sublamination of $\GL G_t$, $\L_t$~is quasiperiodic, each leaf of $\L_t$ is aperiodic, and for each $s<t \in I$, $g_{ts}$ maps $\L_s$ to $\L_t$ inducing a bijection of leaves to leaves and taking each leaf to its corresponding leaf by an isometry.
\end{lemma}

\begin{proof} For any compact lamination minimality implies quasiperiodicity; for convenience we prove this in the present context. Suppose that $\L_t$ is not quasiperiodic. Fix a leaf segment $\alpha \subset \L_t$ that defies the definition of quasiperiodicity. There is a sequence of leaf segments $\alpha'_n$ with $\Length(\alpha'_n) \to +\infinity$ so that $\alpha'_n$ does not traverse $\alpha$ (we are here confusing a leaf segment with its corresponding path in $G_t$). Passing to a limit, we obtain a nonempty sublamination of $\L_t$ in which no leaf traverses the segment $\alpha$, and so this is a sublamination is proper, contradicting minimality.

To prove that no leaf is periodic, if $\ell_t \subset \L_t$ is a periodic leaf then the corresponding leaf $\ell_s \subset \L_s$ is periodic for all $s<t$ and $\Length(\ell_s) = \Length(\ell_t)$. But this implies that $\inj(G_s) \le \frac{1}{2} \Length(\ell_t)$ for $s<t$, contradicting Proposition~\ref{PropInjRad}~(3).
\end{proof}

\paragraph{Remark.} In an earlier draft we gave a direct but intricate construction which produces sublaminations $\L_t \subset \GL G_t$ satisfying Lemma~\ref{LemmaLegalLamination} without invoking the Hausdorff Maximal Principle, by studying the inverse limit as $t \to -\infinity$ of the graphs $G_t$ in the category of compact singular laminations. A careful desingularization process produces the required sublaminations $\L_t$.

\paragraph{Identifying $\L$ with $\Lambda_-$.}  To complete Step 2, we shall show that $\L$ has length zero in $T_-$. By applying Lemma~\ref{LemmaMinimalLamination} it will follow that $\L = \Lambda_-$.

We must prove that for any marked graph $G \in \hat\X_r$ and every $F_r$-equivariant morphism $h \from \wt G \to T_-$ there exists $C \ge 0$ such that for every leaf $\ell \subset \wt G$ of $\wt\L_G$ and for any $x,y \in \ell$, $d_{T_-}(hx,hy) \le C$. To prove this it suffices to check it for \emph{some} $G \in \hat\X_r$ and \emph{some} $F_r$-equivariant morphism $h \from \wt G \to T$; this is because there exists an $F_r$-equivariant morphism between $\wt G, \wt G'$ for any two $G,G' \in \X_r$, and because any two $F_r$-equivariant morphisms $\wt G \to T$ are $F_r$-equivariantly homotopic. 
We will do this for $G=G_0$ and an arbitrarily chosen $F_r$-equivariant continuous map $\phi_- \from \wt G_0 \to T_0 \to T_-$ such that $\phi_-$ is a homothety when restricted to each edge of~$T_0$.

Continuing with the notation established earlier in Step 2, denote $T_s^* = \mu_s T_s$, and recall that $T_s^*$ converges to $T_-$ in $\whclosure\X_r$ as $s \to -\infinity$; when the metric is not important we will use $T_s$ or $T_s^*$ without discrimination. Let $\L_s$ denote $\L$ realized in $G_s$, and let $\wt\L_s$ denote the realization in $T_s$ or $T_s^*$. We want to connect the property ``$\L$ has length zero in $T_-$'' to properties of $\wt\L_s$ for $s$ close to $-\infinity$.


We claim that there are $F_r$-equivariant continuous maps $\phi_{s} \from T_0 \to T_s^*$ defined for $s \le 0$ satisfying the following property:
\begin{description}
\item[Uniform Length:] \quad
For any finite geodesic segment $\gamma_0 \subset T_0$, denoting $\gamma^\#_s = (\phi_s)_\#(\gamma_0) \subset T_s^*$ and $\gamma^\#_- = (\phi_-)_\#(\gamma_0) \subset T_-$, we have $\Length(\gamma^\#_s;T_s^*) \to \Length(\gamma^\#_-;T_-)$ as $s \to -\infinity$.
\end{description}
To prove this claim, fix a base point in $T_0$, for $R>0$ let $B_R \subset T_0$ denote the set of vertices of $T_0$ that are within distance $R$ of the base point, let $G_R$ denote the set of elements of $F_r$ that move one point of $B_R$ to another, and consider the set $\phi_-(B_R) \subset T_-$. For each $R>0$ and $\epsilon>0$, if $s$ is sufficiently large then there exists a $G_R$-equivariant $\epsilon$-approximation from $\phi_-(B_R)$ to a subset of $T_s^*$, and by postcomposing $\phi_-$ with this $\epsilon$-approximation we obtain a $G_R$-equivariant map $\phi^{R,\epsilon}_s \from B_R \mapsto T_s^*$. We may choose functions $R(s) \to +\infinity$ and $\epsilon(s) \to 0$ so that $\phi^{R(s),\epsilon(s)}_s$ exists for all $s$ sufficiently close to $-\infinity$. We may then extend $\phi^{R(s),\epsilon(s)}_s$ in a unique manner to an $F_r$-equivariant map $\phi_s \from T_0 \to T_s^*$ which is homothetic on each edge. The Uniform Length property for finite geodesic edge paths $\gamma_0 \subset T_0$ follows immediately because the endpoints of $\gamma_0$ lie on $B_R$ for sufficiently large $R$. For an arbitrary geodesic segment of $T_0$, which can be written as a finite geodesic edge path with partial edges concatenated at the beginning and at the end, the Uniform Length property follows easily. 

As a consequence of the Uniform Length property, observe that there are constants $L > 0$, $s_0 \in I$ such that the map $\phi_-$ and each of the maps $\phi_s$ for $s \le s_0$ are $L$-lipschitz: just apply the Uniform Length property to one edge in each of the finitely many edge orbits of $T_0$. Combining this observation with the bounded cancellation lemma~\ref{bcl}, there is a constant $B > 0$ such that:
\begin{description}
\item[Uniform Bounded Cancellation:] The map $\phi_-$ and each of the maps $\phi_s$ for $s \le s_0$ have bounded cancellation constant $\le B$. 
\end{description}
The constants $s_0$ and $B$ as above shall be fixed for the rest of the argument.

Suppose then that $\L$ does not have length zero in $T_-$. It follows that there is a finite leaf segment $\gamma_0$ of $\wt\L_0$ such that $\Length(\gamma^\#_-;T^{\vphantom{\#}}_-) > 6B$. Applying the Uniform Length property above, for $s$ sufficiently close to $-\infinity$ it follows that $\Length(\gamma^{\#\vphantom{*}}_s;T^{*\vphantom{\#}}_s) > 5B$, let's say for all $s \le s_1$ where $s_1 \le s_0$. Since $\L$ is minimal there exists a constant $K$ such that for every finite leaf segment $\beta_0$ of $\wt\L_0$, if $\Length(\beta_0;T_0) \ge K$ then $\beta_0$ contains a translate of $\gamma_0$. By applying Uniform Bounded Cancellation it follows that $\Length(\beta^{\#\vphantom{*}}_s;T^{*\vphantom{\#}}_s) > 3B$ for all $s \le s_1$ and all leaf segments $\beta_0$ of $\wt\L_0$ with $\Length(\beta_0;T_0) \ge K$. We will show that this leads to a contradiction.

Every geodesic segment in $G_s$ of length equal to $\inj(G_s)$ is embedded and passes over at most $2r-2$ vertices, and therefore intersects some edge in a segment of length at least $C_s = \inj(G_s) / (2r-3)$. It follows that every geodesic segment of length $C'_s = \inj(G_s) + 2\Length(G_s)$ contains an entire edge of length at least $\C_s$. There are at most $2r-1$ edges in $G_s$ and so every geodesic segment of length $C''_s = [2(2r-1)+1]C'_s$ goes at least three times over some edge whose length is at least $C_s$. As $s \to -\infinity$ we have $\inj(G_s) \to +\infinity$ and so $C_s \to +\infinity$. As $s \to -\infinity$ we also have $\mu_s \Length(G_s) \to 0$, and since $C''_s$ is bounded above by a constant times $\Length(G_s)$ it follows that $\mu^{\vphantom{3}}_s C''_s \to 0$ as $s \to -\infinity$.

Let $\wt\lambda_s$ be a leaf of $\wt\L_s$ covering a leaf $\lambda_s$ of $\L_s$. There is a subpath $\wt\alpha_s$ of $\wt\lambda_s$ of length at most $C''_s$ that decomposes as $\wt\alpha_s = \wt E_1 \wt\alpha' \wt E_2 \wt\alpha'' \wt E_3$ where $\wt E_1, \wt E_2, \wt E_3$ are lifts of the same edge $E$ of $G_s$ with length at least $C_s$.  If every crossing of $E$ by the leaf $\lambda_s$ was immediately followed by $\alpha'' E$ then $\lambda_s$ would be the periodic leaf $(\alpha'' E)^\infinity$. Since $\lambda_s$ is not periodic it follows that there is a covering translate $\wt\lambda_s'$ of $\wt\lambda_s$ that contains $\wt E_2$ but not $\wt E_2 \wt\alpha'' \wt E_3$. Similarly, there is a covering translate $\wt\lambda_s''$ that contains $\wt E_2$ but not $\wt E_1 \wt \alpha' \wt E_2$. Thus the arc $\wt\sigma_s = \wt\lambda_s \intersect \wt\lambda_s' \intersect \wt\lambda_s''$ satisfies
$$C_s \le \Length(\wt\sigma_s;T_s) \le C''_s
$$
and so
$$\Length(\wt\sigma_s;T^*_s) \le \mu^{\vphantom{3}}_s C''_s \to 0 \quad\text{as}\quad s \to -\infinity
$$

Let $\beta(s) = f_{0s}(\wt\sigma_s)$. For $s$ sufficiently close to $-\infinity$ we have $$\Length(\beta(s);T_0) = \Length(\wt\sigma_s;T_s) \ge C_s \ge K
$$
and so
$$\Length\left(\left(\beta(s)\right)^\#_s;T^*_s\right) > 3B
$$
On the other hand, the geodesics $\wt\lambda_s,\wt\lambda'_s,\wt\lambda''_s$, when mapped from $T_s$ to $T_0$ via $f_{0s}$, give geodesics in $T_0$ denoted $\wt\lambda_0,\wt\lambda'_0,\wt\lambda''_0$, and by construction we have
$$\beta(s) \subset \wt\lambda_0 \intersect \wt\lambda'_0 \intersect \wt\lambda''_0
$$
Applying Uniform Bounded Cancellation to the map $\phi_s$ it follows that if $s \le s_0$ then $(\beta(s))^\#_s = (\phi_s)_\#(\beta(s))$ is contained in the $B$-neighborhood of $\wt\lambda_s \intersect \wt\lambda'_s \intersect \wt\lambda''_s = \wt\sigma_s$, and so for $s$ sufficiently close to $-\infinity$ we have
$$\Length\left( (\beta(s))^\#_s; T^*_s \right) < 2B + \Length(\wt\sigma_s;T^*_s) < 3B
$$
which is a contradiction. This contradiction shows that $\L$ must have length zero in $T_-$. Since $\L$ is minimal, Lemma~\ref{LemmaMinimalLamination} shows that $\L = \Lambda_-$, completing Step~2.

\subsection{Weak train tracks on fold lines} 
\label{SectionStep3}
Continuing with the notations of Sections~\ref{SectionDirectLimits} and~\ref{SectionLegalLams}, we now carry out Step~3, completing the proof of Theorem~\ref{TheoremFoldLines}. 

By Section~\ref{SectionDirectLimits} the direct limit of the system of $F_r$-equivariant maps $f_{st} \from T_t \to T_s$ is $T_+$, and we have direct limit maps $f_s \from T_s \to T_+$, each of which is $F_r$-equivariant and an isometry on each edge. Fix a leaf $\ell$ of $\Lambda_-$ and $s \in \reals$, and let $\ell_s \subset T_s$ denote the realization of $\ell$ in $T_s$. We must show that $f_s \restrict \ell_s$ is an isometry from $\ell_s$ into $T_+$, that is, for each $x,y \in \ell_s$ we must prove $d_{T_+}(f_s(x),f_s(y)) = d_{T_s}(x,y)$. In Section~\ref{SectionLegalLams} we proved that $\Lambda_- = \L$ is the unique minimal legal sublamination. The geodesic $\ell_{s}$ is therefore legal, and so $f_{ts} \restrict \ell_s$ is an isometry for each $t \ge s$. For any $t \ge s$ it follows that  $d_{T_t} (f_{ts}(x),f_{ts}(y)) = d_{T_s}(x,y)$. Since the direct sequence $(f_{ts}(x))_{t \ge s}$ determines $f_s(x) \in T_+$, and the sequence $(f_{ts}(y))_{t \ge s}$ determines $f_s(y)$, it follows that $d_{T_+}(f_s(x),f_s(y))$ equals the limit of the constant sequence $d_{T_t} (f_{ts}(x),f_{ts}(y))$ which equals $d_{T_s}(x,y)$, proving that $d_{T_+}(f_s(x),f_s(y)) = d_{T_s}(x,y)$. This proves that each map $f_s \from G_s \to T_+$ is a $\Lambda_-$ isometry, and so $G_s$ is a weak train track.

\bibliographystyle{amsalpha} 
\bibliography{mosher}        


\end{document}